\documentclass[11pt]{article}
\usepackage{graphicx}

\usepackage[utf8]{inputenc}
\usepackage[T1]{fontenc}

\usepackage[]{amsfonts}
\usepackage{amssymb}
\usepackage{amsmath}
\def\couleur(#1 #2 #3)
	{
\special{ps::[end]}}

\def\bx#1{\setbox1=\hbox{\kern3pt{#1}\kern3pt}			% Make a box. Close it by "}"
 \dimen1=\ht1 \advance\dimen1 by 3pt \dimen2=\dp1 \advance\dimen2 by 3pt
 \setbox1=\hbox{\vrule height\dimen1 depth\dimen2\box1\vrule}%
 \setbox1=\vbox{\hrule\box1\hrule}%
 \advance\dimen1 by .4pt \ht1=\dimen1
 \advance\dimen2 by .4pt \dp1=\dimen2 \box1\relax}

\def\wbb#1{\kern#1em}
\def\vci{\vrule  width.02em height1.47ex depth-.0ex}		% le 1 en blackboard
\def\11{{\rm\wbb{.2}\vci\wbb{-.37}1}}

\def\underset#1#2{\mathrel{\mathop{\kern0pt #2}\limits_{#1}}}

\def\overset#1#2{\mathrel{\mathop{\kern0pt #2}\limits^{#1}}}

\newtheorem{Thrm}{Theorem}[section]
\newtheorem{Lmm}[Thrm]{Lemma}
\newtheorem{Dfnt}[Thrm]{Definition}
\newtheorem{Prps}[Thrm]{Proposition}
\newtheorem{Crll}[Thrm]{Corollary}
\newtheorem{Rmrq}[Thrm]{Remark}

\begin{document}

\title{Almost strictly pseudo-convex domains. Examples and Application.}

\author{Eric Amar}
\maketitle
 \ \par 

\tableofcontents
\ \par 
\renewcommand{\abstractname}{Abstract}

\begin{abstract}
\quad \quad 	In this work we introduce a class of smoothly bounded domains
  $\displaystyle \Omega $  in  ${\mathbb{C}}^{n}$  with few non
 strictly pseudo-convex points in  $\displaystyle \partial \Omega
 $  with respect to a certain Minkowski dimension. We call them
 almost strictly pseudo-convex, {\bf aspc}. For these domains
 we prove that a canonical measure associated to a separated
 sequence of points in  $\displaystyle \Omega $  which projects
 on the set of weakly pseudo-convex points is automatically 
 a geometric Carleson measure. This class of {\bf aspc} domains
 contains of course strictly pseudo-convex domains but also pseudo-convex
 domains of finite type in  ${\mathbb{C}}^{2},$  domains locally
 diagonalizable, convex domains of finite type in  ${\mathbb{C}}^{n},$
  domains with real analytic boundary and domains like  $\displaystyle
 \ \left\vert{z_{1}}\right\vert ^{2}+\exp  \lbrace 1-\left\vert{z_{2}}\right\vert
 ^{-2}\rbrace <1,$  which are not of finite type.\ \par 
\quad \quad 	As an application we study interpolating sequences for convex
 domains of finite type in  ${\mathbb{C}}^{n}.$   After proving
 an embedding Carleson type theorem, we get that if  $\Omega
 $  is a convex domain of finite type in  ${\mathbb{C}}^{n}$
  and if  $S\subset \Omega $  is a dual bounded sequence of points
 in  $H^{p}(\Omega ),$  if  $p=\infty $  then for any  $q<\infty
 ,\ S$  is  $H^{q}(\Omega )$  interpolating with the linear extension
 property and if  $p<\infty $  then  $S$  is  $H^{q}(\Omega )$
  interpolating with the linear extension property, provided
 that  $q<\min  (p,\ 2).$ \ \par 
\end{abstract}

\section{Introduction.}
\quad 	The aim of this work is to study a classical problem in harmonic
 analysis  and complex variables namely the interpolating sequences
 in some domains in  ${\mathbb{C}}^{n}.$  In order to do this
 we shall develop notions and tools well adapted to this aim,
 and which may be useful in other areas.\ \par 
\quad \quad 	The first notion we shall study is the existence of a good family
 of polydiscs.\ \par 
Troughout this work, domain will mean bounded open connected
 pseudo-convex set with  ${\mathcal{C}}^{\infty }$  smooth boundary.
 The Lebesgue measure on a manifold of real dimension  $k$  will
 be denoted  $\displaystyle \sigma _{k}.$  We also use the notation
 in formula :  $\displaystyle A::B...$  which means  $A$  such
 that  $B...$ \ \par 
\quad \quad 	Let  $\Omega $  be a domain in  ${\mathbb{C}}^{n}$  and  $a\in
 {\mathcal{U}}\subset \Omega ,$  where  ${\mathcal{U}}$  is a
 neighborhood of  $\partial \Omega $  in  $\Omega $  such that
 we have a well defined normal projection  $\pi $  on  $\partial
 \Omega \ ;$  we set  $\alpha :=\pi (a)$  and to this point 
 $\alpha \in \partial \Omega ,$  we associate a multi-index 
 $\displaystyle m(\alpha )=(1,\ m_{2}(\alpha ),\ ...,\ m_{n}(\alpha
 )),\ m_{1}\leq m_{2}\leq \cdot \cdot \cdot \leq m_{m}$  and
 a orthonormal basis  $b(\alpha )=(L_{1},...,L_{n})$  of  ${\mathbb{C}}^{n}$
  such that  $L_{1}$  is in the complex normal at  $\alpha $
  to  $\partial \Omega $  and  $(L_{2},...,L_{n})$  is a basis
 of the complex tangent space at  $\alpha $  to  $\partial \Omega
 .$  We set  $\displaystyle m(a):=m(\pi (a)),\ d(a):=d(a,\partial
 \Omega ).$ \ \par 
Now we define a polydisc  $Q_{a}(\delta )$  centered at  $a$
  of parameter  $\delta >0$  such that it has a radius  $\delta
 d(a)$  in the  $L_{1}$  direction and radii  $\delta d(a)^{1/m_{j}(a)},\
 j=2,...,\ n$  along  the  $L_{j}$  complex direction. We shall
 say that we have a "good family"  ${\mathcal{Q}}$  of polydiscs
 if these polydiscs reflect well the geometry of the domain,
 i.e. there is a parameter  $\delta _{0}>0$  such that  $\forall
 a\in {\mathcal{U}},\ Q_{a}(\delta _{0})\subset \Omega $  and
  $\displaystyle M({\mathcal{Q}}):=\sup _{\alpha \in \partial
 \Omega }m_{n}(\alpha )<\infty .$  For instance if  $\alpha =\pi
 (a)$  is a point of strict pseudo-convexity then we have that
  $m_{2}(a)=\cdot \cdot \cdot =m_{n}(a)=2.$ \ \par 
\quad \quad 	This notion of good family  ${\mathcal{Q}}$  of polydiscs is
 strongly inspired by the work of Catlin~\cite{Catlin84}.\ \par 
\quad \quad 	This good family  ${\mathcal{Q}}$  allows us to define separated
 sequences of points in  $\Omega .$ \ \par 

\begin{Dfnt}
Let  $\Omega $  be a domain in  ${\mathbb{C}}^{n}$  with a good
 family of polydiscs  ${\mathcal{Q}}.$  We shall say that a sequence
 of points  $S\subset \Omega $  is  $\delta $  {\bf separated}
 if any two distinct points in  $S$  are center of disjoint polydiscs
 in the family  ${\mathcal{Q}}$  with parameter  $\delta .$ 
\end{Dfnt}
Associated to this good family  ${\mathcal{Q}}$  we define  $\displaystyle
 \lambda (a):=\sum_{j=2}^{n}{\frac{1}{m_{j}(a)}}.$  If  $S$ 
 is a separated sequence of points in  $\Omega ,$  we define
 its canonical measure to be\ \par 
\quad \quad \quad \quad 	\begin{equation}  \mu _{S}:=\sum_{a\in S}{d(a)^{1+2\lambda (a)}\delta
 _{a}},\label{aspcGlo62}\end{equation}\ \par 
where  $\delta _{a}$  is the Dirac measure at  $a.$  For instance
 if  $\Omega $  is a strictly pseudo-convex domain, then  $m(a)=(1,2,...,2)$
  and  $\displaystyle \mu _{S}=\sum_{a\in S}{d(a)^{n}\delta _{a}}.$ \ \par 
\quad \quad 	We shall see that the sequences of points we are interested
 in are contained in the zero set of holomorphic functions.\ \par 
\quad Let  $u$  be a holomorphic function in a domain  $\displaystyle
 \Omega ,\ u\in {\mathcal{H}}(\Omega ),$  set  $X:=u^{-1}(0)$
  its zero set and  $\displaystyle \Theta :=\partial \bar \partial
 \ln \left\vert{u}\right\vert $   its associated  $(1,1)$  current
 of integration. As usual we have that  $\displaystyle \mathrm{T}\mathrm{r}\Theta
 (z)$  is the trace of the associated matrix and we have (~\cite{LelongGrum86}
 p 55))  $\displaystyle \mathrm{T}\mathrm{r}\Theta (z)=\Delta
 \ln \left\vert{u(z)}\right\vert .$ \ \par 
We shall define a class of such zero sets which contains the
 zero sets of Nevanlinna functions.\ \par 

\begin{Dfnt}
A holomorphic divisor  $X$  in the domain  $\Omega $  is in the
 {\bf Blaschke class} if, with  $\Theta $  its associated  $(1,1)$
  current of integration and  $\displaystyle d(z):=d(z,\Omega
 ^{c})$  the euclidean distance to the boundary,\par 
\quad \quad \quad \quad \quad  $\displaystyle \ {\left\Vert{\Theta }\right\Vert}_{B}:=\int_{\Omega
 }{d(z)\mathrm{T}\mathrm{r}\Theta (z)}<\infty .$ 
\end{Dfnt}
\quad \quad 	We shall need to study sequences of points contained in such
 sets ; let  $\displaystyle \sigma _{k}$  be the  Lebesgue measure
 on manifold of real dimension  $\displaystyle k.$  	Let  $\Omega
 $  be a domain equipped with a good family  ${\mathcal{Q}}$
  of polydiscs and  $X$  a divisor in  $\Omega .$  We set for
  $a\in \Omega ,\ X_{a}:=X\cap Q_{a}(\delta ),\ X_{a}^{j}$  the
 projection of  $\displaystyle X_{a}$  on  $E_{j}:=\lbrace z\in
 {\mathbb{C}}^{n}::z_{j}=0\rbrace $  in the coordinates in the
 basis  $\displaystyle b(\alpha )$  associated to  $\displaystyle
 \alpha =\pi (a),$  and  $\displaystyle A_{j}(X_{a}):=\sigma
 _{2n-2}(X_{a}^{j}).$  As usual  $\displaystyle \sigma _{2n-2}(X_{a})$
  is the measure of the regular points in  $\displaystyle X_{a},$
  as defined in~\cite{LelongGrum86}, proposition 2.48, p55.\ \par 
\quad \quad 	We get\ \par 

\begin{Thrm}
(Discretized Blaschke condition) Let  $u$  be holomorphic in
  $\Omega ,\ X:=u^{-1}(0)$  and  $\Theta :=\partial \bar \partial
 \ln \left\vert{u}\right\vert $  its current of integration ;
 suppose that  $\Theta $  is in the Blaschke class. Let  $S$
  be a  $\delta $  separated sequence in  $X$  with respect to
 a good family  ${\mathcal{Q}}$  of polydiscs with parameter
  $\delta _{0}.$  Then we have
\end{Thrm}
\quad \quad \quad \quad \quad \quad \quad  \[\displaystyle \ \sum_{a\in S}{d(a)\sigma _{2n-2}(X_{a})}\leq
 \frac{2}{\delta _{0}}{\left\Vert{\Theta }\right\Vert}_{B}.\] \ \par 
\quad \quad 	To go further and get the Malliavin discretized condition, {\sl
 with the right control on the constants,} we need to introduce
 {\bf  quasi convex domains} with respect to the good family
  ${\mathcal{Q}},$  i.e.  ${\mathcal{Q}}$  {\bf quasi convex
 domains}. This is a class of domains containing the convex ones
 and the lineally convex ones and adapted to our aim. They will
 be defined precisely by definition~\ref{3_BlasMall0}.\ \par 

\begin{Thrm}
(Discretized Malliavin condition) Let  $\Omega =\lbrace \rho
 <0\rbrace $  be a domain equipped with a good family  ${\mathcal{Q}}$
  of polydiscs with parameter  $\delta _{0}$  and which is  ${\mathcal{Q}}$
  quasi convex. Let  $\Theta $  be a current in the Blaschke
 class and  $S$  a  $\delta $  separated sequence in  $X\cap
 {\mathcal{U}}.$  Then we have\par 
\quad \quad \quad \quad \quad \quad \quad  $\displaystyle \ \sum_{a\in S}{\sum_{j=2}^{n}{A_{j}(X_{a})}}\leq
 C{\left\Vert{\Theta }\right\Vert}_{B},$ \par 
where  $C$  is a constant depending only on the  ${\mathcal{M}}(Q)+1$
  first order derivatives of  $\rho $  and on  $\delta ,\ \delta
 _{0},$  and the constant of quasi convexity.
\end{Thrm}
\quad \quad 	Together these two results gives the following theorem\ \par 

\begin{Thrm}
~\label{aspcGlo79}Let  $\Omega $  be a domain equipped with a
 good family  ${\mathcal{Q}}$  of polydiscs such that  $\Omega
 $  is  ${\mathcal{Q}}$  quasi convex ; let  $S$  be a  $\delta
 $  separated sequence of points which is contained in the Blaschke
 divisor  $X.$  Then\par 
\quad \quad \quad $\displaystyle \delta ^{2n-2}\sum_{a\in S}{d(a)^{n}}\leq \gamma
 (\Omega ){\left\Vert{\Theta }\right\Vert}_{B},$ \par 
where  $\displaystyle \gamma (\Omega )$  depends only on the
  ${\mathcal{C}}^{{\mathcal{M}}(Q)+1}$  norm of  $\rho ,$  on
  $\displaystyle n$  and  $\delta _{0},$  the parameter of the
 family  ${\mathcal{Q}},$  and on the constant of quasi convexity.
\end{Thrm}
\quad 	We have that  $\displaystyle 1+2\lambda (a)\leq n$  and equality
 for a point  $a$  such that  $\displaystyle \pi (a)$  is a strictly
 pseudo-convex point, hence in general this is not enough to
 deal with all types of sequence in  $\displaystyle \Omega .$
  So we are lead to introduce a class of domains with "few" points
 non strictly pseudo-convex, i.e. few "bad" points. If  $\Omega
 $  is a domain in  ${\mathbb{C}}^{n},$  throughout this work
  $W\subset \partial \Omega $  will denote the set of non strictly
 pseudo-convex points of  $\displaystyle \partial \Omega .$ \ \par 
\quad Let  $\alpha \in \partial \Omega \ ;$  by linear change of variables
 we can suppose that  $\alpha =0\in \partial \Omega \subset {\mathbb{C}}^{n},\
 z_{1}=0$  is the equation of the complex tangent space. The
 projection  $\pi $  locally near  $0\in \partial \Omega $  can
 be seen as a  ${\mathcal{C}}^{\infty }$  diffeomorphism  $\displaystyle
 \tilde \pi \ :\ \partial \Omega \rightarrow T_{0}(\partial \Omega
 ),\ \tilde \pi :=(\pi _{\mid T_{0}(\partial \Omega )})^{-1}.$ \ \par 

\begin{Dfnt}
The pseudo-convex domain  $\Omega $  in  ${\mathbb{C}}^{n}$ 
 is said to be {\bf almost stricly pseudo-convex}, {\bf aspc}
 at  $0\in \partial \Omega $  if there is a neighbourhood  $V_{0}$
  of  $0$  and a basis  $\displaystyle b:=\lbrace L_{1},...,L_{n}\rbrace
 $  of  ${\mathbb{C}}^{n},$  still with  $\displaystyle L_{1}$
  a complex normal unit vector,  such that, with  $\displaystyle
 (z_{1},...,z_{n})$  its associated coordinates, the slices\par 
\quad \quad \quad $\tilde \pi (W\cap V_{0})\cap \lbrace z_{1}=0\rbrace \cap \lbrace
 z_{2}=a_{2}\rbrace \cap \cdot \cdot \cdot \cap \lbrace z_{n-1}=a_{n-1}\rbrace
 $ \par 
have homogeneous Minkowki dimension less than  $2-\beta ,\ \beta >0.$ \par 
\quad $\Omega $  is said to be {\bf aspc} if this is true for all points
 in  $\partial \Omega $  with the same  $\beta >0.$ 
\end{Dfnt}
This means that we need only to find a particular coordinate
 system  $\displaystyle z_{1},...,z_{n}$  such that the slices
 of non s.p.c. points along the  $\displaystyle z_{n}$  direction
 of the tangent space to  $\displaystyle \partial \Omega $  have
 small Minkowski dimension.\ \par 
Of course the strictly pseudo-convex domains are {\bf aspc} because
  $W=\emptyset .$  The {\sl homogeneous Minkowski dimension}
 is defined precisely in section~\ref{4_AspcDomain30} and it
 quantifies the fact that bad points are few.\ \par 
\quad \quad 	This class of domains contains a large family of interesting
 domains such as strictly pseudo-convex domains, convex domains
 of finite type, etc ..., as shown in section~\ref{5_examAspc30}.\ \par 
And also non finite type domains as  $\lbrace z\in {\mathbb{C}}^{2}::\
 \left\vert{z_{1}}\right\vert ^{2}+\exp (1-\left\vert{z_{2}}\right\vert
 ^{-2})<1\rbrace .$ \ \par 
\quad \quad 	Usually we think that strictly pseudo-convex points are easier
 to deal with than non strictly pseudo-convex ones but for these
 domains and the properties we are interested in, this is not
 the case. In fact we have a good control on what happen for
 points projecting on weakly pseudo-convex points.\ \par 

\begin{Thrm}
~\label{aspc0}Let  ${\mathcal{Q}}$  be a good family of polydiscs
 on a {\bf aspc} domain  $\Omega $  in  ${\mathbb{C}}^{n},$ 
 and  $S$  be a  $\delta $  separated sequence of points in 
 $\Omega .$  Let  $W$  be the set of non strictly pseudo-convex
 points on  $\partial \Omega .$  If  $\pi (S\cap {\mathcal{U}})\subset
 V\cap W,$  where  $V$  is an open set of  $\partial \Omega ,$
  then we have :\par 
\quad \quad \quad \quad 	\begin{equation}  \ \sum_{a\in S\cap {\mathcal{U}}}{d(a)^{1+2\lambda
 (a)}}=\delta ^{-2n}\sum_{a\in S\cap {\mathcal{U}}}{\sigma _{2n}(Q_{a}(\delta
 ))}\leq C(\Omega )\frac{\sigma _{2n-1}(V)}{\delta ^{2}},\label{aspcGlo65}\end{equation}\par
 
where  $C(\Omega )$  depends only on  $\rho ,\ n,\ {\mathcal{M}}(Q),$
  and the constant  $\beta $  in the Minkowski dimension of 
 $W\subset \partial \Omega .$ 
\end{Thrm}
\quad \quad 	In fact this theorem says that the canonical measure associated
 to such a sequence is a geometric Carleson measure. So for these
 domains it remains to concentrate only on points which project
 on strictly pseudo-convex points on  $\partial \Omega .$  As
 an application we get :\ \par 

\begin{Thrm}
~\label{1_introduction31}Let  $\Omega $  be a {\bf aspc} domain
 in  ${\mathbb{C}}^{n}.$  Let  ${\mathcal{Q}}=\lbrace Q_{a}(\delta
 _{0})\rbrace _{a\in \Omega }$  be a good family of polydiscs
 for  $\Omega $  and suppose that  $\Omega $  is  ${\mathcal{Q}}$
  quasi convex. Let  $S$  a  $\delta $  separated sequence of
 points contained in a divisor  $X$  of the Blaschke class of
  $\Omega $  which projects on the open set  ${\mathcal{V}}\subset
 \partial \Omega .$  Then we have\par 
\quad \quad \quad \quad \quad  $\displaystyle \ \sum_{a\in S}{d(a)^{1+2\lambda (a)}}\leq \gamma
 (\Omega ){\left\Vert{\Theta _{X}}\right\Vert}_{B}+C(\Omega )\sigma
 ({\mathcal{V}})<\infty .$ 
\end{Thrm}
\quad 	The interpolating sequences are defined via the Hardy spaces
 of the domain  $\Omega .$ \ \par 

\begin{Dfnt}
Let  $\Omega $  be a domain in  ${\mathbb{C}}^{n}$  defined by the function\par 
\quad \quad \quad \quad \quad  $\rho \in {\mathcal{C}}^{\infty }({\mathbb{C}}^{n}),\ \Omega
 :=\lbrace z\in {\mathbb{C}}^{n}::\rho (z)<0\rbrace ,\ \forall
 z\in \partial \Omega ,\ \partial \rho (z)\neq 0.$ \par 
Let  $f$  be a holomorphic function in  $\Omega ,$  we say that
  $f$  is in the {\bf Hardy class}  $H^{p}(\Omega )$  if\par 
\quad \quad \quad $\displaystyle \ {\left\Vert{f}\right\Vert}_{p}^{p}:=\sup  _{\epsilon
 >0}\int_{\lbrace \rho (z)=-\epsilon \rbrace }{\left\vert{f(z)}\right\vert
 ^{p}\,d\sigma _{\epsilon }(z)}<\infty .$ \par 
We say that  $f$  is in the {\bf Nevanlinna class}  ${\mathcal{N}}(\Omega
 )$  if\par 
\quad \quad \quad $\displaystyle \ {\left\Vert{f}\right\Vert}_{{\mathcal{N}}}=\sup
 _{\epsilon >0}\int_{\lbrace \rho (z)=-\epsilon \rbrace }{\log
 ^{+}\left\vert{f(z)}\right\vert \,d\sigma _{\epsilon }(z)}<\infty .$ 
\end{Dfnt}
Here  $\,d\sigma _{\epsilon }$  is the Lebesgue measure on the
 smooth manifold  $\lbrace \rho (z)=-\epsilon \rbrace $  for
  $\epsilon $  small enough.\ \par 
These spaces are independent of the choice of the defining function~\cite{Stein72}.\
 \par 
\quad \quad 	As we shall see the study of interpolating sequences is intimately
 linked to  $p$  Carleson measures.\ \par 

\begin{Dfnt}
Let  $\lambda $  be a positive Borel measure on the domain  $\Omega
 $  and  $p\geq 1.$  We shall say that  $\lambda $  is a  $p$
  {\bf Carleson measure} in  $\Omega $  if :\par 
\quad \quad \quad $\displaystyle \exists C_{p}>0,\ \forall f\in H^{p}(\Omega ),\
 \int_{\Omega }{\left\vert{f}\right\vert ^{p}\,d\lambda }\leq
 C_{p}^{p}{\left\Vert{f}\right\Vert}_{H^{p}}^{p}.$ \par 
This means that we have a continuous embedding of  $H^{p}(\Omega
 )$  in  $L^{p}(\lambda ).$ 
\end{Dfnt}
Usually we have only a geometric condition to work with :\ \par 

\begin{Dfnt}
Let  $\lambda $  be a borelian positive measure on the domain
  $\Omega $  equipped with a good family of polydiscs  ${\mathcal{Q}}.$
  We shall say that  $\lambda $  is a {\bf geometric Carleson
 measure} in  $\Omega $  if :\par 
\quad \quad \quad $\displaystyle \exists C>0::\forall a\in \Omega ,\ \lambda (\Omega
 \cap Q_{a}(2))\leq C\sigma (\partial \Omega \cap Q_{a}(2)).$ 
\end{Dfnt}
So we need a way to go from geometric Carleson measures to  $p$
  Carleson measures and this is why we need to restrict to convex
 domains of finite type. For them we have a Carleson embedding theorem.\ \par 

\begin{Thrm}
~\label{1_introduction32}Let  $\Omega $  be a convex domain of
 finite type. If the measure  $\lambda $  is a geometric Carleson
 measure we have\par 
\quad \quad \quad $\displaystyle \forall p>1,\ \exists C_{p}>0,\ \forall f\in H^{p}(\Omega
 ),\ \int_{\Omega }{\left\vert{f}\right\vert ^{p}\,d\lambda }\leq
 C_{p}^{p}{\left\Vert{f}\right\Vert}_{H^{p}}^{p}.$ \par 
Conversely if the positive measure  $\lambda $  is  $p$  Carleson
 for a  $p\in \lbrack 1,\ \infty \lbrack ,$  then it is a geometric
 Carleson measure, hence it is  $q$  Carleson for any  $q\in
 \rbrack 1,\ \infty \lbrack .$ 
\end{Thrm}
\ \par 
\quad \quad 	It remains to see when the canonical measure associated to a
 separated sequence is a geometric Carleson measure. In the unit
 ball  ${\mathbb{B}}$  of  ${\mathbb{C}}^{n}$  this is done by
 an easy generalization of a lemma of Garnett : a measure  $\lambda
 $  is Carleson in the ball  ${\mathbb{B}}$  iff all its images
 under the automorphisms of  ${\mathbb{B}}$  are uniformly bounded
 measures~\cite{AmarWirtBoule07}. In a general domain there is
 only the identity as automorphism, so we have to overcome this issue.\ \par 
\quad \quad 	We do it by building sub domains associated to a point  $a\in
 \Omega $  and which are equivalent to Carleson windows. This
 can be done with the right control of the constants if  $\Omega
 $  is a {\sl well balanced} domain ; this notion will be defined
 later. Convex domains, lineally convex domains are well balanced.\ \par 
\quad \quad 	The space  $H^{2}(\Omega )$  is a sub space of the Hilbert space
  $L^{2}(\partial \Omega )$  hence there is a orthogonal projection
  $S\ :\ L^{2}(\partial \Omega )\rightarrow H^{2}(\Omega ).$
  We shall denote  $k_{a}(z)$  the kernel of this (Szeg\"o) projection,
 it is a reproducing kernel for  $H^{2}(\Omega ).$ \ \par 
Now we have the tools needed to deal with interpolating sequences.\ \par 

\begin{Dfnt}
We say that the sequence  $S$  of points in  $\Omega $  is  $H^{p}(\Omega
 )$  {\bf interpolating} if\par 
\quad \quad \quad $(i)\ \forall a\in S,\ k_{a}\in H^{p'}(\Omega )\ ;$  (this is
 always true if  $p\geq 2.$ )\par 
\quad \quad \quad $\displaystyle (ii)\ \forall \lambda \in \ell ^{p}(S),\ \exists
 f\in H^{p}(\Omega )::\forall a\in S,\ f(a)=\lambda _{a}{\left\Vert{k_{a}}\right\Vert}_{p'},$
 \par 
with  $p'$  the conjugate exponent of  $p,$  i.e.  $\displaystyle
 \ \frac{1}{p}+\frac{1}{p'}=1.$ 
\end{Dfnt}
We have a weaker notion than interpolation :\ \par 

\begin{Dfnt}
We shall say that the sequence  $S$  of points in  $\Omega $
  is {\bf dual bounded} in  $H^{p}(\Omega )$  if there is a bounded
 sequence of elements in  $H^{p}(\Omega ),\ \lbrace \rho _{a}\rbrace
 _{a\in S}\subset H^{p}(\Omega )$  which dualizes the associated
 sequence of reproducing kernels, i.e.\par 
\quad \quad \quad $(i)\ \forall a\in S,\ k_{a}\in H^{p'}(\Omega )\ ;$  (this is
 always true if  $p\geq 2.$ )\par 
\quad \quad \quad $\displaystyle (ii)\ \exists C>0::\forall a\in S,\ {\left\Vert{\rho
 _{a}}\right\Vert}_{p}\leq C,\ \forall a,b\in S,\ {\left\langle{\rho
 _{a},\ k_{b}}\right\rangle}=\delta _{a,b}{\left\Vert{k_{b}}\right\Vert}_{p'}.$ 
\end{Dfnt}
\quad Clearly if  $S$  is  $H^{p}(\Omega )$  interpolating then  $S$
  is dual bounded in  $H^{p}(\Omega )\ :$  just interpolate the
 basic sequence of  $\ell ^{p}(S).$  In the unit disc of  ${\mathbb{C}}$
  the converse is true, here we have a partial converse of this.\ \par 

\begin{Thrm}
~\label{1_introduction30}Let  $\Omega $  be a convex domain of
 finite type in  ${\mathbb{C}}^{n}$  and let  $S\subset \Omega
 $  be a dual bounded sequence of points in  $\displaystyle H^{p}(\Omega
 ),$  if  $p=\infty $  then for any  $q<\infty ,\ S$  is  $\displaystyle
 H^{q}(\Omega )$  interpolating ; if  $p<\infty $  then  $S$
  is  $\displaystyle H^{q}(\Omega )$  interpolating, provided
 that  $q<\min  (p,\ 2).$ 
\end{Thrm}
\quad \quad 	Let me give a rough sketch of the proof.\ \par 
Take a sequence  $S$  in the convex domain  $\displaystyle \Omega
 \ ;$  to apply a general result on interpolating sequences done
 in~\cite{AmarExtInt06} we need the following facts :\ \par 
\quad \quad  $\bullet $  a link between the  $H^{p}(\Omega )$  norm of the
 reproducing kernels  $\displaystyle k_{a}$  and the geometry
 of the boundary of  $\displaystyle \partial \Omega \ :$  the
  $p$  regularity of the domain  $\displaystyle \Omega $  which says\ \par 
\quad \quad \quad $\exists C>0::\forall a\in \Omega ,\ {\left\Vert{k_{a}}\right\Vert}_{p}^{-p'}\leq
 C\sigma (\partial \Omega \cap Q_{a}(2)),$ \ \par 
where  $p'$  is the conjugate exponent of  $p.$  We shall see
 that this is true for convex domain of finite type.\ \par 
\quad \quad  $\bullet $  Structural hypotheses for the Lebesgue measure on
  $\displaystyle \partial \Omega .$  These are reverse H\"older
 inequalities for the norms of the reproducing kernels  $\displaystyle
 k_{a}.$  We shall see that this is also true for convex domain
 of finite type.\ \par 
\quad \quad  $\bullet $  The fact that the canonical measure associated to
  $S,\ \mu _{S}:=\sum_{a\in S}{d(a)^{1+2\lambda (a)}\delta _{a}}$
  is  $q$  Carleson.\ \par 
And this is the main difficulty. To achieve this we use the fact
 that a convex domain of finite type is almost strictly pseudo-convex,
 so, with  $W$  the set of weakly pseudo-convex points in  $\displaystyle
 \partial \Omega ,$  we have that the measure	 	 $\displaystyle
 \mu _{b}:=\sum_{a\in S\cap \pi ^{-1}(W)}{d(a)^{1+2\lambda (a)}\delta
 _{a}}$  is already a geometric Carleson measure in  $\displaystyle
 \Omega $  by theorem~\ref{aspc0}.\ \par 
\quad \quad 	It remains to deals with the points which project on the strictly
 pseudo-convex points in  $\displaystyle \partial \Omega .$ \ \par 
By assumption  $\displaystyle S\backslash \lbrace a\rbrace $
  is contained in the zero set of  $\displaystyle \rho _{a}\in
 H^{p}(\Omega )\subset {\mathcal{N}}(\Omega ).$  So we can use
 theorem~\ref{aspcGlo79} to get, because a convex domain is quasi
 convex, that  $\displaystyle \mu :=\sum_{a\in S}{d(a)^{n}\delta
 _{a}}$  is a bounded measure in  $\displaystyle \Omega .$  To
 prove that  $\mu $  is a geometric Carleson measure we construct
 sub domains  $\displaystyle \Omega _{a}$  associated to points
  $\displaystyle a\in \Omega $  and which are comparable to the
 Carleson windows  $\displaystyle \Omega \cap Q_{a}(2).$  Because
 we have a precise estimate of the bound of  $\displaystyle \
 \sum_{a\in S}{d(a)^{n}\delta _{a}}$  in terms of  $\displaystyle
 \Omega $  and of the holomorphic function  $u$  whose zero set
 contains  $S$  we can apply what we have done to the sub domain
  $\displaystyle \Omega _{a}$  and get that  $\displaystyle \
 \sum_{b\in S\cap \Omega _{b}}{d(b)^{n}\delta _{b}}$  is bounded
 by a uniform constant times  $\displaystyle \sigma _{2n-1}(\partial
 \Omega _{a}\cap \partial \Omega )$  which means that  $\displaystyle
 \mu :=\sum_{a\in S}{d(a)^{n}\delta _{a}}$  is a geometric Carleson
 measure in  $\displaystyle \Omega .$ \ \par 
\quad \quad 	Now we use the Carleson embedding theorem~\ref{1_introduction32}
 to get that the measure  $\displaystyle \mu :=\sum_{a\in S}{d(a)^{n}\delta
 _{a}}$  is a  $q$  Carleson measure for any  $\displaystyle
 q\in \rbrack 1,\infty \lbrack .$ \ \par 
\quad \quad 	For "good points", i.e. those which project on strictly pseudo-convex
 ones we have that  $\displaystyle 1+2\lambda (a)=n,$  hence
 gluing with the estimate coming from the {\bf aspc} side, we
 get theorem~\ref{1_introduction30} as an application of the
 notion of {\bf aspc} domains.\ \par 
\ \par 
\quad \quad 	I am deeply grateful to the referee who not only had to deal
 with the mathematics in this paper but also gave me a lot of
 valuable suggestions on the presentation of it. Hence even if
 the results here are essentially the same as in the preprint~\cite{AmAspc09},
 the presentation is completely rewritten, the statements are
 precised and the proofs are detailed.\ \par 
\ \par 
\quad \quad 	The general organization is as follow.\ \par 
\ \par 
In section~\ref{2_BonFamille42} we define the good family  ${\mathcal{Q}}$
  of polydiscs in a domain  $\Omega $  and we give two characterisations
 of them :\ \par 
 $\bullet $ 	 an analytic one in term of finite linear type\ \par 
 $\bullet $  a geometric one in term of complex tangentially
 ellipsoid at every point  $\alpha \in \partial \Omega .$ \ \par 
\ \par 
In section~\ref{aspcGlo63} we define precisely the Blaschke class
 of divisors  $X$  in  $\Omega ,$  the notion of  ${\mathcal{Q}}$
  quasi convexity, and we prove the discretized Blaschke and
 Malliavin conditions.\ \par 
\ \par 
In section~\ref{4_AspcDomain30} we introduce the notion of almost
 strictly pseudo-convex domains and we use a nice theorem of
 Ostrowski to get theorem~\ref{1_introduction31}.\ \par 
\ \par 
In section~\ref{5_examAspc30} we prove that domains of finite
 type in  ${\mathbb{C}}^{2},$  locally diagonalizable domains,
 convex domains of finite type, domains with real analytic boundary,
 are all {\bf aspc} domains, together of course with the strictly
 pseudo-convex domains.\ \par 
\ \par 
In section~\ref{7_ConvFini40} we set the geometric properties
 we need for convex domain of finite type and in section~\ref{8_CarlesonMeas30}
 we study Carleson measures in such domains and state and prove
 the Carleson embedding theorem~\ref{1_introduction32}.\ \par 
\ \par 
In section~\ref{6_CarlDomain33} we construct the sub domain associated
 to a point  $a\in \Omega $  which is equivalent to the Carleson
 window  $\displaystyle Q_{a}(2)\cap \Omega $  and which allows
 us to overcome the lack of automorphisms.\ \par 
\ \par 
In section~\ref{9_InterSeq40} we define the notion of  $p$ {\sl
  regularity} making a link between the  $H^{p}(\Omega )$  norm
 of the reproducing kernels and the geometry of  $\displaystyle
 \partial \Omega .$  Then we prove Theorem~\ref{1_introduction30}
 via a tour around properties of reproducing kernels.\ \par 
\ \par 
Finally in the section~\ref{AG1} we state and prove the facts
 we need from potential theory.\ \par 
\vfill\eject\ \par 

\section{Good family of polydiscs.~\label{2_BonFamille42}}
\quad \quad 	In this section we shall study domains with a good family of
 polydiscs and get some properties of these domains we shall use later.\ \par 
\quad Let  $\Omega $  be a domain in  ${\mathbb{C}}^{n},$  recall that
 here, this means a bounded open connected set with a  ${\mathcal{C}}^{\infty
 }$  smooth boundary. Let  ${\mathcal{U}}$  be a neighbourhood
 of  $\partial \Omega $  in  $\Omega $  such that the normal
 projection  $\pi $  onto  $\partial \Omega $  is a smooth well
 defined application. For  $a\in \Omega $  set  $d(a):=d(a,\
 \Omega ^{c})$  the distance from  $a$  to the boundary of  $\Omega .$ \ \par 
\ \par 
We shall need the notion of a "good" family of polydiscs, directly
 inspired by the work of Catlin~\cite{Catlin84}.\ \par 
\quad Let  $\alpha \in \partial \Omega $  and let  $b(\alpha )=(L_{1},\
 L_{2},...,\ L_{n})$  be an orthonormal basis of  ${\mathbb{C}}^{n}$
  such that  $(L_{2},...,\ L_{n})$  is a basis of the tangent
 complex space  $T_{\alpha }^{{\mathbb{C}}}$  of  $\partial \Omega
 $  at  $\alpha \ ;$  hence  $L_{1}$  is the complex normal at
  $\alpha $  to  $\partial \Omega .$ \ \par 
Let  $m(\alpha )=(m_{1},\ m_{2},...,\ m_{n})\in {\mathbb{R}}^{n}$
  be a multi-index at  $\alpha $  with  $m_{1}=1,\ \forall j\geq
 2,\ m_{j}\geq 2.$ \ \par 
\quad For  $a\in {\mathcal{U}},$  let  $\alpha =\pi (a)\in \partial
 \Omega ,\ b(a):=b(\alpha ),\ m(a):=m(\alpha ),$  and  $\delta
 >0$  ; set  $\displaystyle Q_{a}(\delta ):=\prod_{j=1}^{n}{\delta
 D_{j}}$  the polydisc such that  $\delta D_{j}$  is the disc
 centered at  $a,$  parallel to  $L_{j}(\alpha )$  with radius
  $\delta {\times}d(a)^{1/m_{j}(\alpha )},\ j=1,...,n.$ \ \par 
\quad This way we have a family of polydiscs  ${\mathcal{Q}}:=\lbrace
 Q_{a}(\delta )\rbrace _{a\in {\mathcal{U}}}$  defined by the
 family of basis  $\lbrace b(\alpha )\rbrace _{\alpha \in \partial
 \Omega },$  the family of multi-indices  $\lbrace m(\alpha )\rbrace
 _{\alpha \in \partial \Omega }$  and the number  $\delta .$ \ \par 
\quad It will be useful to extend this family to the whole of  $\Omega
 .$  In order to do so let  $(z_{1},\ ...,\ z_{n})$  be the canonical
 coordinates system in  ${\mathbb{C}}^{n}$  and for  $a\in \Omega
 \backslash {\mathcal{U}},$  let  $Q_{a}(\delta )$  be the polydisc
 of center  $a,$  of sides parallel to the axis and radius  $\delta
 d(a)$  in the  $z_{1}$  direction and  $\delta d(a)^{1/2}$ 
 in the other directions. So the points  $a\in \Omega \backslash
 {\mathcal{U}}$  have automatically a "minimal" multi-index 
 $m(a)=(1,\ 2,\ ...,\ 2).$ \ \par 
\quad Now we can set\ \par 

\begin{Dfnt}
~\label{BonFamille39}We say that  ${\mathcal{Q}}$  is a {\bf
 "good family" of polydiscs} for  $\Omega $  if the  $m_{j}(a)$
  are uniformly bounded, i.e.  $\displaystyle M({\mathcal{Q}}):=\sup
 _{j=1,...,n,a\in \Omega }m_{j}(a)<\infty ,$  and if there exists
  $\delta _{0}>0,$  called the parameter of the family  ${\mathcal{Q}},$
  such that all the polydiscs  $\lbrace Q_{a}(\delta _{0})\rbrace
 _{a\in \Omega }$  of  ${\mathcal{Q}}$  are contained in  $\Omega
 .$  In this case we call  $m(a)$  the multi-type at  $a$  of
 the family  ${\mathcal{Q}}.$ 
\end{Dfnt}
\quad We notice that, for a good family  ${\mathcal{Q}},$  by definition
 the multi-type is always finite. Moreover there is no regularity
 assumptions on the way that the basis  $b(\alpha )$  varies
 with respect to  $\alpha \in \partial \Omega .$ \ \par 
\quad We can see easily that there is always good families of polydiscs
 in a domain  $\Omega $  in  ${\mathbb{C}}^{n}\ :$  for a point
  $\alpha \in \partial \Omega ,$  take any orthonormal basis
  $b(\alpha )=(L_{1},\ L_{2},...,\ L_{n}),$  with  $L_{1}$  a
 complex normal vector to  $\partial \Omega ,$  and the "minimal"
 multi-type  $m(\alpha )=(1,\ 2,...,\ 2).$  Then, because  $\Omega
 $  is of class  ${\mathcal{C}}^{2}$  and relatively compact,
 we have the existence of a uniform  $\delta _{0}>0$  such that
 the family  ${\mathcal{Q}}$  is a good one.\ \par 

\subsection{Examples of domains with a good family of polydiscs.}
\quad The stricly pseudo-convex domains in  ${\mathbb{C}}^{n}\ :$ 
 they have a good family of polydiscs associated with the best
 possible multi-type, the one defined by Catlin~\cite{Catlin84},
 which is also the "minimal" one in this case :\ \par 
\quad \quad \quad \quad $\forall a\in {\mathcal{U}},\ m_{1}=1,\ \ \forall j=2,...,n,\
 m_{j}(a)=2.$ \ \par 
Moreover these polydiscs are associated to the pseudo balls of
 a structure of spaces of homogeneous type (Koranyi-Vagi~\cite{KoranyiVagi71},
 Coifman-Weiss ~\cite{CoifWeiss71}).\ \par 
\quad The finite type domains in  ${\mathbb{C}}^{2}:$  also here we
 have the best possible multi-type and a structure of spaces
 of homogeneous type. (Nagel-Rosay-Stein-Wainger~\cite{NaRoStWa89}) \ \par 
\quad The bounded convex finite type domains in  ${\mathbb{C}}^{n}:$
  again we have the best multi-type and a structure of spaces
 of homogeneous type. (McNeal~\cite{McNeal94})\ \par 

\subsection{An analytical characterisation by linear finite
 type.~\label{fortePC74}}
\quad We shall recall precisely the definition of the multi-type~\cite{Catlin84}
 and the linear multi type (McNeal~\cite{McNeal92}, Yu~\cite{Yu92}).
 We shall take the definitions and the notations from J. Yu~\cite{Yu92}.\ \par 
\quad Let  $\Omega $  be a domain in  ${\mathbb{C}}^{n}$  defined by
 the function  $\rho ,$  and let  $p\in \partial \Omega $  be fixed.\ \par 
\quad Let  $\Gamma _{n}$  be the set of the  $n$ -tuples of numbers
  $\Lambda =(m_{1},...,m_{n})$  with  $1\leq m_{j}\leq \infty
 $  and such that\ \par 
(i)  $m_{1}\leq m_{2}\leq \cdot \cdot \cdot \leq m_{n}.$ \ \par 
(ii) for all  $k=1,...,n,$  either  $m_{k}=+\infty $  or there
 are non negative integers  $a_{1},...,a_{k}$  such that  $\displaystyle
 a_{k}>0$  and  $\displaystyle \ \sum_{j=1}^{k}{a_{j}/m_{j}}=1.$ \ \par 
This condition (ii) is automatically fulfilled in the case all
  $\displaystyle m_{j}$  are integers.\ \par 
\quad An element in  $\Gamma _{n}$  will be called a {\sl weight}.
 The set  $\Gamma _{n}$  of weights can be ordered lexicographically :\ \par 
\quad \quad \quad \quad \quad  $\Lambda =(m_{1},...,m_{n})<\Lambda '=(m'_{1},...,m'_{n})$ \ \par 
if there is a  $k$  such that  $\forall j<k,\ m_{j}=m'_{j}$ 
 and  $m_{k}<m'_{k}.$ \ \par 

\begin{Lmm}
~\label{2_BonFamille20}The entries  $\displaystyle m_{j}$  of
 a weight  $\displaystyle m=(m_{1},...,m_{n})$  are rational
 numbers. Given  $\displaystyle M>0$  there is only a finite
 number of weights  $\displaystyle m=(m_{1},...,m_{n})$  such
 that  $\displaystyle m_{n}\leq M.$  Moreover if  $\displaystyle
 m_{1}=1$  then  $\displaystyle m_{2}\in {\mathbb{N}}.$ 
\end{Lmm}
\quad \quad 	Proof.\ \par 
We have by (ii) that  $\displaystyle \exists a_{1}\in {\mathbb{N}}::\frac{a_{1}}{m_{1}}=1$
  hence  $m_{1}=a_{1}\in {\mathbb{N}}.$  Again by (ii)\ \par 
\quad \quad \quad \quad \quad  $\displaystyle \exists a_{1},a_{2}\in {\mathbb{N}}::\frac{a_{1}}{m_{1}}+\frac{a_{2}}{m_{2}}=1\Rightarrow
 a_{1}\leq m_{1}\leq M,\ a_{2}\leq m_{2}\leq M,$ \ \par 
hence we have only a finite number of possible  $\displaystyle
 m_{1},\ a_{1},\ a_{2}.$  For each of such possibility we have\ \par 
\quad \quad \quad \quad \quad  $\displaystyle \ \frac{1}{m_{2}}=\frac{1}{a_{2}}(1-\frac{a_{1}}{m_{1}}),$
 \ \par 
hence only one solution and a rational one.\ \par 
\quad \quad 	So we have only a finite number of solutions for  $\displaystyle
 m_{2}$  and all are rational numbers. We notice that if  $\displaystyle
 m_{1}=1,$  then  $\displaystyle a_{1}=0$  and  $\displaystyle
 m_{2}=a_{2}\in {\mathbb{N}}.$ \ \par 
Suppose now that  $\displaystyle m_{1},...,m_{k}$  are in finite
 number, then, as we just seen,  $\displaystyle a_{1},...,a_{k}$
  are also in finite number and  $\displaystyle a_{k+1}\leq m_{k+1}\leq
 M$  so only a finite number of  $\displaystyle a_{k+1}.$  Now
 as above for each choice of  $\displaystyle a_{1},...,a_{k+1},\
 m_{1},...,m_{k}$  we have only one solution  $\displaystyle
 m_{k+1}$  for\ \par 
\quad \quad \quad \quad \quad  $\displaystyle \ \frac{1}{m_{k+1}}=\frac{1}{a_{k+1}}(1-\frac{a_{1}}{m_{1}}-\cdot
 \cdot \cdot -\frac{a_{k}}{m_{k}})$ \ \par 
which is rational and the lemma is proved by induction.  $\blacksquare $ \ \par 
\quad A weight is said to be {\sl distinguished} if there exist holomorphic
 coordinates  $z_{1},...,z_{n},$  in a neighbourhood of  $p$
  with  $p$  mapped to the origin and such that :\ \par 
\quad \quad \quad \quad 	\begin{equation}  \ \sum_{i=1}^{n}{\frac{\alpha _{i}+\beta _{i}}{m_{i}}}<1\
 \Rightarrow \ \partial ^{\alpha }\bar \partial ^{\beta }\rho
 (p)=0,\label{aspcGlo0}\end{equation}\ \par 
where  $\displaystyle \partial ^{\alpha }:=\frac{\partial ^{\left\vert{\alpha
 }\right\vert }}{\partial z_{1}^{\alpha _{1}}\cdot \cdot \cdot
 \partial z_{n}^{\alpha _{n}}}$  and  $\displaystyle \bar \partial
 ^{\beta }:=\frac{\partial ^{\left\vert{\beta }\right\vert }}{\partial
 \bar z_{1}^{\beta _{1}}\cdot \cdot \cdot \partial \bar z_{n}^{\beta
 _{n}}}.$ \ \par 

\begin{Dfnt}
The {\bf multi-type}  ${\mathcal{M}}(\partial \Omega ,p)$  is
 the smallest weight  ${\mathcal{M}}:=(m_{1},...,m_{n})$  in
  $\Gamma _{n}$  (in lexicographic sense) such that  ${\mathcal{M}}\geq
 \Lambda $  for every distinguished weight  $\Lambda .$ 
\end{Dfnt}
Because  $\partial \Omega $  is smooth at  $p,$  we always have
  $m_{1}=1.$ \ \par 
\quad We call a weight  $\Lambda $  {\sl linearly distinguished} if
 there exists a complex linear change of variables near  $p$
  with  $p$  mapped to the origin and such that~(\ref{aspcGlo0})
 holds in these new coordinates.\ \par 

\begin{Dfnt}
The {\bf linear multi-type}  ${\mathcal{L}}(\partial \Omega ,\
 p)$  is the smallest weight  ${\mathcal{L}}:=(m_{1},...,m_{n})$
  such that  ${\mathcal{L}}\geq \Lambda $  for every linear distinguished
 weight  $\Lambda .$  We shall say that  $\Omega $  is of {\bf
 linear finite type} if\par 
\quad \quad \quad $\exists m\in {\mathbb{N}}::\forall p\in \partial \Omega ,\ {\mathcal{L}}(\partial
 \Omega ,\ p)\leq (m,...,m).$ 
\end{Dfnt}
\quad Clearly we have  ${\mathcal{L}}(\partial \Omega ,\ p)\leq {\mathcal{M}}(\partial
 \Omega ,\ p).$ \ \par 
\ \par 
\quad If, for  $p\in \partial \Omega $  fixed,  $\Omega $  is of linear
 finite type  ${\mathcal{L}}(\partial \Omega ,\ p)=(m_{1},...,m_{n}),$
  then  there is a  ${\mathbb{C}}$ -linear change of variables
 such that~\cite{Yu92} :\ \par 
\quad \quad \quad \quad \quad  $\displaystyle \ \sum_{i=1}^{n}{\frac{\alpha _{i}+\beta _{i}}{m_{i}}}<1\
 \Rightarrow \ \partial ^{\alpha }\bar \partial ^{\beta }\tilde
 \rho (0)=0,$ \ \par 
where  $\tilde \rho $  is the defining function of  $\Omega $
  in these new coordinates  $\zeta =(\zeta _{1},...,\zeta _{n}).$
  Set  $\displaystyle m_{n}':={\left\lceil{m_{n}}\right\rceil}=\min
 _{k\in {\mathbb{N}},k\geq m_{n}}k.$ \ \par 

\begin{Lmm}
We have\par 
\quad \quad \quad \quad \quad  $\tilde \rho (\zeta )=\Re \zeta _{1}+\sum_{2\leq \left\vert{\alpha
 }\right\vert +\left\vert{\beta }\right\vert \leq m_{n}'}{A_{\alpha
 \beta }\zeta ^{\alpha }\bar \zeta ^{\beta }}+o(\left\vert{\zeta
 }\right\vert ^{m_{n}'}),$ \par 
with  $A_{\alpha \beta }\neq 0$  only if  $\displaystyle \ \sum_{i=1}^{n}{\frac{\alpha
 _{i}+\beta _{i}}{m_{i}}}\geq 1.$ 
\end{Lmm}
\quad \quad 	Proof.\ \par 
We expand  $\tilde \rho $  by Taylor formula near  $0$  up to
 order  $m_{n}'$  and we compute\ \par 
\quad \quad \quad \quad \quad  $\displaystyle \partial ^{\alpha }\bar \partial ^{\beta }\tilde
 \rho (0)=\alpha !\beta !A_{\alpha \beta }.$ \ \par 
But if  $\displaystyle \ \sum_{i=1}^{n}{\frac{\alpha _{i}+\beta
 _{i}}{m_{i}}}<1$  then  $\displaystyle \partial ^{\alpha }\bar
 \partial ^{\beta }\tilde \rho (0)=0$  because the linear multi-type
 of  $\Omega $  at  $0$  is  $m=(1,\ m_{2},...,\ m_{n}).$ \ \par 
Because  $\displaystyle j\geq 2\Rightarrow m_{j}\geq 2,$  fixing
  $j$  $\geq 2$  and taking  $\displaystyle \alpha _{j}=1,\alpha
 _{i}=0$  for  $\displaystyle i\neq j$  and  $\displaystyle \beta
 _{i}=0$  for all  $i$  we get\ \par 
\quad \quad \quad \quad \quad  $\displaystyle \ \sum_{i=1}^{n}{\frac{\alpha _{i}+\beta _{i}}{m_{i}}}=\frac{1}{m_{j}}<1$
  hence  $\displaystyle \forall j\geq 2,\ \frac{\partial \tilde
 \rho }{\partial \zeta _{j}}(0)=0.$ \ \par 
Replacing  $\displaystyle \alpha _{j}$  by  $\displaystyle \beta
 _{j}$  we get  $\displaystyle \forall j\geq 2,\ \frac{\partial
 \tilde \rho }{\partial \bar \zeta _{j}}(0)=0$  hence the complex
 tangent plane to  $\displaystyle \partial \Omega $  at  $0$
  is still  $\displaystyle \zeta _{1}=0,$  and the  $\displaystyle
 \zeta _{j},\ j\geq 2$  are coordinates in the complex tangent space.\ \par 
So multiplying  $\zeta _{1}$  by a complex constant of module
  $1$  if necessary, we have\ \par 
\quad \quad \quad \quad \quad  $\displaystyle \tilde \rho (\zeta )=\Re \zeta _{1}+\sum_{2\leq
 \left\vert{\alpha }\right\vert +\left\vert{\beta }\right\vert
 \leq m_{n}'}{A_{\alpha \beta }\zeta ^{\alpha }\bar \zeta ^{\beta
 }}+o(\left\vert{\zeta }\right\vert ^{m_{n}'}),$ \ \par 
with  $\displaystyle \ \sum_{i=1}^{n}{\frac{\alpha _{i}+\beta
 _{i}}{m_{i}}}<1\Rightarrow \partial ^{\alpha }\bar \partial
 ^{\beta }\tilde \rho (0)=0.$   $\blacksquare $ \ \par 
\quad The aim of this sub-section is to show\ \par 

\begin{Thrm}
~\label{linkTfFpc20}If  $\Omega $  is a domain in  ${\mathbb{C}}^{n}$
  of finite linear type, then there is a good family  ${\mathcal{Q}}$
  of polydiscs such that the multi-type associated to  ${\mathcal{Q}}$
  is precisely the linear multi-type of  $\Omega .$ 
\end{Thrm}
\quad \quad 	Proof.\ \par 
Going back to the previous coordinates, this means that there
 are complex directions  $v_{1},\ v_{2},...,v_{n}$  with  $v_{1}$
  the complex normal at  $p,\ v_{2},...,\ v_{n}$  in the complex
 tangent space, such that:\ \par 
\quad \quad \quad \quad \quad  $\displaystyle \ \sum_{i=1}^{n}{\frac{\alpha _{i}+\beta _{i}}{m_{i}}}<1\
 \Rightarrow \ \partial _{v}^{\alpha }\bar \partial _{v}^{\beta
 }\rho (p)=0,$ \ \par 
with now   $\displaystyle \partial _{v}^{\alpha }:=\frac{\partial
 ^{\left\vert{\alpha }\right\vert }}{\partial v_{1}^{\alpha _{1}}\cdot
 \cdot \cdot \partial v_{n}^{\alpha _{n}}}$  and  $\displaystyle
 \bar \partial _{v}^{\beta }:=\frac{\partial ^{\left\vert{\beta
 }\right\vert }}{\partial \bar v_{1}^{\beta _{1}}\cdot \cdot
 \cdot \partial \bar v_{n}^{\beta _{n}}}$  are the derivatives
 in the directions  $v_{j}.$  We can suppose that  $p=0.$ \ \par 
\quad \quad 	To define the polydiscs we need to have an orthonormal basis
 at  $p$  and we shall built it with the vectors  $v_{j},\ j=2,...,v_{n}.$
 \ \par 
\quad \quad 	We have already that  $v_{1}$  is the complex normal direction,
 so choose  $e_{1}$  in the direction  $v_{1}$  and with norm
  $1.$  Now we use the Gram-Schmidt orthogonalisation procedure
 in the complex tangent plane  $\displaystyle \mathrm{S}\mathrm{p}\mathrm{a}\mathrm{n}(v_{2},v_{3},...,\
 v_{n})\ :$ \ \par 
\quad take  $e_{n}$  parallel to  $\displaystyle v_{n}$  and of norm  $1$  ;\ \par 
\quad \quad 	in  $\mathrm{S}\mathrm{p}\mathrm{a}\mathrm{n}(e_{n},\ v_{n-1})$
  take  $e_{n-1}$  of norm  $1$  and orthogonal to  $e_{n}\ ;$ \ \par 
and proceed this way to get an orthonormal basis  $\displaystyle
 (e_{2},...,e_{n})$  of  $\displaystyle T_{0}^{{\mathbb{C}}}(\partial
 \Omega )$  and complete it with  $\displaystyle e_{1}$  to get
 an orthonormal basis  $b(p)=(e_{1},...,e_{n})$  at  $p\ (=0).$ \ \par 
By this construction we have, with  $\displaystyle \zeta _{j}$
  the coordinates associated to the basis  $\displaystyle b(p),$ \ \par 
\quad \quad \quad \quad \quad  $\zeta _{1}=z_{1},\ \zeta _{2}=b_{2}^{2}z_{2},\ ...,\ \zeta
 _{n}=b_{n}^{2}z_{2}+\cdot \cdot \cdot +b_{n}^{n}z_{n},$ \ \par 
i.e. the matrix of change of coordinates is triangular.\ \par 
\quad So the lemma gives, still with  $\displaystyle m_{n}':={\left\lceil{m_{n}}\right\rceil},$
 \ \par 
\quad \quad \quad \quad \quad  $\displaystyle \tilde \rho (\zeta )=\Re \zeta _{1}+\sum_{2\leq
 \left\vert{\alpha }\right\vert +\left\vert{\beta }\right\vert
 \leq m_{n}'}{A_{\alpha \beta }\zeta ^{\alpha }\bar \zeta ^{\beta
 }}+o(\left\vert{\zeta }\right\vert ^{m_{n}'}),$ \ \par 
with  $\displaystyle \ A_{\alpha \beta }\neq 0\Rightarrow \sum_{i=1}^{n}{\frac{\alpha
 _{i}+\beta _{i}}{m_{i}}}\geq 1.$ \ \par 
where now the  $\zeta _{j}=b_{j}\cdot z$  are seen as linear
 forms on  $z.$  Fix  $t>0$  small enough so that  $a:=(-t,0,...,0)\in
 {\mathcal{U}},$  hence  $\displaystyle \pi (a)=0=p.$  Suppose
 that  $z\in Q_{a}(\delta )$  the polydisc based on  $b(p)$ 
 with  $\delta $  to be fixed later ; this means  $t=d(a)$  and\ \par 
\quad \quad \quad \quad  $\forall j=1,...,n,\ \left\vert{z_{j}}\right\vert <\delta d(a)^{1/m_{j}}.$
 \ \par 
This implies, because  $m_{1}=1\leq m_{2}\leq \cdot \cdot \cdot
 \leq m_{n},$  that\ \par 
\quad $\ \left\vert{\zeta _{j}}\right\vert \leq \sum_{k=1}^{j}{\left\vert{b_{j}^{k}}\right\vert
 \left\vert{z_{k}}\right\vert }$  $\leq \sum_{k=1}^{j}{\left\vert{b_{j}^{k}}\right\vert
 \delta d(a)^{1/m_{k}}}\leq \delta d(a)^{1/m_{j}}{\left({\sum_{k=1}^{j}{\left\vert{b_{j}^{k}}\right\vert
 }}\right)}=\delta B_{j}d(a)^{1/m_{j}},$ \ \par 
with  $\displaystyle B_{j}:=\sum_{k=1}^{j}{\left\vert{b_{j}^{k}}\right\vert
 }.$  		So we get\ \par 
\quad \quad \quad \quad \quad  $\displaystyle \ \left\vert{\zeta ^{\alpha }}\right\vert \leq
 \delta ^{\left\vert{\alpha }\right\vert }B^{\left\vert{\alpha
 }\right\vert }\prod_{j=1}^{n}{d(a)^{\frac{\alpha _{j}}{m_{j}}}},$ \ \par 
with  $B:=\max _{j=1,...,n}\left\vert{B_{j}}\right\vert .$  Replacing
  $\alpha _{j}$  by  $\beta _{j}$  in the previous proof, we get\ \par 
\quad \quad  $\displaystyle \ \left\vert{\bar \zeta ^{\beta }}\right\vert
 \leq \delta ^{\left\vert{\beta }\right\vert }B^{\left\vert{\beta
 }\right\vert }\prod_{j=1}^{n}{d(a)^{\frac{\beta _{j}}{m_{j}}}},$ \ \par 
so\ \par 
\quad \quad \quad \quad \quad  $\displaystyle \ \left\vert{\zeta ^{\alpha }\bar \zeta ^{\beta
 }}\right\vert \leq \delta ^{\left\vert{\alpha }\right\vert +\left\vert{\beta
 }\right\vert }B^{\left\vert{\alpha }\right\vert +\left\vert{\beta
 }\right\vert }\prod_{j=1}^{n}{d(a)^{\frac{\alpha _{j}+\beta
 _{j}}{m_{j}}}}=\delta ^{\left\vert{\alpha }\right\vert +\left\vert{\beta
 }\right\vert }B^{\left\vert{\alpha }\right\vert +\left\vert{\beta
 }\right\vert }d(a)^{\sum_{j=1}^{n}{\frac{\alpha _{j}+\beta _{j}}{m_{j}}}}.$
 \ \par 
In the sum, in order to have  $A_{\alpha \beta }\neq 0,$  we
 have  $\displaystyle \ \sum_{j=1}^{n}{\frac{\alpha _{j}+\beta
 _{j}}{m_{j}}}\geq 1,$ \ \par 
hence\ \par 
\quad \quad \quad \quad \quad  $\rho (z)\leq \Re z_{1}+\sum_{2\leq \left\vert{\alpha }\right\vert
 +\left\vert{\beta }\right\vert \leq m_{n}'}{A_{\alpha \beta
 }\delta ^{\left\vert{\alpha }\right\vert +\left\vert{\beta }\right\vert
 }B^{\left\vert{\alpha }\right\vert +\left\vert{\beta }\right\vert
 }d(a)}+o(\left\vert{z}\right\vert ^{m_{n}'})\leq \Re z_{1}+\delta
 d(a)C+o(\left\vert{z}\right\vert ^{m_{n}'})$ \ \par 
with\ \par 
\quad \quad \quad \quad \quad  $C:=\ \sum_{2\leq \left\vert{\alpha }\right\vert +\left\vert{\beta
 }\right\vert \leq m_{n}'}{A_{\alpha \beta }\delta ^{(\left\vert{\alpha
 }\right\vert +\left\vert{\beta }\right\vert -1)}B^{\left\vert{\alpha
 }\right\vert +\left\vert{\beta }\right\vert }}.$ \ \par 
Moreover we have  $\ \left\vert{z_{1}-d(a)}\right\vert <\delta d(a)$  so\ \par 
\quad \quad \quad \quad \quad  $\rho (z)\leq -d(a)+\delta d(a)+\delta d(a)C+o(\left\vert{z}\right\vert
 ^{m_{n}'})=d(a)(-1+\delta (1+C))+o(\left\vert{z}\right\vert
 ^{m_{n}'}).$ \ \par 
The constant  $C$  depends on a finite number of derivatives
 of  $\rho .$  Because the domain is of finite linear type  
 $\exists M({\mathcal{Q}})::\forall p\in \partial \Omega ,\ m_{n}(p)\leq
 M({\mathcal{Q}}),$  by the compactness of  $\partial \Omega
 $  we have  $\exists D>0,\ C=C(p)\leq D$  for any  $p\in \partial
 \Omega .$ \ \par 
Hence if  $\delta _{0}(1+D)\leq 1/2$  we have  $\rho (z)<0$ 
 if  $\ \left\vert{z}\right\vert $  is small enough to absorb
 the  $o(\left\vert{z}\right\vert ^{m_{n}'}).$  This means that
  $\displaystyle Q_{a}(\delta _{0})\subset \Omega .$ \ \par 
So we find a  $\delta _{0}>0$  such that, shrinking  ${\mathcal{U}}$
  if necessary to absorb the  $\displaystyle o(\left\vert{z}\right\vert
 ^{m_{n}'}),$  we get   $\forall a\in {\mathcal{U}},\ Q_{a}(\delta
 _{0})\subset \Omega .$   $\blacksquare $ \ \par 

\begin{Prps}
~\label{strongPC724}If  $\Omega $  is equipped with a good family
  ${\mathcal{Q}}$  with multi-type  $\lbrace m(a)\rbrace _{a\in
 \Omega },$  then it is of linear multi type smaller than  $\displaystyle
 \lbrace (1,{\left\lceil{m_{n}(\alpha )}\right\rceil},...,{\left\lceil{m_{n}(\alpha
 )}\right\rceil})\rbrace _{\alpha \in \partial \Omega }.$ 
\end{Prps}
\quad Proof.\ \par 
Let  $\alpha \in \partial \Omega $  and suppose, by rotation
 and translation that  $\alpha =0,\ \rho (z)=\Re z_{1}+\Gamma
 (\Im z_{1};z')$  with  $z'=(z_{2},...,z_{n}).$ \ \par 
We have that for any points  $a\in {\mathcal{U}}$  such that
  $\pi (a)=\alpha ,$  the polydisc  $Q_{a}(\delta _{0})$  is
 contained in  $\Omega .$ \ \par 
This means that, for  $\displaystyle \delta <\delta _{0},$  the
 point  $A_{\delta }:=(-d(a),\ \delta z_{2},...,\delta z_{j},...,\delta
 z_{n})$  with  $\ \left\vert{z_{j}}\right\vert =d(a)^{1/m_{j}}$
  is in  $\Omega $  hence the real segment  $I:=\lbrace A_{\delta
 }\rbrace _{\delta \in \lbrack 0,\delta _{0}\rbrack }$  centered
 at  $a$  is contained in  $\Omega .$ \ \par 
\quad \quad 	Fix  $\displaystyle (0,z')\in T_{\alpha }^{{\mathbb{C}}}(\partial
 \Omega )$  and set  $\Re z_{1}=\mu (\delta z'),\ \delta \in
 \lbrack 0,\delta _{0}\lbrack $  the graph of  $\displaystyle
 \partial \Omega $  over the segment  $I,$  i.e.  $\mu (\delta
 z')$  is such that  $\displaystyle (\mu (\delta z'),\delta z')\in
 \partial \Omega \ ;$   then we have that 	 $\forall \delta \in
 \lbrack 0,\delta _{0}\lbrack ,\mu (\delta z')<d(a).$  (See the
 following picture, with  $N$  the {\sl inward} normal)\ \par 
\ \par 
\ \par 
\ \par 
\ \par 
\ \par 
\ \par 
\begin{figure}[h]
\begin{center}
\vspace{-4cm}
\resizebox{9cm}{!}{\includegraphics{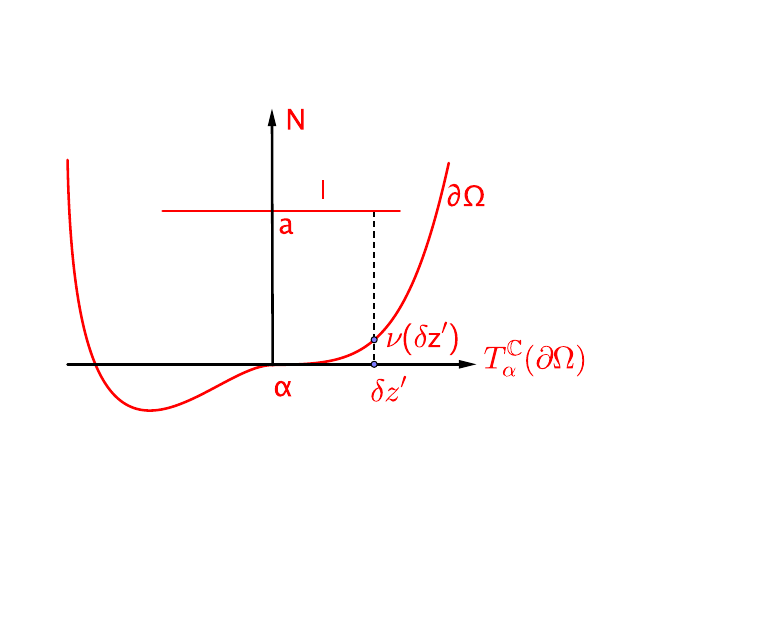}}
\end{center}
\vspace{-4cm}
\end{figure}\ \par 
\ \par 
\ \par 
\ \par 
\ \par 
So the distance from the point  $\displaystyle \delta z'\in T_{\alpha
 }^{{\mathbb{C}}}(\partial \Omega )$  to  $\displaystyle \alpha
 (=0)$  is  $\displaystyle \delta {\sqrt{\sum_{j=2}^{n}{d(a)^{2/m_{j}}}}}.$
 \ \par 
\ \par 
Recall that the order of contact of  $\displaystyle \partial
 \Omega $  with the real direction  $\displaystyle v$  at  $\displaystyle
 \alpha \in \partial \Omega $  is the order of vanishing of 
 $\displaystyle \rho (\alpha +tv)-\rho (\alpha )$  when  $\displaystyle
 t\rightarrow 0.$ \ \par 
\quad \quad 	Setting  $\displaystyle t:=\delta {\sqrt{\sum_{j=2}^{n}{d(a)^{2/m_{j}}}}},\
 x=d(a)$  and  $\displaystyle u=\sum_{j=2}^{n}{x^{2/m_{j}}}$  we have\ \par 
\quad \quad \quad \quad \quad  $\displaystyle \ \frac{dt}{dx}=\delta \frac{\sum_{j=2}^{n}{\frac{2}{m_{j}}x^{-1+2/m_{j}}}}{2{\sqrt{u}}}\neq
 0$  and finite for  $\displaystyle x\neq 0,$ \ \par 
so by the implicit function theorem we have a smooth function
  $f(t)$  such that  $\displaystyle d(a)=f(t)$  for  $\displaystyle
 d(a)\neq 0.$ \ \par 
Now we make the change of variables  $\displaystyle \delta z'=t\zeta
 ',$  still with  $\displaystyle \zeta '\in T_{\alpha }^{{\mathbb{C}}}(\partial
 \Omega ),$  we have  $\displaystyle \mu (\delta z')=\mu (t\zeta
 ')\leq f(t)=d(a)\ ;$  hence the order of contact of  $\displaystyle
 \partial \Omega $  with the direction  $z'$  is bigger than
 the order of contact of  $\displaystyle f(t)$  at  $\displaystyle
 t=0.$  So fix  $\displaystyle \delta <\delta _{0}$  and let
  $\displaystyle d(a)\rightarrow 0\ ;$  because for any  $\displaystyle
 a::\pi (a)=\alpha $  we have  $\displaystyle Q_{a}(\delta _{0})\subset
 \Omega ,$  we get\ \par 
\quad \quad \quad \quad  $\displaystyle t/\delta d(a)^{1/m_{n}(\alpha )}={\sqrt{\sum_{j=2}^{n}{d(a)^{2/m_{j}-2/m_{n}}}}}$
 \ \par 
and because  $\displaystyle d(a)^{2/m_{j}-2/m_{n}}\rightarrow
 0$  if  $\displaystyle m_{j}<m_{n},$  we get\ \par 
\quad \quad \quad \quad \quad  $\displaystyle t/\delta d(a)^{1/m_{n}(\alpha )}\rightarrow {\sqrt{l(\alpha
 )}}$ \ \par 
where  $l(\alpha )$  is the number of  $\displaystyle j::m_{j}=m_{n},$
  hence  $\displaystyle 1<l(\alpha )\leq n-1.$ \ \par 
\quad 	So  $\displaystyle t\simeq {\sqrt{l(\alpha )}}\delta d(a)^{1/m_{n}(\alpha
 )}\Rightarrow f(t)=d(a)\simeq l(\alpha )^{-m_{n}(\alpha )/2}\delta
 ^{-m_{n}}t^{m_{n}},$  hence the order of contact of  $\displaystyle
 f(t)$  at  $0$  is  $\displaystyle m_{n}(\alpha )$  hence the
 order of contact of  $\displaystyle \partial \Omega $  in the
 real direction  $\displaystyle z'$  is at least  $\displaystyle
 m_{n}(\alpha ).$ \ \par 
\quad \quad 	We have proved that for any real direction in  $\displaystyle
 T_{\alpha }^{{\mathbb{C}}}(\partial \Omega )$  the order of
 contact of  $\displaystyle \partial \Omega $  is at least  $\displaystyle
 m_{n}(\alpha ).$ \ \par 
\quad \quad 	Let us make {\sl any} linear change of variables keeping  $\displaystyle
 T_{\alpha }^{{\mathbb{C}}}(\partial \Omega )$  and sending 
 $\alpha $  to  $\displaystyle 0.$  Let us expand  $\rho $  in
 these new coordinates\ \par 
\quad \quad \quad \quad \quad  $\displaystyle \rho (w,\bar w)=\sum_{\left\vert{\alpha }\right\vert
 +\left\vert{\beta }\right\vert =k}{A_{\alpha ,\beta }w^{\alpha
 }\bar w^{\beta }}+{\mathcal{O}}(\left\vert{w}\right\vert ^{k+1}).$ \ \par 
Set  $\displaystyle w_{1}=0,$  fix  $\displaystyle \zeta =(0,\zeta
 _{2},...,\zeta _{n})$  and set  $\displaystyle w=t\zeta \in
 T_{\alpha }^{{\mathbb{C}}}(\partial \Omega ).$  Then we have\ \par 
\quad \quad \quad \quad \quad  $\displaystyle \rho (w,\bar w)=t^{k}\sum_{\left\vert{\alpha
 }\right\vert +\left\vert{\beta }\right\vert =k}{A_{\alpha ,\beta
 }\zeta ^{\alpha }\bar \zeta ^{\beta }}+{\mathcal{O}}(t^{k+1}).$ \ \par 
We already know that the order of contact of  $\displaystyle
 \partial \Omega $  with any real direction of  $\displaystyle
 T_{\alpha }^{{\mathbb{C}}}(\partial \Omega )$  is bigger than
  $\displaystyle m_{n},$  and this is still true if we change
 coordinates {\sl linearly} provided that we keep  $\displaystyle
 T_{\alpha }^{{\mathbb{C}}}(\partial \Omega ).$  So the order
 of vanishing of  $\rho $  along the real line  $\displaystyle
 t\zeta $  is bigger than  $\displaystyle m_{n}(\alpha )$  hence
 in order to have  $\displaystyle A_{\alpha ,\beta }=\frac{\partial
 ^{\alpha +\beta }\rho }{\partial ^{\alpha }w\partial ^{\beta
 }\bar w}(0)$  not all zeros for  $\displaystyle \ \left\vert{\alpha
 }\right\vert +\left\vert{\beta }\right\vert =k,$  we need to
 have  $\displaystyle k\geq m_{n}$  hence  $\displaystyle k\geq
 {\left\lceil{m_{n}}\right\rceil}.$  This implies	 	 $\displaystyle
 \ \sum_{j=2}^{n}{\frac{\alpha _{j}+\beta _{j}}{{\left\lceil{m_{n}}\right\rceil}}}\geq
 \sum_{j=2}^{n}{\frac{\alpha _{j}+\beta _{j}}{k}}=1.$ \ \par 
\quad \quad 	This means that for this change of variables the weight  $\displaystyle
 (1,{\left\lceil{m_{n}}\right\rceil},...,{\left\lceil{m_{n}}\right\rceil})$
  is linearly distinguished and hence  $\displaystyle \partial
 \Omega $  at  $\alpha $  is of finite linear type. Moreover
 the linear multi type  $\displaystyle (1,m'_{2},...,m'_{n})$
  of  $\displaystyle \partial \Omega $  at  $\alpha $  being
 smaller than the weight  $\displaystyle (1,{\left\lceil{m_{n}}\right\rceil},...,{\left\lceil{m_{n}}\right\rceil})$
  by definition, we have\ \par 
\quad \quad \quad \quad \quad  $\displaystyle \forall j=2,...,n,\ m'_{j}(\alpha )\leq {\left\lceil{m_{n}(\alpha
 )}\right\rceil}.$   $\blacksquare $ \ \par 
\ \par 
\quad \quad 	Theorem~\ref{linkTfFpc20} and Proposition~\ref{strongPC724}
 give the characterization :\ \par 

\begin{Crll}
The domain  $\Omega $  has a good family  ${\mathcal{Q}}$  of
 polydiscs associated to  $\lbrace m(\alpha )\rbrace _{\alpha
 \in \partial \Omega }$  iff the linear multi-type of  $\Omega
 $  is  smaller than  $\displaystyle \lbrace (1,{\left\lceil{m_{n}(\alpha
 )}\right\rceil},...,{\left\lceil{m_{n}(\alpha )}\right\rceil})\rbrace
 _{\alpha \in \partial \Omega }.$ 
\end{Crll}
\ \par 

\subsection{A geometrical characterisation by existence of inner
 complex tangential ellipsoids.}
\quad \quad 	First we set tools we shall need.  Recall the standard notations\ \par 
\quad \quad \quad \quad \quad  $\displaystyle \forall \alpha =(\alpha _{1},...,\alpha _{n})\in
 {\mathbb{N}}^{n},\ \partial ^{\alpha }f(x,p):=\frac{\partial
 ^{\left\vert{\alpha }\right\vert }f}{\partial ^{\alpha _{1}}x_{1}\cdot
 \cdot \cdot \partial ^{\alpha _{n}}x_{n}}(x,p),\ x^{\alpha }:=x_{1}^{\alpha
 _{1}}\cdot \cdot \cdot x_{n}^{\alpha _{n}}.$ \ \par 

\begin{Lmm}
~\label{1_goodFam40}Let  $f(x,p)$  be a  ${\mathcal{C}}^{\infty
 }({\mathbb{R}}^{n}{\times}{\mathbb{R}}^{m})$  function ; then
 there exist  ${\mathcal{C}}^{\infty }({\mathbb{R}}^{n}{\times}{\mathbb{R}}^{m})$
  functions  $f_{\alpha },$  for  $\alpha \in {\mathbb{N}}^{n},$
  such that :\par 
\quad \quad \quad \quad \quad  $\displaystyle f(x,p)=f(0,p)+\cdot \cdot \cdot +\frac{1}{k!}\sum_{\left\vert{\beta
 }\right\vert =k}{x^{\beta }\partial ^{\beta }f(0,p)}+\frac{1}{k!}\sum_{\left\vert{\alpha
 }\right\vert =k+1}{x^{\alpha }f_{\alpha }(x,p)}.$ \par 
The  $f_{\alpha }$  are given explicitly by the formulas :\par 
\quad \quad \quad \quad \quad  $\displaystyle f_{\alpha }(x,p):=\int_{0}^{1}{\partial ^{\alpha
 }f(tx,p)(1-t)^{k}dt}.$ 
\end{Lmm}
\quad \quad 	Proof.\ \par 
Set, for  $t\in {\mathbb{R}},\ g(t,p):=f(tx,p)$  then,we have\ \par 
\quad \quad \quad \quad \quad  $\displaystyle g^{(k)}(t,p):=\frac{\partial ^{k}g}{\partial
 t^{k}}(t,p)=\sum_{\left\vert{\alpha }\right\vert =k}{x^{\alpha
 }\partial ^{\alpha }f(tx,p)}.$ \ \par 
Apply to  $g$  the Taylor formula with integral remainder\ \par 
\quad \quad \quad \quad \quad  $\displaystyle g(1,p)=g(0,p)+\cdot \cdot \cdot +\frac{g^{(k)}(0,p)}{k!}+\frac{1}{k!}\int_{0}^{1}{g^{(k+1)}(t,p)(1-t)^{k}dt}.$
 \ \par 
We get\ \par 
\quad \quad \quad \quad \quad  $\displaystyle f(x,p)=f(0,p)+\cdot \cdot \cdot +\sum_{\left\vert{\alpha
 }\right\vert =k}{x^{\alpha }\partial ^{\alpha }f(0,p)}+\frac{1}{k!}\sum_{\left\vert{\beta
 }\right\vert =k+1}{x^{\beta }\int_{0}^{1}{\partial ^{\beta }f(tx,p)(1-t)^{k}dt}}.$
 \ \par 
Now set  $\displaystyle f_{\beta }(x,p):=\int_{0}^{1}{\partial
 ^{\beta }f(tx,p)(1-t)^{k}dt},$  then, deriving under the integral
 sign, we have that  $\displaystyle f_{\beta }$  is  ${\mathcal{C}}^{\infty
 }$  in the two variables  $\displaystyle x,p.$   $\blacksquare $ \ \par 
\ \par 
\quad We suppose we are given a family of orthonormal basis and multitypes
  $\lbrace b_{\alpha },\ m(\alpha )\rbrace _{\alpha \in \partial
 \Omega }.$ \ \par 
First, without loss of generality, we make the assumption that
 the normal derivative of  $\rho $  is  $1$  at any point  $\alpha
 \in \partial \Omega ,$  and  $\rho \in {\mathcal{C}}^{\infty
 }({\mathbb{C}}^{n}).$ \ \par 
\quad \quad 	Fix  $\alpha \in \partial \Omega \ ;$  by translation we can
 suppose  $\alpha =0,$  i.e. with  $\rho _{\alpha }(z):=\rho
 (z+\alpha )$  we have  $\rho _{\alpha }(0)=0.$ \ \par 
Now we make the rotation  $U_{\alpha }$  sending the standard
 basis of  ${\mathbb{C}}^{n}$  to  $b_{\alpha },$  i.e.  $\rho
 _{\alpha }(z):=\rho (U_{\alpha }z+\alpha ).$  In these new coordinates,
 we have that  $z'=(z_{2},...,z_{n})$  are the coordinates in
 the complex tangent space and  $z_{1}=x_{1}+iy_{1}$  is the
 coordinate in the normal complex plane and  $x_{1}$  is the
 coordinate in the real normal at  $\alpha \ (=0).$ \ \par 
\quad \quad 	Set  $h_{\alpha }(z'):=\rho _{\alpha }(0,z')\in {\mathcal{C}}^{\infty
 }({\mathbb{C}}^{n-1})$  then  $\rho _{\alpha }(z)-h_{\alpha
 }(z')=0$  if  $z_{1}=0,\ \forall z'\in {\mathbb{C}}^{n-1}.$ \ \par 
Set\ \par 
\quad \quad \quad \quad \quad  $g_{\alpha }(x_{1},y_{1},z'):=-x_{1}+\rho _{\alpha }(z)-h_{\alpha
 }(z')\in {\mathcal{C}}^{\infty }({\mathbb{C}}^{n})\ ;$ \ \par 
we have  $\rho _{\alpha }(z)=x_{1}+g_{\alpha }(z_{1},z')+h_{\alpha
 }(z').$ \ \par 
Recall that  $x_{1}$  is the coordinate on the real normal, so
  $\displaystyle \ \frac{\partial \rho _{\alpha }}{\partial x_{1}}(0)=1$
  by assumption and, because  $y_{1}$  is a tangent coordinate,
 we have  $\displaystyle \ \frac{\partial \rho _{\alpha }}{\partial
 y_{1}}(0)=0,$  so\ \par 
\quad \quad \quad \quad \quad  $\displaystyle \ \frac{\partial g_{\alpha }}{\partial x_{1}}(0,0)=-1+\frac{\partial
 \rho _{\alpha }}{\partial x_{1}}(0,0)=0,\ \ \frac{\partial g_{\alpha
 }}{\partial y_{1}}(0,0)=\frac{\partial \rho _{\alpha }}{\partial
 y_{1}}(0,0)=0.$ \ \par 

\begin{Lmm}
~\label{1_goodFam71}There is a number  $R>0,$  independent of
  $\displaystyle \alpha \in \partial \Omega ,$  such that, after
 the change of coordinates above, we have the estimate\par 
\quad \quad \quad \quad \quad  $\ \forall \alpha \in \partial \Omega ,\ \forall z\in B(0,R),\
 \left\vert{g_{\alpha }(z)}\right\vert \leq \frac{1}{4}\left\vert{z_{1}}\right\vert
 $ \par 
and the factorization\par 
\quad \quad \quad \quad \quad  $\displaystyle g_{\alpha }(x_{1},y_{1},z')=x_{1}g_{1}(x_{1},y_{1},z')+y_{1}g_{2}(x_{1},y_{1},z'),$
 \par 
with  $\displaystyle \forall z\in B(0,R),\ j=1,2,\ \ \left\vert{g_{j}(z)}\right\vert
 <3/10.$ 
\end{Lmm}
\quad \quad 	Proof.\ \par 
We apply lemma~\ref{1_goodFam40} to  $g:=g_{\alpha }(x_{1},y_{1},z')$
  to order  $\displaystyle 1$  with  $\displaystyle z'$  as the
 parameter  $\displaystyle p.$  We get\ \par 
\quad \quad \quad \quad \quad  $\displaystyle g(x_{1},y_{1},z')=g(0,0,z')+x_{1}\frac{\partial
 g}{\partial x_{1}}(0,0,z')+y_{1}\frac{\partial g}{\partial y_{1}}(0,0,z')+$
 \ \par 
\quad \quad \quad \quad \quad \quad \quad \quad  $\displaystyle +x_{1}^{2}g_{(2,0)}(x_{1},y_{1},z')+x_{1}y_{1}g_{(1,1)}(x_{1},y_{1},z')+y_{1}^{2}g_{(0,2)}(x_{1},y_{1},z').$
 \ \par 
The term  $\displaystyle g(0,0,z')=\rho _{\alpha }(0,z')-h_{\alpha
 }(z')=0,$  it remains the others.\ \par 
Because\ \par 
\quad \quad \quad \quad \quad  $\displaystyle g=g_{\alpha }=-x_{1}+\rho _{\alpha }(z)-h_{\alpha
 }(z')=-x_{1}+\rho (U_{\alpha }z+\alpha )-\rho (U_{\alpha }(0,z')+\alpha
 )$ \ \par 
all its derivatives are controlled by the derivatives of  $\rho
 $  in a neighborhood of  $\bar \Omega ,$  because  $U_{\alpha
 }$  is a rotation independent of  $z,$  so they are controlled
 uniformly  in  $\alpha .$  So are the integrals of them, hence
 the functions  $\displaystyle g_{(j,k)}.$ \ \par 
Because  $\displaystyle \ \frac{\partial g}{\partial x_{1}}(0,0,0)=0,$
  we have that  $\displaystyle \ \frac{\partial g}{\partial x_{1}}(0,0,z')$
  is small when  $\displaystyle \ \left\vert{z'}\right\vert $
  is small and this is uniform with respect to the point  $\displaystyle
 \alpha \in \partial \Omega .$ \ \par 
\quad \quad 	The same for  $\displaystyle \ \frac{\partial g}{\partial y_{1}}(0,0,z').$
  Moreover the functions  $\displaystyle g_{(j,k)}$  are bounded
 again uniformly with respect to the point  $\displaystyle \alpha
 \in \partial \Omega .$  So finally we can choose  $R'>0$   small
 enough and independent of the point  $\displaystyle \alpha \in
 \partial \Omega $   to get\ \par 
\quad \quad \quad \quad \quad  $\displaystyle \ \left\vert{z'}\right\vert <R'\Rightarrow \left\vert{\frac{\partial
 g}{\partial x_{1}}(0,0,z')}\right\vert <1/10,\ \left\vert{\frac{\partial
 g}{\partial y_{1}}(0,0,z')}\right\vert <1/10.$ \ \par 
Take  $\displaystyle R''$  small enough to have\ \par 
\quad \quad \quad \quad \quad  $\displaystyle \forall z\in B(0,R''),\ i,j=0,1,2,\ \ \left\vert{z_{1}}\right\vert
 \left\vert{g_{(i,j)}(z)}\right\vert <1/10\ ;$ \ \par 
then, with  $\displaystyle R:=\min (R',R'',1/6)$  we get\ \par 
\quad \quad \quad \quad \quad  $\displaystyle \forall z\in B(0,R),\ \left\vert{g(x_{1},y_{1},z')}\right\vert
 \leq \frac{1}{10}(\left\vert{x_{1}}\right\vert +\left\vert{y_{1}}\right\vert
 +3R\left\vert{z_{1}}\right\vert )\leq \frac{\left\vert{z_{1}}\right\vert
 }{10}(2+3R)\leq \frac{\left\vert{z_{1}}\right\vert }{4}.$ \ \par 
\quad \quad 	Now setting\ \par 
\quad \quad \quad \quad \quad  $\displaystyle g_{1}(x_{1},y_{1},z'):=\frac{\partial g}{\partial
 x_{1}}(0,0,z')+x_{1}g_{(2,0)}(x_{1},y_{1},z')+y_{1}g_{(1,1)}(x_{1},y_{1},z')$
 \ \par 
and\ \par 
\quad \quad \quad \quad  $\displaystyle g_{2}(x_{1},y_{1},z'):=y_{1}g_{(0,2)}(x_{1},y_{1},z'),$ \ \par 
we have\ \par 
\quad \quad \quad \quad \quad  $\displaystyle \forall z\in B(0,R),\ \left\vert{g_{1}(z)}\right\vert
 \leq \frac{3}{10},\ \left\vert{g_{2}(z)}\right\vert \leq 1/10,$ \ \par 
and the factorization\ \par 
\quad \quad \quad \quad  $\displaystyle g(x_{1},y_{1},z')=g_{\alpha }(x_{1},y_{1},z')=x_{1}g_{1}(x_{1},y_{1},z')+y_{1}g_{2}(x_{1},y_{1},z').$
   $\blacksquare $ \ \par 
\ \par 

\begin{Lmm}
~\label{1_BonFamille20}If  $\Omega $  has a good family of polydiscs
 for  $\displaystyle \alpha \in \partial \Omega ,$  there is
 a complex tangentially elliptic domain  $C=C_{\alpha },$  with
 aperture  $\displaystyle \Gamma >0$  near the point  $\alpha
 \in \partial \Omega ,$  of class  ${\mathcal{C}}^{2}$  and such that\par 
\quad \quad 	-  $C_{\alpha }\subset \Omega ,$  near  $\alpha ,$ \par 
\quad \quad 	-  $\bar C_{\alpha }$  and  $\displaystyle \bar \Omega $  meet
 at  $\alpha ,$ \par 
\quad \quad 	-  $\forall a\in {\mathcal{U}}::\pi (a)=\alpha ,\ Q_{a}(\delta
 )\subset C_{\alpha }$  provided that  $\displaystyle \delta
 ^{2}<\min (\Gamma ,\frac{1}{4(n-1)A},1/4).$ 
\end{Lmm}
\quad \quad 	Proof.\ \par 
We shall build  $C.$  We make the change of variables above then
 we can write\ \par 
\quad \quad \quad \quad \quad  $\rho _{\alpha }(z)=x_{1}+g_{\alpha }(z)+h_{\alpha }(z'),$ \ \par 
where  $z_{1}=x_{1}+iy_{1},\ z':=(z_{2},...,z_{n})$  and  $g_{\alpha
 }\in {\mathcal{C}}^{\infty }({\mathbb{C}}^{n}),\ h_{\alpha }\in
 {\mathcal{C}}^{\infty }({\mathbb{C}}^{n-1}),\ g_{\alpha }(\alpha
 )=h_{\alpha }(\alpha ')=0.$ \ \par 
Let\ \par 
\quad \quad \quad \quad \quad  $\lambda (z):=x_{1}+g_{\alpha }(z)+A{\left({\left\vert{z_{2}}\right\vert
 ^{m_{2}}+\cdot \cdot \cdot +\left\vert{z_{n}}\right\vert ^{m_{n}}}\right)},$
 \ \par 
with  $m=m(\alpha ).$  We fix an aperture  $\Gamma >0$  and we
 shall choose  $A$  in order to have that\ \par 
\quad \quad \quad \quad \quad  $C:=\lbrace \lambda <0\rbrace \cap \lbrace \left\vert{y_{1}}\right\vert
 <-\Gamma x_{1}\rbrace $  \ \par 
fills the requirements of the lemma.\ \par 
This domain  $C=C_{\alpha }$  is what we shall call a "complex
 tangentially elliptic domain with aperture  $\displaystyle \Gamma
 >0$ ". As the referee remarks this can also be seen as the classical
 "approach regions" to the boundary in the strictly pseudo convex case.\ \par 
Fix  $a_{}=(-t,\ 0,...,0),\ t\in {\mathbb{R}},\ (-t,z')\in B(0,R)$
  in order to have  $\ \left\vert{g_{\alpha }(-t,z')}\right\vert
 <t/4,$  by lemma~\ref{1_goodFam71} ; consider the slice  $\displaystyle
 S_{z_{1}}$  of  $\displaystyle C_{\alpha }$ \ \par 
\quad \quad  $\displaystyle \forall z_{1}::(z_{1},z')\in B(0,R),\ S_{z_{1}}:=\lbrace
 z':=(z_{2},...,z_{n})\ ::A{\left({\left\vert{z_{2}}\right\vert
 ^{m_{2}}+\cdot \cdot \cdot +\left\vert{z_{n}}\right\vert ^{m_{n}}}\right)}<-(\Gamma
 +2)x_{1}\rbrace .$ \ \par 
Then\ \par 
\quad \quad \quad \quad \quad  $z_{1}=-t\Rightarrow y_{1}=0\Rightarrow S_{-t}:=\lbrace z':=(z_{2},...,z_{n})\
 ::A{\left({\left\vert{z_{2}}\right\vert ^{m_{2}}+\cdot \cdot
 \cdot +\left\vert{z_{n}}\right\vert ^{m_{n}}}\right)}<(\Gamma
 +2)t\rbrace .$ \ \par 
If  $z'\in S_{-t},$  then  $\displaystyle \forall j\geq 2,\ A\left\vert{z_{j}}\right\vert
 ^{m_{j}}<(\Gamma +2)t\Rightarrow \left\vert{z_{j}}\right\vert
 <\frac{1}{(A/(\Gamma +2))^{1/m_{j}}}t^{1/m_{j}}=\frac{d(a)^{1/m_{j}}}{(A/(\Gamma
 +2))^{1/m_{j}}}.$ \ \par 
\quad \quad 	Hence if  $\displaystyle A\geq \frac{\Gamma +2}{\delta _{0}^{m_{n}(\alpha
 )}}$   then  $\displaystyle \forall j\geq 2,\ (A/(\Gamma +2))^{1/m_{j}}\delta
 _{0}\geq (A/(\Gamma +2))^{1/m_{n}}\delta _{0}\geq 1$  hence\ \par 
\quad \quad \quad \quad \quad  $\displaystyle z'\in S_{-t}\Rightarrow \forall j=2,...,n,\ \
 \left\vert{z_{j}}\right\vert <\frac{d(a)^{1/m_{j}}}{(A/(\Gamma
 +2))^{1/m_{j}}}\leq \delta _{0}d(a)^{1/m_{j}}\Rightarrow (-t,z')\in
 Q_{a}(\delta _{0})\subset \Omega .$ \ \par 
So if  $(-t+iy_{1},z')\in C,$  we have  $\displaystyle \ \left\vert{y_{1}}\right\vert
 <\Gamma t$  and\ \par 
\quad \quad \quad \quad \quad \quad  $\displaystyle \lambda (-t+iy_{1},z')=-t+g_{\alpha }(z)+A{\left({\left\vert{z_{2}}\right\vert
 ^{m_{2}}+\cdot \cdot \cdot +\left\vert{z_{n}}\right\vert ^{m_{n}}}\right)}<0\Rightarrow
 $ \ \par 
\quad \quad \quad \quad \quad \quad \quad \quad \quad \quad  $\displaystyle A{\left({\left\vert{z_{2}}\right\vert ^{m_{2}}+\cdot
 \cdot \cdot +\left\vert{z_{n}}\right\vert ^{m_{n}}}\right)}<t-g_{\alpha
 }(-t+iy_{1},z').$ \ \par 
But, in the ball  $\displaystyle B(0,R),$  we have  $\displaystyle
 \ \left\vert{g_{\alpha }(-t+iy_{1},z')}\right\vert <\frac{\left\vert{y_{1}}\right\vert
 }{4}+\frac{t}{4},$  hence\ \par 
\quad \quad \quad \quad \quad  $\displaystyle A{\left({\left\vert{z_{2}}\right\vert ^{m_{2}}+\cdot
 \cdot \cdot +\left\vert{z_{n}}\right\vert ^{m_{n}}}\right)}<t-g_{\alpha
 }(-t+iy_{1},z')<t+t/4+\left\vert{y_{1}}\right\vert /4<(\Gamma +2)t,$ \ \par 
then  $z'\in S_{-t+iy_{1}}$  and  $(-t+iy_{1},z')\in Q_{a}(\delta
 _{0})\subset \Omega \Rightarrow \rho _{\alpha }(-t+iy_{1},z')<0$
  provided that  $\displaystyle \ \left\vert{y_{1}}\right\vert
 <\delta _{0}t,$  hence we need to take the aperture  $\displaystyle
 \Gamma \leq \delta _{0}.$  To have the same  $A$  for all the
 boundary points, we take\ \par 
\quad \quad \quad \quad \quad  $\displaystyle A=\frac{\Gamma +2}{\delta _{0}^{M({\mathcal{Q}})}}$
  with  $\displaystyle M({\mathcal{Q}})=\sup _{\alpha \in \partial
 \Omega }m_{n}(\alpha ),$ \ \par 
which is bounded because  ${\mathcal{Q}}$  is a good family.\ \par 
With this choice of  $A,$  we have that\ \par 
\quad \quad \quad \quad \quad  $(-t+iy_{1},z')\in C\cap B(0,R)\Rightarrow (-t+iy_{1},z')\in
 Q_{a}(\delta _{0})\subset \Omega $ \ \par 
so  $\lambda (-t+iy_{1},z')<0\Rightarrow \rho (-t+iy_{1},z')<0,$  i.e.\ \par 
\quad \quad \quad \quad \quad  $\displaystyle \lambda (z)=-t+g_{\alpha }(-t+iy_{1},z')+A{\left({\left\vert{z_{2}}\right\vert
 ^{m_{2}}+\cdot \cdot \cdot +\left\vert{z_{n}}\right\vert ^{m_{n}}}\right)}<0,$
 \ \par 
hence\ \par 
\quad \quad \quad \quad \quad  $\displaystyle A{\left({\left\vert{z_{2}}\right\vert ^{m_{2}}+\cdot
 \cdot \cdot +\left\vert{z_{n}}\right\vert ^{m_{n}}}\right)}<t-g_{\alpha
 }(-t+iy_{1},z'),$ \ \par 
and this implies  $\rho _{\alpha }(z)<0,$  i.e.\ \par 
\quad \quad \quad \quad \quad  $\displaystyle h_{\alpha }(z')<t-g_{\alpha }(-t+iy_{1},z'),$ \ \par 
so necessarily  $\displaystyle h_{\alpha }(z')\leq A{\left({\left\vert{z_{2}}\right\vert
 ^{m_{2}}+\cdot \cdot \cdot +\left\vert{z_{n}}\right\vert ^{m_{n}}}\right)},$
  because if not suppose there is a  $z'$  such that\ \par 
\quad \quad \quad \quad \quad  $\displaystyle A{\left({\left\vert{z_{2}}\right\vert ^{m_{2}}+\cdot
 \cdot \cdot +\left\vert{z_{n}}\right\vert ^{m_{n}}}\right)}<h(z'),$  \ \par 
take a  $t>0$  with\ \par 
\quad \quad \quad \quad \quad  $\displaystyle A{\left({\left\vert{z_{2}}\right\vert ^{m_{2}}+\cdot
 \cdot \cdot +\left\vert{z_{n}}\right\vert ^{m_{n}}}\right)}<t-g(-t,z')<h(z'),$
 \ \par 
then the point  $\displaystyle (-t,z')$  belongs to  $C,$  because
 we have  $\displaystyle y_{1}=0<\Gamma t,$  but not to  $\Omega
 ,$  which is a contradiction.\ \par 
Hence we proved\ \par 
\quad \quad \quad \quad 	\begin{equation}  z\in C_{\alpha }\cap B(0,R)\Rightarrow h_{\alpha
 }(z')\leq A{\left({\left\vert{z_{2}}\right\vert ^{m_{2}}+\cdot
 \cdot \cdot +\left\vert{z_{n}}\right\vert ^{m_{n}}}\right)}.\label{aspcGlo0}\end{equation}\
 \par 
\quad \quad 	Hence over any point of  $\displaystyle B(\alpha ,R)\cap \lbrace
 \left\vert{y_{1}}\right\vert <-\Gamma x_{1}\rbrace $  we have
 a a domain of class  ${\mathcal{C}}^{2},$  because		  $\displaystyle
 m_{1}=1\Rightarrow m_{2}\in {\mathbb{N}}$  by lemma~\ref{2_BonFamille20}
 hence\ \par 
\quad \quad \quad \quad \quad  $\displaystyle \ \frac{\partial ^{2}\left\vert{z_{2}}\right\vert
 ^{m_{2}}}{\partial z_{2}^{2}}=\frac{\partial ^{2}(\bar z_{2}z_{2})^{m_{2}/2}}{\partial
 z_{2}^{2}}=\frac{m_{2}(m_{2}-2)}{4}\left\vert{z_{2}}\right\vert
 ^{m_{2}-4}\bar z_{2}^{2}.$ \ \par 
If  $\displaystyle m_{2}=2$  then  $\displaystyle \ \frac{\partial
 ^{2}\left\vert{w}\right\vert ^{m}}{\partial w^{2}}=0$  and this
 term is  ${\mathcal{C}}^{2}.$  If  $\displaystyle m_{2}\geq 3$  then\ \par 
\quad \quad \quad \quad \quad  $\displaystyle \ \frac{\partial ^{2}\left\vert{z_{2}}\right\vert
 ^{m_{2}}}{\partial z_{2}^{2}}=\frac{m_{2}(m_{2}-2)}{4}\left\vert{z_{2}}\right\vert
 ^{m_{2}-4}\bar z_{2}^{2}$ \ \par 
and this is continuous hence again this term is  ${\mathcal{C}}^{2}.$ \ \par 
Now we have  $\displaystyle m_{j}\geq m_{2}$  for  $\displaystyle
 j\geq 3$  hence all the other terms are also  ${\mathcal{C}}^{2}.$ \ \par 
\ \par 
\quad \quad 	It remains to prove the last item of the lemma.\ \par 
\quad \quad 	Take a point  $\displaystyle a\in \Omega ,\ \pi (a)=\alpha $
  then	  $\displaystyle a=(-t,0,...,0)$  after the usual change
 of variables ; fix a  $\delta >0$  to be precised later ; if
  $(x_{1}+iy_{1},z')\in Q_{a}(\delta ),$  then\ \par 
\quad \quad \quad \quad \quad  $\displaystyle \ \forall j=2,...,n,\ \left\vert{z_{j}}\right\vert
 <\delta t^{1/m_{j}},\ \left\vert{x_{1}+t}\right\vert <\delta
 t\Rightarrow x_{1}<-t(1-\delta ),\ \left\vert{y_{1}}\right\vert
 <\delta t,$ \ \par 
so we already choose  $\displaystyle \delta \leq \Gamma $  to
 have  $\displaystyle \ \left\vert{y_{1}}\right\vert <-\Gamma
 x_{1},$  and\ \par 
\quad \quad \quad \quad \quad \quad \quad  $\displaystyle A{\left({\left\vert{z_{2}}\right\vert ^{m_{2}}+\cdot
 \cdot \cdot +\left\vert{z_{n}}\right\vert ^{m_{n}}}\right)}<tA\sum_{j=2}^{n}{\delta
 ^{m_{j}}},$ \ \par 
hence\ \par 
\quad \quad  $\displaystyle \lambda (x_{1},y_{1},z')=x_{1}+g(z)+A{\left({\left\vert{z_{2}}\right\vert
 ^{m_{2}}+\cdot \cdot \cdot +\left\vert{z_{n}}\right\vert ^{m_{n}}}\right)}<-t(1-\delta
 )+tA\sum_{j=2}^{n}{\delta ^{m_{j}}}+\left\vert{g(z)}\right\vert .$ \ \par 
Because  $m_{j}\geq 2,$ 		 $\ A\delta ^{m_{j}}\leq \delta ^{2}A,$  so\ \par 
\quad \quad \quad \quad \quad  $\displaystyle \lambda (x_{1},y_{1},z')<-t(1-\delta )+(n-1)\delta
 ^{2}tA+\left\vert{g(z)}\right\vert .$ \ \par 
But, by lemma~\ref{1_goodFam71}, the smooth function  $g(z)$
  is bounded on  $B(0,R)$  by  $\ \frac{1}{4}\ \left\vert{z_{1}}\right\vert
 ,$  and we have  $\displaystyle \ \left\vert{y_{1}}\right\vert
 <\Gamma \left\vert{x_{1}}\right\vert $  so\ \par 
\quad \quad \quad \quad \quad  $\displaystyle \lambda (x_{1},y_{1},z')<-t(1-\delta )+(n-1)\delta
 ^{2}tA+\frac{t}{4}(1+\Gamma )=t(\delta -\frac{3}{4}+(n-1)\delta
 ^{2}A+\frac{\Gamma }{4}),$ \ \par 
so if we choose  $\displaystyle \Gamma <1,\ (n-1)\delta ^{2}A<\frac{1}{5},\
 \delta <\frac{1}{5}$  and  $z\in B(0,R)$  we get\ \par 
\quad \quad \quad \quad \quad \quad \quad  $\lambda (x_{1},z')<t(-\frac{3}{4}+\frac{1}{5}+\frac{1}{5}+\frac{1}{4})<0,$
  i.e.  $Q_{a}(\delta )\subset C.$ \ \par 
So it remains to choose  $\delta $  with  $\delta <\min (\Gamma
 ,\frac{1}{2{\sqrt{(n-1)A}}},1/5),\ \Gamma <\min (\delta _{0},1)$
  and  $z\in B(0,R)$  to have  $Q_{a}(\delta )\subset C.$ \ \par 
The family  $\lbrace C_{\alpha }\rbrace $  is determined by 
 $\displaystyle \lbrace b(\alpha ),\ m(\alpha )\rbrace _{\alpha
 \in \partial \Omega },$  the aperture  $\Gamma $  and the number
  $A.$   $\blacksquare $ \ \par 
\ \par 
\quad \quad 	It would be nice to have an actual ellipsoid domain osculating
  $\Omega $  at  $\alpha ,$  instead of a conic domain whose
 slices parallel to the complex tangent space centered on the
 real normal are convex ellipsoids.\ \par 
But this is not true in general as shown by the following simple
 example in  ${\mathbb{C}}^{2}.$ \ \par 
Take  $\Omega =\lbrace \rho <0\rbrace $  near  $0,$  with :\ \par 
\quad \quad \quad \quad \quad  $\rho (z)=x_{1}+ay_{1}^{2}+b\left\vert{z_{2}}\right\vert ^{m}+cy_{1}x_{2},$
 \ \par 
with  $m\geq 3,\ c>0.$  Then there is no way to have that  $C:=\lbrace
 \gamma <0\rbrace \subset \Omega $  near  $0$  with :\ \par 
\quad \quad \quad \quad \quad  $\lambda (z)=x_{1}+Ay_{1}^{2}+B\left\vert{z_{2}}\right\vert
 ^{^{m}}$  for any choice of  $A$  and  $B.$ \ \par 
Just take points  $\displaystyle y_{1}=\frac{x_{2}}{k},\ y_{2}=0$
  then  $\displaystyle z\in \partial C\Rightarrow -x_{1}=A\frac{x_{2}^{2}}{k^{2}}+Bx_{2}^{m}$
  and at this point we have\ \par 
\quad \quad \quad \quad \quad  $\displaystyle \rho (z)=(a-A)\frac{x_{2}^{2}}{k^{2}}+(b-B)x_{2}^{m}+c\frac{x_{2}^{2}}{k}$
 \ \par 
and this is {\sl not} negative for  $k$  big enough and  $x_{2}$
  small enough if  $m\geq 3.$ \ \par 
\ \par 
\quad \quad 	Now we shall see that we have a converse to lemma~\ref{1_BonFamille20}.\ \par 

\begin{Lmm}
If a domain  $\Omega $  contains a family of complex tangentially
 elliptic domains  $\displaystyle \lbrace C_{\alpha }\rbrace
 _{\alpha \in \partial \Omega }$  based on  $\displaystyle \lbrace
 b(\alpha ),\ m(\alpha )\rbrace _{\alpha \in \partial \Omega
 },$  aperture  $\displaystyle \Gamma >0,\ \Gamma <1,$  and number
  $A,$  then  $\Omega $  possesses a good family of polydiscs
 still based on  $\displaystyle \lbrace b(\alpha ),\ m(\alpha
 )\rbrace _{\alpha \in \partial \Omega }$  and with parameter
  $\displaystyle \delta _{0}=\min (\Gamma ,\frac{1}{2{\sqrt{(n-1)A}}},1/5).$ 
\end{Lmm}
\quad \quad 	Proof.\ \par 
This is a remake of the proof of the last item in Lemma~\ref{1_BonFamille20}.\
 \par 
Again we make the canonical change of variables associated to
 the basis  $b_{\alpha }\ ;$  we have\ \par 
\quad \quad \quad \quad \quad  $\rho _{\alpha }(z)=x_{1}+g_{\alpha }(z)+h_{\alpha }(z'),$ \ \par 
where  $z_{1}=x_{1}+iy_{1},\ z':=(z_{2},...,z_{n})$  and  $g=g_{\alpha
 }\in {\mathcal{C}}^{\infty }({\mathbb{C}}^{n}),\ h=h_{\alpha
 }\in {\mathcal{C}}^{\infty }({\mathbb{C}}^{n-1}),\ g(\alpha
 )=h(\alpha ')=0.$ \ \par 
\quad \quad 	Let  $\displaystyle a=(-t,0,...,0)$  and fix a  $\delta >0$
  to be precised later ; if  $(x_{1}+iy_{1},z')\in Q_{a}(\delta )$  then\ \par 
\quad \quad \quad \quad \quad  $\displaystyle \ \forall j=2,...,n,\ \left\vert{z_{j}}\right\vert
 <\delta t^{1/m_{j}},\ \left\vert{x_{1}+t}\right\vert <\delta
 t\Rightarrow x_{1}<-t(1-\delta ),\ \left\vert{y_{1}}\right\vert
 <\delta t,$ \ \par 
so we already choose  $\displaystyle \delta \leq \Gamma $  and\ \par 
\quad \quad \quad \quad \quad \quad \quad  $\displaystyle A{\left({\left\vert{z_{2}}\right\vert ^{m_{2}}+\cdot
 \cdot \cdot +\left\vert{z_{n}}\right\vert ^{m_{n}}}\right)}<tA\sum_{j=2}^{n}{\delta
 ^{m_{j}}},$ \ \par 
hence\ \par 
\quad \quad  $\displaystyle \lambda (x_{1},y_{1},z')=x_{1}+g(z)+tA{\left({\left\vert{z_{2}}\right\vert
 ^{m_{2}}+\cdot \cdot \cdot +\left\vert{z_{n}}\right\vert ^{m_{n}}}\right)}<-t(1-\delta
 )+tA\sum_{j=2}^{n}{\delta ^{m_{j}}}+\left\vert{g(z)}\right\vert .$ \ \par 
We get, because  $m_{j}\geq 2,$ 		 $\ A\delta ^{m_{j}}\leq \delta
 ^{2}A,$  so\ \par 
\quad \quad \quad \quad \quad  $\lambda (x_{1},y_{1},z')<-t(1-\delta )+(n-1)\delta ^{2}A+\left\vert{g(z)}\right\vert
 .$ \ \par 
But, by lemma~\ref{1_goodFam71}, the smooth function  $g(z)$
  is bounded on  $B(0,R)$  by  $\ \frac{1}{4}\ \left\vert{z_{1}}\right\vert
 ,$  with  $R>0$  independent of  $\alpha .$  So as in the proof
 of the last item in Lemma~\ref{1_BonFamille20}, if	 	 $\displaystyle
 \Gamma <1,\ \delta <\frac{1}{5},\ \delta ^{2}(n-1)A<\frac{1}{5},$
  we have  $\displaystyle Q_{a}(\delta )\subset C.$ \ \par 
So it remains to choose  $\displaystyle \Gamma <1,\ \delta $
  with  $\displaystyle \delta <\min (\Gamma ,\frac{1}{2{\sqrt{(n-1)A}}},1/5),$
  which is independent of  $\alpha ,$  and  $z\in B(0,R)$  to
 have     $Q_{a}(\delta )\subset C.$ \ \par 
This means that\ \par 
\quad \quad \quad \quad \quad  $(z_{1},z')\in B(0,R)\cap Q_{a}(\delta )\Rightarrow (z_{1},z')\in C.$ \ \par 
For  $d(a)<\frac{1}{2}R^{M({\mathcal{Q}})}\leq \frac{1}{2}R^{m_{n}(\alpha
 )},$  because we can always choose  $\displaystyle R\leq 1,$
  with  $M({\mathcal{Q}})=\sup _{\alpha \in \partial \Omega }m_{n}(\alpha
 )<\infty ,$  then  $Q_{a}(\delta _{0})\subset B(\alpha ,R)$
  hence in this case  $Q_{a}(\delta _{0})\subset C\subset \Omega
 .$   $\blacksquare $ \ \par 
\ \par 
\quad \quad 	Together these lemmas proved\ \par 

\begin{Thrm}
Let   $\Omega $  be a domain in  ${\mathbb{C}}^{n}$  ; there
 is a good family of polydiscs in  $\Omega $  with multi-type
  $\lbrace b(\alpha ),\ m(\alpha )\rbrace _{\alpha \in \partial
 \Omega }$  iff there is a family of complex tangentially ellipsoids
   $\lbrace C_{\alpha }\rbrace _{\alpha \in \partial \Omega }$
  with parameters  $\lbrace b(\alpha ),\ m(\alpha )\rbrace _{\alpha
 \in \partial \Omega }$  such that  $\forall \alpha \in \partial
 \Omega ,\ C_{\alpha }\cap B(\alpha ,R)\subset \Omega \cap B(\alpha
 ,R),$  where  $R$  is given by lemma~\ref{1_goodFam71}.
\end{Thrm}
\ \par 
\vfill\eject\ \par 

\section{Divisors of the Blaschke class.~\label{aspcGlo63}}
\quad Let  $u$  be a holomorphic function in a domain  $\Omega ,$ 
 set  $X:=u^{-1}(0)$  its zero set and  $\Theta :=\partial \bar
 \partial \ln  \left\vert{u}\right\vert $  its associated  $(1,1)$
  current of integration. We shall define a class of such zero
 sets containing the zero sets of Nevanlinna functions.\ \par 

\begin{Dfnt}
A holomorphic divisor  $X$  in the domain  $\Omega $  is in the
 {\bf Blaschke class} if, with  $\Theta $  its associated  $(1,1)$
  current of integration,\par 
\[\displaystyle \ {\left\Vert{\Theta }\right\Vert}_{B}:=\int_{\Omega
 }{d(z)\mathrm{T}\mathrm{r}\Theta (z)}<\infty .\] 
\end{Dfnt}
\quad Let  $S$  be a separated sequence of points in  $\Omega $  contained
 in the zero set  $X$  in the Blaschke class of  $\Omega .$ 
 The aim of this section is to show that the measure  $\mu :=\sum_{a\in
 S}{d(a)^{n}\delta _{a}}$  is finite.\ \par 
\ \par 
\quad We shall need the easy lemma :\ \par 

\begin{Lmm}
~\label{Divisor3}Let  ${\mathcal{Q}}=\lbrace Q_{a}(\delta ),\
 a\in \Omega \rbrace $  be a good family of polydiscs for  $\Omega
 $  with parameter  $\delta _{0}$  and  $\displaystyle \delta
 <\delta _{0}.$  Then we have\par 
\quad \quad \quad $\displaystyle \forall a\in \Omega ,\ \forall z\in Q_{a}(\delta
 ),\ d(a)\leq \frac{1}{\delta _{0}-\delta }d(z,\partial \Omega ).$ 
\end{Lmm}
\quad Proof.\ \par 
We have by definition  $\displaystyle Q_{a}(\delta _{0})\subset
 \Omega ,$  hence  $\displaystyle \forall z\in Q_{a}(\delta ),\
 d(z,\Omega ^{c})\geq d(z,\ Q_{a}(\delta _{0})^{c}),$  but because
  $m_{j}(a)\geq m_{1}(a)=1,d(a)\leq d(a)^{1/m_{j}(a)}$  by the
 construction of the polydisc  $Q_{a}(\delta )$  we have\ \par 
\quad $\displaystyle \forall z\in Q_{a}(\delta ),\ d(z,\partial \Omega
 )\geq d(z,\ Q_{a}(\delta _{0})^{c})\geq \min _{j=1,...,n}(\delta
 _{0}-\delta )d(a)^{1/m_{j}(a)}\geq (\delta _{0}-\delta )d(a).$
   $\blacksquare $ \ \par 

\subsection{The discretized Blaschke condition.}
\quad Let  $u\in {\mathcal{H}}(\Omega ),$  i.e.  $u$  is holomorphic
 in  $\Omega ,$  and let  $X:=u^{-1}(0)\ ;$  put  $\Theta :=\partial
 \bar \partial \ln   \left\vert{u}\right\vert $  the  $(1,1)$
  current associated to  $X.$  Recall that  $\Delta \ln  \left\vert{u(z)}\right\vert
 =Tr\Theta ,$  the trace of  $\Theta ,$  and  $\Theta $  is a
 positive current, hence its trace controls all its coefficients.\ \par 
\quad We have, for any open set  ${\mathcal{V}}\subset \Omega ,$  the
 equality (see for instance ~\cite{LelongGrum86} p 55)\ \par 
\begin{equation}  \ \int_{{\mathcal{V}}}{Tr\Theta }=\sigma _{2n-2}(X\cap
 {\mathcal{V}}).\label{Divisor6}\end{equation}\ \par 
\ \par 
\quad \quad 	Let  $E_{j}:=\lbrace z\in {\mathbb{C}}^{n}::z_{j}=0\rbrace ,$
  this is the sub space orthogonal to the  $\displaystyle z_{j}$
  complex plane.	 	Let  $\Omega $  be a domain equipped with
 a good family  ${\mathcal{Q}}$  of polydiscs and  $X$  a divisor
 in  $\Omega .$  We set for  $a\in \Omega ,\ X_{a}:=X\cap Q_{a}(\delta
 ),\ X_{a}^{j}$  the projection of  $\displaystyle X_{a}$  on
  $\displaystyle E_{j}$  and  $\displaystyle A_{j}(X_{a}):=\sigma
 _{2n-2}(X_{a}^{j}).$ \ \par 
\quad Now let  $S$  be a  $\delta $  separated sequence of points on
  $X.$  We have the discretized Blaschke condition:\ \par 

\begin{Thrm}
(Discretized Blaschke condition) Let  $u$  be holomorphic in
  $\Omega ,\ X:=u^{-1}(0)$  and  $\Theta :=\partial \bar \partial
 \ln   \left\vert{u}\right\vert $  its current of integration
 ; suppose that  $\Theta $  is in the Blaschke class. Let  $S$
  a  $\delta $  separated sequence in  $X$  with respect to a
 good family  ${\mathcal{Q}}$  of polydiscs with parameter  $\delta
 _{0}.$  Then we have, provided that  $\displaystyle \delta <\delta
 _{0}/2,$ \par 
\quad \quad \quad \quad 	\begin{equation}  \ \sum_{a\in S}{d(a)\sigma _{2n-2}(X_{a})}\leq
 \frac{2}{\delta _{0}}{\left\Vert{\Theta }\right\Vert}_{B}.\label{Divisor7}\end{equation}
\end{Thrm}

      Proof.\ \par 
Let  $a\in S$  then by lemma~\ref{Divisor3} we have  $\displaystyle
 \forall z\in Q_{a}(\delta ),\ d(z)\geq \delta _{0}-\delta \geq
 \frac{\delta _{0}}{2}d(a).$  Now\ \par 
\quad \quad \quad $\displaystyle \ \sum_{a\in S}{\int_{Q_{a}(\delta )}{d(z)Tr\Theta
 }}\leq \int_{\Omega }{d(z)Tr\Theta }={\left\Vert{\Theta }\right\Vert}_{B},$
 \ \par 
because  $S$  is  $\delta $  separated, hence the polydiscs 
 $Q_{a}(\delta )$  are disjoint. Then\ \par 
\quad \quad \quad $\displaystyle \ {\left\Vert{\Theta }\right\Vert}_{B}\geq \frac{\delta
 _{0}}{2}\sum_{a\in S}{d(a)\int_{Q_{a}(\delta )}{Tr\Theta }}=\frac{\delta
 _{0}}{2}\sum_{a\in S}{d(a)\sigma _{2n-2}(X_{a})}.$   $\blacksquare $ \ \par 

\subsection{The discretized Malliavin condition.}
\quad Let us set  $\gamma :=i\sum_{j=1}^{n}{\,dz_{j}\wedge \,d\bar
 z_{j}},$  we have that  $\partial \gamma =\bar \partial \gamma
 =0$  and  $\gamma $  is a positive  $(1,1)$  form. We shall
 follow the proof by H. Skoda~\cite{zeroSkoda}, p 277.\ \par 
Set  $\beta :=\gamma ^{\wedge (n-2)}$  and apply Stokes formula
 to  $\rho \Theta \wedge \bar \partial \rho \wedge \beta \ :$ \ \par 
\quad \quad \quad $\displaystyle 0=\int_{\partial \Omega }{\rho \Theta \wedge \bar
 \partial \rho \wedge \beta }=\int_{\Omega }{\Theta \wedge \partial
 \rho \wedge \bar \partial \rho \wedge \beta }-\int_{\Omega }{\rho
 \Theta \wedge \partial \bar \partial \rho \wedge \beta },$ \ \par 
because  $\Theta $  and  $\beta $  are closed. Hence\ \par 
\quad \quad \quad $\displaystyle \ \left\vert{\int_{\Omega }{\Theta \wedge \partial
 \rho \wedge \bar \partial \rho \wedge \beta }}\right\vert =\left\vert{\int_{\Omega
 }{\rho \Theta \wedge \partial \bar \partial \rho \bigwedge \beta
 }}\right\vert \leq {\left\Vert{\partial \bar \partial \rho \wedge
 \beta }\right\Vert}_{\infty }\int_{\Omega }{(-\rho )Tr\Theta }\leq $ \ \par 
\quad \quad \quad \quad \quad \quad \quad $\displaystyle \leq {\left\Vert{\partial \bar \partial \rho }\right\Vert}_{\infty
 }{\left\Vert{\partial \rho }\right\Vert}_{\infty }{\left\Vert{\beta
 }\right\Vert}_{\infty }\int_{\Omega }{d(z,\ \partial \Omega
 )\ Tr\Theta }\leq {\left\Vert{\partial \bar \partial \rho }\right\Vert}_{\infty
 }{\left\Vert{\partial \rho }\right\Vert}_{\infty }{\left\Vert{\beta
 }\right\Vert}_{\infty }{\left\Vert{\Theta }\right\Vert}_{B}<\infty ,$ \ \par 
because the trace of  $\Theta $  controls all its coefficients
 and  $(-\rho (z))\leq {\left\Vert{\partial \rho }\right\Vert}_{\infty
 }d(z,\ \partial \Omega ).$ \ \par 
The norm  $\ {\left\Vert{\beta }\right\Vert}_{\infty }$  is a
 constant depending only on the dimension  $n,$  hence we can
 set  $C(\rho ):={\left\Vert{\partial \bar \partial \rho }\right\Vert}_{\infty
 }{\left\Vert{\partial \rho }\right\Vert}_{\infty }{\left\Vert{\beta
 }\right\Vert}_{\infty }$  which depends only on the first two
 derivatives of the defining function  $\rho .$ \ \par 
Hence we proved\ \par 

\begin{Lmm}
~\label{Divisor5}We have the estimate :\par 
\quad \quad \quad $\displaystyle \ \left\vert{\int_{\Omega }{\Theta \wedge \partial
 \rho \wedge \bar \partial \rho \wedge \beta }}\right\vert \leq
 C(\rho ){\left\Vert{\Theta }\right\Vert}_{B}.$ 
\end{Lmm}
\quad \quad 	Set  $\displaystyle m_{n}':={\left\lceil{m_{n}}\right\rceil}.$ \ \par 

\begin{Lmm}
~\label{3_BlasMall42} If a real smooth function  $h(z)$  verifies\par 
\quad \quad \quad \quad \quad  $\ \left\vert{h(z)}\right\vert \leq \sum_{j=1}^{n}{\left\vert{z_{j}}\right\vert
 ^{m_{j}}}$ \par 
with  $m_{j}\geq 2$  and  $\ \left\vert{z_{j}}\right\vert ^{m_{j}}\leq
 d(a),$  then 	 $\partial h\wedge \bar \partial h=d(a)\Gamma
 (z),$  where  $\Gamma (z)$  is a positive bounded  $(1,1)$ 
 form with its sup norm controlled by the  $\displaystyle m_{n}'+1$
  derivatives of  $\displaystyle h.$ 
\end{Lmm}
\quad \quad 	Proof.\ \par 
We shall use lemma~\ref{1_goodFam40}, this time using complex
 variables notations, for the function  $h$  with no parameter;
 there are smooth functions  $\displaystyle f_{\alpha ,\beta
 }(z,\bar z)$  for  $\displaystyle \ \left\vert{\alpha }\right\vert
 +\left\vert{\beta }\right\vert =m_{n}'$  such that\ \par 
\quad \quad \quad \quad \quad  $\displaystyle h(z)=\sum_{k=0}^{m_{n}'-1}{\sum_{\alpha ,\beta
 ,\left\vert{\alpha }\right\vert +\left\vert{\beta }\right\vert
 =k}{a_{\alpha ,\beta }z^{\alpha }\bar z^{\beta }}}+\sum_{\alpha
 ,\beta ,\left\vert{\alpha }\right\vert +\left\vert{\beta }\right\vert
 =m_{n}'}{f_{\alpha ,\beta }(z,\bar z)z^{\alpha }\bar z^{\beta }}.$ \ \par 
with  $\displaystyle z^{\alpha }:=z_{1}^{\alpha _{1}}\cdot \cdot
 \cdot z_{n}^{\alpha _{n}}$  and the same for  $\bar z^{\beta }.$ \ \par 
Consider the path  $t\in \lbrack 0,\epsilon \rbrack \rightarrow
 z_{j}(t):=\zeta _{j}t^{1/m_{j}}$  then\ \par 
\quad \quad \quad \quad \quad  $z^{\alpha }=\zeta ^{\alpha }t^{\gamma (\alpha )},$  with  $\displaystyle
 \gamma (\alpha )=\sum_{j=1}^{n}{\frac{\alpha _{j}}{m_{j}}},$ \ \par 
hence\ \par 
\quad \quad \quad \quad \quad  $\displaystyle h(z(t))=\sum_{\alpha ,\beta ,\left\vert{\alpha
 }\right\vert +\left\vert{\beta }\right\vert <m_{n}'}{a_{\alpha
 ,\beta }\zeta ^{\alpha }\bar \zeta ^{\beta }t^{\gamma (\alpha
 )+\gamma (\beta )}}+\sum_{\alpha ,\beta ,\left\vert{\alpha }\right\vert
 +\left\vert{\beta }\right\vert =m_{n}'}{f_{\alpha ,\beta }(z(t),\bar
 z(t))\zeta ^{\alpha }\bar \zeta ^{\beta }t^{\gamma (\alpha )+\gamma
 (\beta )}}.$ \ \par 
We also have\ \par 
\quad \quad \quad \quad \quad  $\displaystyle \ \sum_{j=1}^{n}{\left\vert{z_{j}}\right\vert
 ^{m_{j}}}=t\sum_{j=1}^{n}{\left\vert{\zeta _{j}}\right\vert ^{m_{j}}}$ \ \par 
hence let  $s=\gamma (\alpha )+\gamma (\beta ),$  then for  $\displaystyle
 \ \left\vert{\alpha }\right\vert +\left\vert{\beta }\right\vert
 =m_{n}'$  we have\ \par 
\quad \quad \quad \quad \quad  $\displaystyle 1=\sum_{j=1}^{n}{\frac{\alpha _{j}}{m_{n}'}}+\sum_{j=1}^{n}{\frac{\beta
 _{j}}{m_{n}'}}\leq \sum_{j=1}^{n}{\frac{\alpha _{j}}{m_{j}}}+\sum_{j=1}^{n}{\frac{\beta
 _{j}}{m_{j}}}=s,$ \ \par 
because  $\displaystyle m_{j}\leq m_{n}\leq m_{n}'.$ \ \par 
The function  $\displaystyle s=\gamma (\alpha )+\gamma (\beta
 )$  can take only a finite number of values, say  $\displaystyle
 s_{1}<\cdot \cdot \cdot <s_{k},$  then  because  $\ \left\vert{h(z(t))}\right\vert
 \leq t\sum_{j=1}^{n}{\left\vert{\zeta _{j}}\right\vert ^{m_{j}}},$
  if  $\displaystyle s_{1}<1,$ \ \par 
\quad \quad \quad \quad \quad  $\displaystyle \ \left\vert{\sum_{\alpha ,\beta ,\ \gamma (\alpha
 )+\gamma (\beta )=s_{1}}{a_{\alpha ,\beta }\zeta ^{\alpha }\bar
 \zeta ^{\beta }}}\right\vert \leq t^{1-s_{1}}\sum_{j=1}^{n}{\left\vert{\zeta
 _{j}}\right\vert ^{m_{j}}}+$ \ \par 
\quad \quad \quad \quad \quad \quad \quad \quad  $\displaystyle +\sum_{s=s_{2}}^{s_{k}}{t^{s-s_{1}}\left\vert{\sum_{\alpha
 ,\beta ,\ \gamma (\alpha )+\gamma (\beta )=s}{a_{\alpha ,\beta
 }\zeta ^{\alpha }\bar \zeta ^{\beta }}}\right\vert }+\sum_{\alpha
 ,\beta ,\left\vert{\alpha }\right\vert +\left\vert{\beta }\right\vert
 =m_{n}'}{\left\vert{f_{\alpha ,\beta }(z(t))}\right\vert \left\vert{\zeta
 }\right\vert ^{\left\vert{\alpha }\right\vert +\left\vert{\beta
 }\right\vert }t^{\gamma (\alpha )+\gamma (\beta )-s_{1}}}.$ \ \par 
In the last sum we have  $\displaystyle \gamma (\alpha )+\gamma
 (\beta )\geq 1$  because  $\displaystyle \ \left\vert{\alpha
 }\right\vert +\left\vert{\beta }\right\vert =m_{n}'$  and the
 functions  $\displaystyle f_{\alpha ,\beta }$  are bounded,
 hence letting  $\displaystyle t\rightarrow 0,$  we get\ \par 
\quad \quad \quad \quad \quad  $\displaystyle \ \sum_{\alpha ,\beta ,\ \gamma (\alpha )+\gamma
 (\beta )=s_{1}}{a_{\alpha ,\beta }\zeta ^{\alpha }\bar \zeta
 ^{\beta }}=0.$ \ \par 
\quad \quad 	Then we can repeat the same  for  $\displaystyle s_{2},...,\
 s_{j}$  provided that  $\displaystyle s_{j}<1,$  hence we get\ \par 
\quad \quad \quad \quad \quad  $\displaystyle \ \sum_{\alpha ,\beta ,\ \gamma (\alpha )+\gamma
 (\beta )<1}{a_{\alpha ,\beta }\zeta ^{\alpha }\bar \zeta ^{\beta }}=0.$ \ \par 
\ \par 
So in the expansion of  $h$  it remains only  $\alpha ,\beta
 $  such that  $\gamma (\alpha )+\gamma (\beta )\geq 1.$ \ \par 
\quad 	Now we compute\ \par 
\quad \quad \quad \quad \quad  $\displaystyle \partial z^{\alpha }=\sum_{j=1}^{n}{\alpha _{j}z^{\alpha
 }/z_{j}dz_{j}}=z^{\alpha }\sum_{j=1}^{n}{\frac{\alpha _{j}}{z_{j}}dz_{j}}.$
 \ \par 
And\ \par 
\quad \quad \quad \quad \quad  $\displaystyle \bar \partial \bar z^{\beta }=\sum_{j=1}^{n}{\beta
 _{j}\bar z^{\beta }/\bar z_{j}d\bar z_{j}}=\bar z^{\beta }\sum_{j=1}^{n}{\frac{\beta
 _{j}}{\bar z_{j}}d\bar z_{j}}.$ \ \par 
Set  $\displaystyle \omega (z,\alpha ):=\sum_{j=1}^{n}{\frac{\alpha
 _{j}}{z_{j}}dz_{j}},$  we have\ \par 
\quad \quad \quad \quad \quad  $\displaystyle \partial h=\sum_{\alpha ,\beta ,1\leq \gamma
 (\alpha )+\gamma (\beta )<m_{n}'}{a_{\alpha ,\beta }z^{\alpha
 }\bar z^{\beta }\omega (z,\alpha )}+\sum_{\alpha ,\beta ,\left\vert{\alpha
 }\right\vert +\left\vert{\beta }\right\vert =m_{n}'}{f_{\alpha
 ,\beta }(z,\bar z)z^{\alpha }\bar z^{\beta }\omega (z,\alpha )}+$ \ \par 
\quad \quad \quad \quad \quad \quad \quad \quad \quad \quad  $\displaystyle +\sum_{\alpha ,\beta ,\left\vert{\alpha }\right\vert
 +\left\vert{\beta }\right\vert =m_{n}'}{z^{\alpha }\bar z^{\beta
 }\partial f_{\alpha ,\beta }(z,\bar z)}.$ \ \par 
and\ \par 
\quad \quad \quad \quad \quad  $\displaystyle \bar \partial h=\sum_{\alpha ,\beta ,\gamma (\alpha
 )+\gamma (\beta )\geq 1}{a_{\alpha ,\beta }z^{\alpha }\bar z^{\beta
 }\bar \omega (z,\beta )}+\sum_{\alpha ,\beta ,\left\vert{\alpha
 }\right\vert +\left\vert{\beta }\right\vert =m_{n}'}{f_{\alpha
 ,\beta }(z,\bar z)z^{\alpha }\bar z^{\beta }\bar \omega (z,\beta )}+$ \ \par 
\quad \quad \quad \quad \quad \quad \quad \quad \quad \quad  $\displaystyle +\sum_{\alpha ,\beta ,\left\vert{\alpha }\right\vert
 +\left\vert{\beta }\right\vert =m_{n}'}{z^{\alpha }\bar z^{\beta
 }\bar \partial f_{\alpha ,\beta }(z,\bar z)}.$ \ \par 
So we have as the generic term for  $\partial h\wedge \bar \partial h$ \ \par 
\quad \quad \quad \quad \quad  $\displaystyle Adz_{j}\wedge d\bar z_{k}:=z^{\alpha +\alpha
 '}\bar z^{\beta +\beta '}\frac{dz_{j}}{z_{j}}\wedge \frac{d\bar
 z_{k}}{\bar z_{k}}$ \ \par 
hence\ \par 
\quad \quad \quad \quad \quad  $\ \left\vert{A}\right\vert =\left\vert{z_{1}}\right\vert ^{\alpha
 _{1}+\beta _{1}+\alpha '_{1}+\beta '_{1}}\cdot \cdot \cdot \left\vert{z_{n}}\right\vert
 ^{\alpha _{n}+\beta _{n}+\alpha '_{n}+\beta '_{n}}\left\vert{z_{j}}\right\vert
 ^{-1}\left\vert{z_{k}}\right\vert ^{-1}$ \ \par 
with  $\ \left\vert{z_{l}}\right\vert ^{m_{l}}\leq d(a)$  we get\ \par 
\quad \quad \quad \quad \quad  $\ \left\vert{A}\right\vert \leq d(a)^{\gamma (\alpha )+\gamma
 (\beta )+\gamma (\alpha ')+\gamma (\beta ')}d(a)^{-1/m_{j}}d(a)^{-1/m_{k}}\leq
 d(a),$ \ \par 
because\ \par 
\quad \quad \quad \quad \quad  $\displaystyle \gamma (\alpha )+\gamma (\beta )+\gamma (\alpha
 ')+\gamma (\beta ')\geq 2$  and  $m_{j}\geq 2,m_{k}\geq 2.$ \ \par 
The special terms are of the forms\ \par 
\quad \quad \quad \quad \quad  $\displaystyle Bdz_{j}\wedge d\bar z_{k}:=f_{\alpha ,\beta }(z)z^{\alpha
 +\alpha '}\bar z^{\beta +\beta '}\frac{dz_{j}}{z_{j}}\wedge
 \frac{d\bar z_{k}}{\bar z_{k}}$ \ \par 
and, by the same argument, they verify  $\displaystyle \ \left\vert{B}\right\vert
 \leq {\left\Vert{f_{\alpha ,\beta }}\right\Vert}_{\infty }d(a)\ ;$ \ \par 
or\ \par 
\quad \quad \quad \quad \quad  $\displaystyle Cdz_{j}\wedge d\bar z_{k}:=\frac{\partial f_{\alpha
 ,\beta }}{\partial z_{j}}z^{\alpha +\alpha '}\bar z^{\beta +\beta
 '}dz_{j}\wedge \frac{d\bar z_{k}}{\bar z_{k}}$ \ \par 
and they verify a fortiori  $\displaystyle \ \left\vert{C}\right\vert
 \leq {\left\Vert{\partial f_{\alpha ,\beta }}\right\Vert}_{\infty
 }d(a)\ ;$ \ \par 
or\ \par 
\quad \quad \quad \quad \quad  $\displaystyle Ddz_{j}\wedge d\bar z_{k}:=f_{\alpha ,\beta }(z)f_{\alpha
 ',\beta '}(z)z^{\alpha +\alpha '}\bar z^{\beta +\beta '}\frac{dz_{j}}{z_{j}}\wedge
 \frac{d\bar z_{k}}{\bar z_{k}}$ \ \par 
and they verify  $\displaystyle \ \left\vert{D}\right\vert \leq
 {\left\Vert{f_{\alpha ,\beta }f_{\alpha ',\beta '}}\right\Vert}_{\infty
 }d(a)\ ;$ \ \par 
or\ \par 
\quad \quad \quad \quad \quad  $\displaystyle Edz_{j}\wedge d\bar z_{k}:=\frac{\partial f_{\alpha
 ,\beta }}{\partial z_{j}}\frac{\partial f_{\alpha ',\beta '}}{\partial
 z_{j}}z^{\alpha +\alpha '}\bar z^{\beta +\beta '}dz_{j}\wedge
 d\bar z_{k}$ \ \par 
and they verify  $\ \left\vert{E}\right\vert \leq {\left\Vert{\partial
 f_{\alpha ,\beta }\bar \partial f_{\alpha ',\beta '}}\right\Vert}_{\infty
 }d(a)\ ;$ \ \par 
or\ \par 
\quad \quad \quad \quad \quad  $\displaystyle Fdz_{j}\wedge d\bar z_{k}:=\frac{\partial f_{\alpha
 ,\beta }}{\partial z_{j}}f_{\alpha ,\beta }(z)z^{\alpha +\alpha
 '}\bar z^{\beta +\beta '}dz_{j}\wedge \frac{d\bar z_{k}}{\bar z_{k}}$ \ \par 
and they verify  $\displaystyle \ \left\vert{F}\right\vert \leq
 {\left\Vert{f_{\alpha ',\beta '}\partial f_{\alpha ,\beta }}\right\Vert}_{\infty
 }d(a)\ ;$ \ \par 
and the conjugates of these expressions are also bounded. All
 the bounds are controlled by the  $\displaystyle m_{n}'+1$ 
 derivatives of  $h$  and we have a finite set of such smooth
 coefficients so\ \par 
\quad \quad \quad \quad \quad  $\partial h\wedge \bar \partial h=d(a)\Gamma (z),$ \ \par 
where  $\Gamma (z)$  is a positive bounded  $(1,1)$  form controlled
 by the  $\displaystyle m_{n}'+1$  derivatives of  $h.$   $\blacksquare
 $ \ \par 
\ \par 
\quad \quad 	We shall evaluate the integral  $\displaystyle \ \int_{Q_{a}(\delta
 )}{\Theta \wedge \partial \rho \wedge \bar \partial \rho \wedge
 \beta },$  and we start first with the following  $\lambda ,$
  defining a complex tangential ellipsoid  $C$  as in the previous
 section, lemma 2.9, 2.10, but here we choose twice the previous
 one to ease the computations,\ \par 
\quad \quad \quad \quad \quad  $\lambda (z):=2x_{1}+2x_{1}g_{1}(z)+2y_{1}g_{2}(z)+2A(\left\vert{z_{2}}\right\vert
 ^{m_{2}}+\cdot \cdot \cdot +\left\vert{z_{n}}\right\vert ^{m_{n}}).$ \ \par 
We have, with  $a=(a_{1},0,...,0),\ a_{1}<0,$ \ \par 
\quad  $\displaystyle \partial \lambda (z)=(1+g_{1}(z)-ig_{2}(z)+2x_{1}\frac{\partial
 g_{1}}{\partial z_{1}}+2y_{1}\frac{\partial g_{2}}{\partial
 z_{1}})dz_{1}+2x_{1}\sum_{j\geq 2}{\frac{\partial g_{1}}{\partial
 z_{j}}(z)dz_{j}}+2y_{1}\sum_{j\geq 2}{\frac{\partial g_{2}}{\partial
 z_{j}}(z)dz_{j}}+$ \ \par 
\quad \quad \quad \quad \quad \quad \quad \quad  $\displaystyle +A(m_{2}\left\vert{z_{2}}\right\vert ^{m_{2}-2}\bar
 z_{2}dz_{2}+\cdot \cdot \cdot +m_{n}\left\vert{z_{2}}\right\vert
 ^{m_{n}-2}\bar z_{n}dz_{n})$ \ \par 
and\ \par 
\quad  $\displaystyle \bar \partial \lambda (z)=(1+g_{1}(z)+ig_{2}(z)+2x_{1}\frac{\partial
 g_{1}}{\partial \bar z_{1}}+2y_{1}\frac{\partial g_{2}}{\partial
 \bar z_{1}})d\bar z_{1}$  $\displaystyle +2x_{1}\sum_{j\geq
 2}{\frac{\partial g_{1}}{\partial z_{j}}d\bar z_{j}}+2y_{1}\sum_{j\geq
 2}{\frac{\partial g_{2}}{\partial z_{j}}d\bar z_{j}}+$ \ \par 
\quad \quad \quad \quad \quad \quad \quad \quad  $\displaystyle +A(m_{2}\left\vert{z_{2}}\right\vert ^{m_{2}-2}z_{2}d\bar
 z_{2}+\cdot \cdot \cdot +m_{n}\left\vert{z_{n}}\right\vert ^{m_{n}-2}z_{n}d\bar
 z_{n}),$ \ \par 
because\ \par 
\quad \quad \quad \quad \quad  $\displaystyle \ \frac{\partial \left\vert{w}\right\vert ^{m}}{\partial
 w}=\partial _{w}((\bar ww)^{m/2})=\frac{m}{2}(\bar ww)^{m/2-1}\bar
 w=\frac{m}{2}\left\vert{w}\right\vert ^{m-2}{\times}\bar w.$ \ \par 

\begin{Lmm}
~\label{3_BlasMall64}We have\par 
\quad \quad \quad \quad \quad  $\displaystyle \forall z\in Q_{a}(\delta ),\ \partial \lambda
 \wedge \bar \partial \lambda =$ \par 
\quad \quad \quad \quad \quad  $\displaystyle B(z)dz_{1}\wedge d\bar z_{1}+\sum_{j=2}^{n}{C_{j}(z)\left\vert{z_{j}}\right\vert
 ^{m_{j}-1}dz_{1}\wedge d\bar z_{j}}+\sum_{j=2}^{n}{D_{j}(z)\left\vert{z_{j}}\right\vert
 ^{m_{j}-1}dz_{j}\wedge d\bar z_{1}}+d(a)\Gamma (z),$ \par 
where  $\displaystyle B,C_{j},D_{j}$  are bounded with bounds
 depending only on the the  ${\mathcal{C}}^{1}$  norms of  $\displaystyle
 g_{1},g_{2}$  and  $\Gamma $  is a  $\displaystyle (1,1)$  form
 still with bounded coefficients depending only on the the  ${\mathcal{C}}^{1}$
  norms of  $\displaystyle g_{1},g_{2}.$ 
\end{Lmm}
\quad \quad 	Proof.\ \par 
Because  $\displaystyle \forall z\in Q_{a}(\delta )$  we have
  $\displaystyle \ \left\vert{z_{1}}\right\vert \leq \delta d(a)\Rightarrow
 \left\vert{x_{1}}\right\vert \leq \delta d(a),\ \left\vert{y_{1}}\right\vert
 \leq \delta d(a)$  and 		 $\displaystyle \forall j\geq 2,\ \
 \left\vert{z_{j}}\right\vert \leq \delta d(a)^{1/m_{j}},$  so
 the terms in  $\partial \lambda \wedge \bar \partial \lambda
 $  containing  $\displaystyle \ \frac{\partial g_{1}}{\partial
 z_{j}}$  or  $\displaystyle \ \frac{\partial g_{2}}{\partial
 z_{j}}$  or  $\displaystyle \ \frac{\partial g_{1}}{\partial
 \bar z_{j}}$  or  $\displaystyle \ \frac{\partial g_{2}}{\partial
 \bar z_{j}}$  can be put in  $\displaystyle \Gamma .$ \ \par 
\quad \quad 	For the terms in\ \par 
\quad \quad \quad \quad \quad  $\displaystyle A^{2}\sum_{j,k=1,...,n}{m_{j}m_{k}\left\vert{z_{j}}\right\vert
 ^{m_{j}-2}\left\vert{z_{k}}\right\vert ^{m_{k}-2}\bar z_{j}z_{k}dz_{j}\wedge
 d\bar z_{k}}\ $ \ \par 
we have\ \par 
\quad \quad \quad \quad \quad \quad \quad \quad \quad \quad  $\displaystyle \forall j,k\geq 2,\ \left\vert{z_{j}}\right\vert
 ^{m_{j}-1}\left\vert{z_{k}}\right\vert ^{m_{k}-1}\leq \delta
 ^{2}d(a)^{\frac{m_{j}-1}{m_{j}}+\frac{m_{k}-1}{m_{k}}},$ \ \par 
suppose that  $m_{j}\geq m_{k},$  then\ \par 
\quad \quad \quad \quad \quad  $\displaystyle \ \frac{m_{j}-1}{m_{j}}+\frac{m_{k}-1}{m_{k}}\geq
 \frac{m_{j}-1+m_{k}-1}{m_{j}}\geq \frac{m_{j}}{m_{j}}=1,$ \ \par 
because  $\forall k\geq 1,\ m_{k}\geq 2.$  Hence, they also can
 be put in  $\displaystyle \Gamma .$ \ \par 
\quad Hence it remains\ \par 
\quad \quad \quad \quad \quad  $\displaystyle B(z)dz_{1}\wedge d\bar z_{1}$  with  $\displaystyle
 B(z):=(1+g_{1})^{2}+g_{2}^{2},$ \ \par 
\quad \quad \quad \quad \quad  $\displaystyle \ \sum_{j=2}^{n}{C_{j}(z)\left\vert{z_{j}}\right\vert
 ^{m_{j}-2}z_{j}dz_{1}\wedge d\bar z_{j}}$  with  $\displaystyle
 C_{j}(z):=(1+g_{1}(z)-ig_{2}(z))Am_{j},$ \ \par 
and\ \par 
\quad \quad \quad \quad \quad  $\displaystyle \ \sum_{j=2}^{n}{D_{j}(z)\left\vert{z_{j}}\right\vert
 ^{m_{j}-2}\bar z_{j}dz_{j}\wedge d\bar z_{1}}$  with  $\displaystyle
 D_{j}(z):=(1+g_{1}(z)+ig_{2}(z))Am_{j}.$ \ \par 
Clearly the bounds on those terms and in  $\displaystyle \Gamma
 $  depend only on the  ${\mathcal{C}}^{1}$  norms of  $\displaystyle
 g_{1},g_{2}.$   $\blacksquare $ \ \par 

\begin{Lmm}
~\label{3_BlasMall65}Let  $\Theta $  be a positive  $\displaystyle
 (1,1)$  current and  $\displaystyle F(z_{j})$  a function ;
 then for all  $\displaystyle \eta >0$  we have\par 
\quad \quad \quad \quad \quad  $\ 2\left\vert{\Theta \wedge F(z_{j})dz_{1}\wedge d\bar z_{j}\wedge
 \beta }\right\vert \leq \eta \Theta \wedge dz_{1}\wedge d\bar
 z_{1}\wedge \beta +\frac{1}{\eta }\Theta \wedge \left\vert{F(z_{j})}\right\vert
 ^{2}dz_{j}\wedge d\bar z_{j}\wedge \beta .$ 
\end{Lmm}
\quad \quad 	Proof.\ \par 
by Cauchy-Schwarz, because  $\Theta \wedge \beta $  is positive, we get\ \par 
\quad \quad \quad \quad \quad  $\displaystyle \ \left\vert{\Theta \wedge dz_{1}\wedge F(z_{j})d\bar
 z_{j}\wedge \beta }\right\vert ^{2}\leq \left\vert{\Theta \wedge
 dz_{1}\wedge d\bar z_{1}\wedge \beta }\right\vert \left\vert{\Theta
 \wedge F(z_{j})dz_{j}\wedge \bar F(z_{j})d\bar z_{j}\wedge \beta
 }\right\vert $ \ \par 
hence, because  $\displaystyle 2ab\leq \eta a^{2}+\frac{1}{\eta }b^{2},$ \ \par 
\quad \quad \quad \quad \quad  $\displaystyle 2\ \left\vert{\Theta \wedge dz_{1}\wedge F(z_{j})d\bar
 z_{j}\wedge \beta }\right\vert \leq \eta \Theta \wedge dz_{1}\wedge
 d\bar z_{1}\wedge \beta +\frac{1}{\eta }\Theta \wedge \left\vert{F(z_{j})}\right\vert
 ^{2}dz_{j}\wedge d\bar z_{j}\wedge \beta .$   $\blacksquare $ \ \par 
\quad \quad 	Now let us go back to the general case. Fix  $a\in \Omega ,\
 \alpha =\pi (a),$  we know by lemma~\ref{1_BonFamille20} that
 there is a complex tangential ellipsoid  $C=C_{\alpha }$  with
 exponents  $\lbrace m_{j}(\alpha )\rbrace $  meeting  $\displaystyle
 \partial \Omega $  at  $\alpha $  and contained in  $\Omega
 .$  Moreover we have, after the canonical change of variables
 of lemma~\ref{1_goodFam71}, and multiplying by  $2$  the functions
 for making the following computations slightly easier,\ \par 
\quad \quad \quad \quad \quad  $\rho (z)=2x_{1}+2x_{1}g_{1}(z)+2y_{1}g_{2}(z)+h_{\alpha }(z')=2x_{1}+2x_{1}g_{1}(z)+2y_{1}g_{2}(z)+\rho
 (0,z'),$ \ \par 
as the defining function for  $\Omega $  and\ \par 
\quad \quad \quad \quad \quad  $\displaystyle \lambda (z):=2x_{1}+2x_{1}g_{1}(z)+2y_{1}g_{2}(z)+2A(\left\vert{z_{2}}\right\vert
 ^{m_{2}}+\cdot \cdot \cdot +\left\vert{z_{n}}\right\vert ^{m_{n}}),\
 \left\vert{y_{1}}\right\vert <-\Gamma x_{1},$ \ \par 
as the defining functions for  $C_{\alpha }.$  We notice that
 the functions  $g_{j}$  in  $\lambda $  are the same as the
 functions  $g_{j}$  in  $\rho $  and depend only on  $\rho .$
  In particular the  ${\mathcal{C}}^{1}$  norms of the  $\displaystyle
 g_{j}$  are controlled by the  ${\mathcal{C}}^{2}$  norm of
  $\rho $  hence they are uniformly bounded with respect to   $\alpha .$ \ \par 

\begin{Lmm}
~\label{AG0}We have, with  $\Theta $  a positive  $\displaystyle
 (1,1)$  current,\par 
\quad \quad  $\displaystyle \ \int_{Q_{a}(\delta )}{\Theta \wedge dz_{1}\wedge
 d\bar z_{1}\wedge \beta }\leq 5\int_{Q_{a}(\delta )}{\Theta
 \wedge \partial \lambda \wedge \bar \partial \lambda \wedge
 \beta }+\Gamma d(a)\int_{Q_{a}(\delta )}{\mathrm{T}\mathrm{r}\Theta }.$ \par 
with the constant  $\Gamma $  depending only on the  ${\mathcal{C}}^{2}$
  norm of  $\rho ,$  on  $\displaystyle n$  and  $\delta _{0}.$ 
\end{Lmm}
\quad \quad 	Proof.\ \par 
Using lemma~\ref{3_BlasMall64}, we get\ \par 
\quad \quad \quad \quad  $\displaystyle \Theta \wedge \partial \lambda \wedge \bar \partial
 \lambda \wedge \beta -B(z)\Theta \wedge dz_{1}\wedge d\bar z_{1}\wedge
 \beta =$ \ \par 
\quad \quad \quad \quad \quad \quad \quad  $\displaystyle \ \sum_{j=2}^{n}{C_{j}(z)\left\vert{z_{j}}\right\vert
 ^{m_{j}-2}\bar z_{j}\Theta \wedge dz_{1}\wedge d\bar z_{j}\wedge
 \beta }+$ \ \par 
\quad \quad \quad \quad \quad \quad \quad \quad \quad \quad  $\displaystyle +\sum_{j=2}^{n}{D_{j}(z)\left\vert{z_{j}}\right\vert
 ^{m_{j}-2}z_{j}\Theta \wedge dz_{j}\wedge d\bar z_{1}\wedge
 \beta }+d(a)\Theta \wedge \Gamma \wedge \beta .$ \ \par 
Hence\ \par 
\quad \quad \quad \quad \quad  $\displaystyle B(z)\Theta \wedge dz_{1}\wedge d\bar z_{1}\wedge
 \beta =\Theta \wedge \partial \lambda \wedge \bar \partial \lambda
 \wedge \beta -U-d(a)\Theta \wedge \Gamma \wedge \beta ,$ \ \par 
with\ \par 
\quad \quad \quad \quad \quad  $\displaystyle U:=\sum_{j=2}^{n}{C_{j}(z)\left\vert{z_{j}}\right\vert
 ^{m_{j}-2}\bar z_{j}\Theta \wedge dz_{1}\wedge d\bar z_{j}\wedge
 \beta }+\sum_{j=2}^{n}{D_{j}(z)\left\vert{z_{j}}\right\vert
 ^{m_{j}-2}z_{j}\Theta \wedge dz_{j}\wedge d\bar z_{1}\wedge \beta }$ \ \par 
By lemma~\ref{3_BlasMall65} we get, with  $\eta >0$  to be fixed later,\ \par 
\quad \quad \quad \quad \quad  $\displaystyle 2\left\vert{C_{j}(z)\left\vert{z_{j}}\right\vert
 ^{m_{j}-1}\Theta \wedge dz_{1}\wedge d\bar z_{j}\wedge \beta
 }\right\vert \leq $ \ \par 
\quad \quad \quad \quad \quad \quad \quad  $\displaystyle \eta \Theta \wedge dz_{1}\wedge d\bar z_{1}\wedge
 \beta +\frac{1}{\eta }\left\vert{C_{j}}\right\vert ^{2}\left\vert{z_{j}}\right\vert
 ^{2m_{j}-2}\Theta \wedge dz_{j}\wedge d\bar z_{j}\wedge \beta .$ \ \par 
But for  $\displaystyle z\in Q_{a}(\delta )$  we have  $\ \left\vert{z_{j}}\right\vert
 ^{2m_{j}-2}\leq \left\vert{z_{j}}\right\vert ^{m_{j}-2}\delta
 d(a)$  because  $\displaystyle \ \left\vert{z_{j}}\right\vert
 ^{m_{j}}\leq \delta d(a)$   hence\ \par 
\quad \quad \quad \quad \quad  $2\left\vert{C_{j}(z)\left\vert{z_{j}}\right\vert ^{m_{j}-1}\Theta
 \wedge dz_{1}\wedge d\bar z_{j}\wedge \beta }\right\vert \leq $ \ \par 
\quad \quad \quad \quad \quad \quad \quad  $\displaystyle \eta \Theta \wedge dz_{1}\wedge d\bar z_{1}\wedge
 \beta +\frac{1}{\eta }\ \left\vert{C_{j}}\right\vert ^{2}\ \left\vert{z_{j}}\right\vert
 ^{m_{j}-2}\delta d(a)\Theta \wedge dz_{j}\wedge d\bar z_{j}\wedge
 \beta .$ \ \par 
Set  $\displaystyle C'_{j}:=2\left\vert{C_{j}}\right\vert ^{2}\left\vert{z_{j}}\right\vert
 ^{m_{j}-2}$  whose bound depend on the  ${\mathcal{C}}^{1}$
  norm of the  $\displaystyle g_{j}$  we get\ \par 
\quad \quad \quad \quad \quad  $\displaystyle \ \left\vert{\sum_{j=2}^{n}{C_{j}(z)\left\vert{z_{j}}\right\vert
 ^{m_{j}-2}\bar z_{j}\Theta \wedge dz_{1}\wedge d\bar z_{j}\wedge
 \beta }}\right\vert \leq $ \ \par 
\quad \quad \quad \quad \quad \quad \quad  $\displaystyle \ \frac{1}{2}(n-1)\eta \Theta \wedge dz_{1}\wedge
 d\bar z_{1}\wedge \beta +\frac{\delta }{\eta }d(a)\sum_{j=2}^{n}{C'_{j}\Theta
 \wedge dz_{j}\wedge d\bar z_{j}\wedge \beta }.$ \ \par 
Doing exactly the same proof, with  $\displaystyle D'_{j}:=2\left\vert{D_{j}}\right\vert
 ^{2}\left\vert{z_{j}}\right\vert ^{m_{j}-2}$  we get\ \par 
\quad \quad \quad \quad \quad  $\displaystyle \ \left\vert{\sum_{j=2}^{n}{D_{j}(z)\left\vert{z_{j}}\right\vert
 ^{m_{j}-2}z_{j}\Theta \wedge dz_{j}\wedge d\bar z_{1}\wedge
 \beta }}\right\vert \leq $ \ \par 
\quad \quad \quad \quad \quad \quad \quad  $\displaystyle \ \frac{1}{2}(n-1)\eta \Theta \wedge dz_{1}\wedge
 d\bar z_{1}\wedge \beta +\frac{\delta }{\eta }d(a)\sum_{j=2}^{n}{D'_{j}\Theta
 \wedge dz_{j}\wedge d\bar z_{j}\wedge \beta }.$ \ \par 
\ \par 
So we get\ \par 
\quad \quad \quad \quad \quad \quad \quad  $\displaystyle \ \left\vert{U}\right\vert \leq (n-1)\eta \Theta
 \wedge dz_{1}\wedge d\bar z_{1}\wedge \beta +\frac{\delta }{\eta
 }d(a)\sum_{j=2}^{n}{(C'_{j}+D'_{j})\Theta \wedge dz_{j}\wedge
 d\bar z_{j}\wedge \beta }.$ \ \par 
\quad \quad 	Now we choose		  $\displaystyle \eta :=\frac{1}{4(n-1)}$  and
 we get, because  $\displaystyle \Theta \wedge \partial \lambda
 \wedge \bar \partial \lambda \wedge \beta $  and  $\displaystyle
 B(z)\Theta \wedge dz_{1}\wedge d\bar z_{1}\wedge \beta $  are positive,\ \par 
\quad \quad \quad \quad \quad  $\displaystyle B(z)\Theta \wedge dz_{1}\wedge d\bar z_{1}\wedge
 \beta \leq \Theta \wedge \partial \lambda \wedge \bar \partial
 \lambda \wedge \beta +\left\vert{U}\right\vert +d(a)\left\vert{\Theta
 \wedge \Gamma \wedge \beta }\right\vert .$ \ \par 
Hence\ \par 
\quad \quad \quad  $\displaystyle B(z)\Theta \wedge dz_{1}\wedge d\bar z_{1}\wedge
 \beta \leq \Theta \wedge \partial \lambda \wedge \bar \partial
 \lambda \wedge \beta +\frac{1}{4}\Theta \wedge dz_{1}\wedge
 d\bar z_{1}\wedge \beta +d(a)(4(n-1)\delta (a)\left\vert{\Theta
 \wedge \Gamma '\wedge \beta }\right\vert +\left\vert{\Theta
 \wedge \Gamma \wedge \beta }\right\vert ).$ \ \par 
with\ \par 
\quad \quad \quad \quad \quad  $\displaystyle \Gamma ':=\sum_{j=2}^{n}{C'_{j}dz_{j}\wedge d\bar
 z_{j}}+\sum_{j=2}^{n}{D'_{j}dz_{j}\wedge d\bar z_{j}}.$ \ \par 
So finally\ \par 
\quad \quad \quad \quad \quad  $\displaystyle (B(z)-\frac{1}{4})\Theta \wedge dz_{1}\wedge
 d\bar z_{1}\wedge \beta \leq \Theta \wedge \partial \lambda
 \wedge \bar \partial \lambda \wedge \beta +d(a)(4(n-1)\delta
 (a)\left\vert{\Theta \wedge \Gamma '\wedge \beta }\right\vert
 +\left\vert{\Theta \wedge \Gamma \wedge \beta }\right\vert ).$ \ \par 
\ \par 
\quad Recall that  $\displaystyle B(z)=(1+g_{1}(z))^{2}+g_{2}(z)^{2}\geq
 0$  and we know by lemma~\ref{1_goodFam71} that  $\displaystyle
 \forall z\in B(0,R),\ \left\vert{g_{1}(z)}\right\vert \leq \frac{3}{10}$
  hence, provided that  $\displaystyle Q_{a}(\delta )\subset
 B(0,R),$  i.e.  $d(a)<R^{1/M({\mathcal{Q}})},$ with  $M({\mathcal{Q}})=\sup
 _{\alpha \in \partial \Omega }m_{n}(\alpha )<\infty ,$  because
  ${\mathcal{Q}}$  is a good family, we get\ \par 
\quad \quad \quad \quad \quad  $\displaystyle (1+g_{1}(z))^{2}+g_{2}(z)^{2}-1/4\geq (\frac{7}{10})^{2}-\frac{1}{4}=\frac{24}{100}.$
 \ \par 
So dividing by  $\displaystyle B(z)-\frac{1}{4}$  we get\ \par 
\quad \quad \quad \quad \quad  $\displaystyle \Theta \wedge dz_{1}\wedge d\bar z_{1}\wedge
 \beta \leq 5\Theta \wedge \partial \lambda \wedge \bar \partial
 \lambda \wedge \beta +5d(a)(4(n-1)\delta (a)\left\vert{\Theta
 \wedge \Gamma '\wedge \beta }\right\vert +5\left\vert{\Theta
 \wedge \Gamma \wedge \beta }\right\vert ).$ \ \par 
\quad Integrating, we get the lemma because the trace of  $\Theta $
  controls all its coefficients.  $\blacksquare $ \ \par 
\ \par 
\quad \quad 	We shall need the following definition.\ \par 

\begin{Dfnt}
~\label{3_BlasMall0}The domain  $\Omega ,$  equipped with a good
 family  ${\mathcal{Q}},$  will be said {\bf quasi convex} at
  $a\in \Omega $  if, with  $\alpha =\pi (a),\ m=m(\alpha ),$
  taking the coordinates associated to the basis  $b(\alpha ),$
  centered at  $\alpha ,$  we have with  $\rho _{\alpha }$  a
 defining function for  $\Omega ,$ \par 
\quad \quad \quad \quad \quad  $\forall z\in Q_{a}(2)::\rho _{\alpha }(0,z')<0,\ -\rho _{\alpha
 }(0,z')\leq \gamma (\left\vert{z_{2}}\right\vert ^{m_{2}}+\cdot
 \cdot \cdot +\left\vert{z_{n}}\right\vert ^{m_{n}}).$ \par 
The domain will be said {\bf quasi convex} if  $\Omega $  is
 quasi convex at  $a$  for all  $a\in {\mathcal{U}}\cap \Omega
 $  with the same constant  $\gamma .$ 
\end{Dfnt}
\quad \quad 	A convex  $\Omega $  or a lineally convex  $\Omega $  are quasi
 convex because for them  $\displaystyle \Omega \cap T_{\alpha
 }^{{\mathbb{C}}}(\partial \Omega )=\emptyset $  hence  $\rho
 (0,z')\geq 0.$ \ \par 
\ \par 
\quad We have  $\Theta =\sum_{i,j=1}^{n}{\Theta _{ij}\,dz_{i}\wedge
 \partial \bar z_{j}}\ $  and\ \par 
\quad \quad \quad $\displaystyle \Theta \wedge dz_{1}\wedge d\bar z_{1}\wedge \beta
 =\sum_{i,j=2}^{n}{\Theta _{ij}\,dz_{i}\wedge \partial \bar z_{j}\wedge
 \,dz_{1}\wedge \,d\bar z_{1}}\wedge \beta .$ \ \par 
In the integral  $\displaystyle \ \int_{Q_{a}(\delta )}{\Theta
 \wedge dz_{1}\wedge d\bar z_{1}\wedge \beta },$  it remains
 precisely the sum of the  $\sigma _{2n-1}$  areas of the projections
 of  $X_{a}$  on the  $\displaystyle E_{j},\ j\geq 2,$  see~\cite{LelongGrum86},
 Proposition 2.48, p 55. So recall that for  $a\in \Omega ,\
 X_{a}:=X\cap Q_{a}(\delta ),\ X_{a}^{j}$  the projection of
  $\displaystyle X_{a}$  on  $\displaystyle E_{j}$  is denoted
  $\displaystyle X_{a}^{j}$  and  $\displaystyle A_{j}(X_{a}):=\sigma
 _{2n-2}(X_{a}^{j})\ ;$  we get\ \par 
\quad \quad \quad $\displaystyle \ \int_{Q_{a}(\delta )}{\Theta \wedge dz_{1}\wedge
 d\bar z_{1}\wedge \beta }=\sum_{j=2}^{n}{A_{j}(X_{a})}.$ \ \par 
\quad \quad 	So by lemma~\ref{AG0} we have\ \par 
\quad  $\displaystyle \ \sum_{j=2}^{n}{A_{j}(X_{a})}=\ \int_{Q_{a}(\delta
 )}{\Theta \wedge dz_{1}\wedge d\bar z_{1}\wedge \beta }\leq
 \ 5\int_{Q_{a}(\delta )}{\Theta \wedge \partial \lambda \wedge
 \bar \partial \lambda \wedge \beta }+\Gamma d(a)\ \int_{Q_{a}(\delta
 )}{\mathrm{T}\mathrm{r}\Theta },$ \ \par 
hence, because if  $\displaystyle z\in Q_{a}(\delta ),$  by lemma~\ref{Divisor3},
 we have  $\displaystyle d(a)\leq \frac{1}{\delta _{0}-\delta
 }d(z),$  then\ \par 
\quad \quad \quad \quad \quad  $\displaystyle \ \sum_{j=2}^{n}{A_{j}(X_{a})}\leq 5\int_{Q_{a}(\delta
 )}{\Theta \wedge \partial \lambda \wedge \bar \partial \lambda
 \wedge \beta }+\frac{\Gamma }{\delta _{0}-\delta }\ \int_{Q_{a}(\delta
 )}{d(z)\mathrm{T}\mathrm{r}\Theta }.$ \ \par 
\ \par 
\quad \quad 	At this point we shall use equation~\ref{aspcGlo0} which says,
 recall we multiply by  $2$  the defining function  $\rho $ 
 of the domain and the defining function  $\lambda $  of the cone,\ \par 
\quad \quad \quad \quad \quad  $\forall z\in C_{\alpha }\cap B(0,R),\ \rho (0;z')\leq 2A(\left\vert{z_{2}}\right\vert
 ^{m_{2}}+\cdot \cdot \cdot +\left\vert{z_{n}}\right\vert ^{m_{n}}).$ \ \par 
So either  $\rho (0;z')\geq 0,$  then we have\ \par 
\quad \quad \quad \quad \quad  $\displaystyle 0\leq \rho (0;z')\leq 2A(\left\vert{z_{2}}\right\vert
 ^{m_{2}}+\cdot \cdot \cdot +\left\vert{z_{n}}\right\vert ^{m_{n}}),$ \ \par 
or  $\rho (0,z')<0,$  then we use that  $\Omega $  is  $m(\alpha
 )$  quasi convex at  $\alpha $  to get\ \par 
\quad \quad \quad \quad \quad  $-\rho (0,z')\leq \gamma (\left\vert{z_{2}}\right\vert ^{m_{2}}+\cdot
 \cdot \cdot +\left\vert{z_{n}}\right\vert ^{m_{n}}),$ \ \par 
in any case we can apply lemma~\ref{3_BlasMall42} with  $\displaystyle
 z'$  instead of  $z$  to\ \par 
\quad \quad  $h(z'):=\lambda -\rho =2A(\left\vert{z_{2}}\right\vert ^{m_{2}}+\cdot
 \cdot \cdot +\left\vert{z_{n}}\right\vert ^{m_{n}})-\rho (0;z'),$ \ \par 
to get\ \par 
\quad \quad \quad \quad \quad  $\partial h(z')\wedge \bar \partial h(z')=d(a)\Gamma (z'),$ \ \par 
with the sup norm of  $\Gamma $  controlled by the  $\displaystyle
 m_{n}(a)+1$  derivatives of  $\displaystyle h.$ \ \par 
So we have  $\lambda =\rho +h,$  with  $\partial h\wedge \bar
 \partial h=d(a)\Gamma (z'),$  hence\ \par 
\quad \quad \quad \quad \quad  $\Theta \wedge \partial \lambda \wedge \bar \partial \lambda
 \wedge \beta =\Theta \wedge \partial \rho \wedge \bar \partial
 \rho \wedge \beta +\Theta \wedge \partial h\wedge \bar \partial
 h\wedge \beta +\Theta \wedge \partial \rho \wedge \bar \partial
 \lambda \wedge \beta +\Theta \wedge \partial h\wedge \bar \partial
 \rho \wedge \beta ,$ \ \par 
by Cauchy-Schwartz, because  $\Theta \wedge \beta $  is positive, we get\ \par 
\quad \quad \quad \quad \quad  $\displaystyle \ \left\vert{\Theta \wedge \partial h\wedge \bar
 \partial \rho \wedge \beta }\right\vert ^{2}\leq \left\vert{\Theta
 \wedge \partial \rho \wedge \bar \partial \rho \wedge \beta
 }\right\vert \left\vert{\Theta \wedge \partial h\wedge \bar
 \partial h\wedge \beta }\right\vert $ \ \par 
hence, because  $2ab\leq a^{2}+b^{2},$ \ \par 
\quad \quad \quad \quad \quad  $\displaystyle 2\ \left\vert{\Theta \wedge \partial h\wedge
 \bar \partial \rho \wedge \beta }\right\vert \leq \Theta \wedge
 \partial \rho \wedge \bar \partial \rho \wedge \beta +\Theta
 \wedge \partial h\wedge \bar \partial h\wedge \beta $ \ \par 
and\ \par 
\quad \quad \quad \quad \quad  $\displaystyle \Theta \wedge \partial \lambda \wedge \bar \partial
 \lambda \wedge \beta \leq 2\Theta \wedge \partial \rho \wedge
 \bar \partial \rho \wedge \beta +2\Theta \wedge \partial h\wedge
 \bar \partial h\wedge \beta .$ \ \par 
Finally\ \par 
\quad \quad \quad \quad \quad  $\displaystyle \Theta \wedge \partial \lambda \wedge \bar \partial
 \lambda \wedge \beta \leq 2\Theta \wedge \partial \rho \wedge
 \bar \partial \rho \wedge \beta +2d(a)\Theta \wedge \Gamma (z')\wedge
 \beta .$ \ \par 
\ \par 
\quad \quad 	So we have, with  $\displaystyle \Gamma '=\frac{1}{\delta _{0}-\delta
 }\sup _{z\in Q_{a}(\delta )}(\Gamma (z'),\gamma ),$  still controlled
 by the  $\displaystyle m_{n}(a)+1$  derivatives of  $h.$ \ \par 
\quad \quad \quad \quad \quad  $\displaystyle \ \int_{Q_{a}(\delta )}{\Theta \wedge \partial
 \lambda \wedge \bar \partial \lambda \wedge \beta }\leq \int_{Q_{a}(\delta
 )}{\Theta \wedge \partial \rho \wedge \bar \partial \rho \wedge
 \beta }+\Gamma '\int_{Q_{a}(\delta )}{d(z)\mathrm{T}\mathrm{r}\Theta
 }.$ \ \par 
Hence\ \par 
\quad \quad \quad \quad \quad  $\displaystyle \ \sum_{j=2}^{n}{A_{j}(X_{a})}\leq 5\int_{Q_{a}(\delta
 )}{\Theta \wedge \partial \rho \wedge \bar \partial \rho \wedge
 \beta }+5\Gamma '\ \int_{Q_{a}(\delta )}{d(z)\mathrm{T}\mathrm{r}\Theta
 }.$ \ \par 
By use of lemma~\ref{Divisor5} and setting  $S':=S\cap {\mathcal{U}},$
  we have, because the polydiscs  $Q_{a}(\delta ),\ a\in S'$
  are disjoint\ \par 
\quad \quad \quad \quad \quad  \begin{equation}  \ \sum_{a\in S'}{\sum_{j=2}^{n}{A_{j}(X_{a})}}\leq
 C{\left\Vert{\Theta _{X}}\right\Vert}_{B(\Omega )},\label{Divisor9}\end{equation}\
 \par 
where  $C=5C(\rho )\ +5\Gamma '\ .$  	Notice that the constant
  $\displaystyle C$  does not depend on  $\alpha $  and depends
 only on the derivatives of  $\rho $  up to order  $M({\mathcal{Q}})+1,$
  with  $\displaystyle M({\mathcal{Q}})=\sup _{a\in \Omega }m_{n}(a)<\infty
 ,$  because  ${\mathcal{Q}}$  is a good family.\ \par 
Hence we proved the discretized Malliavin condition :\ \par 

\begin{Thrm}
(Discretized Malliavin condition) ~\label{3_BM0}Let  $\Omega
 $  be a domain in  ${\mathbb{C}}^{n}$  equipped with a good
 family  ${\mathcal{Q}}$  of polydiscs and which is  ${\mathcal{Q}}$
  quasi convex. Let  $\Theta $  be a current in the Blaschke
 class and  $S$  a  $\delta $  separated sequence in  $X\cap
 {\mathcal{U}}$  with respect to a good family  ${\mathcal{Q}}$
  of polydiscs with parameter  $\delta _{0}.$  Then we have,\par 
\quad \quad \begin{equation}  \ \sum_{a\in S}{\sum_{j=2}^{n}{A_{j}(X_{a})}}\leq
 C{\left\Vert{\Theta }\right\Vert}_{B},\label{Divisor9}\end{equation}\par 
where  $C$  is a constant depending only on the derivatives of
  $\rho $  up to order  $M({\mathcal{Q}})+1,$  on  $\delta ,\
 \delta _{0}$  and on the constant of quasi convexity.
\end{Thrm}

\subsection{A geometrical lemma.}
\quad Let  $\Omega $  be a domain in  ${\mathbb{C}}^{n}.$  Let  $a\in
 {\mathcal{U}},\ \alpha =\pi (a)$  and  $Q_{a}(\delta )$  the
 polydisc of a good family   ${\mathcal{Q}}$  associated to  $\Omega .$ \ \par 
\quad Let  ${\mathbb{D}}^{n}$  be the unit polydisc in  ${\mathbb{C}}^{n},$
  and let  $\Phi _{a}$  be the bi-holomorphic application from
  ${\mathbb{D}}^{n}$  onto  $Q_{a}(\delta )\ :$ \ \par 
\quad \quad \quad $\displaystyle \forall z=(z_{1},...,z_{n})\in {\mathbb{D}}^{n},\
 1\leq j\leq n,\ Z_{j}=a_{j}+\delta d(a)^{1/m_{j}(a)}z_{j}L_{j}.$ \ \par 
If  $X$  is the zero set of a holomorphic function in  $\Omega
 $  with  $a\in X,$  we can lift  $X_{a}:=X\cap Q_{a}(\delta
 )$  in  ${\mathbb{D}}^{n}$  by  $\Phi _{a}^{-1}.$  Set  $Y_{a}:=\Phi
 _{a}^{-1}(X_{a}),$  and recall that the multi-type is such that
  $m_{1}=1$  and  $m_{n}$  is always bounded,  $m_{n}(a)\leq
 M({\mathcal{Q}}).$  We have :\ \par 

\begin{Lmm}
~\label{Divisor8} (i)  $\displaystyle \sigma _{2n-2}(X_{a})=\sum_{j=1}^{n}{A_{j}(X_{a})}.$
 \par 
\quad (ii)  $\displaystyle \forall j=1,...,n,\ A_{j}(X_{a})=\delta
 ^{2n-2}d(a)^{2\lambda _{j}(a)}A_{j}(Y_{a}),$  with  $\displaystyle
 \lambda _{j}(a)=\sum_{k\neq j}{\frac{1}{m_{k}(a)}}.$ \par 
\quad (iii) with  $m_{1}<m_{2}\leq \cdot \cdot \cdot \leq m_{n},$ \par 
\quad \quad \quad \quad \quad  $\displaystyle \ \frac{n-1}{M({\mathcal{Q}})}\leq \lambda _{1}(a)\leq
 \lambda _{2}(a)\leq \cdot \cdot \cdot \leq \lambda _{n}(a)\leq n/2\ ;$ \par 
\quad \quad \quad \quad \quad  $\displaystyle \lambda _{1}(a)=\sum_{k=2}^{n}{\frac{1}{m_{k}(a)}}\leq
 \frac{n-1}{2}\ ;$ \par 
\quad (iv)  $\displaystyle c_{n}\leq \sigma _{2n-2}(Y_{a}).$ 
\end{Lmm}
\quad Proof.\ \par 
The {\sl (i)} is classical( ~\cite{LelongGrum86}, Proposition
 2.48, p 55).\ \par 
\quad The application  $\Phi _{a}$  sends  $E_{k}=\lbrace z_{k}=0\rbrace
 $  in  $F_{k}:=$  \{ the orthogonal to  $L_{k}$  axis\}  and
 the jacobian of this restriction at the point  $a,$   $J_{k}\Phi
 ,$  is  $J_{k}\Phi =\delta ^{n-1}d(a)^{\lambda _{k}(a)}.$  Because
 the application is holomorphic, we get that the jacobian for
 the change of real variables is\ \par 
\quad \quad \quad $\ \left\vert{J_{k}}\right\vert ^{2}=\delta ^{2n-2}d(a)^{2\lambda
 _{k}(a)},$ \ \par 
which gives the {\sl (ii)}.\ \par 
\quad For the {\sl (iii)} we notice that\ \par 
\quad \quad \quad $\displaystyle 2\leq j,k\leq n,\ m_{k}(a)\geq 2\Rightarrow \frac{1}{m_{k}(a)}\leq
 \frac{1}{2}\Rightarrow \sum_{k\neq j,\ 2\leq k\leq n}{\frac{1}{m_{k}(a)}}\leq
 \frac{n-2}{2}.$ \ \par 
Hence if  $\displaystyle 2\leq j\leq n,\ \lambda _{j}(a)=\sum_{k\neq
 j,2\leq k\leq n}{\frac{1}{m_{k}(a)}}+1\leq n/2\ ;$ \ \par 
if  $\displaystyle j=1,\ \lambda _{1}(a)\leq \frac{n-1}{2}\leq
 n/2.$  Hence {\sl (iii).}\ \par 
\quad The {\sl (iv)} is just the Wirtinger inequality~\cite{Federer69},
 adapted from the ball to the polycube as follows :	 	 $\displaystyle
 Y_{a}\cap B(0,1)\subset Y_{a}$  hence by Wirtinger inequality we get\ \par 
\quad \quad \quad \quad \quad  $\displaystyle c_{n}\leq \sigma _{2n-2}(Y_{a}\cap B(0,1))\leq
 \sigma _{2n-2}(Y_{a}),$ \ \par 
hence the lemma is proved.  $\blacksquare $ \ \par 

\subsection{The result.}

\begin{Thrm}
~\label{3_AspcDomain0}Let  $\Omega $  be a domain in  ${\mathbb{C}}^{n}$
  equipped with a good family  ${\mathcal{Q}}$  of polydiscs
 and which is  ${\mathcal{Q}}$  quasi convex. Let  $S$  be a
  $\delta $  separated sequence of points which is contained
 in the Blaschke divisor  $X.$  Then\par 
\quad \quad \quad $\delta ^{2n-2}\sum_{a\in S}{d(a)^{n}}\leq \gamma (\Omega ){\left\Vert{\Theta
 _{X}}\right\Vert}_{B},$ \par 
where  $\displaystyle \gamma (\Omega )$  depends only on the
 derivatives of  $\rho $  up to order  $M({\mathcal{Q}})+1,$
  on  $\displaystyle \ n$  and  $\delta _{0},$  the parameter
 of the family  ${\mathcal{Q}},$  and on the constant of quasi convexity.
\end{Thrm}
\quad Proof.\ \par 
We have by lemma~\ref{Divisor8}\ \par 
\quad \quad \quad \quad \quad  $\displaystyle \forall j=2,...,n,\ A_{j}(X_{a})=\delta ^{2n-2}d(a)^{2\lambda
 _{j}(a)}A_{j}(Y_{a}),$ \ \par 
but for  $\displaystyle j\geq 2,\ 2\lambda _{j}(a)\leq n,$  hence\ \par 
\quad \quad \quad \quad \quad  $\displaystyle \forall j=2,...,n,\ A_{j}(X_{a})\geq \delta ^{2n-2}d(a)^{n}A_{j}(Y_{a}).$
 \ \par 
For  $\displaystyle j=1$  we have	 	 $\displaystyle 1+2\lambda
 _{1}(a)\leq n,$  by lemma~\ref{Divisor8}, {\sl (iii)}, hence\ \par 
\quad \quad \quad \quad \quad  $\displaystyle A_{1}(X_{a})\geq \delta ^{2n-2}d(a)^{n}A_{1}(Y_{a}).$ \ \par 
Then theorem~\ref{3_BM0} gives\ \par 
\quad \quad \quad \quad \quad  $\displaystyle C{\left\Vert{\Theta _{X}}\right\Vert}_{B}\geq
 \ \sum_{a\in S}{\sum_{j=2}^{n}{A_{j}(X_{a})}}\geq \delta ^{2n-2}\sum_{a\in
 S}{d(a)^{n}\sum_{j=2}^{n}{A_{j}(X_{a})}}.$ \ \par 
And the Blaschke condition gives\ \par 
\quad \quad \quad \quad \quad  $\displaystyle \ \frac{2}{\delta _{0}}{\left\Vert{\Theta }\right\Vert}_{B}\geq
 \sum_{a\in S}{d(a)\mathrm{A}\mathrm{r}\mathrm{e}\mathrm{a}(X_{a})}\geq
 \sum_{a\in S}{d(a)A_{1}(X_{a})}\geq \delta ^{2n-2}\sum_{a\in
 S}{d(a)^{1+2\lambda _{1}(a)}A_{1}(Y_{a})},$ \ \par 
hence\ \par 
\quad \quad \quad \quad \quad  $\displaystyle \ \frac{2}{\delta _{0}}{\left\Vert{\Theta }\right\Vert}_{B}\geq
 \delta ^{2n-2}\sum_{a\in S}{d(a)^{n}A_{1}(Y_{a})}.$ \ \par 
So\ \par 
\quad \quad \quad \quad \quad  $\displaystyle \delta ^{2n-2}\sum_{a\in S}{d(a)^{n}{\left({A_{1}(Y_{a})+\sum_{j\geq
 2}{A_{j}(Y_{a})}}\right)}}\leq (C+\frac{2}{\delta _{0}}){\left\Vert{\Theta
 _{X}}\right\Vert}_{B}.$ \ \par 
Now with  $\displaystyle A_{1}(Y_{a})+\sum_{j\geq 2}{A_{j}(Y_{a})}\geq
 c_{n}$  by Wirtinger inequality, we get the theorem.  $\blacksquare $ \ \par 
\ \par 
\quad We already have defined in the introduction~(\ref{aspcGlo62})
 the canonical measure associated to a sequence  $S$ \ \par 
\quad \quad \quad \quad \quad  $\displaystyle \mu _{S}:=\sum_{a\in S\cap {\mathcal{U}}}{d(a)^{1+2\lambda
 (a)}\delta _{a}}$  with  $\displaystyle \lambda (a):=\sum_{j=2}^{n}{\frac{1}{m_{j}(a)}}\
 .$ \ \par 
\quad The theorem says that the measure  $\displaystyle \ \sum_{a\in
 S\cap {\mathcal{U}}}{d(a)^{n}\delta _{a}}$  is bounded which
 is weaker than the fact the measure  $\displaystyle \mu _{S}$
  is bounded unless  $S$  is a separated sequence projecting
 on points of strict pseudo-convexity, because there we have
  $1+2\lambda (a)=n.$  In the next section we shall introduce
 domains for which we can control the right measure  $\displaystyle
 \mu _{S}.$ \ \par 
\ \par 
\vfill\eject\ \par 
\ \par 

\section{Almost strictly pseudo-convex domains.~\label{4_AspcDomain30}}
\quad \quad 	We shall introduce a family of domains  $\Omega $  with "few"
 non strictly pseudo-convex points on  $\partial \Omega .$ \ \par 
\quad This family will be big enough to contain interesting cases,
 as convex domains of finite type for instance and will allow
 us to manage these "bad" points.\ \par 

\subsection{Minkowski dimension.}

\subsubsection{Definitions and first properties.}

\begin{Lmm}
~\label{4_AspcDomain42}Let  $f$  be a function Lipschitz  $\alpha
 >0,\ \alpha \leq 1,$  on the closed interval   $I=\lbrack 0,\
 h\rbrack $  of  ${\mathbb{R}}.$  Then the graph\par 
\quad \quad \quad \quad \quad  $G:=\lbrace (x,y)::x\in I,\ y=f(x)\rbrace \subset {\mathbb{R}}^{2}$ \par 
of  $f$  can be covered by  $\displaystyle N_{r}(h)\leq Chr^{\alpha
 -2}$  disjoint discs  $\displaystyle D(a,r)$  centered at  $\displaystyle
 a\in G$  and of radius  $\displaystyle r,$  provided that  $\displaystyle
 r\leq h.$ 
\end{Lmm}
\quad Proof.\ \par 
This is Corollary 11.2 p 147 in~\cite{Falconer90}. The proof
 is as follows. Let  $\displaystyle 0<r<1$  and  $m$  the least
 integer greater than or equal to  $\displaystyle h/r.$  We have
 by Proposition 11.1 p. 146 in~\cite{Falconer90} :\ \par 
\quad \quad 	\begin{equation}  r^{-1}\sum_{j=0}^{m-1}{R_{f}(jr,(j+1)r)}\leq
 N_{r}(h)\leq 2m+r^{-1}\sum_{j=0}^{m-1}{R_{f}(jr,(j+1)r)},\label{4_AspcDomain41}\end{equation}\
 \par 
with  $\displaystyle R_{f}(t_{1},t_{2}):=\sup _{t_{1}<t,u<t_{2}}\left\vert{f(t)-f(u)}\right\vert
 .$ \ \par 
Because  $f$  is Lipschitz  $\alpha $  we have  $\displaystyle
 R_{f}(t_{1},t_{2})\leq C\left\vert{t_{1}-t_{2}}\right\vert ^{\alpha
 }$  hence  $\displaystyle R_{f}(jr,(j+1)r)\leq Cr^{\alpha }.$
  Putting this in~(\ref{4_AspcDomain41}) we get\ \par 
\quad \quad \quad \quad \quad  $\displaystyle N_{r}(h)\leq 2m+mCr^{\alpha -1}.$ \ \par 
But provided that  $\displaystyle m>0,$  i.e.  $\displaystyle
 h\geq r,$  we have  $\displaystyle m\leq 2h/r$  so\ \par 
\quad \quad \quad \quad \quad  $\displaystyle N_{r}(h)\leq 4\frac{h}{r}+2Chr^{\alpha -2}\leq
 C'hr^{\alpha -2}.$   $\blacksquare $ \ \par 
\ \par 
\quad \quad 	We shall define an homogeneous Minkowski dimension. Denote 
 $\displaystyle \# A$  the number of points in the set  $\displaystyle
 A.$ \ \par 

\begin{Dfnt}
Let  $W\subset {\mathbb{R}}^{2}$  be a bounded set and  $\displaystyle
 \alpha >0\ ;$  let  $\displaystyle D(a,h)$  be a disc centered
 at  $a$  and of radius  $h$  ; let  ${\mathcal{R}}_{r}(W\cap
 D(a,h))$  be a covering of  $\displaystyle W\cap D(a,h)$  by
 discs of radius  $r\ ;$  we shall say that  $W$  has {\bf homogeneous
 Minkowski} dimension  $\alpha $  if :\par 
\quad \quad  $\displaystyle \ \exists C>0,\ \forall a\in W,\ \forall h>0,\
 \forall r>0,\ r\leq h,\ \exists {\mathcal{R}}_{r}(W\cap D(a,h))::\#
 {\mathcal{R}}_{r}(W\cap D(a,h))\leq \max (1,Chr^{-\alpha }).$ 
\end{Dfnt}
The number  $C$  will be called the {\sl constant} of  $W$  with
 respect to the homogeneous Minkowski dimension  $\alpha .$ \ \par 
\quad \quad 	Clearly if  $W$  has homogeneous Minkowski dimension  $\alpha
 $  with constant  $\displaystyle C$  then it has upper Minkowski
 dimension  $\alpha ,$  see~\cite{Falconer90}, but the converse
 is false as can be seen with the canonical example of  $\displaystyle
 W=\lbrace 0,1,1/2,1/3,...,1/n,...\rbrace $  which has Minkowski
 dimension  $\displaystyle 1/2$  but it is not homogeneous, i.e.
  $W\cap D(0,h)$  has no covering with the property above for
 any  $\displaystyle h>0.$ \ \par 
\quad \quad 	On the other hand, lemma~\ref{4_AspcDomain42} gives examples
 of such sets.\ \par 
\quad The following result is a corollary of a nice theorem of Ostrowski~\cite{Ostrowski40}.\
 \par 

\begin{Crll}
~\label{strongPC13}Let  $P(y)$  be a monic polynomial of degree
  $d$  in the real variable  $y$  whose coefficients are  ${\mathcal{C}}^{\infty
 }$  functions of  $x\in {\mathbb{R}}.$  Then the graph of the
 zero set of  $P$  has homogeneous Minkowski dimension less than
  $\displaystyle 2-\frac{1}{d}.$ 
\end{Crll}
\quad Proof.\ \par 
by a theorem of Ostrowski~\cite{Ostrowski40} we have that locally
 the roots  $y$  of the equation\ \par 
\quad \quad \quad $\displaystyle P(y)=y^{d}+a_{1}y^{d-1}+\cdot \cdot \cdot +a_{d}=0,$ \ \par 
are Lipschitz  $\displaystyle \ \frac{1}{d}$  functions of the
 coefficients  $a_{j}.$  Composing with the  ${\mathcal{C}}^{\infty
 }$  function\ \par 
\quad \quad \quad \quad $x\rightarrow a(x):=\lbrace a_{j}(x),\ j=1,...,d\rbrace ,$ \ \par 
we get that the roots  $y_{k}(x),\ k\leq d,$  are still Lipschitz
  $\displaystyle \ \frac{1}{d}$  and we can apply Lemma~\ref{4_AspcDomain42}
 to the graph of each root. Because there is at most  $d$  such
 graphs, the corollary is proved.  $\blacksquare $ \ \par 

\subsubsection{Domains in  ${\mathbb{C}}^{n}$ }
\quad Let  $D(\rho )$  be the disc in  ${\mathbb{C}}$  of center  $0$
  and radius  $\rho $  and denote  $\sigma _{2n}$  the Lebesgue
 measure in  ${\mathbb{C}}^{n}={\mathbb{R}}^{2n}.$  We have the lemma:\ \par 

\begin{Lmm}
~\label{suiteFaible63}Let  $W\subset {\mathcal{D}}:=D(d){\times}D(R)^{n-2}{\times}D(h)\subset
 {\mathbb{C}}^{n}$  and  $\alpha >0$  such that the homogeneous
 Minkowski dimension of\par 
\quad \quad \quad \quad $W\cap \lbrace z_{1}=a_{1},\ ...,\ z_{n-1}=a_{n-1}\rbrace $ \par 
is  $2-\alpha $  for all  $a'=(a_{1},...,\ a_{n-1})\in D(d){\times}D(R)^{n-2}.$
  Let  $S\subset W$  and let another orthonormal basis  $\displaystyle
 b=\lbrace L_{1},..,L_{n}\rbrace $  with  $\displaystyle w=(w_{1},...,w_{n})$
  as coordinates ;  let  $P_{a}$  be a polydisc with respect
 to the basis  $b(a)$  centered on  $a\in S,\ b(a)$  varying
 with  $a\in S,\ P_{a}$  with fixed radii  $(r,\ l_{2}r...,\
 l_{n}r),$  and such that these polydiscs are disjoint. Let 
 $l=\max  _{j=2,...,n}l_{j}.$  Then\par 
\quad \quad $\exists C::\sum_{a\in S}{\sigma _{2n}(P_{a})}\leq Chd^{2}R^{2(n-2)}l^{\alpha
 }r^{\alpha }=Ch^{-1}\sigma _{2n}({\mathcal{D}})l^{\alpha }r^{\alpha }.$ 
\end{Lmm}
\quad Proof.\ \par 
Denote by  $C$  the canonical basis of  ${\mathbb{C}}^{n}$  with
 the  $z_{j}$  as coordinates.\ \par 
\quad \quad 	First set  $C(b,\ lr)$  a polycube, i.e. a polydisc with all
 its radii are equal,  with respect to the canonical basis  $C$
  in  ${\mathcal{D}},$  centered at  $b$  and of radii  $(lr,...,\
 lr).$  Any polydisc  $P_{a}$  with  $a$  in  $C(b,\ lr)$  is
 contained in the "double" polycube  $C(b,\ 2(2)^{n}lr),$  the
  $\displaystyle 2^{n}$  because of the "angle" between the two
 bases ; hence the measure of the union of all those polydiscs
  $P_{a}$  is bounded by the measure of  $C(b,\ 2^{n+1}lr).$
  These polydiscs being disjoint we get\ \par 
\quad \quad \quad $\displaystyle \ \sum_{a\in S\cap C(b,\ lr)}{\sigma _{2n}(P_{a})}\leq
 \sigma _{2n}(C(b,\ 2^{n+1}lr))=2^{2n+1}\pi ^{n}l^{2n}r^{2n}.$ \ \par 
Each polydisc verifies  $\displaystyle \sigma _{2n}(P_{a})=\pi
 ^{n}l_{2}^{2}\cdot \cdot \cdot l_{n}^{2}r^{2n},$  hence the
 number of points  $N_{C}$  of  $S$  in  $C(b,\ lr)$  can be
 estimated by :\ \par 
\quad \quad \quad \quad \quad  $\displaystyle N_{C}\leq 2^{2n+1}\pi ^{n}l^{2n}r^{2n}/\pi ^{n}l_{2}^{2}\cdot
 \cdot \cdot l_{n}^{2}r^{2n}=2^{2n+1}\frac{l^{2n}}{l_{2}^{2}\cdot
 \cdot \cdot l_{n}^{2}}.$ \ \par 
Let  $b'=(b_{1},...,b_{n-1})$  be fixed, then the set  $\displaystyle
 C((b',\ b_{n}),\ lr)\cap \lbrace z'=b'\rbrace \subset D(h)$
  is a disc centered at  $\displaystyle b_{n}\in D(h)$  and of
 radius  $\displaystyle lr.$  The homogeneous Minkowski assumption
 gives that there is a subfamily of these discs which covers
  $S$  whose number  $n_{B}$  of elements verifies\ \par 
\quad \quad \quad $\displaystyle n_{B}\leq Ch(lr)^{\alpha -2}.$ \ \par 
Define the slice of depth  $\displaystyle lr$  to be  $\displaystyle
 B(b',lr):=\bigcup_{b_{n}\in D(h)}{C((b',b_{n}),lr)},$  then
 the number  $N_{B}$  of points of  $S$  in this slice verifies\ \par 
\quad \quad \quad $\displaystyle N_{B}\leq n_{B}{\times}N_{C}\leq Ch(lr)^{\alpha
 -2}{\times}2^{2n+1}\frac{l^{2n}}{l_{2}^{2}\cdot \cdot \cdot
 l_{n}^{2}}.$ \ \par 
The number of such slices, when  $b'$  varies, is bounded by
  $\displaystyle \ \frac{d^{2}R^{2(n-2)}}{(lr)^{2(n-1)}},$  hence
 the total number  $N$  of points in  $S$  can be estimated by :\ \par 
\quad \quad \quad $\displaystyle N\leq \frac{N_{B}d^{2}R^{2(n-2)}}{l^{2(n-1)}r^{2(n-1)}}\leq
 2^{2n+1}d^{2}R^{2(n-2)}Chl^{\alpha }r^{\alpha }\frac{1}{l_{2}^{2}\cdot
 \cdot \cdot l_{n}^{2}r^{2n}}.$ \ \par 
Hence the total measure of the polydiscs  $P_{a}$  is\ \par 
\quad $A:=\sum_{a\in S}{\sigma _{2n}(P_{a})}=N{\times}\pi ^{n}l_{2^{2}}\cdot
 \cdot \cdot l_{n^{2}}r^{2n}\leq 2^{2n+1}\pi ^{n}d^{2}R^{2(n-2)}Chl^{\alpha
 }r^{\alpha }=C'h^{-1}\sigma _{2n}({\mathcal{D}})l^{\alpha }r^{\alpha
 },$ \ \par 
with  $\displaystyle C'=2^{2n+1}\pi ^{n}C$  which depends only
 on  $C,$  the Minkowski constant of  $\displaystyle W.$   $\blacksquare
 $ \ \par 

\subsection{Domains almost strictly pseudo-convex.}
\quad Let  $W$  be the set of weakly pseudo-convex points of  $\partial
 \Omega ,$  i.e.  $W$  is the zero set of the determinant of
 the Levi form  ${\mathcal{L}}$  of  $\partial \Omega .$  Let
  $\pi $  be the normal projection from  $\Omega $  onto  $\partial
 \Omega ,$  defined in a neighbourhood  ${\mathcal{U}}$  of 
 $\partial \Omega $  in  $\Omega .$ \ \par 
\quad Let  $\alpha \in \partial \Omega \ ;$  by linear change of variables
 we can suppose that  $\alpha =0\in \partial \Omega \subset {\mathbb{C}}^{n},\
 z_{1}=0$  is the equation of the complex tangent space. The
 projection  $\pi $  locally near  $0\in \partial \Omega $  gives
 a  ${\mathcal{C}}^{\infty }$  diffeomorphism  $\displaystyle
 \tilde \pi \ :\ \partial \Omega \rightarrow T_{0}(\partial \Omega
 ),\ \tilde \pi :=(\pi _{\mid T_{0}(\partial \Omega )})^{-1}.$ \ \par 

\begin{Dfnt}
The pseudo-convex domain  $\Omega $  in  ${\mathbb{C}}^{n}$ 
 is said to be {\bf almost stricly pseudo-convex}, {\bf aspc},
 at  $0$  if there is a neighbourhood  $V_{0}$  of  $0,$  a positive
 number  $\displaystyle \beta ,$  and a basis  $\displaystyle
 b:=\lbrace L_{1},...,L_{n}\rbrace $  of  ${\mathbb{C}}^{n},$
  still with  $\displaystyle L_{1}$  a complex normal unit vector,
 such that the slices in the associated coordinates for the basis   $b,$ \par 
\quad \quad \quad $\tilde \pi (W\cap V_{0})\cap \lbrace z_{1}=0\rbrace \cap \lbrace
 z_{2}=a_{2}\rbrace \cap \cdot \cdot \cdot \cap \lbrace z_{n-1}=a_{n-1}\rbrace
 $ \par 
have homogeneous Minkowki dimension less than  $2-\beta ,\ \beta >0.$ \par 
\quad $\Omega $  is said to be {\bf aspc} if this is true for all points
 in  $\partial \Omega $  with the same  $\displaystyle \beta
 >0$  and the same underlying constant.
\end{Dfnt}
The basis  $b$  is in general different from the basis  $\displaystyle
 b(\alpha )$  used in the definition of the good family  ${\mathcal{Q}}.$
 \ \par 
\quad Of course the strictly pseudo-convex domains are {\bf aspc} because
  $W=\emptyset .$  \ \par 

\subsection{Sequences projecting on weak pseudo-convex points.}
\quad We still shall use the notations :\ \par 
\quad \quad \quad $\displaystyle \forall a\in {\mathcal{U}},\ \alpha :=\pi (a),\
 \ m(a):=m(\alpha )=(m_{1}(\alpha ),\ ...,\ m_{n}(\alpha ))$
  is the multi-type of a point ;\ \par 
\quad \quad \quad \quad \quad  $W$  is the set of non strictly pseudo-convex points on  $\partial
 \Omega \ ;$ \ \par 
\quad \quad \quad \quad $\displaystyle \forall a\in {\mathcal{U}},\ \lambda (a):=\sum_{j=2}^{n}{\frac{1}{m_{j}(a)}}$
  is the weight exponent.\ \par 

\begin{Thrm}
~\label{aspc0}Let  ${\mathcal{Q}}$  be a good family of polydiscs
 on a {\bf aspc} domain  $\Omega $  in  ${\mathbb{C}}^{n},$ 
 and  $S$  be a  $\delta $  separated sequence of points in 
 $\Omega .$  If  $\pi (S\cap {\mathcal{U}})\subset V\cap W,$
  where  $V$  is an open set of  $\partial \Omega ,$  then we have:\par 
\quad \quad \quad \begin{equation}  \ \sum_{a\in S\cap {\mathcal{U}}}{d(a)^{1+2\lambda
 (a)}}=\delta ^{-2n}\sum_{a\in S\cap {\mathcal{U}}}{\sigma _{2n}(Q_{a}(\delta
 ))}\leq C(\Omega )\sigma _{2n-1}(V),\label{aspcGlo65}\end{equation}\par 
where  $C(\Omega )$  depends only on  $\rho ,\ n,$  the good
 family  ${\mathcal{Q}}$  and the constant  $\beta $  in the
 Minkowski dimension of  $W\subset \partial \Omega .$ 
\end{Thrm}
\quad Proof.\ \par 
The polydisc  $Q_{a}(\delta )$  has radius  $\gamma :=\delta
 d(a)$  in the normal direction and in its conjugate and has radii\ \par 
\quad \quad \quad \quad \quad  $(\delta d(a)^{1/m_{2}(a)},...,\delta d(a)^{1/m_{n}(a)})$ \ \par 
in the complex tangent directions. Let us denote  $L_{2},...,\
 L_{n}$  the complex tangent directions in the basis  $\displaystyle
 b(\alpha )$  associated to  $\pi (a)$  with multi-type  $(m_{2}(a),...,m_{n}(a)).$
 \ \par 
\quad Now fix  $\zeta \in W\subset \partial \Omega $  and let  $V_{\zeta
 }:=B(\zeta ,\epsilon )\cap \partial \Omega $  be a neighbourhood
 of  $\zeta $  in  $\partial \Omega $  such that  $\tilde \pi
 $  is a diffeomorphism from  $V_{\zeta }$  on a neighbourhood
 of  $\zeta $  on the (real) tangent space  $T_{\zeta }.$  One
 can choose the radius of the euclidean ball  $B(\zeta ,\epsilon
 ),\ \epsilon >0$  to be fixed independently of  $\zeta ,$  because
  $\partial \Omega $  is of class  ${\mathcal{C}}^{2}$  and compact.\ \par 
\quad Because  $\Omega $  is {\bf aspc}, we know that there is a basis
  $b=\lbrace v_{1},...,v_{n}\rbrace $  of  ${\mathbb{C}}^{n}$
  such that  $\displaystyle v_{1}=L_{1}$  is still in the complex
 normal space, and a complex direction in the complex tangent
 space at  $\zeta ,$  say  $v_{n},$  along which  $W$  is of
 homogeneous Minkowski dimension  $2-\beta ,\ \beta >0.$  I.e.
 these two basis are different in the complex tangent space only.\ \par 
\ \par 
\quad Let  $S\subset \Omega ::\pi (S)\subset W$  be the  $\delta $
  separated given sequence. First we shall prove the theorem
 with  $V=V_{\zeta }$  and then complete it.\ \par 
\quad The proof will follow from several reductions.\ \par 

\subsubsection{Reduction to a layer parallel to the complex tangent space.}
\quad \quad 	As usual we suppose that  $\displaystyle \zeta =0,\ \Re z_{1}=0$
  is the tangent space  $T_{0}(\partial \Omega ).$ \ \par 
\quad By use of the  ${\mathcal{C}}^{\infty }$  diffeomorphism  $\tilde
 \pi ,$  we can suppose that  $\partial \Omega \simeq T_{0}(\partial
 \Omega )$  in a ball  $\displaystyle B(0,\epsilon )$  with a
 uniform  $\displaystyle \epsilon >0$  which depends only on
  $\Omega $  via its defining function  $\displaystyle \rho .$ \ \par 
Consider the polydisc, {\sl in the basis}  $\displaystyle b,\
 P_{0}(R,h,d):=D(d){\times}D(R)^{n-2}{\times}D(h)\subset B(0,\epsilon
 )$  where  $\displaystyle D(r)$  is a disc centered at  $0$
  and of radius  $\displaystyle r.$  We can manage it to have
  $\displaystyle \epsilon /2{\sqrt{n}}\leq d\leq h$  and still
  $\displaystyle P_{0}(R,h,d)\subset B(0,\epsilon ).$ \ \par 
\quad \quad 	In this ball  $\displaystyle B(0,\epsilon )$  we consider  $\Omega
 $  as a half space  $T_{0}(\partial \Omega ){\times}\rbrack
 0,\epsilon \rbrack $  by use of the diffeomorphism  $\tilde \pi .$ \ \par 
\quad \quad 	From now on we shall restrict everything to  $\displaystyle
 P_{0}(R,h,d),$  which means, in particular, that  $\displaystyle
 z\in P_{0}(R,h,d)\Rightarrow \left\vert{z_{1}}\right\vert \leq d.$ \ \par 
\quad Let  $C_{\gamma }\subset P_{0}(R,h,d)$  be a layer parallel to
  $T_{0}(\partial \Omega )$  at a distance  $\gamma \leq d$ 
 from the boundary, i.e.\ \par 
\quad \quad \quad \quad $a=(a_{1},...,\ a_{n})\in C_{\gamma }\iff \Re a_{1}\simeq d(a)\in
 \lbrack (1-\delta )\gamma ,\ (1+\delta )\gamma \rbrack ,$  with
  $\delta $  the separating constant.\ \par 
Now let  $S_{\gamma }:=S\cap C_{\gamma }\cap P_{0}(R,h,d).$ \ \par 

\subsubsection{Reduction to a fixed multi-type.}
\quad There is only a finite set of possible multi-types for the points
 of  $S$  because we have a good family of polydiscs and the
 multi-type is uniformly bounded by lemma~\ref{2_BonFamille20}.
 Hence it is enough to show the inequality~(\ref{aspcGlo65})
 for the points  $a\in S$  with a fixed multi-type,  $m(a)=(1,\
 m_{2},...,\ m_{_{n}}).$  Of course the axes of the polydisc
  $Q_{a}(\delta )$  are still {\sl dependent of } $a.$ \ \par 
\quad \quad 	We can apply Lemma~\ref{suiteFaible63} to the sequence  $S_{\gamma
 }\ ;$  because  $m_{2}\leq \cdot \cdot \cdot \leq m_{n},$  we set :\ \par 
\quad \quad \quad $\displaystyle r:=\gamma ^{1/m_{2}},\ l:=\gamma ^{\frac{1}{m_{n}}-\frac{1}{m_{2}}}.$
 \ \par 
The lemma gives :\ \par 
\quad \quad \quad \quad \quad  $\displaystyle \ \sum_{a\in S_{\gamma }}{\sigma _{2n}(Q_{a})}\leq
 ChR^{2(n-2)}d^{2}l^{\beta }r^{\beta }=ChR^{2(n-2)}d^{2}\gamma
 ^{\beta /m_{n}}.$ \ \par 
The measure of the trace of  $\displaystyle P_{0}(R,h,d)$  on
 the {\sl real} tangent space  $\displaystyle T_{0}(\partial
 \Omega )$  is\ \par 
\quad \quad \quad \quad \quad  $\displaystyle \sigma _{2n-1}(P_{0}(R,h,d)\cap T_{0}(\partial
 \Omega ))=R^{2(n-2)}h^{2}d,$ \ \par 
because the disc  $\displaystyle D(d)$  is in the complex normal.\ \par 
\quad \quad 	So we get\ \par 
\quad \quad \quad \quad 	\begin{equation}  \ \sum_{a\in S_{\gamma }}{\sigma _{2n}(Q_{a})}\leq
 ChR^{2(n-2)}d^{2}\gamma ^{\beta /m_{n}}=C\frac{d}{h}\gamma ^{\beta
 /m_{n}}\sigma _{2n-1}(P_{0}(R,h,d)\cap T_{0}(\partial \Omega
 )).\label{M0}\end{equation}\ \par 

\subsubsection{Adding the layers.}
\quad Because the sequence is separated, the layers can be ordered
 this way  $\gamma _{k}=\mu ^{k}\gamma _{0},\ k\in {\mathbb{N}},$
  where  $\gamma _{0}\leq d$  is the farthest point from the
 boundary and  $\displaystyle \mu =\frac{1-\delta }{1+\delta }<1.$ \ \par 
We have to add them and, because of inequality~(\ref{M0}), we get\ \par 
\quad \quad \quad \quad \quad  $\displaystyle \ \sum_{k\in {\mathbb{N}}}{\sum_{a\in S_{\gamma
 _{k}}}{\sigma _{2n}(Q_{a})}}\leq C\frac{d}{h}\sigma _{2n-1}(P_{0}(R,h,d)\cap
 T_{0}(\partial \Omega ))\sum_{k\in {\mathbb{N}}}{\gamma _{k}^{\beta
 /m_{n}}}.$ \ \par 
But  $\displaystyle \gamma _{k}=\mu ^{k}\gamma _{0},\ k\in {\mathbb{N}},$
  so\ \par 
\quad \quad \quad \quad \quad  $\displaystyle \ \sum_{k\in {\mathbb{N}}}{\gamma _{k}^{\beta
 /m_{n}}}=\gamma _{0}^{\beta /m_{n}}\sum_{k\in {\mathbb{N}}}{\mu
 ^{k\beta /m_{n}}}=\frac{\gamma _{0}^{\beta /m_{n}}}{1-\mu ^{\beta
 /m_{n}}}\leq \frac{d^{\beta /m_{n}}}{1-\mu ^{\beta /m_{n}}}.$ \ \par 
Hence we get\ \par 
\quad \quad \quad \quad 	\begin{equation}  \ \sum_{k\in {\mathbb{N}}}{\sum_{a\in S_{\gamma
 _{k}}}{\sigma _{2n}(Q_{a})}}\leq C'\frac{d^{1+\beta /m_{n}}}{h}\sigma
 _{2n-1}(P_{0}(R,h,d)\cap T_{0}(\partial \Omega )),\label{M1}\end{equation}\
 \par 
with  $\displaystyle C':=C\frac{1}{1-\mu ^{\beta /m_{n}}}.$ \ \par 

\subsubsection{Adding for all the multi-types.}
\quad \quad 	Because we have a good family of polydiscs the multi type is
 bounded, hence  $\forall a\in \Omega ,\ m_{n}(a)\leq M({\mathcal{Q}}),$
  so we have that, for any multi type,\ \par 
\quad \quad \quad \quad \quad  $\displaystyle C':=C\frac{1}{1-\mu ^{\beta /m_{n}}}\leq C\frac{1}{1-\mu
 ^{\beta /M({\mathcal{Q}})}}=:D,$ \ \par 
hence the inequality~(\ref{M1}) implies\ \par 
\quad \quad \quad \quad 	\begin{equation}  \ \sum_{k\in {\mathbb{N}}}{\sum_{a\in S_{\gamma
 _{k}}}{\sigma _{2n}(Q_{a})}}\leq D\frac{d^{1+\beta /m_{n}}}{h}\sigma
 _{2n-1}(P_{0}(R,h,d)\cap T_{0}(\partial \Omega )).\label{M2}\end{equation}\
 \par 
Because  $\displaystyle d\leq h$  we have  $\displaystyle \ \frac{d^{1+\beta
 /m_{n}}}{h}\leq d^{\beta /m_{n}}.$ \ \par 
\quad \quad 	Recall that  $\displaystyle \sigma _{2n-1}(Q_{a}(\delta ))=\delta
 ^{-2n}d(a)^{1+2\lambda (a)},$  then we get\ \par 
\quad \quad \quad \quad \quad  $\displaystyle \ \sum_{k\in {\mathbb{N}}}{\sum_{a\in S_{\gamma
 _{k}}}{d(a)^{1+2\lambda (a)}}}\leq D2^{n}\delta ^{-2n}d^{\beta
 /m_{n}}\sigma _{2n-1}(P_{0}(R,h,d)\cap T_{0}(\partial \Omega )).$ \ \par 
\quad \quad 	Now set  $\displaystyle V_{0}:=P_{0}(R,h,d)\cap T_{0}(\partial
 \Omega ),\ d\leq h\ ;$  the number of possible multi types being
 finite, we have a finite sum of finite numbers so  $\displaystyle
 \ \sum_{a\in S}{d(a)^{1+2\lambda (a)}}$  is finite, for  $S\cap
 {\mathcal{U}}\cap \lbrace d(a)\leq d\rbrace \subset \pi ^{-1}(V_{0}),$
  with constant  $\displaystyle \ C(\Omega )\sigma _{2n-1}(V_{0}),$
  where  $C(\Omega )$  depends only on the defining function
  $\rho $  of  $\Omega ,$  the Minkowski constants of  $W$  and
 of the good family  ${\mathcal{Q}}.$ \ \par 
\ \par 
\quad Now let  $V$  be an open set in  $\partial \Omega \ ;$  because
  $\partial \Omega $  is a bounded smooth manifold in  ${\mathbb{R}}^{2n}$
  we can cover it by a finite number of sets  $\lbrace V_{\zeta
 }\rbrace _{\zeta \in {\mathcal{R}}}$  "almost" disjoint, i.e. such that\ \par 
\quad $\bullet $  the union  $\displaystyle \ \bigcup_{\zeta \in {\mathcal{R}}}{V_{\zeta
 }}$  covers  $\partial \Omega .$ \ \par 
\quad $\bullet $  any point of  $\partial \Omega $  belongs to at most
  $N$  of the  $\displaystyle V_{\zeta }.$ \ \par 
This gives\ \par 
\quad \quad \quad $\displaystyle V\subset \bigcup_{\zeta \in {\mathcal{R}}}{V_{\zeta
 }\cap V}.$ \ \par 
Hence\ \par 
\quad \quad \quad $\displaystyle \sigma _{2n-1}(V)\leq \sum_{\zeta \in {\mathcal{R}}}{\sigma
 _{2n-1}(V_{\zeta }\cap V)}.$ \ \par 
On the other hand we just proved, shrinking  ${\mathcal{U}}$
  to  ${\mathcal{U}}\cap \lbrace d(a)\leq d\rbrace $  if necessary,\ \par 
\quad \quad \quad $\displaystyle \ \sum_{a\in S\cap {\mathcal{U}}\cap \pi ^{-1}(V_{\zeta
 }\cap V)}{d(a)^{1+2\lambda (a)}}\leq C(\Omega )\sigma _{2n-1}(V_{\zeta
 }),$ \ \par 
so\ \par 
\quad $\displaystyle \ \sum_{a\in S\cap {\mathcal{U}}\cap \pi ^{-1}(V)}{d(a)^{1+2\lambda
 (a)}}\leq \sum_{\zeta \in {\mathcal{R}}}{\sum_{a\in S\cap \pi
 ^{-1}(V_{\zeta }\cap V)}{d(a)^{1+2\lambda (a)}}}\leq C(\Omega
 )\sum_{\zeta \in {\mathcal{R}}}{\sigma _{2n-1}(V_{\zeta })}\leq $ \ \par 
\quad \quad \quad \quad \quad \quad \quad \quad \quad \quad \quad \quad \quad \quad  $\displaystyle \leq CN{\times}\sigma _{2n-1}(V),$ \ \par 
the last inequality because any point of  $V$  is covered at
 most  $N$  times.  $\blacksquare $ \ \par 

\begin{Rmrq}
In fact this theorem says that the measure\par 
\quad \quad \quad \quad \quad  $\displaystyle \lambda _{S}:=\sum_{a\in S\cap {\mathcal{U}}}{d(a)^{1+2\lambda
 (a)}\delta _{a}}$ \par 
associated to the a separated sequence  $S$  of points projecting
 on the weakly pseudo-convex points in  $\partial \Omega $  is
 a geometric Carleson measure, as we shall see later.
\end{Rmrq}

\subsection{Sequence of points in a Blaschke divisor.}
\quad We shall glue the previous result with the one we got in theorem~\ref{3_AspcDomain0}
 to have the control of the canonical measure  $\displaystyle
 \mu _{S}$  associated to a separated sequence  $\displaystyle S.$ \ \par 
\ \par 

\begin{Thrm}
Let  $\Omega $  be a {\bf aspc} domain in  ${\mathbb{C}}^{n}$
  equipped with a good family  ${\mathcal{Q}}$  of polydiscs
 and which is  ${\mathcal{Q}}$  quasi convex. Let  $S$  a  $\delta
 $  separated sequence of points contained in a divisor  $X$
  of the Blaschke class of  $\Omega ,$  with  $\Theta $  as its
 current of integration, which projects on the open set  ${\mathcal{V}}\subset
 \partial \Omega .$  Then we have \par 
\quad \quad \quad \quad \quad  $\displaystyle \ \sum_{a\in S}{d(a)^{1+2\lambda (a)}}\leq \gamma
 (\Omega ){\left\Vert{\Theta }\right\Vert}_{B}+C(\Omega )\sigma
 _{2n-1}({\mathcal{V}})<\infty ,$ \par 
where  $d(a)$  is the distance from  $a$  to the boundary of
  $\Omega $  and  $\displaystyle \lambda (a):=\sum_{j=2}^{n}{\frac{1}{m_{j}(a)}},$
  with  $(1,\ m_{2}(a),...,m_{n}(a))$  is the multi-type associated
 to the family  ${\mathcal{Q}}.$ \par 
Moreover the constants  $C(\Omega ),\ \gamma (\Omega ),$  depend
 only on the  ${\mathcal{C}}^{2}$  norm of the defining function
  $\rho ,\ n,\ \delta $  and  $\delta _{0}$  the parameter of
 the good family  ${\mathcal{Q}},$  the Minkowski constants of
 the {\bf aspc} domain  $\Omega $  and the constant of quasi convexity.
\end{Thrm}
\quad Proof.\ \par 
Let  $B_{S}$  be the set of (bad) points in the sequence  $S,$
  i.e. which project on the weakly pseudo-convex points in  ${\mathcal{V}}\subset
 \partial \Omega \ ;$ \ \par 
let  $G_{S}$  be the set of (good) points in the sequence  $S,$
  i.e. which project on the strictly pseudo-convex points in
  ${\mathcal{V}}\subset \partial \Omega \ ;$ \ \par 
then  $S=B_{S}\cup G_{S}$  and we have by theorem~\ref{3_AspcDomain0}\ \par 
\quad \quad \quad $\ \sum_{a\in G_{S}}{d(a)^{n}}\leq \sum_{a\in S}{d(a)^{n}}\lesssim
 {\left\Vert{\Theta }\right\Vert}_{B}\Rightarrow \sum_{a\in G_{S}}{d(a)^{1+2\lambda
 (a)}}\leq \gamma (\Omega ){\left\Vert{\Theta }\right\Vert}_{B},$ \ \par 
because\!\!\!\! , for these points we have  $m_{1}=1,\ m_{2}=2,...,\
 m_{n}=2,$  hence  $n=1+2\lambda (a).$ \ \par 
By theorem~\ref{aspc0} we have\ \par 
\quad \quad \quad $\ \sum_{a\in B_{S}}{d(a)^{1+2\lambda (a)}}\leq C(\Omega )\sigma
 _{2n-1}({\mathcal{V}})<\infty ,$ \ \par 
so adding these two inequalities, we get\ \par 
\quad \quad \quad $\ \sum_{a\in S}{d(a)^{1+2\lambda (a)}}\leq \gamma (\Omega ){\left\Vert{\Theta
 _{X}}\right\Vert}_{B}+C(\Omega )\sigma _{2n-1}({\mathcal{V}})<\infty
 .$   $\blacksquare $ \ \par 
\vfill\eject\ \par 

\section{Examples of almost strongly pseudo-convex domains.~\label{5_examAspc30}}

      	An example of {\bf aspc} domain {\sl not} of finite type
 is the following\ \par 
\quad \quad \quad \quad \quad  $\displaystyle \ \left\vert{z_{1}}\right\vert ^{2}+\exp (1-\left\vert{z_{2}}\right\vert
 ^{-2})<1,$ \ \par 
because the set  $W$  of its weakly pseudo-convex points is the
 circle  $\displaystyle \ \left\vert{z_{1}}\right\vert =1,\ z_{2}=0,$
  hence it has Minkowski dimension  $1.$ \ \par 
\quad The other examples are mainly based on the following theorem.\ \par 
\quad Let  $\Omega $  be a domain in  ${\mathbb{C}}^{n}$  and  ${\mathcal{L}}$
  its Levi form. Set  ${\mathcal{D}}:=det{\mathcal{L}}\ ;$  then
 the set  $W$  of points of weak pseudo-convexity is  $W:=\lbrace
 z\in \partial \Omega ::{\mathcal{D}}(z)=0\rbrace .$ \ \par 

\begin{Thrm}
~\label{aspc3}Let  $\Omega $  be a domain in  ${\mathbb{C}}^{n}$
  of finite linear type, and  ${\mathcal{D}}$  the determinant
 of its Levi form. Suppose that:\par 
\quad \quad \quad $\displaystyle \forall \alpha \in \partial \Omega ,\ \exists
 v\in T_{\alpha }^{{\mathbb{C}}}(\partial \Omega )::\exists k\in
 {\mathbb{N}},\ \frac{\partial ^{k}{\mathcal{D}}}{\partial v^{k}}(\alpha
 )\neq 0,$ \par 
then  $\Omega $  is aspc and can be equipped with a family of
 polydiscs whose multi-type is the given linear multi-type.
\end{Thrm}
\quad Proof.\ \par 
The fact that there is a good family of polydiscs associated
 to the linear type is given by Theorem~\ref{linkTfFoc20}.\ \par 
\quad It remains to verify the condition on the smallness of the set
  $W$  of weakly pseudo-convex points.\ \par 
\quad Let  $\alpha \in \partial \Omega ,$  we may suppose that  $\alpha
 =0$  and that the complex normal is the  $z_{1}$  axis.\ \par 
\quad Because  $\Omega $  fullfills the hypothesis of the theorem,
 there is a  $j::1<j\leq n,$  a real direction, for instance
 the  $y_{j}$  axis, with  $z_{j}=x_{j}+iy_{j},$  and an integer
  $m,$  such that, with  $\tilde {\mathcal{D}}$  being the restriction
 of  ${\mathcal{D}}$  to the  $z_{j}$  complex plane via the
 diffeomorphism  $\pi ,\ \tilde {\mathcal{D}}:={\mathcal{D}}\circ
 \pi \ :$ \ \par 
\quad \quad \quad $\displaystyle \ \frac{\partial ^{m}\tilde {\mathcal{D}}}{\partial
 y_{j}^{m}}(0)=\frac{\partial ^{m}{\mathcal{D}}}{\partial y_{j}^{m}}(0)\neq
 0.$ \ \par 
\quad The differentiable preparation theorem of Malgrange gives that
 there is a polynomial with  ${\mathcal{C}}^{\infty }$  coefficients,\ \par 
\quad \quad \quad $P(x_{j},\ y_{j})=y_{j}^{m}+\sum_{k=1}^{m}{a_{k}(x_{j})y_{j}^{m-k}}$ \ \par 
and a  ${\mathcal{C}}^{\infty }$  function  $Q(x_{j},\ y_{j}),\
 Q(0)\neq 0$  such that\ \par 
\quad \quad \quad $\tilde {\mathcal{D}}(x_{j},\ y_{j})=Q(x_{j},\ y_{j})P(x_{j},\ y_{j}).$ \ \par 
\quad Hence the zero set of  $\tilde {\mathcal{D}}$  is the same as
 the one of  $P$  and we know, by corollary~\ref{strongPC13},
 that the homogeneous Minkowski dimension of it is less or egal
 to  $\displaystyle 2-\frac{1}{m}.$ \ \par 
\quad Because  ${\mathcal{D}}$  and  $\tilde {\mathcal{D}}$  are  ${\mathcal{C}}^{\infty
 }$  functions,  $\displaystyle \ \frac{\partial ^{m}\tilde {\mathcal{D}}}{\partial
 y_{j}^{m}}\neq 0$  in a neighbourhood of  $0$  with the same
 number  $m,$  hence we have that the  homogeneous Minkowski
 dimension of  $\lbrace \tilde {\mathcal{D}}=0\rbrace $  is less
 or egal to  $\displaystyle 2-\frac{1}{m}$  in all the slices
 parallel to the  $z_{j}$  axis in a neighbourhood of  $0,$ 
 and we are done.  $\blacksquare $ \ \par 
\quad \quad 	A natural question, asked by the referee, is :\ \par 
{\sl Is the condition } $\displaystyle \forall \alpha \in \partial
 \Omega ,\ \exists v\in T_{\alpha }^{{\mathbb{C}}}(\partial \Omega
 )::\exists k\in {\mathbb{N}},\ \frac{\partial ^{k}{\mathcal{D}}}{\partial
 v^{k}}(\alpha )\neq 0$ {\sl  actually necessary ?}\ \par 
\quad \quad 	I have no answer to it, but we shall see that for convex domains
 this condition is a consequence of the linear finite type of
  $\displaystyle \Omega .$ \ \par 
\quad \quad 	We shall need the definition.\ \par 

\begin{Dfnt}
Let  $f$  be a function defined on an open set  ${\mathcal{V}}\subset
 {\mathbb{R}}^{n},\ f\in {\mathcal{C}}^{\infty }({\mathcal{V}})\
 ;$  we shall say that  $f$  is {\bf flat} at  $a\in {\mathcal{V}}$
  if  $\displaystyle \forall \alpha \in {\mathbb{N}}^{n},\ \frac{\partial
 ^{\left\vert{\alpha }\right\vert }f}{\partial x_{1}^{\alpha
 _{1}}\cdot \cdot \cdot \partial x_{n}^{\alpha _{n}}}(a)=0.$ 
\end{Dfnt}

\subsection{Pseudo-convex domains of finite type in  ${\mathbb{C}}^{2}$
 \!\!\!\! .}

\begin{Lmm}
~\label{aspcGlo712} Let  $\displaystyle f(z)$  be a real valued
 smooth function of  $z\in {\mathbb{D}},$  the unit disc in 
 ${\mathbb{C}}\ ;$  if  $\Delta f$  is flat at  $\displaystyle
 0$  then for any  $m\in {\mathbb{N}}$  there is a harmonic function
  $h$  in  ${\mathbb{D}}$  such that  $f-h={\mathcal{O}}(\left\vert{z}\right\vert
 ^{m})$  at the origin.
\end{Lmm}
\quad \quad 	Proof.\ \par 
Take the Taylor expansion of  $f$  at  $0$  :\ \par 
\quad \quad \quad \quad \quad  $f(x+iy)=\sum_{k,l=0}^{m+2}{a_{kl}x^{k}y^{l}}+{\mathcal{O}}(\left\vert{z}\right\vert
 ^{m+3}).$ \ \par 
We get the expansion of  $\Delta f$  near  $0$  :\ \par 
\quad \quad \quad \quad \quad  $\Delta f(x+iy)=\sum_{k=2,l=0}^{m}{k(k-1)a_{kl}x^{k-2}y^{l}}+\sum_{k=0,l=2}^{m}{l(l-1)a_{kl}x^{k}y^{l-2}}+{\mathcal{O}}(\left\vert{z}\right\vert
 ^{m+1}).$ \ \par 
Hence\ \par 
\quad \quad \quad \quad \quad  $\displaystyle \Delta f(x+iy)=\sum_{k,l=0}^{m}{\lbrack (k+1)(k+2)a_{k+2,l}+(l+1)(l+2)a_{k,l+2}\rbrack
 x^{k}y^{l}}+{\mathcal{O}}(\left\vert{z}\right\vert ^{m+1}).$ \ \par 
But  $\Delta f$  flat at  $0$  means that  $\displaystyle \lbrack
 (k+1)(k+2)a_{k+2,l}+(l+1)(l+2)a_{k,l+2}\rbrack =0,$  hence setting\ \par 
\quad \quad \quad \quad \quad  $\displaystyle h:=\sum_{k,l=0}^{m+2}{a_{kl}x^{k}y^{l}},$ \ \par 
we have that \ \par 
\quad \quad \quad \quad \quad  $\displaystyle \Delta h(x+iy)=\sum_{k,l=0}^{m}{\lbrack (k+1)(k+2)a_{k+2,l}+(l+1)(l+2)a_{k,l+2}\rbrack
 x^{k}y^{l}}=0$ \ \par 
because all the coefficients are zero. So we get that  $h$  is
 harmonic and  $f-h={\mathcal{O}}(\left\vert{z}\right\vert ^{m+3}).$
   $\blacksquare $ \ \par 

\begin{Thrm}
~\label{aspcGlo814} Let  $\displaystyle \Omega $  be a domain
 of finite type in  ${\mathbb{C}}^{2}$  then  $\displaystyle \Omega $  is aspc.
\end{Thrm}
\quad \quad 	For the proof of this theorem we shall use the following lemma.\ \par 

\begin{Lmm}
~\label{DC0} Let  $h$  be a real valued harmonic function in
 a disc  $D(0,R)\subset {\mathbb{C}}\ ;$  then  $h$  cannot have
 isolated zeroes.
\end{Lmm}
\quad \quad 	Proof.\ \par 
Suppose that  $\displaystyle h(0)=0,\ h\not\in 0,$  then by the
 mean formula we have for any  $\displaystyle 0\leq r<R,$ \ \par 
\quad \quad \quad \quad \quad  $\displaystyle 0=h(0)=\frac{1}{2\pi }\int_{0}^{2\pi }{h(re^{i\theta
 })d\theta }.$ \ \par 
But  $h$  being real valued on the circle  $\displaystyle C(r):=\lbrace
 \left\vert{z}\right\vert =r\rbrace $  cannot be always positive
 or always negative, hence it must change sign on  $\displaystyle
 C(r)$  so it must be zero at least twice, because  $h$  is continuous.
 This is true for any  $\displaystyle 0<r<R,$  hence the lemma
 is proved.  $\blacksquare $ \ \par 
\quad \quad 	Proof of the theorem.\ \par 
Let  $\Omega \subset {\mathbb{C}}^{2}$  be defined near the origin by\ \par 
\quad \quad \quad \quad \quad  $\displaystyle \rho (z)=\Re z_{1}+f(\Im z_{1},z_{2}).$ \ \par 
We have that\ \par 
\quad \quad \quad \quad \quad  $\displaystyle \rho (z)=\Re z_{1}+f(0,z_{2})+(f(\Im z_{1},z_{2})-f(0,z_{2})).$
 \ \par 
Suppose that  ${\mathcal{D}}:=\Delta f(0,z_{2})$  is flat at
  $0$  then by lemma~\ref{aspcGlo712} for any  $m\in {\mathbb{N}}$
  there is  $h(z_{2})$  harmonic near  $\displaystyle z_{2}=0$
  and such that\ \par 
\quad \quad \quad \quad \quad  $f(0,z_{2})=h(z_{2})+{\mathcal{O}}(\left\vert{z_{2}}\right\vert ^{m}).$ \ \par 
There is a conjugate  $\displaystyle \tilde h$  to  $h$  such
 that  $\displaystyle u:=h+i\tilde h$  is holomorphic in  $\displaystyle
 z_{2}$  near  $0$  and  $\displaystyle \tilde h(0)=0\Rightarrow
 u(0)=0.$  We have  $\displaystyle f(\Im z_{1},z_{2})-f(0,z_{2})=\Im
 z_{1}{\times}g(\Im z_{1},z_{2}),$  with  $g$  smooth as we seen
 in lemma~\ref{1_goodFam40} ; hence we have\ \par 
\quad \quad \quad \quad \quad  $\displaystyle \rho (z)=\Re z_{1}+h(z_{2})+\Im z_{1}{\times}g(\Im
 z_{1},z_{2})+{\mathcal{O}}(\left\vert{z_{2}}\right\vert ^{m}).$ \ \par 
Let  $\displaystyle X:=\lbrace z_{1}=-u(z_{2})\rbrace $  be this
 holomorphic variety. By lemma~\ref{DC0} there is a sequence
  $Z:=\lbrace w_{n}\rbrace _{n\in {\mathbb{N}}}\subset \lbrace
 z_{1}=0\rbrace $  such that  $\displaystyle \tilde h(w_{n})=0$
  and  $\displaystyle w_{n}\rightarrow 0.$ \ \par 
Now take a point  $\displaystyle a_{n}=(a_{n}^{1},w_{n})\in X\Rightarrow
 \Re a_{n}^{1}=-h(w_{n}),\ \Im a_{n}^{1}=-\tilde h(w_{n}).$ 
 We have that  $\displaystyle \Im a_{n}^{1}=-\tilde h(w_{n})=0$  hence\ \par 
\quad \quad \quad \quad \quad  $\displaystyle \rho (a_{n})=\Re a_{n}^{1}+h(w_{n})+{\mathcal{O}}(\left\vert{z_{2}}\right\vert
 ^{m})=-h(w_{n})+h(w_{n})+{\mathcal{O}}(\left\vert{w_{n}}\right\vert
 ^{m})={\mathcal{O}}(\left\vert{w_{n}}\right\vert ^{m}),$ \ \par 
because  $\displaystyle \Im z_{1}{\times}g(\Im z_{1},z_{2})=0$
  on  $\displaystyle a_{n}.$ \ \par 
\quad \quad 	Hence the distance from  $\displaystyle \partial \Omega $  to
 the holomorphic variety  $X$  is  $\displaystyle {\mathcal{O}}(\left\vert{z_{2}}\right\vert
 ^{m})$  near  $0$  along the sequence  $Z$  going to  $0,$ 
 so the type of  $\displaystyle \partial \Omega $  is bigger
 than  $m$  at  $\displaystyle 0.$ \ \par 
\quad \quad 	This being true for any  $m\in {\mathbb{N}}$  we have a contradiction
 with the fact that  $\displaystyle \Omega $  is of finite type
 in D'Angelo sense~\cite{Angelo82}.\ \par 
Hence  $\displaystyle \Delta f(0,z_{2})$  is not flat at  $0$
  and we can apply directly theorem~\ref{aspc3} to get that 
 $\displaystyle \Omega $  is {\bf aspc}.  $\blacksquare $ \ \par 
\ \par 

\subsection{Locally diagonalizable domains.}
\quad In this context, the domains with a locally diagonalizable Levi
 form where introduced by C. Fefferman, J. Kohn and M. Machedon~\cite{FefKohMach90}
 in order to obtain H\"older estimates for the  $\bar \partial
 _{b}$  operator.\ \par 
\quad \quad 	Recall that  $\Omega $  locally diagonalizable means that there
 is a neighbourhood  $\displaystyle V_{\alpha }\subset \partial
 \Omega $  of  $\displaystyle \alpha \in \partial \Omega $  and
  $(L_{1},...,\ L_{n})$  a basis of  ${\mathbb{C}}^{n}$  depending
 smoothly on  $\zeta \in V_{\alpha }$  and diagonalizing the
 Levi form  ${\mathcal{L}}.$ \ \par 
\quad We shall need the following lemma.\ \par 

\begin{Lmm}
~\label{strongPC16}Let  $\Omega $  be a domain locally diagonalizable
 in  ${\mathbb{C}}^{n}$  and of finite linear type. Then the
 determinant of its Levi form is not flat on the complex tangent
 space of  $\partial \Omega .$ 
\end{Lmm}
\quad Proof.\ \par 
Let  $\alpha \in \partial \Omega ,$  then there is a neighbourhood
  $V_{\alpha }$  of  $\alpha $  and  $(L_{1},...,\ L_{n})$  a
 basis of  ${\mathbb{C}}^{n}$  depending smoothly on  $z\in V_{\alpha
 },$  and diagonalizing the Levi form  ${\mathcal{L}},$  with
  $L_{1}$  the complex normal direction, so we have, restricting
  ${\mathcal{L}}$  to the complex tangent space :\ \par 
\ \par 
\[ {\mathcal{L}}(z)={\left({
\begin{matrix}
{\lambda _{2}}&{0}&{\cdots }&{0}\cr 
{\vdots }&{\vdots }&{\vdots }&{\vdots }\cr 
{0}&{\cdots }&{0}&{\lambda _{n}}\cr 
\end{matrix}
}\right)}.\] \ \par 
\ \par 
Hence  ${\mathcal{D}}:=det{\mathcal{L}}=\lambda _{2}\cdots \lambda
 _{n}.$  Now suppose that, for any complex direction  $L_{j},\
 j=2,...,\ n,$  at  $\alpha ,$  there is a real direction  $v_{j},\
 v_{j}\in L_{j},$  such that  $\displaystyle \exists k=k_{j}\in
 {\mathbb{N}},\ \frac{\partial ^{k}\lambda _{j}}{\partial v_{j}^{k}}(\alpha
 )\neq 0,$  then with  $k:=(k_{2},...,\ k_{n})\ :$ \ \par 
\quad \quad \quad $\displaystyle \ \frac{\partial ^{\left\vert{k}\right\vert }{\mathcal{D}}}{\partial
 v_{2}^{k_{2}}\cdots \partial v_{n}^{k_{n}}}(\alpha )\neq 0,$ \ \par 
and  ${\mathcal{D}}$  is not flat at  $\alpha .$  Hence if  ${\mathcal{D}}$
  is flat at  $\alpha ,$  we must have\ \par 
\quad \quad \quad $\displaystyle \exists L_{j},\ \forall v_{j}\in L_{j},\ \forall
 k\in {\mathbb{N}},\ \frac{\partial ^{k}\lambda _{j}}{\partial
 v_{j}^{k}}(\alpha )=0.$ \ \par 
Now this  $j$  is fixed and we slice  $\displaystyle \Omega $ \ \par 
\quad \quad  $\displaystyle \Omega _{j}:=\lbrace z_{2}=\cdot \cdot \cdot
 =z_{j-1}=z_{j+1}=\cdot \cdot \cdot =0\rbrace \cap \Omega .$ \ \par 
We are exactly in the situation of a domain in  ${\mathbb{C}}^{2}$
  and we can use the proof of theorem~\ref{aspcGlo814} to get
 a contradiction with the fact that  $\rho $  has a finite order
 of contact with a real direction in  $\displaystyle L_{j}$ 
 because  $\displaystyle \Omega $  is of finite linear type.\ \par 
\quad Hence we proved\ \par 

\begin{Thrm}
Let  $\Omega $  be a domain locally diagonalizable in  ${\mathbb{C}}^{n}$
  and of finite linear type. Then  $\Omega $  is aspc.
\end{Thrm}

\subsection{Convex domains.}

\begin{Thrm}
~\label{strongPC14}Let  $\Omega $  be convex in a neighborhood
 of  $0\in {\mathbb{R}}^{n+1}.$  Suppose that the tangent space
 at  $0$  is  $x_{n+1}=0$  and  $\partial \Omega =\lbrace x_{n+1}=f(x_{1},...,x_{n})\rbrace
 ,$  with  $f$  convex. If the determinant of the hessian of
  $f$  is flat at  $0$  then  $f$  is flat in a direction  $x=(x_{1},\
 ...,\ x_{n})\in {\mathbb{R}}^{n}$  of the tangent space at  $0.$ 
\end{Thrm}
\quad Proof.\ \par 
If  $f$  is not flat in any direction, we can find  $\alpha >0$
  and  $m\in {\mathbb{N}}$  such that  $f(x)\geq \alpha \left\vert{x}\right\vert
 ^{2m}$  in a ball  $B(0,\ R)\subset {\mathbb{R}}^{n}.$  Let
 us see the functions\ \par 
\quad \quad \quad \quad $\displaystyle h(x):=\frac{\alpha }{2}\left\vert{x}\right\vert
 ^{2m},\ g(x):=\frac{\alpha }{2}\left\vert{x}\right\vert ^{2m}+\epsilon
 \left\vert{x}\right\vert ^{2}+\delta ,$ \ \par 
with  $\epsilon >0$  and  $\delta >0.$  Denote  $H_{f}$  the
 hessian of the function  $f.$ \ \par 
\quad Because  $detH_{f}$  is flat at  $0,$  there is a ball  $B(0,\
 r)\subset {\mathbb{R}}^{n}$  such that:\ \par 
\quad \quad \quad \quad $\forall x\in B(0,r),\ detH_{f}(x)\leq detH_{h}(x)$ \ \par 
and\ \par 
\begin{equation}  \forall \epsilon >0,\ \forall \delta >0,\ detH_{h}<detH_{g}.\label{conveFini32}\end{equation}\
 \par 
We choose  $\epsilon $  and  $\delta $  so small that there is
 a real  $t$  such that\ \par 
\quad \quad \quad \quad \quad  $\displaystyle r>t>(2\frac{\delta +\epsilon r^{2}}{\alpha })^{1/2m}$ \ \par 
then\ \par 
\quad \quad \quad $\displaystyle \ \frac{\alpha }{2}t^{2m}>\epsilon t^{2}+\delta
 \Rightarrow \alpha t^{2m}>\frac{\alpha }{2}t^{2m}+\epsilon t^{2}+\delta
 ,$ \ \par 
hence,\ \par 
\quad \quad \quad \quad 	\begin{equation}  \forall x::\left\vert{x}\right\vert =t,\ f(x)\geq
 \alpha \left\vert{x}\right\vert ^{2m}>g(x).\label{conveFini33}\end{equation}\
 \par 
On the other hand, because  $g(0)=\delta >f(0)=0,$  and  $f$
  and  $g$  are continuous, we get\ \par 
\quad \quad \quad $\exists s>0,\ s<t::\forall x,\ \left\vert{x}\right\vert <s,\
 f(x)<g(x).$ \ \par 
The maximum principle for the Monge-Amp\`ere operator says~\cite{Zuily95}
 :\ \par 

\begin{Lmm}
Let  $v$  be a convex function (i.e.  $H_{v}\geq 0$ ) defined
 in a bounded open set  $V$  and a regular function  $\rho $  such that\par 
\quad \quad \quad $\displaystyle \mathrm{d}\mathrm{e}\mathrm{t}H_{v}(x)>\mathrm{d}\mathrm{e}\mathrm{t}H_{\rho
 }(x),\ v\leq \rho $  on  $\partial V,$ \par 
then  $v\leq \rho $  on  $V.$ 
\end{Lmm}
\ \par 
Because  $detH_{g}>detH_{f}$  in  $B(0,r)$  by~(\ref{conveFini32})
 and  $g<f$  on  $\partial B(0,t)$  by~(\ref{conveFini33}), we
 can apply this principle, i.e.  $g\leq f$  everywhere in  $B(0,\
 t)$  which is a contradiction in the ball  $B(0,s).$  Hence
  $f$  has to be flat in some direction.  $\blacksquare $ \ \par 

\begin{Crll}
~\label{strongPC15}Let  $\Omega $  be a convex domain in a neighbourhood
 of  $0\in \partial \Omega \subset {\mathbb{R}}^{n}.$  If  $\partial
 \Omega $  is flat in no direction of its tangent space at  $0,$
  then the determinant of the hessian of   $\Omega $  is not flat at  $0.$ 
\end{Crll}
\quad Proof.\ \par 
If not we have a contraction with theorem~\ref{strongPC14}. 
 $\blacksquare $ \ \par 
\ \par 
\quad Let us see now the case of a convex domain of finite type in
  ${\mathbb{C}}^{n}.$  We shall need the following lemma.\ \par 

\begin{Lmm}
~\label{5_examAspc41} Let  $\Omega $  be a convex domain of finite
 type in  ${\mathbb{C}}^{n}$  then for any complex line  $L$
  in the tangent complex space at  $0\in \partial \Omega $  there
 is at most one real direction  $v$  in  $L$  such that  $\partial
 \Omega $  is flat in this direction at  $\displaystyle 0.$ 
\end{Lmm}
\quad \quad 	Proof.\ \par 
We can choose  $\rho (z)=\Re z_{1}-f(\Im z_{1},\ z_{2},...,z_{n})$
  as defining function for  $\Omega $  with  $f$  a positive
 real valued convex function and with the  $\displaystyle z_{n}$
  axis  $\displaystyle L_{n}$  as the given  $L.$  (Here the
 complex normal direction is  $L_{1}$  as usual).\ \par 
Suppose there are two such directions  $\displaystyle v_{1},v_{2}$
  in  $\displaystyle L_{n}$  this means\ \par 
\quad \quad \quad \quad \quad  $\displaystyle \forall k\in {\mathbb{N}},\ \frac{\partial ^{k}\rho
 }{\partial v_{j}^{k}}(0)=0,\ j=1,2.$ \ \par 
The vector  $\displaystyle v_{1}$  can be seen as a point  $\displaystyle
 a_{1}$  in the complex plane  $\displaystyle P_{n}=\lbrace z_{1}=z_{2}=\cdot
 \cdot \cdot =z_{n-1}=0\rbrace $  and also  $\displaystyle v_{2}$
  corresponds to the point  $\displaystyle a_{2}\in P_{n}.$ 
 Let  $\displaystyle t\in \lbrack 0,1\rbrack ,\ a_{t}:=ta_{1}+(1-t)a_{2}\in
 P_{n},$  because  $f$  is convex this implies that  $\displaystyle
 0\leq f(a_{t})\leq tf(a_{1})+(1-t)f(a_{2})$  and this means
 that the order of contact in the direction  $\displaystyle v=tv_{1}+(1-t)v_{2}$
  is bigger than the minimum of the order of contact in the directions
  $\displaystyle v_{1}$  and  $\displaystyle v_{2},$  hence\ \par 
\quad \quad \quad \quad \quad  $\displaystyle \forall k\in {\mathbb{N}},\ \frac{\partial ^{k}\rho
 }{\partial v^{k}}(0)=0,\ j=1,2,$  with  $\displaystyle v=tv_{1}+(1-t)v_{2}.$
 \ \par 
This being true for any  $\displaystyle t\in \lbrack 0,1\rbrack
 $  we have that  $f$  is flat in the sector of  $\displaystyle
 P_{n}$  between  $\displaystyle v_{1}$  and  $\displaystyle
 v_{2},$  but  $f$  being  ${\mathcal{C}}^{\infty }$  this implies
 that  $f$  is flat at  $\displaystyle 0.$ \ \par 
\quad By a result of Boas and Straube~\cite{BoasStraube92} we have
 that for a convex domain the multi type or the order of contact
 with real lines is the same, the multi type of  $\partial \Omega
 $  being finite, this means that there is a real direction in
  $\displaystyle L$  which is not flat, hence a contradiction
 which gives the lemma.  $\blacksquare $ \ \par 

\begin{Crll}
~\label{strongPC10}Let  $\Omega $  be a convex domain of finite
 type in  ${\mathbb{C}}^{n}$  near  $0\in \partial \Omega .$
  There is a complex line  $L$  in the tangent complex space
 at  $0$  and a real vector  $v\in {\mathbb{C}}L,$  such that
 the determinant of the Levi form of a defining function for
  $\Omega $  near  $0$  is not flat in the direction  $v.$ 
\end{Crll}
\quad Proof.\ \par 
Let  $L_{2},...,\ L_{n}$  be an orthonormal basis of  $T_{0}^{{\mathbb{C}}}(\partial
 \Omega ).$  Because  $\Omega $  is of finite type, we know by
 lemma~\ref{5_examAspc41} that in any complex direction  $L_{j},\
 2\leq j\leq n,$  there is at most one real direction in which
  $\partial \Omega $  is flat ; we can always take that direction
 to be the  $\displaystyle y_{j}$  axis without changing the
 ambiant complex structure. If such a direction does not exist
 we still take the  $\displaystyle y_{j}$  axis in the following.\ \par 
We set  $E$  to be the subspace  $E:=\lbrace y_{2}=\cdot \cdot
 \cdot =y_{n}=0\rbrace \cap \lbrace z_{1}=0\rbrace .$ \ \par 
We write the defining function as usual\ \par 
\quad \quad \quad \quad \quad  $\displaystyle \rho (z)=\Re z_{1}-f(\Im z_{1},z_{2},...,z_{n}),$ \ \par 
hence the domain  $\Omega \cap E$  has defining function\ \par 
\quad \quad \quad $\tilde \rho (x):=-f(0,x_{2},...,x_{n}).$ \ \par 
Let  ${\mathcal{L}}(z_{1},...,z_{n}):=\partial \bar \partial
 \rho (z)$  be the Levi form of  $\Omega ,$  we have\ \par 
\quad \quad \quad \quad 	\begin{equation}  \partial \bar \partial f(x,0)=-{\mathcal{L}}(x,\
 0)=\lbrace \frac{\partial ^{2}f}{\partial x_{j}\partial x_{k}}(x,\
 0)\rbrace _{j,k=2,...,n}=H_{\tilde f}(x)\label{conveFini45}\end{equation}\
 \par 
with  $\displaystyle \tilde f(x_{2},...,x_{n}):=f(0,x),$  and
 the new convex set  $\Omega _{1}:=\Omega \cap E$  still verifies
 the conditions of corollary~\ref{strongPC15} :  $\tilde {\mathcal{D}}(x):=detH_{\tilde
 f}(x)$  is not flat because we get rid of the flat directions
 ; hence there is a real vector  $v$  in the tangent space at
  $0$  for  $\partial \Omega _{1}$  such that  $\tilde {\mathcal{D}}$
  is not flat in the direction  $v.$  This means\ \par 
\quad \quad \quad \quad $\displaystyle \exists k\in {\mathbb{N}}::\frac{\partial ^{k}\tilde
 {\mathcal{D}}}{\partial v^{k}}(0)\neq 0\ ;$ \ \par 
but, using~(\ref{conveFini45}), we get\ \par 
\quad \quad \quad $\displaystyle \ \frac{\partial ^{k}{\mathcal{D}}}{\partial v^{k}}(0)=\frac{\partial
 ^{k}\tilde {\mathcal{D}}}{\partial v^{k}}(0)\neq 0.\ \blacksquare $ \ \par 

\begin{Thrm}
Let  $\Omega $  be a convex domain of finite type in  ${\mathbb{C}}^{n};$
  then  $\Omega $  is aspc.
\end{Thrm}
\quad Proof.\ \par 
By use of corollary~\ref{strongPC10}, it remains to apply theorem~\ref{aspc3}.
  $\blacksquare $ \ \par 

\subsection{Domains with real analytic boundary.}

\begin{Lmm}
Let  $\Omega $  be a bounded domain with real analytic boundary,
 then  $\Omega $  is of finite linear type.
\end{Lmm}
\quad Proof.\ \par 
Take a point  $\alpha \in \partial \Omega $  and suppose that
 a real line through  $\alpha $  has a contact of infinite order
 with  $\partial \Omega ,$  then, using Lojasiewicz~\cite{Tougeron72}
 we get that the line, which is real analytic, and  $\partial
 \Omega $  are regularly situated, hence the line must be contained
 in  $\partial \Omega .$  But this cannot happen because  $\partial
 \Omega $  is bounded.  $\blacksquare $ \ \par 
\quad In fact we have a better result because we know, by the work
 of K. Diederich and J-E. Fornaess~\cite{DiedForn78}, that  $\Omega
 $  is of finite type.\ \par 
\ \par 
\quad The function  ${\mathcal{D}}=det{\mathcal{L}}$  is also real
 analytic, hence if  ${\mathcal{D}}$  is flat at a point  $\alpha
 \in \partial \Omega ,$  this means in particular that  $\displaystyle
 \forall v\in T_{\alpha }(\partial \Omega ),\ \forall k\in {\mathbb{N}},\
 \frac{\partial ^{k}{\mathcal{D}}}{\partial v^{k}}(\alpha )=0,$
  hence  ${\mathcal{D}}$  is identically zero on  $\partial \Omega
 .$  This says that all the points of  $\partial \Omega $  are
 non stricly pseudo-convex points. But this is impossible because
  $\partial \Omega $  is compact, hence contains at least a strictly
 pseudo-convex point, because of the following simple and well
 known lemma~\cite{Krantz01} :\ \par 

\begin{Lmm}
Let  $\Omega $  be a bounded domain in  ${\mathbb{R}}^{n},$ 
 with a smooth boundary of class  ${\mathcal{C}}^{3}.$  Then
  $\partial \Omega $  contains a point of strict convexity.
\end{Lmm}
\quad Now let  $\alpha \in \partial \Omega $  and suppose that  ${\mathcal{D}}$
  is flat in all the complex tangent directions of  $T_{\alpha
 }^{{\mathbb{C}}}(\partial \Omega ),$  then, because  $\partial
 \Omega $  is of finite type, we can recover the derivatives
 in the "missing direction", namely the real direction conjugate
 to the normal one, by brackets of derivatives in the complex
 tangent directions.\ \par 
Hence we have that  ${\mathcal{D}}$  is also flat in the direction
 conjugate to the normal one, but this will implies that  ${\mathcal{D}}$
  is flat at the point  $\alpha ,$  and this is forbidden by
 the lemma. So we can apply theorem~\ref{aspc3} to conclude:\ \par 

\begin{Thrm}
Let  $\Omega $  be a domain in  ${\mathbb{C}}^{n}$  with real
 analytic boundary, then  $\Omega $  is aspc and of finite linear type.
\end{Thrm}
\vfill\eject\ \par 
\ \par 

\section{Convex domains of finite type.~\label{7_ConvFini40}}
\quad McNeal~\cite{McNeal94}, introduced tools for studying the geometry
 of convex domains of finite type : a family of polydiscs and
 a related pseudo-distance which are well suited to these domains.\ \par 
\quad These tools were used and a little bit modified by different authors :\ \par 
McNeal and Stein~\cite{McNealStein97}, J. Bruna, P. Charpentier
 and Y. Dupain~\cite{BruChaDup98}, K. Diederich and E. Mazzilli~\cite{DiedMazz01},
 A. Cumenge~\cite{Cum01} and also  T. Hefer~\cite{Hefer04}, among others.\ \par 
\quad We start first with notations and definitions taken from Hefer~\cite{Hefer04}
 in order for the reader to follow easily the citations we use.
 This means that the polydiscs in the family seem different but
 we shall show, in section~\ref{8_CarlesonMeas30}, that they
 are the same than the ones defined in section~\ref{2_BonFamille42}.\ \par 
\quad Let  $r$  be a defining function for  $\Omega ,\ \Omega :=\lbrace
 z\in {\mathbb{C}}^{n}::r(z)<0\rbrace ,$  where  $\Omega $  is
 a convex domain of finite type.\ \par 
\quad Hefer uses the  $\epsilon $  distance in the direction  $v$  :\ \par 
\quad \quad \quad \quad 	\begin{equation}  \tau (\zeta ,\ v,\ \epsilon ):=\sup  \lbrace
 c\ :\ \left\vert{r(\zeta +\lambda v)-r(\zeta )}\right\vert \leq
 \epsilon ,\ \forall \lambda \in {\mathbb{C}}::\left\vert{\lambda
 }\right\vert \leq c\rbrace ,\label{suPr67}\end{equation}\ \par 
and built two  $\epsilon $  extremal bases (introduced in~\cite{BruChaDup98}),
 a variant of the original one of McNeal~\cite{McNeal94}, which
 are equivalent, from which we keep one :\ \par 
\quad \quad \quad \quad \quad  $b_{\epsilon }(\zeta )=(v_{1}(\zeta ,\ \epsilon ),...,v_{n}(\zeta
 ,\ \epsilon )),$ \ \par 
and the  $\epsilon $  distance in the direction  $v_{k}\ :$ \ \par 
\quad \quad \quad $\displaystyle \forall k=1,...,\ n,\ \tau _{k}(\zeta ,\ \epsilon
 ):=\tau (\zeta ,\ v_{k},\ \epsilon ).$ \ \par 
\quad This allows him to define a family of polydiscs\ \par 
\quad \quad \quad \quad 	\begin{equation}  \forall t>0,\ tP_{\epsilon }(\zeta ):=\lbrace
 z=\zeta +\sum_{k=1}^{n}{w_{_{k}}v_{k}(\zeta ,\ \epsilon )}\in
 {\mathbb{C}}^{^{n}}::\ \forall k=1,...,\ n,\ \left\vert{w_{k}}\right\vert
 <t\tau _{k}(\zeta ,\ \epsilon )\rbrace ,\label{suPr65}\end{equation}\ \par 
and the {\bf pseudo-distance}  $d(z,\ \zeta ):=\inf  \lbrace
 \epsilon ::z\in P_{\epsilon }(\zeta )\rbrace $  associated to
 it. See the nice introduction in~\cite{Hefer04} to see why this
 definition is relevant.\ \par 
\quad From his theorem 1.7 of~\cite{Hefer04} I just keep the "geometrical"
 part.\ \par 

\begin{Thrm}
~\label{subPrinciple66}Let  $\Omega \subset {\mathbb{C}}^{n}$
  be a smooth convex domain of finite type and let  $(m_{1},...,\
 m_{n})$  be its multitype.\par 
If  ${\mathcal{U}}$  is a sufficiently small compact neighborhood
 of  $\partial \Omega $ \!\!\!\! , if  $\zeta \in {\mathcal{U}}$
  and if   $(m_{1}(\zeta ),...,\ m_{n}(\zeta ))$  is the multitype
 of  $\partial \Omega _{\zeta }:=\lbrace z\in {\mathbb{C}}^{n}::r(z)=r(\zeta
 )\rbrace $  at the point  $\zeta ,$  then there are constants
  $c,C>0$  depending only on  ${\mathcal{U}}$  (and on the fixed
 defining function  $r$  of  $\Omega $ ) such that\par 
\quad \quad \quad \quad 	\begin{equation}  c\epsilon ^{1/m_{j}(\zeta )}\leq \tau _{j}(\zeta
 ,\epsilon )\leq C\epsilon ^{1/m_{j}(\zeta )}.\label{subPrinciple616}\end{equation}
\end{Thrm}

      Hence we have a family of polydiscs\ \par 
\quad \quad \quad \quad \quad \quad \quad  \begin{equation}  {\mathcal{P}}:=\lbrace P_{\epsilon }(\zeta
 )\rbrace _{\zeta \in {\mathcal{U}},\epsilon >0}\label{subPrinciple68}\end{equation}\
 \par 
which is equivalent to the family used by McNeal and Stein~\cite{McNealStein97}.\
 \par 
\quad We shall extract from proposition 2.7 of~\cite{Hefer04} the following
 facts we shall need later.\ \par 
\quad \quad \quad $\forall t>0,\ \exists c_{t},\ \exists C_{t}$  depending only
 on  $t$  such that\ \par 
\quad \quad \quad \quad \quad \quad \quad  \begin{equation}  \forall \zeta \in {\mathcal{U}},\ P_{c_{t}\epsilon
 }(\zeta )\subset tP_{\epsilon }(\zeta )\subset P_{C_{t}\epsilon
 }(\zeta ).\label{convexInterpol526}\end{equation}\ \par 
\quad There are constants  $C_{1}>1,\ c_{2}<1$  and  $c_{3}>0,$  independant
 of  $\zeta $  and  $\epsilon ,$  such that\ \par 
\quad \quad \quad \quad \quad  $\displaystyle \forall \zeta \in {\mathcal{U}},\forall \epsilon
 >0,\ \frac{1}{2}P_{\epsilon }(\zeta )\subset C_{1}P_{\epsilon
 /2}(\zeta )\ ;$ \ \par 
\quad \quad \quad \quad $\displaystyle \forall \epsilon >0,\ \forall t<c_{2}\epsilon
 ,\ \forall \zeta ,\ C_{1}P_{t}(\zeta )\subset P_{\epsilon }(\zeta )\ ;$ \ \par 
\quad \quad \quad \quad \quad \quad \quad 	\begin{equation}  \forall \zeta \in \Omega ,\ c_{3}P_{\left\vert{r(\zeta
 )}\right\vert }(\zeta )\subset \Omega .\label{suPr610}\end{equation}\ \par 
There is a constant  $C_{3}$  independent of  $z,\ \zeta \in
 {\mathcal{U}},$  and independent of  $s>0$  such that\ \par 
\quad \quad \quad \quad \quad \quad \quad  \begin{equation}  P_{s}(z)\cap P_{s}(\zeta )\neq \emptyset \Rightarrow
 P_{s}(z)\subset C_{3}P_{s}(\zeta ).\label{suPr611}\end{equation}\ \par 
This implies, with  $\sigma _{2n}(P)$  the euclidean volume of  $P,$ \ \par 
\quad \quad \quad \quad $\displaystyle \ \frac{1}{C_{3}^{2n}}\sigma _{2n}(P_{s}(\zeta
 ))\leq \sigma _{2n}(P_{s}(z))\leq C_{3}^{2n}\sigma _{2n}(P_{s}(\zeta
 )).$ \ \par 
But  $\displaystyle \sigma _{2n}(P_{s}(\zeta ))=\tau _{1}(\zeta
 ,s)^{2}\prod_{j=2}^{n}{\tau _{j}(\zeta ,\ s)^{2}}$  and  $\tau
 _{1}(\zeta ,s)\simeq s,$  so we have\ \par 
\quad \quad \quad \quad \quad \quad \quad 	\begin{equation}  P_{s}(z)\cap P_{s}(\zeta )\neq \emptyset \Rightarrow
 \prod_{j=2}^{n}{\tau _{j}(\zeta ,\ s)^{2}}\simeq \prod_{j=2}^{n}{\tau
 _{j}(z,\ s)^{2}}.\label{subPrinciple618}\end{equation}\ \par 
If  $\pi (z)$  is the projection of  $z$  to  $\partial \Omega
 ,$  then we have the estimate\ \par 
\quad \quad $\displaystyle d(z,\ \pi (z))\simeq \left\vert{r(z)}\right\vert
 \ ;\ z\in P_{\epsilon }(\zeta )\Rightarrow d(z,\ \zeta )\leq
 \epsilon \ ;\ z\notin P_{\epsilon }(\zeta )\Rightarrow d(z,\
 \zeta )\gtrsim \epsilon \ ;$ \ \par 
and\ \par 
\quad \quad \quad \quad $d(z,\ \zeta )\leq \epsilon \Rightarrow z\in P_{t}(\zeta )$ 
 for all  $t\gtrsim \epsilon $  and  $d(z,\ \zeta )\geq \epsilon
 \Rightarrow z\notin P_{t}(\zeta )$  for all  $t\lesssim \epsilon .$ \ \par 

\subsection{Szeg\"o and Poisson-Szeg\"o kernels.}
\quad We shall continue with notions introduced by McNeal and Stein~\cite{McNealStein97}
 ; we modify slightly the previous notations :\ \par 
\quad $\forall x,y\in \partial \Omega ,\ \rho (x,y):=d(x,y)$  is the
 pseudo-distance which, proved by McNeal~\cite{McNeal92}, gives
 a structure of space of homogeneous type to  $\partial \Omega .$ \ \par 
\quad The "distance" in  $\bar \Omega ,\ \rho ^{*}(z,w)$  is defined by:\ \par 
\quad \quad \quad $\displaystyle \rho ^{*}(z,w):=\left\vert{r(z)}\right\vert +\left\vert{r(w)}\right\vert
 +\rho (\pi (z),\ \pi (w)),$ \ \par 
where  $\pi $  is the normal projection on the boundary  $\partial
 \Omega $  of  $\Omega ,$  well defined in  ${\mathcal{U}},$
  (shrinking  ${\mathcal{U}}$  if necessary).\ \par 
\ \par 
\quad We have, still following McNeal and Stein~\cite{McNealStein97} :\ \par 
$\bullet $  the pseudo-balls on  $\partial \Omega $  are defined by\ \par 
\quad \quad \quad $\forall \alpha \in \partial \Omega ,\ B(\alpha ,\ \epsilon ):=P_{\epsilon
 }(\alpha )\cap \partial \Omega \ ;$  \ \par 
$\bullet $  the "tents" are defined in  $\displaystyle {\mathcal{U}}\cap
 \bar \Omega ,$  where  ${\mathcal{U}}$  is a sufficiently small
 compact neighborhood of  $\partial \Omega $  defined in theorem~\ref{subPrinciple66},
 by  $\displaystyle \forall a\in {\mathcal{U}}\cap \bar \Omega
 ,\ T_{a}(\ \epsilon )=P_{\epsilon }(a)\cap \bar \Omega .$ \ \par 
We shall also need this notation :\ \par 
\quad $\displaystyle \forall z\in \Omega ,\ \forall w\in \bar \Omega
 ,\ T(z,w)$  is the smallest "tent" containing the two points
  $z$  and  $w.$  The  $\epsilon $  underlying this tent is equivalent
 to  $\displaystyle \rho ^{*}(z,w)$  as done in~\cite{McNealStein97}.\ \par 
\ \par 
\quad Let  $S(z,w)$  the Szeg\"o kernel of  $\Omega ,$  i.e. the kernel
 associated to the orthogonal projection from  $L^{2}(\partial
 \Omega )$  onto the Hardy space  $H^{2}(\Omega ).$ \ \par 
\quad We have(~\cite{McNealStein97}, p 521) :\ \par 
\quad \quad \quad \quad \quad  $\displaystyle \forall (z,w)\in \Omega {\times}\Omega \backslash
 \Delta ,\ \left\vert{S(z,w)}\right\vert \lesssim \frac{\delta
 }{\sigma _{2n}(T(z,w))},\ \delta :=\rho ^{*}(z,w).$ \ \par 
Keeping  $z\in \Omega $  and pushing  $w$  to  $\partial \Omega
 ,$  we still have\ \par 
\quad \quad \quad \quad \quad \quad \quad  \begin{equation}  \forall (z,y)\in \Omega {\times}\partial \Omega
 ,\ \left\vert{S(z,y)}\right\vert \lesssim \frac{\delta }{\sigma
 _{2n}(T(z,y))},\ \delta :=\rho ^{*}(z,y).\label{subPrinciple30}\end{equation}\
 \par 
\quad We also have the following estimates (~\cite{McNealStein97}, p 525)\ \par 
\quad \quad  \begin{equation}  \sigma _{2n-1}(B(x,\ \epsilon ))\simeq \epsilon
 \prod_{j=2}^{n}{\tau _{j}(x,\ \epsilon )^{2}}\ ;\ \sigma _{2n}(T_{z}(\
 \epsilon ))\simeq \epsilon ^{2}\prod_{j=2}^{n}{\tau _{j}(x,\
 \epsilon )^{2}}\simeq \epsilon \sigma _{2n-1}(B(x,\epsilon )).\label{strongPC930}\end{equation}\
 \par 
Now we have, by its very definition~(\ref{suPr67}), that  $\displaystyle
 \forall k\in {\mathbb{N}},\ \tau _{j}(x,\ 2^{k}\epsilon )\geq
 \tau _{j}(x,\ \epsilon ),$  hence using~(\ref{strongPC930})\ \par 
\quad \quad \quad \quad \quad \quad \quad  \begin{equation}  \forall k\in {\mathbb{N}},\ \sigma _{2n-1}(B(x,\
 2^{k}\epsilon ))\gtrsim 2^{k}\sigma _{2n-1}(B(x,\ \epsilon )).\label{strongPC933}\end{equation}\
 \par 
\ \par 
\quad Let  $z\in \Omega ,\ x=\pi (z)\in \partial \Omega $  be fixed
 and cover  $\partial \Omega $  by annuli\ \par 
\quad \quad $C_{k}:=B(x,\ 2^{k}\delta )\backslash B(x,\ 2^{k-1}\delta ),\
 k\geq 1$  and  $C_{0}:=B(x,\ \delta )$  with  $\delta =\rho
 ^{*}(z,\ z)=2\left\vert{r(z)}\right\vert .$ \ \par 

\begin{Lmm}
~\label{subPrinciple40}With  $z\in \Omega ,\ x=\pi (z)\in \partial
 \Omega ,\ \delta :=\rho ^{*}(z,z)=2\left\vert{r(z)}\right\vert
 ,$  we have :\par 
\quad \quad 	\begin{equation}  \forall z\in \Omega ,\ \forall y\in \partial
 \Omega ,\ \left\vert{S(z,y)}\right\vert \lesssim \frac{1}{\sigma
 _{2n-1}(B(x,\delta /2))}{\11}_{B(x,\delta )}(y)+\sum_{k\in {\mathbb{N}}}{\frac{1}{\sigma
 _{2n-1}(B(x,\ 2^{k}\delta ))}{\11}_{C_{k}}(y)}.\label{subP31}\end{equation}
\end{Lmm}
\quad Proof.\ \par 
This is a well known technique of harmonic analysis (we already
 used it in~\cite{AmarBonami} for the same goal, for instance).\ \par 
By inequality~(\ref{subPrinciple30}) we get, with  $y\in \partial
 \Omega \cap C_{k},$  hence  $\rho (x,y)\leq 2^{k}\delta ,$ \ \par 
\quad \quad \quad $\displaystyle \ \left\vert{S(z,y)}\right\vert \lesssim \frac{\rho
 ^{*}(z,y)}{\sigma _{2n}(T(z,y))}=\frac{\left\vert{r(z)}\right\vert
 +\rho (x,y)}{\sigma _{2n}(T(z,y))}\leq \frac{(1+2^{k})\delta
 }{\sigma _{2n}(T(z,y))}.$ \ \par 
So\ \par 
\quad \quad \quad $\displaystyle \ \left\vert{S(z,y)}\right\vert \lesssim \frac{\delta
 }{\sigma _{2n}(T(z,y))}{\11}_{B(x,\delta )}(y)+\sum_{k\geq 1}{\frac{(1+2^{k})\delta
 }{\sigma _{2n}(T(z,y))}{\11}_{C_{k}}(y)}.$ \ \par 
If  $y\in C_{k},\ k\geq 1,$  we have  $\sigma _{2n}(T(z,y))\gtrsim
 \sigma _{2n}(T_{z}(\ 2^{k-1}\delta ))$  and for  $y\in B(x,\delta
 )$  we have\ \par 
\quad \quad \quad \quad \quad  $\sigma _{2n}(T(z,y))\gtrsim \sigma _{2n}(T_{z}(\left\vert{r(z)}\right\vert
 ))=\sigma _{2n}(T_{z}(\delta /2)),$  so\ \par 
\quad \quad \quad \quad \quad  $\displaystyle \ \left\vert{S(z,y)}\right\vert \lesssim \frac{\delta
 }{\sigma _{2n}(T_{z}(\delta /2))}{\11}_{B(x,\delta )}(y)+\sum_{k\geq
 1}{\frac{(1+2^{k})\delta }{\sigma _{2n}(T_{z}(2^{k}\delta ))}{\11}_{C_{k}}(y)}.$
 \ \par 
Now by the equivalences~(\ref{strongPC930}) we have  $\displaystyle
 \sigma _{2n}(T_{z}(h))\simeq h\sigma _{2n-1}(B(x,h)),$  so we get\ \par 
\quad \quad \quad \quad \quad  $\displaystyle \ \left\vert{S(z,y)}\right\vert \lesssim \frac{\delta
 }{\delta \sigma _{2n-1}(B(x,\delta /2))}{\11}_{B(x,\delta )}(y)+\sum_{k\geq
 1}{\frac{(1+2^{k})\delta }{\delta (1+2^{k})\sigma _{2n-1}(B(x,\
 N^{k}\delta ))}{\11}_{C_{k}}(y)},$ \ \par 
hence\ \par 
\quad \quad \quad \quad \quad  $\displaystyle \ \left\vert{S(z,y)}\right\vert \lesssim \frac{{\11}_{B(x,\delta
 )}(y)}{\sigma _{2n-1}(B(x,\delta /2))}+\sum_{k\geq 1}{\frac{{\11}_{C_{k}}(y)}{\sigma
 _{2n-1}(B(x,\ 2^{k-1}\delta ))}},$   $\blacksquare $ \ \par 

\begin{Lmm}
~\label{subPrinciple41}We have, with   $z\in \Omega ,\ x=\pi
 (z)\in \partial \Omega ,\ \delta :=\rho ^{*}(z,z)=2\left\vert{r(z)}\right\vert
 ,$ \par 
\quad \quad \quad \quad \quad  $\displaystyle \ {\left\Vert{S(z,\cdot )}\right\Vert}_{p}\lesssim
 \frac{1}{\sigma _{2n-1}(B(x,\delta ))^{1/p'}},$ \par 
where  $p'$  is the conjugate exponent of  $p.$ 
\end{Lmm}
\quad Proof:\ \par 
Lemma~\ref{subPrinciple40} gives us\ \par 
\quad \quad \quad \quad \quad  $\displaystyle \ \left\vert{S(z,y)}\right\vert \lesssim \frac{{\11}_{B(x,\delta
 )}(y)}{\sigma _{2n-1}(B(x,\delta /2))}+\sum_{k\geq 1}{\frac{{\11}_{C_{k}}(y)}{\sigma
 _{2n-1}(B(x,\ 2^{k-1}\delta ))}},$ \ \par 
hence integrating on  $\partial \Omega ,$  we get\ \par 
\quad \quad \quad \quad \quad  $\displaystyle \ {\left\Vert{S(z,\cdot )}\right\Vert}_{p}^{p}\lesssim
 \frac{\sigma _{2n-1}(B(x,\delta ))}{\sigma _{2n-1}(B(x,\delta
 /2))^{p}}+\sum_{k\geq 1}{\frac{\sigma _{2n-1}(C_{k})}{\sigma
 _{2n-1}(B(x,\ 2^{k-1}\delta ))^{p}}}.$ \ \par 
From  $\displaystyle C_{k}\subset B(x,\ 2^{k}\delta ),$  we get
  $\sigma _{2n-1}(C_{k})\leq \sigma _{2n-1}(B(x,\ 2^{k}\delta )),$  hence\ \par 
\quad \quad \quad \quad \quad  $\displaystyle \ {\left\Vert{S(z,\cdot )}\right\Vert}_{p}^{p}\lesssim
 \frac{\sigma _{2n-1}(B(x,\delta ))}{\sigma _{2n-1}(B(x,\delta
 /2))^{p}}+\sum_{k\geq 1}{\frac{1}{\sigma _{2n-1}(B(x,\ 2^{k-1}\delta
 ))^{p-1}}}.$ \ \par 
\quad Because these pseudo-balls are associated to a space of homogeneous
 type, there is a constant  $K$  such that  $\sigma _{2n-1}(B(x,\
 2h))\leq K\sigma _{2n-1}(B(x,h)).$  Using also inequality~(\ref{strongPC933})
 we get, with  $t=\delta /2=\left\vert{r(z)}\right\vert \ :$ \ \par 
\quad \quad \quad $\displaystyle \ {\left\Vert{S(z,\cdot )}\right\Vert}_{p}^{p}\lesssim
 \frac{1}{\sigma _{2n-1}(B(x,t))^{p-1}}+\sum_{k\geq 1}{\frac{1}{2^{(p-1)k}}\frac{1}{\sigma
 _{2n-1}(B(x,\ t))^{p-1}}}$ \ \par 
\quad \quad \quad \quad \quad \quad \quad \quad $\displaystyle \lesssim \frac{1}{\sigma _{2n-1}(B(x,\ t))^{p-1}}{\left({1+\sum_{k\geq
 1}{\frac{1}{2^{(p-1)k}}}}\right)}\lesssim \frac{1}{\sigma _{2n-1}(B(x,\
 t))^{p-1}}.$ \ \par 
for  $p>1,$  we get the estimate:\ \par 
\quad \quad \quad \quad \quad  $\displaystyle \ {\left\Vert{S(z,\cdot )}\right\Vert}_{p}^{p}\lesssim
 \sigma _{2n-1}(B(x,t))^{1-p}\Rightarrow {\left\Vert{S(z,\cdot
 )}\right\Vert}_{p}\lesssim \frac{1}{\sigma _{2n-1}(B(x,\delta
 ))^{1/p'}}.$   $\blacksquare $ \ \par 
\ \par 
\quad Now let  $\displaystyle K_{\Omega }(z,w)$  be the Bergman kernel
 of  $\Omega ,$  i.e. the kernel associated to the orthogonal
 projection  $L^{2}(\Omega )\rightarrow A^{2}(\Omega ),$  where
  $A^{2}$  is the Bergman space of square summable holomorphic
 functions in  $\Omega .$ \ \par 
We have a lower bound(~\cite{McNeal94}, theorem 3.4):\ \par 
\quad \quad \quad \quad \quad \quad \quad 	\begin{equation}  K_{\Omega }(a,a)\gtrsim \prod_{j=1}^{n}{\tau
 _{j}(a,\ \delta )^{-2}}\simeq \frac{1}{\delta \sigma _{2n-1}(B(\alpha
 ,\delta ))},\label{aspcGlo710}\end{equation}\ \par 
here with  $\delta =\left\vert{r(a)}\right\vert $  and  $a$ 
 in a neighbourhood  ${\mathcal{V}}_{p}$  of the point  $p\in
 \partial \Omega $  and  $\alpha =\pi (a).$  We also have an
 upper bound(~\cite{McNeal94}, theorem 5.2): \ \par 
\quad \quad \quad \quad \quad \quad \quad 	\begin{equation}  K_{\Omega }(a,z)\lesssim \prod_{j=1}^{n}{\tau
 _{j}(a,\ \delta )^{-2}}\simeq \frac{1}{\sigma _{2n}(T(z,a))},\label{aspcGlo711}\end{equation}\
 \par 
always in a neighbourhood of uniform size of  $p\in \partial
 \Omega ,$  and here with\ \par 
\quad \quad \quad $\delta =\left\vert{r(a)}\right\vert +\left\vert{r(z)}\right\vert
 +\rho (\pi (a),\ \pi (z))=\rho ^{*}(a,z).$ \ \par 
\quad So, with  $\alpha \in \partial \Omega $  fixed,  $\pi (a)=\alpha
 ,$  and  ${\mathcal{V}}$  a neighbourhood of  $\alpha $  valid
 for these two estimates, we have\ \par 

\begin{Lmm}
~\label{subPrinciple620}We have, with  $\alpha =\pi (a),\ \delta
 =\left\vert{r(a)}\right\vert ,$ \par 
\quad \quad \quad $\displaystyle \ {\left\Vert{K_{\Omega }(a,\cdot )}\right\Vert}_{H^{p}}^{p}\lesssim
 \frac{1}{\delta ^{p}\sigma _{2n-1}(B(\alpha ,\ \delta ))^{p-1}}$ \par 
and\par 
\quad \quad \quad \quad 	\begin{equation}  \ {\left\Vert{S(a,\cdot )}\right\Vert}_{H^{p}(\Omega
 )}\geq \frac{1}{\sigma _{2n-1}(B(\alpha ,\delta ))^{1/p'}}.\label{strongPC934}\end{equation}
\end{Lmm}

      Proof of the lemma.\ \par 
From the inequality~(\ref{aspcGlo711}) and using the annuli\ \par 
\quad \quad \quad \quad \quad  $C_{k}:=B(x,\ 2^{k}\delta )\backslash B(x,\ 2^{k-1}\delta ),\
 k\geq 1,\ C_{0}=B(x,\delta ),$ \ \par 
we already used in the proof of lemma~\ref{subPrinciple40} we
 get, exactly as before, with  $x=\pi (z),$  and  $\displaystyle
 \delta =2\left\vert{r(z)}\right\vert ,$ \ \par 
\quad \quad \quad \quad \quad  $\displaystyle \ \left\vert{K_{\Omega }(z,w)}\right\vert \lesssim
 \frac{{\11}_{B(x,\delta )}(w)}{\delta \sigma _{2n-1}(B(x,\delta
 /2))}+\sum_{k\geq 1}{\frac{{\11}_{C_{k}}(w)}{2^{k-1}\delta \sigma
 _{2n-1}(B(x,\ 2^{k-1}\delta ))}}.$ \ \par 
Hence\!\!\!\! , with  $\alpha =\pi (a),$ \ \par 
\quad \quad \quad $\displaystyle \ \int_{{\mathcal{V}}\cap \lbrace r(z)=-\delta
 /2\rbrace }{\left\vert{K_{\Omega }(a,z)}\right\vert ^{p}\,d\sigma
 (z)}\lesssim $ \ \par 
\quad \quad \quad \quad \quad \quad \quad  $\displaystyle \lesssim \frac{\sigma _{2n-1}(B(\alpha ,\delta
 ))}{\delta ^{p}\sigma _{2n-1}(B(\alpha ,\delta /2))^{p}}+\sum_{k\geq
 1}{\frac{\sigma _{2n-1}(C_{k})}{2^{p(k-1)}\delta ^{p}\sigma
 _{2n-1}(B(\alpha ,\ 2^{k}\delta )^{p}}}.$ \ \par 
Hence, again as before,\ \par 
\quad \quad \quad $\displaystyle \ \int_{{\mathcal{V}}\cap \lbrace r(z)=-\delta
 /2\rbrace }{\left\vert{K_{\Omega }(a,z)}\right\vert ^{p}\,d\sigma
 _{2n-1}(z)}\leq \frac{1}{\delta ^{p}\sigma _{2n-1}(B(\alpha
 ,\delta ))^{p-1}}.$ \ \par 
Ouside of  ${\mathcal{V}},\ K_{\Omega }(a,\ \cdot )$  is bounded
 because by~\cite{McNealStein94}, p 178 :\ \par 
\quad \quad \quad $\displaystyle \ \left\vert{K_{\Omega }(a,\ z)}\right\vert \lesssim
 \frac{1}{\sigma _{2n}(T(a,z))},$ \ \par 
and if  $z\notin {\mathcal{V}}$  then  $\displaystyle 1\lesssim
 \sigma _{2n}(T(a,z))$  uniformly in  $a\in \Omega .$ \ \par 
Hence\ \par 
\quad \quad \quad $\displaystyle \ {\left\Vert{K_{\Omega }(a,\cdot )}\right\Vert}_{H^{p}}^{p}=\int_{{\mathcal{U}}\cap
 \lbrace r(z)=-\delta /2\rbrace }{\left\vert{K_{\Omega }(a,z)}\right\vert
 ^{p}\,d\sigma (z)}+\int_{(\partial \Omega \backslash {\mathcal{U}})\cap
 \lbrace r(z)=-\delta /2\rbrace }{\left\vert{K_{\Omega }(a,z)}\right\vert
 ^{p}\,d\sigma (z)}\leq $ \ \par 
\quad \quad \quad \quad \quad \quad \quad \quad \quad \quad $\displaystyle \lesssim \frac{1}{\delta ^{p}\sigma _{2n-1}(B(\alpha
 ,\ \delta ))^{p-1}}+c\lesssim \frac{1}{\delta ^{p}\sigma _{2n-1}(B(\alpha
 ,\ \delta ))^{p-1}},$ \ \par 
because  $c$  is uniformly bounded, hence\ \par 
\quad \quad \quad $\displaystyle \ {\left\Vert{K_{\Omega }(a,\cdot )}\right\Vert}_{H^{p}}^{p}\lesssim
 \frac{1}{\delta ^{p}\sigma _{2n-1}(B(\alpha ,\ \delta ))^{p-1}},$ \ \par 
which proves the first part of the lemma.\ \par 
\quad \quad 	Notice that even if  $K_{\Omega }$  is linked to Bergman space,
 we have an estimate of its {\sl Hardy}  $H^{p}(\Omega )$  norm.\ \par 
\quad Using the lower bound~(\ref{aspcGlo710}) and the previous inequality,
 we get\ \par 
\quad \quad \quad \quad \quad  $\displaystyle \ \frac{K_{\Omega }(a,a)}{{\left\Vert{K_{\Omega
 }(a,\cdot )}\right\Vert}_{H^{p'}(\Omega )}}\geq \frac{1}{\delta
 \sigma _{2n-1}(B(\alpha ,\delta ))}{\times}\delta \sigma _{2n-1}(B(\alpha
 ,\delta ))^{1/p}\geq \frac{1}{\sigma _{2n-1}(B(\alpha ,\delta
 ))^{1/p'}}.$ \ \par 
\quad Hence, because\ \par 
\quad \quad \quad \quad \quad  $\displaystyle \ {\left\Vert{S(a,\ \cdot )}\right\Vert}_{H^{p}}=\sup
  \lbrace \left\vert{f(a)}\right\vert =\left\vert{{\left\langle{f,\
 S(a,\ \cdot )}\right\rangle}}\right\vert ::f\in H^{p'}(\Omega
 ),\ {\left\Vert{f}\right\Vert}_{p'}=1\rbrace $ \ \par 
we get\ \par 
\quad \quad \quad \quad 	\begin{equation}  \ {\left\Vert{S(a,\cdot )}\right\Vert}_{H^{p}(\Omega
 )}\geq \frac{1}{\sigma _{2n-1}(B(\alpha ,\delta ))^{1/p'}}\label{strongPC934}\end{equation}\
 \par 
by the choice of  $\displaystyle f(z):=\frac{K_{\Omega }(a,z)}{{\left\Vert{K_{\Omega
 }(a,\cdot )}\right\Vert}_{H^{p'}}}$   $\blacksquare $ \ \par 
\ \par 
\quad Recall that the Poisson Szeg\"o kernel is\ \par 
\quad \quad \quad $\displaystyle \forall z\in \Omega ,\ y\in \partial \Omega ,\
 P(z,y):=\frac{\left\vert{S(z,y)}\right\vert ^{2}}{{\left\Vert{S(z,\
 \cdot )}\right\Vert}_{H^{2}}^{2}}=\frac{\left\vert{S(z,y)}\right\vert
 ^{2}}{S(z,z)}.$ \ \par 
We have that this kernel reproduces the holomorphic functions:\ \par 
\quad \quad \quad $\displaystyle \forall f\in A(\Omega ),\ \int_{\partial \Omega
 }{f(y)P(z,y)\,d\sigma (y)}=\frac{1}{S(z,z)}{\left\langle{fS(z,\
 \cdot ),\ S(z,\ \cdot )}\right\rangle}=f(z),$ \ \par 
because of the reproducing property of the Szeg\"o kernel. The
 kernel  $P(z,y)$  is positive and has a  $\displaystyle L^{1}(\partial
 \Omega ,\,d\sigma _{2n-1})$  norm equals to one.\ \par 
\quad Also recall the Hardy-Littlewood kernel\ \par 
\quad \quad \quad $\displaystyle \forall x,\ y\in \partial \Omega ,\ P^{0}_{t}(x,\
 y):=\frac{1}{\sigma _{2n-1}(B(x,\ t))}{\11}_{B(x,t)}(y).$ \ \par 
We have\ \par 

\begin{Lmm}
The Poisson Szeg\"o kernel  $P(z,y)$  is dominated by the Hardy-Littlewood
 one : this means precisely that we have, with  $x=\pi (z),\
 t=\left\vert{r(z)}\right\vert ,$ \par 
\quad \quad \quad $\displaystyle \forall z\in \Omega ,\ \forall y\in \partial \Omega
 ,\ P(z,y)\lesssim P^{0}_{2t}(x,\ y)+\sum_{k\in {\mathbb{N}}}{\frac{1}{2^{k+1}}P^{0}_{2^{k+1}t}(x,\
 y)}.$ 
\end{Lmm}
\quad Proof:\ \par 
\quad Using~(\ref{subP31}) we get, still with  $x=\pi (z),\ t=\left\vert{r(z)}\right\vert
 ,$ \ \par 
\quad \quad \quad \quad \quad  $\displaystyle \ \left\vert{S(z,y)}\right\vert ^{2}\lesssim
 \frac{{\11}_{B(x,2t)(y)}}{\sigma _{2n-1}(B(x,t))^{2}}+\sum_{k\geq
 1}{\frac{{\11}_{C_{k}}(y)}{\sigma _{2n-1}(B(x,\ 2^{k}t))^{2}}}.$ \ \par 
Because  $C_{k}\subset B(x,\ 2^{k+1}t)$  we have  $\displaystyle
 {\11}_{C_{k}}\leq {\11}_{B_{k}},$  with  $\displaystyle B_{k}:=B(x,\
 2^{k+1}t),$  hence setting\ \par 
\quad \quad \quad $\displaystyle \forall x,\ y\in \partial \Omega ,\ P^{0}_{t}(x,\
 y):=\frac{1}{\sigma _{2n-1}(B(x,\ t))}{\11}_{B(x,t)}(y),$ \ \par 
the Hardy-Littlewood kernel, we have\ \par 
\quad \quad \quad $\displaystyle \ \left\vert{S(z,y)}\right\vert ^{2}\lesssim \frac{1}{\sigma
 _{2n-1}(B(x,\ t))}P^{0}_{2t}(x,\ y)+\sum_{k\geq 1}{\frac{1}{\sigma
 _{2n-1}(B(x,\ 2^{k}t))}P^{0}_{2^{k}t}(x,\ y)}.$ \ \par 
But by~(\ref{strongPC933}) we have  $\displaystyle \sigma _{2n-1}(B(x,\
 2^{k}t))\gtrsim 2^{k}\sigma _{2n-1}(B(x,\ t)),$  hence\ \par 
\quad \quad \quad \[\displaystyle \ \left\vert{S(z,y)}\right\vert ^{2}\lesssim
 \frac{1}{\sigma _{2n-1}(B(x,\ t))}(P^{0}_{2t}(x,\ y)+\sum_{k\geq
 1}{\frac{1}{2^{k}}P^{0}_{2^{k}t}(x,\ y)}).\] \ \par 
\quad By~(\ref{strongPC934}) we have, with  $p=2,$ \ \par 
\quad \quad \quad \[\displaystyle \ {\left\Vert{S(a,\cdot )}\right\Vert}_{H^{2}}^{2}=S(a,a)\geq
 \frac{1}{\sigma _{2n-1}(B(\alpha ,\ 2t))},\] \ \par 
Hence we get for the Poisson-Szeg\"o kernel, still with  $t=\left\vert{r(z)}\right\vert
 ,\ x=\pi (z),$ \ \par 
\quad \quad \quad $\displaystyle P(z,y)\lesssim P^{0}_{2t}(x,\ y)+\sum_{k\geq 1}{\frac{1}{2^{k}}P^{0}_{2^{k}t}(x,\
 y)}.$   $\blacksquare $ \ \par 
\ \par 
\quad Combining the previous results we have\ \par 

\begin{Thrm}
~\label{subPrinciple42}Let  $\Omega $  be a convex domain of
 finite type in  ${\mathbb{C}}^{n},$  then, with  $S(z,\ y)$
  its Szeg\"o kernel we have, setting\par 
$x=\pi (z),\ t=\left\vert{r(z)}\right\vert ,\ C_{0}:=B(x,t),\
 \forall k\geq 1,\ C_{k}:=B(x,\ 2^{k}t)\backslash B(x,\ 2^{k-1}t),$ \par 
$\displaystyle \bullet \ \forall z\in \Omega ,\ \forall y\in
 \partial \Omega ,\ \left\vert{S(z,y)}\right\vert \lesssim \frac{1}{\sigma
 _{2n-1}(B(x,t))}{\11}_{B(x,\ 2t)}(y)+\sum_{k\geq 1}{\frac{1}{\sigma
 _{2n-1}(B(x,\ 2^{k}t))}{\11}_{C_{k}}(y)};$ \par 
 $\bullet $   $\displaystyle \ {\left\Vert{S(z,\cdot )}\right\Vert}_{H^{p}(\Omega
 )}\simeq \frac{1}{\sigma _{2n-1}(B(x,t))^{1/p'}}.$ \par 
And with  $P^{0}_{t}(x,y)$  the Hardy-Littlewood kernel and 
 $P(z,y)$  the Poisson-Szeg\"o kernel\par 
$\displaystyle \bullet \ P(z,y)\lesssim P^{0}_{2t}(x,\ y)+\sum_{k\geq
 1}{\frac{1}{2^{k}}P^{0}_{2^{k}t}(x,\ y)}.$ 
\end{Thrm}
\vfill\eject\ \par 

\section{Carleson measures.~\label{8_CarlesonMeas30}}

\subsection{Harmonic analysis.}
\quad \quad 	We start by a "copy and paste" from~\cite{AmarBonami}, where
 we introduced the notion of Carleson measures of order  $\alpha .$ \ \par 
\quad Let  $(X,\ \rho ,\ \,d\sigma )$  be a homogeneous type space~\cite{CoifWeiss71}.
 Denote  $B(x,t):=\lbrace y\in X::\rho (x,y)<t\rbrace $  the
 pseudo-ball centered at  $x$  and of radius  $t>0.$ \ \par 
We define the Carleson windows (or "tents") on  ${\mathbb{R}}^{+}{\times}X$
  the following way: let  $A$  be an open set in  $X,$  then\ \par 
\quad \quad \quad $W(A):=\lbrace (t,x)\in {\mathbb{R}}^{+}{\times}X::B(x,t)\subset
 A\rbrace .$ \ \par 
\quad We set  $W(A)$  instead of  $T(A)$  to differentiate notations
 from the case of the convex domains of finite type we seen in
 the previous section.\ \par 

\begin{Dfnt}
We say that the mesure  $\lambda $  on  ${\mathbb{R}}^{+}{\times}X$
  is a {\bf homogeneous geometric Carleson measure} of  order
  $\alpha $  if, for any open set  $A\subset X,$ \par 
\quad \quad \quad $\ \left\vert{\lambda }\right\vert (W(A))\leq C\sigma (A)^{\alpha }.$ \par 
The usual homogeneous geometric Carleson measures are those with  $\alpha =1.$ 
\end{Dfnt}
We shall abbreviate homogeneous geometric Carleson measure by
 h.g. Carleson measure.\ \par 
\quad In the case  $\alpha =1$  it is enough to test on the sets  $A=B(x,t)$
  because the pseudo-balls generate all open sets in a homogeneous
 type space~\cite{CoifWeiss71}. In the case  $\alpha =1$  we
 shall speak simply of h.g. Carleson measure.\ \par 
\quad \quad 	The action of a kernel  $\displaystyle P_{t}$  on a function
  $f$  will be denoted  $\displaystyle P_{t}f,$  precisely\ \par 
\quad \quad \quad \quad \quad  $\displaystyle P_{t}f(y):=\int_{X}{P_{t}(x,y)f(x)d\sigma (x)}.$ \ \par 
\quad Now we have the abstract Carleson embedding theorem.\ \par 

\begin{Thrm}
~\label{subPrinciple69}If the kernel  $P_{t}$  is dominated by
 the Hardy-Littlewood kernel, and if  $\lambda $  is a h.g. Carleson
 measure on  $\displaystyle {\mathbb{R}}^{+}{\times}X,$  we have\par 
\quad \quad \quad $\displaystyle \forall f\in L^{p}(X,\ \sigma ),\ \int_{X{\times}{\mathbb{R}}^{+}}{\left\vert{P_{t}f(x)}\right\vert
 ^{p}\,d\left\vert{\lambda }\right\vert (x,t)}\lesssim {\left\Vert{f}\right\Vert}_{L^{p}(\sigma
 )}^{p}.$ 
\end{Thrm}
\quad Proof.\ \par 
This is quite well known and implicitly contained in H\"ormander~\cite{HormPSH67},
 Theorem 2.4. But I shall give a proof taken from~\cite{AmarBonami}
 where the same notations as here are used and which uses h.g.
 Carleson measures of order  $\alpha .$ \ \par 
\quad Let  $V^{0}$  the space of finite measure on  $\displaystyle
 {\mathbb{R}}^{+}{\times}X,\ V^{1}$  the space of h.g. Carleson
 ones and, with  $\displaystyle \ \alpha :=1-1/p,\ W^{\alpha
 }:=(V^{0},\ V^{1})_{(\alpha ,p)}$  the intermediate class by
 the real interpolating method. We proved in~\cite{AmarBonami},
 Proposition 1, p 30, that\ \par 
\quad \quad \quad \quad \quad  \begin{equation}  w\in W^{\alpha }\iff \exists \lambda \in V^{1},\
 \exists h\in L^{p}(\lambda )::\,dw=h\,d\lambda .\label{convexInterpol324}\end{equation}\
 \par 
Moreover the norm of  $w$  in  $W^{\alpha }$  is equivalent to
 the norm of  $h$  in  $L^{p}(\lambda ).$ \ \par 
Because  $P_{t},$  being dominated by the Hardy-Littlewood kernel,
 verifies the  $(H1)$  hypothesis of~\cite{AmarBonami}, by theorem
 2, p 27, we have that\ \par 
\begin{equation}  \forall w\in W^{\alpha },\ \forall g\in L^{p'}(\sigma
 ),\ \int_{{\mathbb{R}}^{+}{\times}X}{\left\vert{P_{t}g(x)}\right\vert
 \,d\left\vert{w}\right\vert (t,x)}\leq C_{w}{\left\Vert{g}\right\Vert}_{L^{p'}(\sigma
 )}.\label{convexInterpol323}\end{equation}\ \par 
\quad Now let  $\lambda $  be a geometric Carleson measure and  $\displaystyle
 f\in L^{p'}(X,\ \sigma )\ ;$  we want to prove that  $\displaystyle
 P_{t}f(x)\in L^{p'}(\lambda ).$  Let  $h\in L^{p}(\lambda )$
  and set  $\,dw:=h\,d\lambda $  then  $w\in W^{\alpha }$  by~(\ref{convexInterpol324}).
 We have by~(\ref{convexInterpol323})\ \par 
\quad \quad \quad \quad \quad \quad \quad  $\displaystyle \ \int_{{\mathbb{R}}^{+}{\times}X}{\left\vert{P_{t}f(x)}\right\vert
 \,d\left\vert{w}\right\vert (t,x)}=\int_{{\mathbb{R}}^{+}{\times}X}{\left\vert{P_{t}f(x)}\right\vert
 \left\vert{h}\right\vert \,d\left\vert{\lambda }\right\vert
 (t,x)}\leq C{\left\Vert{h}\right\Vert}_{L^{p}(\lambda )}{\left\Vert{f}\right\Vert}_{L^{p'}(\sigma
 )}\ ;$ \ \par 
but this being true for all functions  $h$  in  $L^{p}(\lambda
 ),$  we have that  $\displaystyle P_{t}f(x)\in L^{p'}(\lambda
 )$  and the theorem is proved by exchanging  $p'$  and  $p.$
   $\blacksquare $ \ \par 

\subsection{Carleson measures in convex domain of finite type.}
\quad Now to define the geometric Carleson measures in our domains
 we have 2 possibilities for a positive Borel measure on  $\Omega \ :$  \ \par 
\quad \quad \quad \quad \quad  $\displaystyle \bullet \ \exists C>0::\forall a\in \Omega ,\
 \epsilon :=2\left\vert{r(a)}\right\vert ,\ \lambda (T_{a}(\epsilon
 ))\leq C\sigma (\partial \Omega \cap P_{\epsilon }(a)),$ \ \par 
with  $P_{\epsilon }(a)\in {\mathcal{P}}$  is the family defined
 in~(\ref{subPrinciple68}).\ \par 
\quad \quad \quad \quad \quad  $\displaystyle \bullet \ \exists C>0::\forall a\in \Omega ,\
 \alpha =\pi (a),\ \lambda (\Omega \cap W(B(\alpha ,\ \left\vert{r(a)}\right\vert
 ))\leq C\sigma (B(\alpha ,\ \left\vert{r(a)}\right\vert )),$ \ \par 
where  $\displaystyle B(\alpha ,\ \left\vert{r(a)}\right\vert
 )$  is the pseudo-ball on  $\partial \Omega $  of center  $\alpha
 $  and radius  $\ \left\vert{r(a)}\right\vert ,$  and  $\displaystyle
 W(B(\alpha ,\ \left\vert{r(a)}\right\vert ))$  is the Carleson
 window defined in the previous subsection. For this section
 we set  $\displaystyle \sigma =\sigma _{2n-1}.$ \ \par 
\ \par 
\quad We shall show that they are equivalent. We have that\ \par 
\quad \quad \quad \quad $\forall a\in {\mathcal{U}}\cap \Omega ,\ \epsilon :=\left\vert{r(a)}\right\vert
 ,\ B(\alpha ,\epsilon ):=\partial \Omega \cap P_{\epsilon }(\alpha ),$ \ \par 
by definition of the family  ${\mathcal{P}}.$  Then we want to show the\ \par 

\begin{Lmm}
There is a constant  $\gamma ,$  independent of  $a,$  such that\par 
\quad \quad \quad $W(B(\alpha ,\ \left\vert{r(a)}\right\vert ))\subset T_{a}(\
 \gamma \left\vert{r(a)}\right\vert ).$ 
\end{Lmm}
\quad Proof.\ \par 
We have, by definition of the Carleson window :\ \par 
\quad \quad \quad $z\in W(B(\alpha ,\ \left\vert{r(a)}\right\vert ))\iff B(x,\
 \left\vert{r(z)}\right\vert )\subset B(\alpha ,\ \left\vert{r(a)}\right\vert
 )),$ \ \par 
where  $x=\pi (z).$  This implies, because  $\partial \Omega
 $  is a space of homogeneous type, that we have  $\ \left\vert{r(z)}\right\vert
 \leq c\left\vert{r(a)}\right\vert ,$  with a uniform constant
  $c\geq 1.$ \ \par 
\quad But with  $\delta :=\left\vert{r(z)}\right\vert ,\ P_{\delta
 }(z)\cap B(x,\ \delta )\neq \emptyset $  hence with  $\epsilon
 =\left\vert{r(a)}\right\vert ,$ \ \par 
\quad \quad \quad \quad $P_{\delta }(z)\cap B(\alpha ,\ \epsilon )\neq \emptyset \Rightarrow
 P_{\delta }(z)\cap P_{\epsilon }(\alpha )\neq \emptyset .$ \ \par 
Let  $s=\max  (\delta ,\ \epsilon )$  then\ \par 
\quad \quad \quad \quad $P_{\delta }(z)\cap P_{\epsilon }(\alpha )\neq \emptyset \Rightarrow
 P_{s}(z)\cap P_{s}(\alpha )\neq \emptyset \Rightarrow P_{s}(z)\subset
 C_{3}P_{s}(\alpha )$ \ \par 
by~(\ref{suPr611}).\ \par 
\quad But if  $s=\delta $  then  $s\leq c\epsilon $  and if  $s=\epsilon
 $  then again  $s\leq c\epsilon $  because  $c\geq 1$  ; so
 in any case  $s\leq c\epsilon $  and this implies\ \par 
\quad \quad \quad $P_{\delta }(z)\subset P_{s}(z)\subset C_{3}P_{s}(\alpha )\subset
 C_{4}P_{\epsilon }(\alpha )\subset P_{\gamma \epsilon }(\alpha ),$ \ \par 
by~(\ref{convexInterpol526}) with  $t=C_{4},\ \gamma =C_{t}.$ \ \par 
And again because  $P_{\epsilon }(a)\cap B(\alpha ,\ \epsilon
 )\neq \emptyset ,$  we get  $P_{\gamma \epsilon }(\alpha )\subset
 P_{\gamma \epsilon }(a)$  by~(\ref{suPr611}) and~(\ref{convexInterpol526})
 ; and finally  $P_{\delta }(z)\subset P_{\gamma \epsilon }(a).$
  Cutting with  $\Omega $  we get\ \par 
\quad \quad \quad $z\in P_{\delta }(z)\cap \Omega \subset P_{\gamma \epsilon }(a)\cap
 \Omega =T_{a}(\ \gamma \left\vert{r(a)}\right\vert ).\ \blacksquare $ \ \par 
\quad \quad 	We shall use the following definition for geometric Carleson
 measure in a convex domain of finite type to continue with the
 same notations.\ \par 

\begin{Dfnt}
Let  $\lambda $  be a positive Borel measure on the bounded convex
 domain of finite type  $\Omega .$  We shall say that  $\lambda
 $  is a {\bf geometric Carleson measure} in  $\Omega $  if :\par 
\quad \quad \quad $\displaystyle \exists C>0::\forall a\in \Omega ,\ \epsilon =2\left\vert{r(a)}\right\vert
 ,\ \lambda (T_{a}(\epsilon ))\leq C\sigma (\partial \Omega \cap
 P_{\epsilon }(a)).$ 
\end{Dfnt}

\subsection{Carleson embedding.}
\quad \quad 	We are in position to prove a Carleson embedding theorem for
 convex domains of finite type.\ \par 
\quad To prove it we shall need the lemma :\ \par 

\begin{Lmm}
~\label{ConvInt2}Let  $a\in \Omega ,\ \alpha =\pi (a),\ \delta
 =\left\vert{r(a)}\right\vert $  ; there is a uniform constant
  $\gamma >0$  such that\par 
\quad \quad \quad \quad \quad \quad \quad 	\begin{equation}  \forall z\in \gamma P_{\delta }(a),\ \left\vert{K_{\Omega
 }(z,\ a)}\right\vert \geq \frac{c}{\delta \sigma _{2n-1}(B(\alpha
 ,\delta ))}.\label{subPrinciple619}\end{equation}
\end{Lmm}
\quad Proof.\ \par 
We have the lower bound~(\ref{aspcGlo710}) of the Bergman kernel\ \par 
\quad \quad \quad \quad $\displaystyle K_{\Omega }(a,a)\gtrsim \prod_{j=1}^{n}{\tau _{j}(a,\
 \delta )^{-2}}\simeq \frac{1}{\delta \sigma _{2n-1}(B(\alpha
 ,\delta ))},$ \ \par 
the last equivalence by equations~(\ref{aspcGlo710}) and~(\ref{aspcGlo711})
 and a upper bound of its derivatives(~\cite{McNeal94}, theorem
 5.2 and ~\cite{McNeal02})\ \par 
\quad \quad \quad \quad \quad \quad 	\begin{equation}  \ \left\vert{\partial _{z}^{\lambda }\bar
 \partial _{a}^{\mu }K_{\Omega }(z,a)}\right\vert \leq C_{\lambda
 \mu }\prod_{j=1}^{n}{\tau _{j}(a,\ \beta )^{-2-\lambda _{j}-\mu
 _{j}}},\label{ConvInt3}\end{equation}\ \par 
with  $\beta =\rho ^{*}(a,\ z).$ \ \par 
\quad \quad 	Set for  $\displaystyle t\in \lbrack 0,1\rbrack ,\ f(t):=K_{\Omega
 }(a+t(z-a),a),$  then  $f$  being complex valued, we have  $\displaystyle
 f(t)=(f_{1}+if_{2}).$ \ \par 
Apply the mean value theorem :\ \par 
\quad \quad \quad \quad \quad  $\displaystyle \exists t_{1},t_{2}\in \lbrack 0,1\rbrack ::f(1)-f(0)=(f_{1}'(t_{1})+if_{2}'(t_{2}))\Rightarrow
 \left\vert{f(1)-f(0)}\right\vert \leq 2\sup _{t\in \lbrack 0,1\rbrack
 }\left\vert{f'(t)}\right\vert .$ \ \par 
Hence\ \par 
\quad \quad \quad \quad \quad  $\displaystyle \ \left\vert{K_{\Omega }(z,a)-K_{\Omega }(a,a)}\right\vert
 \leq 2\sup _{t\in \lbrack 0,1\rbrack }\left\vert{\sum_{j=1}^{n}{(z_{j}-a_{j})\frac{\partial
 K_{\Omega }(a+t(z-a),\ a)}{\partial \zeta _{j}}}}\right\vert
 \lesssim \sum_{j=1}^{n}{\frac{\left\vert{z_{j}-a_{j}}\right\vert
 }{\tau _{j}(a,\ \beta )}}\prod_{k=1}^{n}{\tau _{k}(a,\ \beta )^{-2}},$ \ \par 
by inequality~\ref{ConvInt3} ; so\ \par 
\quad \quad \quad \quad \quad  $\ \left\vert{K_{\Omega }(z,a)-K_{\Omega }(a,a)}\right\vert
 \lesssim \frac{1}{\beta \sigma _{2n-1}(B(\alpha ,\beta ))}\sum_{j=1}^{n}{\frac{\left\vert{z_{j}-a_{j}}\right\vert
 }{\tau _{j}(a,\ \beta )}},$ \ \par 
by equations~(\ref{aspcGlo710}) and~(\ref{aspcGlo711}).\ \par 
Now choose  $z$  such that  $\ \left\vert{z_{j}-a_{j}}\right\vert
 \leq \gamma \tau _{j}(a,\ \delta )\Rightarrow \beta \lesssim
 \delta $  and the homogeneous nature of  $\Omega $  gives that
  $\tau _{j}(a,\ \beta )\simeq \tau _{j}(a,\ \delta )$  hence\ \par 
\quad \quad \quad $\displaystyle \ \left\vert{K_{\Omega }(z,a)-K_{\Omega }(a,a)}\right\vert
 \lesssim \frac{1}{\beta \sigma _{2}n-1(B(\alpha ,\delta ))}\sum_{j=1}^{n}{\frac{\left\vert{z_{j}-a_{j}}\right\vert
 }{\tau _{j}(a,\ \delta )}}\lesssim \frac{n\gamma }{\delta \sigma
 _{2n-1}(B(\alpha ,\delta ))}.$ \ \par 
Take  $\gamma $  uniformly small enough to compensate the constant
 in the last inequality above to get\ \par 
\quad \quad \quad $\displaystyle \ \left\vert{K_{\Omega }(z,a)-K_{\Omega }(a,a)}\right\vert
 \leq \frac{1}{2}{\times}\frac{1}{\delta \sigma _{2n-1}(B(\alpha
 ,\delta ))},$ \ \par 
this means that, for  $z$  in the polydisc  $\gamma P_{\delta
 }(a),$  we have  $\displaystyle \ \left\vert{K_{\Omega }(z,\
 a)}\right\vert \geq \frac{c}{\delta \sigma _{2n-1}(B(\alpha
 ,\delta ))},$  the positive constants  $c,\ \gamma $  being
 uniform.  $\blacksquare $ \ \par 
\quad \quad 	We shall need the definition.\ \par 

\begin{Dfnt}
Let  $\lambda $  be a positive Borel measure on the domain  $\Omega
 $  and  $p\geq 1.$  We shall say that  $\lambda $  is a  $p$
  {\bf Carleson measure} in  $\Omega $  if :\par 
\quad \quad \quad $\displaystyle \exists C_{p}>0,\ \forall f\in H^{p}(\Omega ),\
 \int_{\Omega }{\left\vert{f}\right\vert ^{p}\,d\lambda }\leq
 C_{p}^{p}{\left\Vert{f}\right\Vert}_{H^{p}}^{p}.$ \par 
This means that we have a continuous embedding of  $H^{p}(\Omega
 )$  in  $L^{p}(\lambda ).$ 
\end{Dfnt}
\quad \quad 	Now we have.\ \par 

\begin{Thrm}
~\label{aspcGlo50}If the measure  $\lambda $  is a geometric
 Carleson measure we have\par 
\quad \quad \quad $\displaystyle \forall p>1,\ \exists C_{p}>0,\ \forall f\in H^{p}(\Omega
 ),\ \int_{\Omega }{\left\vert{f}\right\vert ^{p}\,d\lambda }\leq
 C_{p}^{p}{\left\Vert{f}\right\Vert}_{H^{p}}^{p}.$ \par 
Conversely if the positive measure  $\lambda $  is  $p$  Carleson
 for a  $p\in \lbrack 1,\ \infty \lbrack ,$  then it is a geometric
 Carleson measure, hence it is  $q$ -Carleson for any  $q\in
 \rbrack 1,\ \infty \lbrack .$ 
\end{Thrm}
\quad Proof.\ \par 
We apply theorem~\ref{subPrinciple69} to the Poisson-Szeg\"o
 kernel  $P(z,y)$  which is dominated by the Hardy-Littlewood
 kernel. Because a function in  $A(\Omega ),$  the algebra of
 holomorphic function in  $\Omega $  continuous up to  $\partial
 \Omega ,$  is reproduced by  $P(z,y)$  and because this algebra
 is dense in  $H^{p}(\Omega ),$  the first part of the theorem is proved.\ \par 
\quad Suppose now that  $\lambda $  is  $p$ -Carleson for a  $p\in
 \lbrack 1,\ \infty \lbrack ,$  then we have\ \par 
\quad \quad \quad $\displaystyle \exists C>0,\ \forall a\in \Omega ,\ \int_{\Omega
 }{\left\vert{K_{\Omega }(z,\ a)}\right\vert ^{p}\,d\lambda (z)}\leq
 C{\left\Vert{K_{\Omega }(\cdot ,\ a)}\right\Vert}_{H^{p}}^{p},$ \ \par 
with  $\displaystyle K_{\Omega }(z,\ a)$  the Bergman kernel
 at  $a.$  Using the inequality~(\ref{subPrinciple619}) of the
 lemma, we get\ \par 
\quad \quad \quad $\displaystyle \forall a\in \Omega ,\ \int_{\Omega \cap \gamma
 P_{\delta }(a)}{{\left({\frac{1}{\delta \sigma (B(\alpha ,\delta
 ))}}\right)}^{p}\,d\lambda (z)}\leq \int_{\Omega }{\left\vert{K_{\Omega
 }(z,\ a)}\right\vert ^{p}\,d\lambda (z)}\leq C{\left\Vert{K_{\Omega
 }(\cdot ,\ a)}\right\Vert}_{H^{p}}^{p},$ \ \par 
hence\ \par 
\quad \quad \quad $\displaystyle \forall a\in \Omega ,\ {\left({\frac{1}{\delta
 \sigma (B(\alpha ,\delta ))}}\right)}^{p}\lambda (\Omega \cap
 \gamma P_{\delta }(a))\leq C{\left\Vert{K_{\Omega }(\cdot ,\
 a)}\right\Vert}_{H^{p}}^{p}.$ \ \par 
We can use the estimate of  $\displaystyle \ {\left\Vert{K_{\Omega
 }(\cdot ,\ a)}\right\Vert}_{H^{p}}$  done in lemma~(\ref{subPrinciple620})\
 \par 
\quad \quad \quad $\displaystyle \ {\left\Vert{K_{\Omega }(\cdot ,\ a)}\right\Vert}_{H^{p}}^{p}\lesssim
 \frac{1}{\delta ^{p}\sigma (B(\alpha ,\ \delta ))^{p-1}},$ \ \par 
to get\ \par 
\quad \quad \quad $\displaystyle \forall a\in \Omega ,\ {\left({\frac{1}{\delta
 \sigma (B(\alpha ,\delta ))}}\right)}^{p}\lambda (\Omega \cap
 \gamma P_{\delta }(a))\leq C\frac{1}{\delta ^{p}\sigma (B(\alpha
 ,\ \delta ))^{p-1}},$ \ \par 
hence\ \par 
\quad \quad \quad $\displaystyle \forall a\in \Omega ,\ \lambda (\Omega \cap \gamma
 P_{\delta }(a))\leq C\sigma (B(\alpha ,\ \delta )).$ \ \par 
Still by homogeneity we have  $\displaystyle \gamma P_{\delta
 }(a)\supset P_{c\delta }(\alpha )$  and\ \par 
\quad \quad \quad \quad $\displaystyle B(\alpha ,\ \delta )\subset CB(\alpha ,\ c\delta
 )\Rightarrow \sigma (B(\alpha ,\ \delta ))\leq C'\sigma (B(\alpha
 ,\ c\delta )),$ \ \par 
so\ \par 
\quad \quad \quad $\displaystyle \forall a\in \Omega ,\ \lambda (\Omega \cap P_{c\delta
 }(\alpha ))\leq CC'\sigma (B(\alpha ,\ c\delta )),$ \ \par 
and the measure  $\lambda $  is a geometric Carleson measure,
 hence it is a  $q$  Carleson measure by the first part of the
 theorem.  $\blacksquare $ \ \par 
\ \par 
\quad If  $\Omega $  is a convex domain of finite type, with the family
  ${\mathcal{P}}$  of polydiscs of McNeal, we define a related
 family  ${\mathcal{Q}}$  of polydiscs :\ \par 
\quad \quad \quad $\forall a\in {\mathcal{U}}\backslash \partial \Omega ,\ \forall
 t>0,\ \epsilon :=\left\vert{r(a)}\right\vert ,\ Q_{a}(t):=tP_{\epsilon
 }(a),$ \ \par 
where  $tP_{\epsilon }(a)$  si the dilated polydisc as defined
 in~(\ref{suPr65}).\ \par 

\begin{Lmm}
The family  ${\mathcal{Q}}:=\lbrace Q_{a}(t),\ t>0,\ a\in {\mathcal{U}}\rbrace
 $  is a good family of polydiscs in  $\Omega .$ 
\end{Lmm}
\quad \quad 	Proof.\ \par 
By~(\ref{suPr610})((3) of proposition 2.7 in~\cite{Hefer04}),
 we get that  $\exists \delta _{0}>0,$  such that\ \par 
\quad \quad \quad \quad \quad  $a\in \Omega \cap {\mathcal{U}}\Rightarrow \delta _{0}P_{\left\vert{r(a)}\right\vert
 }(a)\subset \Omega ,$ \ \par 
because  $d(a):=d(a,\Omega ^{c})\simeq \left\vert{r(a)}\right\vert
 ,$  the constants being independent of  $a\in \Omega ,$  we have with\ \par 
\quad \quad \quad \quad \quad  $\displaystyle Q_{a}(t):=tP_{\epsilon }(a),\ t=\delta _{0},\
 \epsilon =\left\vert{r(a)}\right\vert \simeq d(a),$ \ \par 
that\ \par 
\quad \quad \quad \quad \quad  $a\in \Omega \cap {\mathcal{U}}\Rightarrow Q_{a}(\delta _{0})\subset
 \Omega ,$ \ \par 
which means precisely that the family  ${\mathcal{Q}}=\lbrace
 Q_{a}(t)\rbrace _{a\in {\mathcal{U}}\cap \Omega ,\ t>0}$  is
 a good family of polydiscs in the sense of section 1. Moreover
 the Hefer's theorem~\ref{subPrinciple66} gives that the size
 of the sides of  $Q_{a}(t)$  are precisely equivalent to\ \par 
\quad \quad \quad \quad \quad  $\displaystyle \ \left\vert{r(a)}\right\vert ^{1/m_{j}}\simeq
 d(a)^{1/m_{j}},$ \ \par 
which means that the multi-type for this family in the sense
 of definition~\ref{BonFamille39} is precisely  $m_{j}(a),\ j=2,\
 ...,\ n.$   $\blacksquare $ \ \par 
\ \par 
\quad \quad 	So we can give a general definition for geometric Carleson measures
 equivalent to the one we gave in the case of convex domains
 of finite type.\ \par 

\begin{Dfnt}
Let  $\lambda $  be a positive Borel measure on the domain  $\Omega
 $  equipped with a good family of polydiscs  ${\mathcal{Q}}.$
  We shall say that  $\lambda $  is a {\bf geometric Carleson
 measure} in  $\Omega $  if :\par 
\quad \quad \quad $\displaystyle \exists C>0::\forall a\in \Omega ,\ \lambda (\Omega
 \cap Q_{a}(2))\leq C\sigma (\partial \Omega \cap Q_{a}(2)).$ 
\end{Dfnt}
\ \par 
\vfill\eject\ \par 

\section{Construction of balanced sub domains.~\label{6_CarlDomain33}}
\quad In the unit ball of  ${\mathbb{C}}^{n}$  a measure whose images
 by all automorphisms of the ball is uniformly bounded is a geometric
 Carleson measure, and this is a fact we used for instance in~\cite{intBall09}.
 Unfortunately in general domain, even convex ones or strictly
 pseudo-convex ones, there is just the identity as automorphism,
 so we have to overcome this issue.\ \par 
\quad \quad 	The aim now is to build a sub domain  $\Omega _{a}$  associated
 to a point  $a\in \Omega $  near the boundary such that the
 restriction to it of the measure we want to study is bounded
 by the right bound. If the domain  $\Omega _{a}$  is equivalent
 to a Carleson window, as defined at the beginning of section~\ref{8_CarlesonMeas30},
 then it will work.\ \par 
\quad \quad 	The main difficulty here is to get bounds {\sl independent}
 of  $a\in \Omega .$  We shall start with convex domains and
 define later a more general kind of domains for which our methods work.\ \par 
\ \par 
\quad Let  $\Omega $  be a  ${\mathcal{C}}^{\infty }$  smooth convex
 domain in  ${\mathbb{C}}^{n},\ a\in \Omega .$  By translation
 and rotation we can suppose that  $a=0,\ \alpha =\pi (a)=(d(a),0,....,0)$
  and the defining function  $\rho =d(a)+\Re z_{1}+\Gamma (z),$
  with  $\Gamma (z)={\mathcal{O}}(\left\vert{z}\right\vert ^{2}).$
  Let  ${\mathcal{E}}_{a},\ {\mathcal{E}}_{a}'$  be smooth complex
 ellipsoids centered at  $a=0\ :$ \ \par 
\quad \quad \quad \quad \quad  $\displaystyle {\mathcal{E}}_{a}:=\lbrace z\in {\mathbb{C}}^{n}::\sum_{j=1}^{n}{\frac{\left\vert{z_{j}}\right\vert
 ^{2}}{d(a)^{2/m_{j}}}}<4n\rbrace ,\ {\mathcal{E}}_{a}':=\lbrace
 z\in {\mathbb{C}}^{n}::\sum_{j=1}^{n}{\frac{\left\vert{z_{j}}\right\vert
 ^{2}}{d(a)^{2/m_{j}}}}<5n\rbrace \ ;$ \ \par 
consider the convex domain  ${\mathcal{E}}_{a}'\cap \Omega $
  and smooth it to get a smoothly bounded convex domain  $\Omega
 _{a}$  such that  $\displaystyle {\mathcal{E}}_{a}\cap \Omega
 \subset \Omega _{a}\subset {\mathcal{E}}_{a}'\cap \Omega .$
  This can be done as in~\cite{AmCoh84}, p. 129 : Suppose that
  $\displaystyle \alpha =\pi (a)=0$  and, as usual,  $\displaystyle
 \rho (z)=\Re z_{1}+f(\Im z_{1},z'),$  with  $\displaystyle z'=(z_{2},...,z_{n})$
  then there is a function  $S(x,y),$  convex and  ${\mathcal{C}}^{\infty
 }({\mathbb{R}}^{2})$  such that a defining function  $\rho _{a}$
  for  $\Omega _{a}$  is given by  $\displaystyle \rho _{a}:=S(2\left\vert{\Im
 z_{1}}\right\vert ^{2}+\left\vert{z'}\right\vert ^{2},\ \rho
 )\ ;$  hence any the  ${\mathcal{C}}^{k}$  norm of  $\displaystyle
 \rho _{a}$  is controlled by the  ${\mathcal{C}}^{k}$  norm
 of the defining function  $\rho $  of  $\Omega ,$  i.e.  $\displaystyle
 \ \forall k\in {\mathbb{N}},\ {\left\Vert{\rho _{a}}\right\Vert}_{{\mathcal{C}}^{k}}\leq
 C_{k}{\left\Vert{\rho }\right\Vert}_{{\mathcal{C}}^{k}}.$  Moreover
 we have that the outward normal derivative  $\displaystyle \
 \frac{\partial \rho }{\partial \eta }$  is uniformly bounded
 below because of the compactness of  $\displaystyle \partial
 \Omega $  and we have also  $\displaystyle \ \frac{\partial
 \rho _{a}}{\partial \eta }\geq \delta \frac{\partial \rho }{\partial
 \eta }>0$  independently of  $\displaystyle a,$  by the construction
 of  $\displaystyle \Omega _{a}.$  We shall need this last fact
 when we shall apply theorem~\ref{8_CarlDomain3} to interpolating
 sequences in section~\ref{9_InterSeq40} ; see remark~\ref{1_P1}.\ \par 
\ \par 
\quad \quad 	Let  ${\mathbb{S}}$  be the unit sphere in  ${\mathbb{C}}^{n}$
  and because  $\Omega _{a}$  is convex it is starlike with respect
 to  $a\ (=0),$   $\partial \Omega _{a}$  admits a spherical
 parametrization, i.e. there is a function  $R(\zeta )\in {\mathcal{C}}^{1}({\mathbb{S}}),\
 R(\zeta )>0,$  such that :\ \par 
\quad \quad \quad \quad \quad  $\partial \Omega _{a}=\lbrace z\in {\mathbb{C}}^{n}::\exists
 \zeta \in {\mathbb{S}},\ z=R(\zeta )\zeta \rbrace .$ \ \par 
Let  $\zeta \in {\mathbb{S}}$  and define  $D_{\zeta }$  to be
 the complex plane slice through  $\zeta \ :$  \ \par 
\quad \quad \quad \quad  $D_{\zeta }:=\lbrace tR(e^{i\theta }\zeta )e^{i\theta }\zeta
 ,\ \theta \in \lbrack 0,2\pi \rbrack ,\ t\in \lbrack 0,1\lbrack
 \rbrace .$ \ \par 
We shall use the notations\ \par 
\quad \quad \quad \quad \quad  $\forall \zeta \in {\mathbb{S}},\ d_{\zeta }(0)=\inf _{\theta
 \in \lbrack 0,2\pi \rbrack }R(e^{i\theta }\zeta )\ ;\ d_{\zeta
 max}(0)=\sup _{\theta \in \lbrack 0,2\pi \rbrack }R(e^{i\theta
 }\zeta ).$ \ \par 

\begin{Lmm}
~\label{CarlWindow72} We have\par 
\quad \quad \quad \quad \quad  $Q_{a}(2)\cap \Omega \subset \Omega _{a}\subset Q_{a}({\sqrt{5n}}).$ \par 
and\par 
\quad \quad \quad \quad \quad  $\displaystyle \forall \zeta \in \partial \Omega ,\ d_{\zeta
 max}(0)\leq \frac{{\sqrt{5n}}}{\delta _{0}}d_{\zeta }(0).$ 
\end{Lmm}
\quad \quad \quad 	Proof.\ \par 
 $\displaystyle z\in Q_{a}(2)\cap \Omega \Rightarrow \forall
 j=1,...,n,\ \left\vert{z_{j}}\right\vert <2d(a)^{1/m_{j}}\Rightarrow
 \sum_{j=1}^{n}{\frac{\left\vert{z_{j}}\right\vert ^{2}}{d(a)^{2/m_{j}}}<4n}\Rightarrow
 z\in {\mathcal{E}}_{a}\cap \Omega \subset \Omega _{a}.$ \ \par 
If  $z\in \Omega _{a}\subset {\mathcal{E}}_{a}'$  then\ \par 
\quad \quad \quad \quad \quad  $\displaystyle \ \sum_{j=1}^{n}{\frac{\left\vert{z_{j}}\right\vert
 ^{2}}{d(a)^{2/m_{j}}}<5n}\Rightarrow \forall j=1,...,n,\ \left\vert{z_{j}}\right\vert
 <{\sqrt{5n}}d(a)^{1/m_{j}}\Rightarrow z\in Q_{a}({\sqrt{5n}}),$ \ \par 
and the first assertion.\ \par 
\quad \quad 	Let us see that  $a$  is "in the middle" of the slices  $D_{\zeta }.$  \ \par 
Choose  $\theta $  such that  $\displaystyle d_{\zeta }(0)=R(e^{i\theta
 }\zeta ),$  then the real segment from  $0$  to  $\displaystyle
 R(e^{i\theta }\zeta )e^{i\theta }\zeta \in \partial \Omega _{a}$
  cross the boundary of  $\displaystyle Q_{a}(\delta _{0})$ 
 at a point  $\displaystyle tR(e^{i\theta }\zeta )e^{i\theta
 }\zeta $  with  $\displaystyle 0<t\leq 1$  because  $\displaystyle
 Q_{a}(\delta _{0})\subset \Omega _{a}.$ \ \par 
But if  $\displaystyle z=(z_{1},...,z_{n})\in \partial Q_{a}(\delta
 _{0})$  then		  $\displaystyle \exists j::\left\vert{z_{j}}\right\vert
 =\delta _{0}d(a)^{1/m_{j}}$  so we have here\ \par 
\quad \quad \quad \quad \quad  $\displaystyle \exists j::tR(e^{i\theta }\zeta )\left\vert{\zeta
 _{j}}\right\vert =\delta _{0}d(a)^{1/m_{j}}$ \ \par 
and because  $\displaystyle 0<t\leq 1$  we get\ \par 
\quad \quad \quad \quad \quad  $\displaystyle \delta _{0}d(a)^{1/m_{j}}\leq R(e^{i\theta }\zeta
 )\left\vert{\zeta _{j}}\right\vert =d_{\zeta }(0)\left\vert{\zeta
 _{j}}\right\vert .$ \ \par 
On the other hand, because  $\displaystyle \Omega _{a}\subset
 Q_{a}({\sqrt{5n}})$  which is a polydisc with sides parallel
 to the axis, we have\ \par 
\quad \quad \quad \quad \quad  $\displaystyle \forall k=1,...,n,\ \forall \varphi \in \lbrack
 0,2\pi \rbrack ,\ R(e^{i\varphi }\zeta )\left\vert{\zeta _{k}}\right\vert
 \leq {\sqrt{5n}}d(a)^{1/m_{k}}\Rightarrow d_{\zeta max}(0)\left\vert{\zeta
 _{k}}\right\vert \leq {\sqrt{5n}}d(a)^{1/m_{k}}\ ;$ \ \par 
in particular for  $\displaystyle \varphi =\theta $  and  $\displaystyle
 k=j$  we get\ \par 
\quad \quad \quad \quad \quad  $\displaystyle \delta _{0}d(a)^{1/m_{j}}\leq d_{\zeta }(0)\left\vert{\zeta
 _{j}}\right\vert \leq d_{\zeta max}(0)\left\vert{\zeta _{j}}\right\vert
 \leq {\sqrt{5n}}d(a)^{1/m_{j}}.$ \ \par 
This implies  $\displaystyle \ \frac{1}{\left\vert{\zeta _{j}}\right\vert
 }\leq \frac{d_{\zeta }(0)}{\delta _{0}d(a)^{1/m_{j}}}$  and
  $\displaystyle d_{\zeta max}(0)\leq \frac{{\sqrt{5n}}d(a)^{1/m_{j}}}{\left\vert{\zeta
 _{j}}\right\vert }\leq \frac{{\sqrt{5n}}}{\delta _{0}}d_{\zeta
 }(0).$   $\blacksquare $ \ \par 
\ \par 
\quad \quad 	Let  $D$  be a bounded convex domain in  ${\mathbb{C}}$  ; take
 a biggest disc contained in  $D,$  say  $\displaystyle D(0,r)$
  with  $\displaystyle 0\in D$  being its center and  $\displaystyle
 D(0,R)$  the smallest disc containing  $D$  with the same center
  $\displaystyle 0.$ \ \par 
Now parametrize the boundary  $\partial D$  of the convex  $D$
  by polar coordinates  $s(\theta )e^{i\theta }$  and set  $\displaystyle
 \gamma :=\frac{R}{r}.$ \ \par 

\begin{Lmm}
~\label{CarlWindow71} Let  $D$  be a convex domain in  ${\mathbb{C}},\
 0\in D$  with the previous notations ; let  $\displaystyle s'$
  be the derivative of  $\displaystyle s,$  then we have\par 
\quad \quad \quad \quad  $\displaystyle \ \left\vert{\frac{s'}{s}}\right\vert \leq {\sqrt{\gamma
 ^{2}-1}}.$ 
\end{Lmm}
\quad \quad 	Proof.\ \par 
We have that  $D(0,r)\subset D\subset D(0,R).$  Let  $z\in \partial
 D$  such that  $\tan V$  is minimal, where  $V$  is the angle
 between  $(0,z)$  and the tangent at  $z$  to  $\partial D.$
  Take the segment tangent  $T$  from  $z$  to  $t$  on the circle
  $\partial D(0,r)\ ;$  because  $D$  is convex we have  $T\subset
 D$  and the points  $w\in \partial D$  near  $z$  are such that
 the angle between  $(w,z)$  and  $(0,z)$  is bigger than the
 angle  $\alpha $  between  $(t,z)$  and  $(0,z)$  hence the
 angle  $V$  is bigger than  $\alpha .$ \ \par 
\ \par 
\ \par 
\ \par 
\ \par 
\ \par 
\ \par 
\ \par 
\ \par 
\begin{figure}[h]
\begin{center}
\vspace{-4cm}
\resizebox{8cm}{!}{\includegraphics{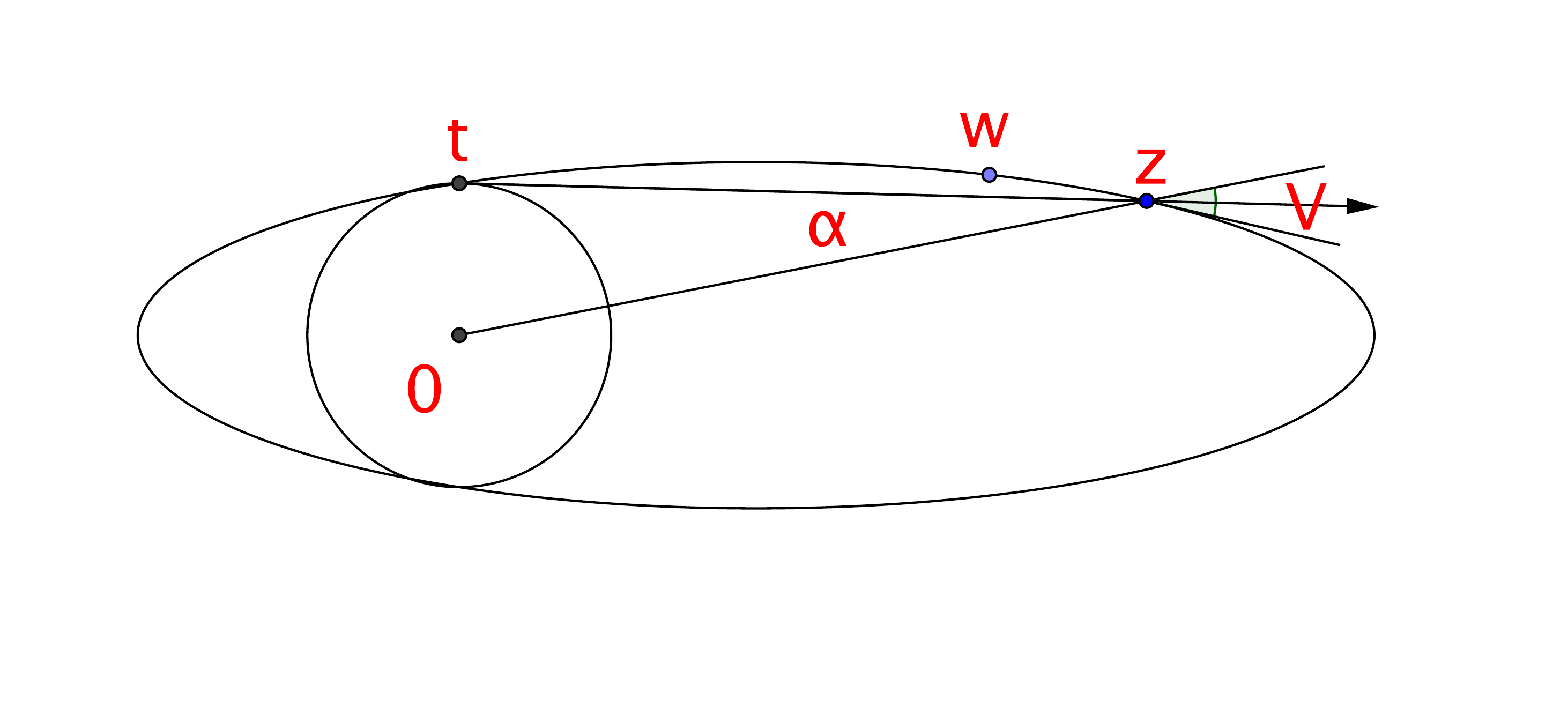}}
\end{center}
\vspace{-4cm}
\end{figure}\ \par 
\ \par 
\ \par 
\ \par 
\ \par 
Now we have that  $\displaystyle \ \left\vert{\sin \alpha }\right\vert
 =\frac{r}{\left\vert{z}\right\vert }\geq \frac{r}{R},$  hence
  $\displaystyle \ \left\vert{\tan \alpha }\right\vert \geq \frac{1}{{\sqrt{\gamma
 ^{2}-1}}},$  where  $\displaystyle \gamma :=\frac{R}{r}.$ \ \par 
So, because  $\displaystyle \ \left\vert{\frac{s'}{s}}\right\vert
 =\frac{1}{\left\vert{\tan V}\right\vert },$  we have\ \par 
\quad \quad \quad \quad \quad  $\displaystyle \ \left\vert{\frac{s'}{s}}\right\vert =\frac{1}{\left\vert{\tan
 V}\right\vert }\leq {\sqrt{\gamma ^{2}-1}}.$   $\blacksquare $ \ \par 
\ \par 
We shall apply this lemma to the slices  $\displaystyle D_{\zeta
 }$  of  $\Omega _{a}.$ \ \par 
Recall that  $U_{\zeta }(\theta )=R(e^{i\theta }\zeta )$  is
 precisely the polar coordinates parametrization of  $\partial
 D_{\zeta }$  in the coordinates of  ${\mathbb{C}}\zeta $  and
  $d_{\zeta }(0)$  is the distance from  $a(=0)$  to  $\partial
 D_{\zeta },$  hence here we have  $r=d_{\zeta }(0),\ R=d_{\zeta
 max}(0).$ \ \par 
We shall say that  $\displaystyle \Omega _{a}$  is  $\gamma $
  balanced with respect to  $a$  (definition~\ref{CarlDomain31})
 if  $\displaystyle \forall \zeta \in {\mathbb{S}},\ d_{\zeta
 max}(a)\leq \gamma d_{\zeta }(a)$  and  $\ \left\vert{U_{\zeta
 }'(\theta )}\right\vert \leq \gamma d_{\zeta max}(a)$  ; with
 this we have\ \par 

\begin{Lmm}
~\label{5_CarlDomain32}Because  $\Omega _{a}$  is such that all
 its slices  $D_{\zeta }=\Omega _{a}\cap \lbrace z::z=a+\lambda
 \zeta ,\ \lambda \in {\mathbb{C}}\rbrace $  are convex we have
 that  $\Omega _{a}$  is  $\gamma $  balanced with  $\displaystyle
 \gamma =\frac{{\sqrt{5n}}}{\delta _{0}}.$ 
\end{Lmm}
\quad \quad 	Proof.\ \par 
By lemma~\ref{CarlWindow71} with  $\displaystyle s(\theta ):=U_{\zeta
 }(\theta ),$  we have\ \par 
\quad \quad \quad \quad \quad  $\displaystyle \ \left\vert{\frac{U_{\zeta }'}{U_{\zeta }}}\right\vert
 \leq {\sqrt{\gamma ^{2}-1}}\leq \gamma \Rightarrow \left\vert{U_{\zeta
 }'}\right\vert \leq \gamma \left\vert{U_{\zeta }}\right\vert
 \leq \gamma d_{\zeta max}(0).$ \ \par 
Now using lemma~\ref{CarlWindow72}  $\displaystyle \gamma =\frac{R}{r}$
  hence we have that  $\Omega _{a}$  is  $\gamma $  balanced
 with  $\displaystyle \gamma =\frac{5n}{\delta _{0}^{2}}.$  
 $\blacksquare $ \ \par 
\quad \quad 	All we have done work as soon as the domain  $\Omega $  verifies
 the following definition.\ \par 

\begin{Dfnt}
A smoothly  ${\mathcal{C}}^{m},\ m\geq 2$  bounded domain  $\Omega
 $  with a good family of polydiscs is {\bf well balanced} if\par 
\quad \quad  $\exists R>2,\ \exists \epsilon >0,\ \exists \gamma >0::\forall
 a\in \Omega ,\ d(a)<\epsilon ,\ \exists \Omega _{a}\ \gamma
 $  balanced such that  $\displaystyle Q_{a}(2)\cap \Omega \subset
 \Omega _{a}\subset Q_{a}(R).$ 
\end{Dfnt}
And we have the theorem\ \par 

\begin{Thrm}
If  $\Omega $  is a smoothly  ${\mathcal{C}}^{m},\ m\geq 2,$
  bounded convex domain in  ${\mathbb{C}}^{n},$  with a good
 family of polydiscs then  $\Omega $  is well balanced.
\end{Thrm}
\quad \quad 	Proof.\ \par 
This is lemma~\ref{5_CarlDomain32}.  $\blacksquare $ \ \par 

\begin{Thrm}
~\label{5_CarlDomain20}If  $\Omega $  is well balanced, then
 for any  $a\in \Omega ,\ d(a)<\epsilon ,$  there is a  $\gamma
 $  balanced sub domain  $\Omega _{a}::Q_{a}(2)\cap \Omega \subset
 \Omega _{a}\subset Q_{a}({\sqrt{5n}})$  with the property\par 
\quad \quad \quad \quad \quad  $\forall u\in {\mathcal{N}}(\Omega _{a}),\ \ln \left\vert{u(a)}\right\vert
 =0$  then, with  $\displaystyle \Theta :=\partial \bar \partial
 \ln \left\vert{u}\right\vert ,\ \int_{\Omega _{a}}{d(z)\mathrm{T}\mathrm{r}}\Theta
 \leq C{\left\Vert{u}\right\Vert}_{{\mathcal{N}}(\Omega _{a})},$ \par 
where the constant depends only on  $\Omega $  and not on  $a.$ 
\end{Thrm}
\quad \quad 	Proof.\ \par 
We apply theorem~\ref{8_CarlDomain3} to  $\Omega _{a},$  then
 we have that\ \par 
\quad \quad \quad \quad \quad  $\displaystyle \ \int_{\Omega _{a}}{d(z)\mathrm{T}\mathrm{r}}\Theta
 \leq C{\left\Vert{u}\right\Vert}_{{\mathcal{N}}(\Omega _{a})},$ \ \par 
where the constant  $C$  depends only on  $\Omega .$    $\blacksquare $ \ \par 

\begin{Rmrq}
If  $\Omega $  is locally biholomorphic to a well balanced domain,
 then we have an analogous result by constructing the  $\Omega
 _{a}$  via the biholomorphism. Precisely let  $p\in \partial
 \Omega $  and  $\Phi $  a biholomorphism of  $\Omega \cap B(p,R)$
  on a well balanced domain  $\Omega '\cap \Phi (B(p,R)).$  Then
 we build the sub domains  $\Omega '_{\Phi (a)}$  and consider
  $\Omega _{a}:=\Phi ^{-1}(\Omega '_{\Phi (a)}).$  Because  $\Phi
 $  is biholomorphic in a neighborhood of  $\bar \Omega \cap
 B(p,R)$  we get easily that theorem~\ref{5_CarlDomain20} is still valid.\par 
In particular if  $\Omega $  is strictly pseudo-convex, then it works.
\end{Rmrq}
\vfill\eject\ \par 

\section{Interpolating and dual bounded sequences in  $H^{p}(\Omega
 ).$ ~\label{9_InterSeq40}}
\quad \quad 	Let  $\Omega $  be a domain in  ${\mathbb{C}}^{n}$  equipped
 with a good family of polydiscs. We shall study interpolating
 sequences in  $\Omega $  and generalise previous results we
 got for the unit ball to convex domains of finite type.\ \par 

\subsection{Reproducing kernels.}
\quad \quad 	Let  $S(z,\zeta )$  be the Szeg\"o kernel of  $\Omega ,$  i.e.
 the kernel of the orthogonal projection from  $L^{2}(\partial
 \Omega )$  onto  $H^{2}(\Omega ).$ \ \par 
\quad \quad 	To any point  $a\in \Omega $  we associate the vector  $k_{a}(\cdot
 ):=S(\cdot ,a)=\bar S(a,\cdot )\in H^{2}(\Omega ).$  This is
 a reproducing kernel for  $a$  because\ \par 
\quad \quad \quad \quad \quad  $\displaystyle \forall f\in H^{2}(\Omega ),\ f(a)=\int_{\partial
 \Omega }{f(\zeta )S(a,\zeta )d\sigma (\zeta )}$ \ \par 
by the definition of the Szeg\"o kernel, but\ \par 
\quad \quad \quad \quad \quad  $\displaystyle \ \int_{\partial \Omega }{f(\zeta )S(a,\zeta
 )d\sigma (\zeta )}=\int_{\partial \Omega }{f(\zeta )\bar k_{a}(\zeta
 )d\sigma (\zeta )}={\left\langle{f,k_{a}}\right\rangle},$ \ \par 
by the definition of  $k_{a}.$ \ \par 

\begin{Dfnt}
We say that the sequence  $S$  of points in  $\Omega $  is  $H^{p}(\Omega
 )$  {\bf interpolating} if\par 
\quad \quad \quad $(i)\ \forall a\in S,\ k_{a}\in H^{p'}(\Omega )\ ;$  (this is
 always true if  $p\geq 2.$ )\par 
\quad \quad \quad $\displaystyle (ii)\ \forall \lambda \in \ell ^{p}(S),\ \exists
 f\in H^{p}(\Omega )::\forall a\in S,\ f(a)=\lambda _{a}{\left\Vert{k_{a}}\right\Vert}_{p'},$
 \par 
with  $p'$  the conjugate exponent of  $p,$  i.e.  $\displaystyle
 \ \frac{1}{p}+\frac{1}{p'}=1.$ 
\end{Dfnt}
\quad A weaker notion is\ \par 

\begin{Dfnt}
We shall say that the sequence  $S$  of points in  $\Omega $
  is {\bf dual bounded} in  $H^{p}(\Omega )$  if there is a bounded
 sequence of elements in  $H^{p}(\Omega ),\ \lbrace \rho _{a}\rbrace
 _{a\in S}\subset H^{p}(\Omega )$  which dualizes the associated
 sequence of reproducing kernels, i.e.\par 
\quad \quad \quad $(i)\ \forall a\in S,\ k_{a}\in H^{p'}(\Omega )\ ;$  (this is
 always true if  $p\geq 2.$ )\par 
\quad \quad \quad $\displaystyle (ii)\ \exists C>0::\forall a\in S,\ {\left\Vert{\rho
 _{a}}\right\Vert}_{p}\leq C,\ \forall a,b\in S,\ {\left\langle{\rho
 _{a},\ k_{b}}\right\rangle}=\delta _{a,b}{\left\Vert{k_{b}}\right\Vert}_{p'}.$ 
\end{Dfnt}
\quad Clearly if  $S$  is  $H^{p}(\Omega )$  interpolating then  $S$
  is dual bounded in  $H^{p}(\Omega )\ :$  just interpolate the
 basic sequence of  $\ell ^{p}(S).$ \ \par 

\begin{Dfnt}
We say that  $S$  has the {\bf linear extension property} if
  $S$  is  $H^{p}(\Omega )$  interpolating and if moreover there
 is a bounded linear operator  $E\ :\ \ell ^{p}(S)\rightarrow
 H^{p}(\Omega )$  making the interpolation, i.e.\par 
\quad \quad \quad $\exists C>0,\ \forall \lambda \in \ell ^{p}(S),\ \forall a\in
 S,\ E(\lambda )(a)=\lambda _{a}{\left\Vert{k_{a}}\right\Vert}_{p'}$
  and  $\ {\left\Vert{E(\lambda )}\right\Vert}_{H^{p}(\Omega
 )}\leq C{\left\Vert{\lambda }\right\Vert}_{p}.$ 
\end{Dfnt}

\subsection{The  $p$  regularity.}
\quad \quad 	Let us introduce a link between the  $H^{p}$  norm of the reproducing
 kernels and the geometry of the boundary of  $\Omega ,$  with
 respect to the good family  ${\mathcal{Q}}.$ \ \par 

\begin{Dfnt}
~\label{BonFamille6}We shall say that  $\Omega $  is  $p$  {\bf
 regular} with respect to the family  ${\mathcal{Q}}$  if :\par 
\quad \quad \quad $\exists C>0::\forall a\in \Omega ,\ {\left\Vert{k_{a}}\right\Vert}_{p}^{-p'}\leq
 C\sigma (\partial \Omega \cap Q_{a}(2)),$ \par 
where  $p'$  is the conjugate exponent of  $p.$  Here we use
 the convention that if  $k_{a}\notin H^{p}(\Omega ),$  then
  $\displaystyle \ {\left\Vert{k_{a}}\right\Vert}_{p}=+\infty
 \Rightarrow {\left\Vert{k_{a}}\right\Vert}_{p}^{-p'}=0,$  so
 the inequality is true in this case.
\end{Dfnt}

\begin{Lmm}
If  $\displaystyle \Omega $  is a convex domain of finite type
 in  ${\mathbb{C}}^{n},$  then  $\displaystyle \Omega $  is 
 $p$  regular for any  $\displaystyle p>1.$ 
\end{Lmm}
\quad \quad 	Proof.\ \par 
Theorem~\ref{subPrinciple42} gives\ \par 
\quad \quad \quad \quad \quad  $\displaystyle \ {\left\Vert{k_{a}}\right\Vert}_{H^{p}(\Omega
 )}={\left\Vert{S(a,\cdot )}\right\Vert}_{H^{p}(\Omega )}\simeq
 \frac{1}{\sigma _{2n-1}(B(\alpha ,d(a)))^{1/p'}}\simeq \frac{1}{\sigma
 _{2n-1}(\partial \Omega \cap Q_{a}(2))^{1/p'}},$ \ \par 
which\!\!\!\! , by the definition~\ref{BonFamille6} of  $p$ 
 regularity, implies this lemma.  $\blacksquare $ \ \par 

\begin{Prps}
~\label{aspcGlo51}Let  $\Omega $  be a convex domain of finite
 type in  ${\mathbb{C}}^{n},\ a\in \Omega $  and  $\Omega _{a}$
  the sub domain  associated to  $a.$  The measure  $\displaystyle
 d\sigma _{2n-1\mid \partial \Omega _{a}\backslash \partial \Omega
 }$  is a geometric Carleson measure in  $\Omega .$ 
\end{Prps}
\quad \quad 	To prove this proposition we shall use the following lemmas.\ \par 

\begin{Lmm}
~\label{9_I0}Let  $U$  be an open set in  ${\mathbb{R}}^{k}$
  and  $V$  a graph in  ${\mathbb{R}}^{k+1}$  over  $U,$  i.e.\par 
\quad \quad \quad \quad \quad  $V:=\lbrace (x,y)\in {\mathbb{R}}^{k+1}::y=f(x),\ x=(x_{1},...,x_{k})\in
 U\rbrace ,$ \par 
with  $f$  of class  ${\mathcal{C}}^{1}(U).$  Then  $\displaystyle
 \sigma _{k}(V)\geq \sigma _{k}(U).$ 
\end{Lmm}
\quad \quad 	Proof.\ \par 
We shall use the formula for the Lebesgue measure for such a
 graph given in~\cite{GeoDiff87}, p. 203, formula 6.4.1.1 : let
  $\displaystyle (U,g)$  be a parametrization of  $V,$  then
 we have that :\ \par 
\quad \quad \quad \quad \quad  $\displaystyle g^{*}\omega ={\sqrt{\mathrm{d}\mathrm{e}\mathrm{t}(\frac{\partial
 g}{\partial x_{i}}\mid \frac{\partial g}{\partial x_{j}})}}dx_{1}\wedge
 \cdot \cdot \cdot \wedge dx_{k},$ \ \par 
where  $\displaystyle M:=(\frac{\partial g}{\partial x_{i}}\mid
 \frac{\partial g}{\partial x_{j}})$  is the matrix of the scalar
 product of the vectors  $\displaystyle \ \frac{\partial g}{\partial
 x_{i}}$  and  $\displaystyle \ \frac{\partial g}{\partial x_{j}}.$ \ \par 
Here we have that  $\displaystyle g(x)=(x_{1},...,x_{k},f(x))$  hence\ \par 
\quad \quad \quad \quad \quad  $\displaystyle \ \frac{\partial g}{\partial x_{j}}=(0,...,0,1,0,...,0,f_{j}'(x)),$
 \ \par 
with the  $1$  at the  $\displaystyle j^{th}$  position and 
 $\displaystyle f_{j}':=\frac{\partial f}{\partial x_{j}}.$  So we get\ \par 
\quad \quad \quad \quad \quad  $\displaystyle \ {\left\langle{\frac{\partial g}{\partial x_{i}},\frac{\partial
 g}{\partial x_{j}}}\right\rangle}=f_{i}'f_{j}'$  if  $\displaystyle
 i\neq j$  and  $\displaystyle \ {\left\langle{\frac{\partial
 g}{\partial x_{j}},\frac{\partial g}{\partial x_{j}}}\right\rangle}=1+(f_{j}')^{2}$
  if  $\displaystyle i=j.$ \ \par 
Hence the matrix  $M$  can be written	 	 $\displaystyle M=I+FF^{t}$
  where  $\displaystyle F$  is the column  vector  $\displaystyle
 F:=(f_{1}',...,f_{k}'),$  and  $\displaystyle F^{t}$  is the
 transpose line matrix.\ \par 
\quad \quad 	Clearly the matrix  $\displaystyle FF^{t}$  is positive, because
 for any vector  $\displaystyle v=(v_{1},...,v_{k})$  we have\ \par 
\quad \quad \quad \quad \quad  $\displaystyle v^{t}FF^{t}v=(v^{t}F)(F^{t}v)=(\sum_{j=1}^{k}{f_{j}'v_{j}})^{2}.$
 \ \par 
The eigenvalues  $\displaystyle \lambda _{j}$  of  $\displaystyle
 FF^{t}$  are all  $0$  but one because  $\displaystyle FF^{t}v=0$
  as soon as  $\displaystyle \ \sum_{j=1}^{k}{f_{j}'v_{j}}=0$
  which is a hyperplane and hence the only non zero eigenvalue,
  $\displaystyle \lambda _{k},$  is such that  $\displaystyle
 \lambda _{k}=\mathrm{T}\mathrm{r}FF^{t}=\sum_{j=1}^{k}{(f_{j}')^{2}}$
  because the sum of the eigenvalues is the trace of the matrix.\ \par 
\quad Now we have that the eigenvalues of  $M$  are  $\displaystyle
 1+\lambda _{j},$  hence the determinant of  $M$  is their product, so\ \par 
\quad \quad \quad \quad \quad  $\displaystyle \mathrm{d}\mathrm{e}\mathrm{t}M=1+\lambda _{k}=1+\sum_{j=1}^{k}{(f_{j}')^{2}}\geq
 1.$ \ \par 
The case  $\displaystyle k=2$  was already done in~\cite{GeoDiff87},
 p. 204, and here we provide the generalisation.\ \par 
Now we have\ \par 
\quad \quad \quad \quad \quad  $\displaystyle \sigma _{k}(V)=\int_{U}{{\sqrt{\mathrm{d}\mathrm{e}\mathrm{t}M}}dx_{1}\cdot
 \cdot \cdot dx_{k}}\geq \int_{U}{dx_{1}\cdot \cdot \cdot dx_{k}}=\sigma
 _{k}(U).$   $\blacksquare $ \ \par 

\begin{Rmrq}
In fact this lemma just says that the measure of the orthogonal
 projection  $U$  of  $V$  on  ${\mathbb{R}}^{k}$  has a Lebesgue
 measure smaller than the measure of  $\displaystyle V.$  I.e.
 the orthogonal projection is contracting for the Lebesgue measure,
 which seems quite natural.
\end{Rmrq}

\begin{Lmm}
~\label{c0} Let  $\displaystyle b\in \Omega ,\ \beta =\pi (b)$
  and  $\displaystyle T_{\beta }(\partial \Omega )$  be the real
 tangent space to  $\displaystyle \partial \Omega $  at  $\beta
 .$  Let  $\displaystyle F_{b}:=T_{\beta }(\partial \Omega )+d(b)n_{\beta
 },$  where  $\displaystyle n_{\beta }$  is the outward real
 normal at  $\beta $  to  $\displaystyle \partial \Omega $  and
  $\displaystyle B_{b}:=F_{b}\cap \bar Q_{b}(2),$  the "bottom"
 of   $\displaystyle Q_{b}(2)\ ;$  a sufficient condition to
 have  $\displaystyle \pi (\Omega \cap Q_{b}(2))\subset \pi (\partial
 \Omega \cap Q_{b}(2))$  is that  $\displaystyle \Omega \cap
 B_{b}=\emptyset .$ \par 
As a consequence if  $\displaystyle \Omega $  is convex then
  $\displaystyle \pi (\Omega \cap Q_{b}(2))\subset \pi (\partial
 \Omega \cap Q_{b}(2)).$ 
\end{Lmm}
\quad \quad 	Proof.\ \par 
Suppose that  $\displaystyle \Omega \cap B_{b}=\emptyset $  and
 take  $\displaystyle z\in \Omega ,$  take  $\zeta =\pi _{b}(z)$
  where  $\displaystyle \pi _{b}$  is the orthogonal projection
 on  $\displaystyle F_{b}\ ;$  then we have  $\displaystyle \zeta
 \in B_{b}$  hence  $\displaystyle \zeta \notin \Omega $  so
  $\displaystyle \rho (\zeta )\geq 0.$  On the other hand  $\displaystyle
 z\in \Omega \Rightarrow \rho (z)<0$  hence  $\rho $  being continuous
 on the real segment  $\displaystyle \lbrack z,\zeta \rbrack
 $  there is a  $\displaystyle w\in \rbrack z,\zeta \rbrack $
  such that  $\displaystyle \rho (w)=0,$  so  $\displaystyle
 w\in \partial \Omega .$  Now  $\displaystyle T_{\beta }(\partial
 \Omega )$  and  $\displaystyle F_{b}$  being parallel the segment
  $\displaystyle \lbrack z,\zeta \rbrack $  is orthogonal to
  $\displaystyle T_{\beta }(\partial \Omega )$  hence  $\displaystyle
 \pi (z)=T_{\beta }(\partial \Omega )\cap \lbrack z,\zeta \rbrack
 =\pi (w).$  Hence, because  $\displaystyle \pi (w)\in \pi (\partial
 \Omega )$  we have  $\displaystyle \pi (\Omega \cap Q_{b}(2))\subset
 \pi (\partial \Omega \cap Q_{b}(2)).$ \ \par 
\quad \quad 	If  $\displaystyle \Omega $  is convex then it lies on the same
 side of  $\displaystyle T_{\beta }(\partial \Omega )$  hence
 we have  $\displaystyle \Omega \cap B_{b}=\emptyset .$   $\blacksquare
 $ \ \par 
\ \par 
\quad \quad 	Proof of the proposition.\ \par 
We have to see that\ \par 
\quad \quad \quad \quad \quad  $\displaystyle \exists C>0::\forall b\in \Omega ,\ \sigma _{2n-1}(\partial
 \Omega _{a}\cap Q_{b}(2))\leq C\sigma _{2n-1}(\partial \Omega
 \cap Q_{b}(2)).$ \ \par 
Of course we take   $b$  such that  $\displaystyle (\partial
 \Omega _{a}\backslash \partial \Omega )\cap Q_{b}(2)\neq \emptyset
 $  and take the convex hull  $E$  of  $\displaystyle (\partial
 \Omega _{a}\backslash \partial \Omega )\cap Q_{b}(2)\ ;$  because
 the domain  $\displaystyle Q_{b}(2)$  is convex,  $\displaystyle
 E\subset Q_{b}(2)$  hence, by~\cite{ConvexMosynska06}, Corollary
 7.2.9 p 82, we have that\ \par 
\quad \quad \quad \quad \quad  $\sigma _{2n-1}(\partial E)\leq \sigma _{2n-1}(\partial Q_{b}(2))\ ;$ \ \par 
but, because  $\displaystyle \Omega _{a}$  $\displaystyle \subset
 \Omega $  is convex,\ \par 
\quad \quad \quad \quad \quad  $\displaystyle (\partial \Omega _{a}\backslash \partial \Omega
 )\cap Q_{b}(2)\subset \partial E\Rightarrow \sigma _{2n-1}((\partial
 \Omega _{a}\backslash \partial \Omega )\cap Q_{b}(2))\leq \sigma
 _{2n-1}(\partial E).$ \ \par 
We have\ \par 
\quad \quad \quad \quad \quad  $\displaystyle \partial Q_{b}(2))=\bigcup_{j=1}^{n}{(\partial
 D_{j}(b,d(b)^{1/m_{j}(b)})\prod_{k\neq j,k=1}^{n}{D_{k}(b,d(b)^{1/m_{k}(b)})})},$
 \ \par 
where  $\displaystyle D_{j}(b,r_{j})$  is the disc of center
  $b$  in the direction  $\displaystyle L_{j}$  given by the
 basis of the good family at the point  $\displaystyle \pi (b).$ So\ \par 
\quad \quad \quad \quad \quad  $\displaystyle \sigma _{2n-1}(\partial Q_{b}(2))=\sum_{j=1}^{n}{2\pi
 d(b)^{1/m_{j}(b)}\prod_{k\neq j,k=1}^{n}{\pi d(b)^{2/m_{k}(b)}}}\leq
 2\pi ^{n}\sum_{j=1}^{n}{d(b)^{2\lambda (b)+2-1/m_{j}(b)}}$ \ \par 
because  $\displaystyle 2\lambda (b)=\sum_{k=2}^{n}{\frac{2}{m_{k}}}\
 ;$  but  $\displaystyle \forall j=1,...,n,\ 2-1/m_{j}\geq 1$
  and we can restrict ourself to  $b$  such that  $\displaystyle
 d(b)\leq 1$  because we need to test only with the  $b$  near
 the boundary. Hence\ \par 
\quad \quad \quad \quad \quad  $\displaystyle \sigma _{2n-1}(\partial Q_{b}(2))\leq 2\pi ^{n}nd(b)^{1+2\lambda
 (b)}.$ \ \par 
\quad \quad 	So far we have\ \par 
\quad \quad \quad \quad 	\begin{equation}  \sigma _{2n-1}((\partial \Omega _{a}\backslash
 \partial \Omega )\cap Q_{b}(2))\leq \sigma _{2n-1}(\partial
 E)\leq \sigma _{2n-1}(\partial Q_{b}(2))\leq 2\pi ^{n}nd(b)^{1+2\lambda
 (b)}.\label{9_I1}\end{equation}\ \par 
\ \par 
To get   $d(b)^{1+2\lambda (b)}\lesssim \sigma _{2n-1}(\partial
 \Omega \cap Q_{b}(2))$  we shall use lemma~\ref{9_I0}. Set 
 $\displaystyle k=2n-1,\ U=Q_{b}(2)\cap T_{\beta }(\partial \Omega
 )$  where  $\displaystyle \beta =\pi (b)\in \partial \Omega
 $  and  $\displaystyle V=\partial \Omega \cap Q_{b}(2).$  For
  $b$  uniformly near  $\displaystyle \partial \Omega ,\ V$ 
 is a graph over  $U':=\pi (V)\subset U,$  with  $\pi $  the
 orthogonal  projection on the real tangent space  $\displaystyle
 T_{\beta }(\partial \Omega ),$  and we have by lemma~\ref{9_I0}
 that  $\displaystyle \sigma _{2n-1}(V)\geq \sigma _{2n-1}(U')$
  so it remains to estimate  $\displaystyle \sigma _{2n-1}(U').$ \ \par 
\quad \quad 	Recall that, by the definition of a good family, we have  $\displaystyle
 Q_{b}(\delta _{0})\subset \Omega $  hence  $\pi (Q_{b}(\delta
 _{0}))\subset \pi (\Omega \cap Q_{b}(2)).$ \ \par 
We apply lemma~\ref{c0} to  $\displaystyle \Omega $  convex to
 get  $\displaystyle \pi (Q_{b}(\delta _{0}))\subset \pi (V)=U'.$  So\ \par 
\quad \quad \quad \quad \quad  $\displaystyle \sigma _{2n-1}(U')\geq \sigma _{2n-1}(\pi (Q_{b}(\delta
 _{0}))).$ \ \par 
\quad \quad 	Because the basis for  $\displaystyle Q_{b}$  is the basis at
  $\displaystyle \beta =\pi (b)$  and  $\displaystyle T_{\beta
 }(\partial \Omega )$  is the real tangent space, the only missing
 direction is the real normal at  $\beta ,$  hence we have\ \par 
\quad \quad \quad \quad \quad  $\displaystyle \sigma _{2n-1}(\pi (Q_{b}(\delta _{0}))=\delta
 _{0}d(b){\times}\prod_{j=2}^{n}{\delta _{0}^{2}d(b)^{2/m_{j}(b)}}=\delta
 _{0}^{2n+1}d(b)^{1+2\lambda (b)}.$ \ \par 
\quad \quad 	Finally we get\ \par 
\quad \quad \quad \quad \quad  $\displaystyle \delta _{0}^{2n+1}d(b)^{1+2\lambda (b)}\leq \sigma
 _{2n-1}(V)=\sigma _{2n-1}(\partial \Omega \cap Q_{b}(2))$ \ \par 
and by~(\ref{9_I1})\ \par 
\quad \quad \quad \quad \quad  $\displaystyle \sigma _{2n-1}((\partial \Omega _{a}\backslash
 \partial \Omega )\cap Q_{b}(2))\leq 2\pi ^{n}nd(b)^{1+2\lambda (b)}$ \ \par 
hence\ \par 
\quad \quad \quad \quad \quad \quad  $\displaystyle \sigma _{2n-1}((\partial \Omega _{a}\backslash
 \partial \Omega )\cap Q_{b}(2))\leq \frac{2\pi ^{n}}{\delta
 _{0}^{2n+1}}n\sigma _{2n-1}(\partial \Omega \cap Q_{b}(2)),$ \ \par 
which says precisely that the measure  $\displaystyle d\sigma
 _{2n-1\mid \partial \Omega _{a}\backslash \partial \Omega }$
  is a geometric Carleson measure in  $\Omega .$   $\blacksquare $ \ \par 
\ \par 
\quad In order to continue we shall need the easy remark :\ \par 

\begin{Rmrq}
~\label{aspcGlo815}For any smoothly bounded domain  $\Omega $
  we have the inequality\par 
\quad \quad \quad $\displaystyle \forall f\in H^{p}(\Omega ),\ {\left\Vert{f}\right\Vert}_{{\mathcal{N}}(\Omega
 )}\leq \root{p'}\of{\sigma _{2n-1}(\partial \Omega )}{\left\Vert{f}\right\Vert}_{H^{p}(\Omega
 )}.$ 
\end{Rmrq}
\quad Proof.\ \par 
We have  $\ln ^{+}\left\vert{f}\right\vert \leq \left\vert{f}\right\vert
 ,$  so, with  $\displaystyle \sigma _{\epsilon }$  the  $\displaystyle
 \sigma _{2n-1}$  Lebesgue measure on the manifold  $\displaystyle
 r(z)=-\epsilon ,$ \ \par 
$\displaystyle \ {\left\Vert{f}\right\Vert}_{{\mathcal{N}}(\Omega
 )}:=\sup _{\epsilon >0}\int_{r(z)=-\epsilon }{\ln  ^{+}\left\vert{f(z)}\right\vert
 \,d\sigma _{\epsilon }(z)}\leq \sup _{\epsilon >0}\int_{r(z)=-\epsilon
 }{\left\vert{f(z)}\right\vert \,d\sigma _{\epsilon }(z)},$ \ \par 
hence  $\ {\left\Vert{f}\right\Vert}_{{\mathcal{N}}(\Omega )}\leq
 {\left\Vert{f}\right\Vert}_{H^{1}(\Omega )}.$  But  $\sigma
 (\partial \Omega )$  being finite, we have  $\displaystyle \
 {\left\Vert{f}\right\Vert}_{H^{1}(\Omega )}\leq \root{p'}\of{\sigma
 _{2n-1}(\partial \Omega )}{\left\Vert{f}\right\Vert}_{H^{p}(\Omega
 )}.$   $\blacksquare $ \ \par 

\begin{Prps}
If  $S$  is a  $H^{p}$  dual bounded sequence in a convex domain
 of finite type, then  $S$  is separated.
\end{Prps}
\quad Proof.\ \par 
The hypothesis on the sequence  $S$  implies that\ \par 
\quad \quad \quad \quad $\exists C>0,\ \exists \rho _{a}\in H^{p}(\Omega )::{\left\Vert{\rho
 _{a}}\right\Vert}_{p}\leq C,\ {\left\langle{\rho _{a},\ k_{b}}\right\rangle}=0,\
 \left\vert{{\left\langle{\rho _{a},\ k_{a}}\right\rangle}}\right\vert
 \gtrsim {\left\Vert{k_{a}}\right\Vert}_{p'}\ ;$ \ \par 
then\ \par 
\begin{equation}  \ {\left\Vert{k_{a}-k_{b}}\right\Vert}_{p'}\geq
 \left\vert{{\left\langle{\frac{\rho _{a}}{{\left\Vert{\rho _{a}}\right\Vert}_{p}},\
 k_{a}-k_{b}}\right\rangle}}\right\vert \geq \frac{1}{C}\left\vert{{\left\langle{\rho
 _{a},\ k_{a}}\right\rangle}}\right\vert \gtrsim {\left\Vert{k_{a}}\right\Vert}_{p'}.\label{ConvInt16}\end{equation}\
 \par 
\quad Now for  $\epsilon >0$  we get the existence of  $\gamma $  such
 that, if we suppose that  $b\in Q_{a}(t)$  for a  $t<\gamma ,$ \ \par 
\quad \quad \quad $\displaystyle \ \left\vert{{\left\langle{\rho _{a},\ k_{a}-k_{b}}\right\rangle}}\right\vert
 \leq C{\left\Vert{k_{a}-k_{b}}\right\Vert}_{p'}\leq C\epsilon
 {\left\Vert{k_{a}}\right\Vert}_{p'}$ \ \par 
and a contradiction with ~(\ref{ConvInt16}) if we choose  $\epsilon
 $  small enough.  $\blacksquare $ \ \par 

\begin{Thrm}
~\label{BonFamille5}Let  $\Omega $  be a convex domain of finite
 type in  ${\mathbb{C}}^{n}.$  If the sequence of points  $S\subset
 \Omega $  is dual bounded in  $H^{p}(\Omega ),$  then the measure
  $\displaystyle \mu :=\sum_{a\in S}{d(a)^{n}\delta _{a}},$ 
 is a geometric Carleson measure in  $\Omega .$ 
\end{Thrm}
\quad Proof.\ \par 
We have to show that\ \par 
\quad \quad \quad $\displaystyle \forall a\in \Omega ,\ \mu (\Omega \cap Q_{a}(2))=\sum_{b\in
 S\cap Q_{a}(2)}{d(b)^{n}}\leq C\sigma _{2n-1}(\partial \Omega
 \cap Q_{a}(2)).$ \ \par 
Dual boundedness means that we have a sequence  $\lbrace \rho
 _{a}\rbrace _{a\in S}\subset H^{p}(\Omega )$  such that\ \par 
\quad \quad \quad $\forall a,b\in S,\ \rho _{a}(b)=\delta _{ab}{\left\Vert{k_{a}}\right\Vert}_{p'},\
 {\left\Vert{\rho _{a}}\right\Vert}_{p}\leq C.$ \ \par 
This implies that\ \par 
\quad \quad \quad \quad \quad  $\displaystyle \forall a\in S,\ {\left\Vert{\rho _{a}/\rho _{a}(a)}\right\Vert}_{H^{p}(\Omega
 )}\leq C{\left\Vert{k_{a}}\right\Vert}_{p'}^{-1}.$ \ \par 
\ \par 
\quad In the case of the unit ball~\cite{AmarWirtBoule07} we used the
 automorphisms and a classical lemma by Garnett to pass from
 bounded measures to geometric Carleson ones. Here of course
 we have to overcome the lack of automorphisms.\ \par 
\quad Because  $S$  is dual bounded it is a separated sequence of points
 in  $\Omega .$  Consider the sub-domain  $\Omega _{a}$  associated
 to the point  $a,$  built in section~\ref{6_CarlDomain33} and
 the sequence  $S_{a}:=S\cap \Omega _{a}\subset \Omega _{a}.$ \ \par 
\quad 	Let  $a\in S$  and  $u:=\rho _{a}/\rho _{a}(a)$  ; we have 
 $u\in H^{p}(\Omega )$  by hypothesis. We notice that  $S\backslash
 \lbrace a\rbrace \subset u^{-1}(0),$  and that  $u(a)=1,$  so
 we get by theorem~\ref{3_AspcDomain0}, with  $X=u^{-1}(0)\cap
 \Omega _{a}$  and  $\Theta $  its  $\displaystyle (1,1)$  current
 of integration,\ \par 
\quad \quad \quad \quad \quad  $\displaystyle \ \sum_{c\in S_{a}}{d(c)^{n}}\leq \Gamma (\Omega
 _{a}){\left\Vert{\Theta }\right\Vert}_{B},$ \ \par 
where  $\Gamma (\Omega _{a})$  depends on the  ${\mathcal{C}}^{M({\mathcal{Q}})+1}$
  norm of the defining function of  $\Omega _{a}$  which by construction
 of  $\Omega _{a}$  is controlled by the  ${\mathcal{C}}^{M({\mathcal{Q}})+1}$
  norm of the defining function of  $\Omega .$ \ \par 
\quad \quad 	Now because  $\Omega $  is convex of finite type,  $\Omega _{a}$
  is  $\displaystyle \ \frac{4n}{\delta _{0}^{2}}$  balanced
 by lemma~\ref{5_CarlDomain32} with respect to  $a,$  hence by
 theorem~\ref{8_CarlDomain3} we get  $\ {\left\Vert{\Theta }\right\Vert}_{B}\leq
 C{\left\Vert{u}\right\Vert}_{{\mathcal{N}}(\Omega _{a})},$ 
 the constant  $C$  depending only on  $\Omega $  and not on  $a.$  So\ \par 
\quad \quad \quad \quad \quad  $\displaystyle \ \sum_{c\in S_{a}}{d(c)^{n}}\leq C\Gamma {\left\Vert{u}\right\Vert}_{{\mathcal{N}}(\Omega
 _{a})},$ \ \par 
and again the constant being independent of  $a.$ \ \par 
\quad \quad 	By the remark~\ref{aspcGlo815} we get\ \par 
\quad \quad \quad \quad \quad  $\displaystyle \ \sum_{c\in S_{a}}{d(c)^{n}}\leq C\Gamma {\left\Vert{u}\right\Vert}_{{\mathcal{N}}(\Omega
 _{a})}\leq C\Gamma \root{p'}\of{\sigma _{2n-1}(\partial \Omega
 _{a})}{\left\Vert{u}\right\Vert}_{H^{p}(\Omega _{a})}.$ \ \par 
Set  $C(\Omega ):=C\Gamma $  which depends only on  $\Omega ,$  we get\ \par 
\quad \quad \quad \quad \quad  $\displaystyle \ \sum_{c\in S_{a}}{d(c)^{n}}\leq C(\Omega )\
 \root{p'}\of{\sigma _{2n-1}(\partial \Omega _{a})}{\left\Vert{u}\right\Vert}_{H^{p}(\Omega
 _{a})}.$ \ \par 
The measure  $\displaystyle d\sigma _{\mid \partial \Omega _{a}\backslash
 \partial \Omega }$  is a geometric Carleson measure in  $\Omega
 $  by lemma~\ref{aspcGlo51} hence by the embedding Carleson
 theorem~\ref{aspcGlo50} we have\ \par 
\quad \quad \quad \quad $u\in H^{p}(\Omega _{a})$  and  $\displaystyle \ {\left\Vert{u}\right\Vert}_{H^{p}(\Omega
 _{a})}\leq C{\left\Vert{u}\right\Vert}_{H^{p}(\Omega )},$ \ \par 
with the constant  $C$  independent of  $a.$ \ \par 
Hence\ \par 
\quad \quad \quad $\displaystyle \forall a\ \in S,\ \sum_{c\in S_{a}}{d(c)^{n}}\leq
 C(\Omega )\ \root{p'}\of{\sigma _{2n-1}(\partial \Omega _{a})}{\left\Vert{u}\right\Vert}_{H^{p}(\Omega
 _{a})}.$ \ \par 
The dual boundedness then gives, because  $\displaystyle u:=\rho
 _{a}/\rho _{a}(a),$ \ \par 
\quad \quad \quad $\displaystyle \forall a\ \in S,\ {\left\Vert{u}\right\Vert}_{H^{p}(\Omega
 )}\lesssim {\left\Vert{k_{a}}\right\Vert}_{p'}^{-1},$ \ \par 
hence\ \par 
\quad \quad \quad $\displaystyle \forall a\ \in S,\ \sum_{c\in S_{a}}{d(c)^{n}}\leq
 C(\Omega ){\left\Vert{u}\right\Vert}_{H^{p}(\Omega )}\leq C(\Omega
 ){\left\Vert{k_{a}}\right\Vert}_{p'}^{-1},$ \ \par 
where the (new) constant  $C(\Omega )$  is still independent of  $a.$ \ \par 
Finally the  $p$ -regularity of  $\Omega $  gives\ \par 
\quad \quad \quad $\displaystyle \ {\left\Vert{k_{a}}\right\Vert}_{p'}^{-1}\lesssim
 (\sigma _{2n-1}(\partial \Omega \cap Q_{a}(2))^{1/p},$ \ \par 
hence\ \par 
\quad \quad \quad  \begin{equation}  \forall a\ \in S,\ \sum_{c\in S_{a}}{d(c)^{n}}\leq
 C(\Omega )(\sigma _{2n-1}(\partial \Omega _{a}))^{1/p'}{\left\Vert{k_{a}}\right\Vert}_{p'}^{-1}\leq
 C(\Omega )\sigma _{2n-1}(\partial \Omega \cap Q_{a}(2)),\label{BonFamille4}\end{equation}\
 \par 
because we have that  $\sigma _{2n-1}(\partial \Omega _{a})\simeq
 \sigma _{2n-1}(\partial \Omega \cap Q_{a}(2)),$  still with
 the constant independent of  $a.$ \ \par 
\ \par 
\quad So we have proved the right inequality for a point  $a\in S.$
  It remains to have it for any point in  $\Omega .$ \ \par 
\quad Fix  $b\in \Omega $  ; take a point  $a_{1}\in S\cap Q_{b}(2)$
  such that  $d(a_{1}):=d(a_{1},\ \partial \Omega )$  is as big
 as possible. Now set  $E_{1}:=Q_{b}(2)\backslash Q_{a_{1}}(2)$
  and take a point  $a_{2}\in S\cap E_{1}$  such that  $d(a_{2}):=d(a_{2},\
 \partial \Omega )$  is as big as possible ; set  $E_{2}:=E_{1}\backslash
 Q_{a_{2}}(2)$  and take a point  $a_{3}\in S\cap E_{2}$  such
 that  $d(a_{3}):=d(a_{3},\ \partial \Omega )$  is as big as
 possible etc... This way we have a sequence  $G:=\lbrace a_{j}\rbrace
 $  of points in  $Q_{b}(2)\cap S$  with  $d(a_{j})$  decreasing.
 Moreover we have\ \par 
\quad \quad \quad $S\cap Q_{b}(2)=\bigcup_{j=1}^{\infty }{S\cap Q_{a_{j}}(2)}.$ \ \par 
For any  $j=1,\ ...$  we have, by~\ref{BonFamille4}\ \par 
\quad \quad \quad $\ \sum_{c\in S_{a_{j}}}{d(c)^{n}}\lesssim {\left\Vert{k_{a_{j}}}\right\Vert}_{p'}^{-1}\leq
 C(\Omega )\sigma _{2n-1}(\partial \Omega \cap Q_{a_{j}}(2)).$ \ \par 
\quad Now define  $P_{a}(\delta ):=Q_{a}(\delta )\cap T_{a}(\partial
 \Omega ),$  where  $\displaystyle T_{a}(\partial \Omega )$ 
 is the parallel hyperplane to the tangent to  $\partial \Omega
 $  at  $\alpha $  passing through  $a.$  We have that if  $a\in
 \Omega \cap {\mathcal{U}}$  then  $\sigma _{2n-1}(P_{a}(\delta
 ))\simeq \sigma _{2n-1}(\pi (P_{a}(\delta ))),$  where the constants
 behind the sign  $\simeq $  are independent of  $a,$  because
 the projection  $\pi $  is a diffeomorphism from  $\displaystyle
 P_{a}(\delta )$  onto its image in  $\partial \Omega .$  Its
 jacobian  $J$  is still a smooth function, hence we have that
  $C={\left\Vert{J}\right\Vert}_{\infty }$  is uniformly bounded
 by the compactness of  $\displaystyle \partial \Omega ,$  and
 so is  $\ {\left\Vert{J^{-1}}\right\Vert}_{\infty }.$ \ \par 
\ \par 
\quad Because the sets  $Q_{a_{j}}(\delta )$  are disjoint we get that\ \par 
\quad $\bullet \ \pi (P_{a_{j}}(\delta ))$  are disjoint and  $\sigma
 _{2n-1}(P_{a_{j}}(\delta ))\simeq \sigma _{2n-1}(\pi (P_{a_{j}}(\delta
 ))),$ \ \par 
\quad $\displaystyle \bullet \ \sigma _{2n-1}(P_{a_{j}}(2))\simeq {\left({\frac{2}{\delta
 }}\right)}^{2n-1}\sigma _{2n-1}(P_{a_{j}}(\delta )).$ \ \par 
So\ \par 
\quad $\displaystyle \sigma _{2n-1}(P_{a_{j}}(2))\simeq {\left({\frac{2}{\delta
 }}\right)}^{2n-1}\sigma _{2n-1}(\pi (P_{a_{j}}(\delta ))).$ \ \par 
We want to estimate\ \par 
\quad \quad \quad $\displaystyle \ \sum_{c\in S_{b}}{d(c)^{n}}=\sum_{j=1}^{\infty
 }{\sum_{c\in S_{a_{j}}}{d(c)^{n}}}\lesssim \sum_{j=1}^{\infty
 }{\sigma _{2n-1}(\partial \Omega \cap Q_{a_{j}}(2))},$ \ \par 
but\ \par 
\quad \quad \quad $\displaystyle \ \sum_{j=1}^{\infty }{\sigma _{2n-1}(\pi (P_{a_{j}}(\delta
 )))}\leq \sigma _{2n-1}(\partial \Omega \cap Q_{b}(2)),$ \ \par 
because the  $\displaystyle \pi (P_{a_{j}}(\delta ))$  are disjoint
 and contained in  $\displaystyle \partial \Omega \cap Q_{b}(2)$  and\ \par 
\quad $\displaystyle \ \sum_{j=1}^{\infty }{\sigma _{2n-1}(\pi (P_{a_{j}}(\delta
 )))}\gtrsim \sum_{j=1}^{\infty }{\sigma _{2n-1}(P_{a_{j}}(\delta
 )))}\gtrsim \sum_{j=1}^{\infty }{\sigma _{2n-1}(P_{a_{j}}(2)))}\gtrsim
 \sum_{j=1}^{\infty }{\sigma _{2n-1}(\partial \Omega \cap Q_{a_{j}}(2))}.$
 \ \par 
So\ \par 
\quad \quad \quad $\displaystyle \ \sum_{c\in S_{b}}{d(c)^{n}}\lesssim \sum_{j=1}^{\infty
 }{\sigma _{2n-1}(\partial \Omega \cap Q_{a_{j}}(2))}\lesssim
 \sum_{j=1}^{\infty }{\sigma _{2n-1}(\pi (P_{a_{j}}(\delta )))}\lesssim
 \sigma _{2n-1}(\partial \Omega \cap Q_{b}(2)).\ \blacksquare $ \ \par 

\begin{Thrm}
~\label{aspc2}Let  $\Omega $  be a convex domain of finite type
 in  ${\mathbb{C}}^{n}.$  If the sequence of points  $S\subset
 \Omega $  is dual bounded in  $H^{p}(\Omega ),$  then the canonical
 measure  $\displaystyle \mu :=\sum_{a\in S}{d(a)^{1+2\lambda
 (a)}\delta _{a}}$  is a geometric Carleson measure in  $\Omega .$ 
\end{Thrm}
\quad Proof.\ \par 
We take advantage of the fact that a convex domain of finite
 type is {\bf aspc} to separate the sequence  $S$  in two parts
  $S=B_{S}\cup G_{S}.$  For the bad points we need not the hypothesis
 of dual boundedness because theorem~\ref{aspc0} gives\ \par 
\quad \quad \quad \quad \quad  $\displaystyle \ \sum_{c\in B_{S}\cap Q_{a}(2)}{d(c)^{1+2\lambda
 (c)}}\leq C(\Omega )\frac{\sigma _{2n-1}(\partial \Omega \cap
 Q_{a}(2))}{\delta ^{2}},$ \ \par 
which is true for any  $a\in \Omega ,$  and this is precisely
 the definition of a geometric Carleson measure, so we get that
 the measure  $\displaystyle \mu _{b}:=\sum_{a\in B_{S}}{d(a)^{1+2\lambda
 (a)}\delta _{a}}$  is a geometric Carleson measure.\ \par 
We have  $\displaystyle \mu _{g}:=\sum_{a\in G_{S}}{d(a)^{1+2\lambda
 (a)}\delta _{a}}\leq \lambda :=\sum_{a\in S}{d(a)^{n}\delta
 _{a}},$  and  $\lambda $  is a geometric Carleson measure by
 theorem~\ref{BonFamille5}. So adding  $\mu _{b}$  and  $\mu
 _{g}$  we get that  $\mu $  is a geometric Carleson measure.
  $\blacksquare $ \ \par 

\subsection{Interpolating sequences.}
\quad \quad 	We shall need the definition\ \par 

\begin{Dfnt}
The sequence  $S$  is a  $q$  Carleson sequence if\par 
\quad \quad \quad $\displaystyle \exists D>0,\ \forall \lambda \in \ell ^{q}(S),\
 {\left\Vert{\sum_{a\in S}{\lambda _{a}\frac{k_{a}}{{\left\Vert{k_{a}}\right\Vert}_{q}}}}\right\Vert}_{q}\leq
 D{\left\Vert{\lambda }\right\Vert}_{\ell ^{q}(S)}.$ 
\end{Dfnt}
\quad \quad 	In~\cite{AmarExtInt06}, we proved by duality that if the canonical
 measure  $\displaystyle \mu :=\sum_{a\in S}{d(a)^{1+2\lambda
 (a)}\delta _{a}}$  is a  $\displaystyle q'$  Carleson measure
 and if  $\displaystyle \ {\left\Vert{k_{a}}\right\Vert}_{q}^{-q'}\simeq
 d(a)^{1+2\lambda (a)}$  then  $S$  is a  $q$  Carleson sequence.
 We shall do it again in this setting.\ \par 

\begin{Lmm}
If  $\Omega $  is a convex domain of finite type in  ${\mathbb{C}}^{n}$
  and if  $S\subset \Omega $  is a dual bounded sequence of points
 in  $\displaystyle H^{p}(\Omega ),$  then  $S$  is a  $q$  Carleson
 sequence for any  $q\in \rbrack 1,\ \infty \lbrack .$ 
\end{Lmm}
\quad Proof.\ \par 
Because the Szeg\"o projection is bounded on  $L^{p}(\partial
 \Omega )$  for  $1<p<\infty ,$  (~\cite{McNealStein97}, theorem
 5.1) we have that the dual of  $H^{p}(\Omega )$  is  $H^{p'}(\Omega
 ),$  with  $p'$  the conjugate exponent of  $p.$  Hence we can
 evaluate the norm this way\ \par 
\quad \quad \quad $\displaystyle \ {\left\Vert{\sum_{a\in S}{\lambda _{a}\frac{k_{a}}{{\left\Vert{k_{a}}\right\Vert}_{q}}}}\right\Vert}_{q}\simeq
 \sup  _{f\in H^{q'}(\Omega ),\ {\left\Vert{f}\right\Vert}_{q'}\leq
 1}\left\vert{\sum_{a\in S}{\lambda _{a}\frac{f(a)}{{\left\Vert{k_{a}}\right\Vert}_{q}}}}\right\vert
 \lesssim $ \ \par 
\quad \quad \quad \quad \quad \quad \quad \quad \quad \quad \quad \quad \quad $\displaystyle \lesssim {\left\Vert{\lambda }\right\Vert}_{\ell
 ^{q}(S)}\sup _{f\in H^{q'}(\Omega ),\ {\left\Vert{f}\right\Vert}_{q'}\leq
 1}\left\vert{\sum_{a\in S}{\frac{\left\vert{f(a)}\right\vert
 ^{q'}}{{\left\Vert{k_{a}}\right\Vert}_{q}^{q'}}}}\right\vert ^{1/q'},$ \ \par 
by H\"older. But  \[\displaystyle \ {\left\Vert{k_{a}}\right\Vert}_{p}={\left\Vert{S(\cdot
 ,\ a)}\right\Vert}_{H^{p}}\simeq \frac{1}{\sigma _{2n-1}(B(\alpha
 ,d(a)))^{1/p'}}\]  by theorem~\ref{subPrinciple42} and we have,
 by~(\ref{strongPC930}),  $\displaystyle \sigma _{2n-1}(B(\alpha
 ,\ \epsilon ))\simeq \epsilon \prod_{j=2}^{n}{\tau _{j}(\alpha
 ,\ \epsilon )^{2}}\ $  and by~(\ref{subPrinciple616}) in Hefer's
 theorem~\ref{subPrinciple66}, we have  $\displaystyle \tau _{j}(\zeta
 ,\ \epsilon )\simeq \epsilon ^{1/m_{j}(\zeta )},$  hence, \ \par 
\quad \quad \quad $\displaystyle \sigma (B(\alpha ,\ \epsilon ))\simeq \epsilon
 \prod_{j=2}^{n}{\tau _{j}(\alpha ,\ \epsilon )^{2}}\ \simeq
 \epsilon ^{1+2\lambda (\alpha )},\ \lambda (\alpha ):=\sum_{j=2}^{n}{\frac{1}{m_{j}(\alpha
 )}}.$ \ \par 
We shall apply this with  $\alpha =\pi (a),\ \epsilon =d(a)$
  and, because  $P_{\epsilon }(a)\cap P_{\epsilon }(\alpha )\neq
 \emptyset $  we have by~(\ref{subPrinciple618})\ \par 
\quad \quad \quad $\displaystyle \ \prod_{j=2}^{n}{\tau _{j}(a,d(a))^{2}}\simeq
 \prod_{j=2}^{n}{\tau _{j}(\alpha ,d(a))^{2}}\Rightarrow d(a)^{2\lambda
 (a)}\simeq d(a)^{2\lambda (\alpha )}.$ \ \par 
Putting this in  $\displaystyle \ {\left\Vert{k_{a}}\right\Vert}_{q}$
  we get\ \par 
\quad \quad \quad $\displaystyle \ {\left\Vert{k_{a}}\right\Vert}_{q}^{-1}\simeq
 d(a)^{\frac{1+2\lambda (\alpha )}{q'}}\Rightarrow {\left\Vert{k_{a}}\right\Vert}_{q}^{-q'}\simeq
 d(a)^{1+2\lambda (\alpha )}\simeq d(a)^{1+2\lambda (a)}.$ \ \par 
\quad Hence\ \par 
\quad \quad \quad $\displaystyle \forall f\in H^{q'}(\Omega ),\ \sum_{a\in S}{\frac{\left\vert{f(a)}\right\vert
 ^{q'}}{{\left\Vert{k_{a}}\right\Vert}_{q}^{q'}}}\simeq \sum_{a\in
 S}{d(a)^{1+2\lambda (a)}\left\vert{f(a)}\right\vert ^{q'}}.$ \ \par 
\ \par 
But theorem~\ref{aspc2} gives that the measure :  $\displaystyle
 \mu :=\sum_{a\in S}{d(a)^{1+2\lambda (a)}\delta _{a}}$  is a
 geometric Carleson measure. We apply the embedding Carleson
 theorem~\ref{aspcGlo50} to  $\mu $  to get\ \par 
\quad \quad \quad $\displaystyle \forall q'>1,\ \exists C_{q'}>0,\ \forall f\in
 H^{q'}(\Omega ),\ \int_{\Omega }{\left\vert{f}\right\vert ^{q'}\,d\mu
 }\leq C_{q'}^{q'}{\left\Vert{f}\right\Vert}_{H^{q'}(\Omega )}^{q'},$ \ \par 
explicitly\ \par 
\quad \quad \quad $\displaystyle \forall f\in H^{q'}(\Omega ),\ \sum_{a\in S}{\left\vert{d(a)}\right\vert
 ^{1+2\lambda (a)}\left\vert{f(a)}\right\vert ^{q'}}\leq C_{q'}^{q'}{\left\Vert{f}\right\Vert}_{H^{q'}(\Omega
 )}^{q'},$ \ \par 
hence\ \par 
\quad \quad \quad $\displaystyle \ \sum_{a\in S}{\frac{\left\vert{f(a)}\right\vert
 ^{q'}}{{\left\Vert{k_{a}}\right\Vert}_{q}^{q'}}}\lesssim {\left\Vert{f}\right\Vert}_{H^{q'}(\Omega
 )}^{q'}.$   $\blacksquare $ \ \par 

\subsection{Structural hypotheses.}
\quad We get easily the structural hypotheses~\cite{AmarExtInt06} for
 the domain  $\Omega .$ \ \par 

\begin{Crll}
~\label{ConvInt0}If  $\Omega $  is a convex domain of finite
 type in  ${\mathbb{C}}^{n},$  then the structural hypotheses
  $SH(q)$  and  $SH(p,s)$  are true for the Lebesgue measure
  $\displaystyle \sigma _{2n-1}$  on  $\displaystyle \partial
 \Omega ,$  i.e.  $\forall q\in \rbrack 1,\infty \lbrack ,$ \par 
\quad \quad \quad $\displaystyle SH(q):\ \ \ \ \ {\left\Vert{k_{a}}\right\Vert}_{q}{\left\Vert{k_{a}}\right\Vert}_{q'}\lesssim
 {\left\Vert{k_{a}}\right\Vert}_{2}^{2},$ \par 
and\!\!\!\! , for  $\displaystyle \forall p,s\in \lbrack 1,\infty
 \rbrack ,\ \frac{1}{s}=\frac{1}{p}+\frac{1}{q},$ \par 
\quad \quad \quad $\displaystyle SH(p,s):\ \ \ \ \ {\left\Vert{k_{a}}\right\Vert}_{s'}\lesssim
 {\left\Vert{k_{a}}\right\Vert}_{p'}{\left\Vert{k_{a}}\right\Vert}_{q'}.$ 
\end{Crll}
\quad Proof.\ \par 
Theorem~\ref{subPrinciple42} gives again\ \par 
\quad \quad \quad \quad \quad  $\displaystyle \ {\left\Vert{k_{a}}\right\Vert}_{H^{p}(\Omega
 )}={\left\Vert{S(a,\cdot )}\right\Vert}_{H^{p}(\Omega )}\simeq
 \frac{1}{\sigma _{2n-1}(B(\alpha ,d(a)))^{1/p'}}$ \ \par 
hence, just replacing,\ \par 
\quad \quad \quad \quad $\displaystyle \ {\left\Vert{k_{a}}\right\Vert}_{q}{\left\Vert{k_{a}}\right\Vert}_{q'}\simeq
 {\left\Vert{k_{a}}\right\Vert}_{2}^{2},\ {\left\Vert{k_{a}}\right\Vert}_{s'}\simeq
 {\left\Vert{k_{a}}\right\Vert}_{p'}{\left\Vert{k_{a}}\right\Vert}_{q'}.$
   $\blacksquare $ \ \par 
\ \par 
Now we are in position to prove theorem~\ref{1_introduction30}:\ \par 

\begin{Thrm}
If  $\Omega $  is a convex domain of finite type in  ${\mathbb{C}}^{n}$
  and if  $S\subset \Omega $  is a dual bounded sequence of points
 in  $\displaystyle H^{p}(\Omega ),$  if  $p=\infty $  then for
 any  $q<\infty ,\ S$  is  $\displaystyle H^{q}(\Omega )$  interpolating
 with the linear extension property ; if  $p<\infty $  then 
 $S$  is  $\displaystyle H^{q}(\Omega )$  interpolating with
 the linear extension property, provided that  $q<\min  (p,\ 2).$ 
\end{Thrm}
\quad Proof\ \par 
we shall apply the main theorem from~\cite{AmarExtInt06} : we
 state it in the special case of a domain  $\Omega \subset {\mathbb{C}}^{n}$
  and of the uniform algebra  $\displaystyle A(\Omega )$  of
 holomorphic functions in  $\Omega ,$  continuous up to  $\displaystyle
 \partial \Omega \ :$ \ \par 

\begin{Thrm}
~\label{ConvInt1}Let  $\displaystyle \Omega $  be a domain in
  ${\mathbb{C}}^{n}$  with  $\sigma $  the Lebesgue measure on
  $\displaystyle \partial \Omega \ ;$  if we have, with  $\displaystyle
 \ \frac{1}{s}=\frac{1}{p}+\frac{1}{q},$  that the measure  $\sigma
 $  verifies the structural hypotheses  $SH(q),\ SH(p,\ s)$  ;\par 
$\bullet \ S$  is dual bounded in  $H^{p}(\Omega )\ ;$ \par 
$\bullet \ S$  is a  $q$ -Carleson sequence ;\par 
then  $S$  is  $H^{s}(\Omega )$  interpolating and has the linear
 extension property, provided that either  $p=\infty $  or  $p\leq 2.$ 
\end{Thrm}
\quad All the requirements of Theorem~\ref{ConvInt1} are by now verified
 so we have that for any  $q<p,\ S$  is  $\displaystyle H^{q}(\Omega
 )$  interpolating with the linear extension property, provided
 that  $p=\infty $  or  $p\leq 2.$  So if  $p=\infty $  or  $p\leq
 2,$  the theorem is proved.\ \par 
\quad If  $2\leq p<\infty ,$  then  $S$  dual bounded in  $H^{p}(\Omega
 )$  means  $\exists \lbrace \rho _{a}\rbrace _{a\in S}\subset
 H^{p}(\Omega )$  with :\ \par 
\quad \quad \quad $\displaystyle \exists C>0::\forall a\in S,\ {\left\Vert{\rho
 _{a}}\right\Vert}_{p}\leq C,\ \forall a,b\in S,\ {\left\langle{\rho
 _{a},\ k_{b}}\right\rangle}=\delta _{a,b}{\left\Vert{k_{b}}\right\Vert}_{p'}.$
 \ \par 
\quad Let  $\displaystyle s::\frac{1}{2}=\frac{1}{p}+\frac{1}{s}$  then we set\ \par 
\quad \quad \quad $\displaystyle \forall a\in S,\ \tilde \rho _{a}:=\rho _{a}{\times}\frac{k_{a}}{{\left\Vert{k_{a}}\right\Vert}_{s}}\Rightarrow
 {\left\Vert{\tilde \rho _{a}}\right\Vert}_{2}\leq C$ \ \par 
and\!\!\!\! , by the reproducing property of  $k_{a}\ :$ \ \par 
\quad \quad \quad $\displaystyle \ {\left\langle{\tilde \rho _{a},\ \frac{k_{a}}{{\left\Vert{k_{a}}\right\Vert}_{2}}}\right\rangle}=\rho
 _{a}(a){\times}\frac{k_{a}(a)}{{\left\Vert{k_{a}}\right\Vert}_{s}}{\times}\frac{1}{{\left\Vert{k_{a}}\right\Vert}_{2}},$
 \ \par 
but  $k_{a}(a)={\left\Vert{k_{a}}\right\Vert}_{2}^{2},$  and
  $\rho _{a}(a)={\left\Vert{k_{a}}\right\Vert}_{p'}$  by definition,
 hence\ \par 
\quad \quad \quad $\displaystyle \ {\left\langle{\tilde \rho _{a},\ \frac{k_{a}}{{\left\Vert{k_{a}}\right\Vert}_{2}}}\right\rangle}=\frac{{\left\Vert{k_{a}}\right\Vert}_{p'}{\times}{\left\Vert{k_{a}}\right\Vert}_{2}}{{\left\Vert{k_{a}}\right\Vert}_{s}}.$
 \ \par 
The structural hypotheses, by corollary~\ref{ConvInt0}, gives\ \par 
\quad \quad \quad $\displaystyle SH(s):\ \ \ \ \ {\left\Vert{k_{a}}\right\Vert}_{s}{\left\Vert{k_{a}}\right\Vert}_{s'}\lesssim
 {\left\Vert{k_{a}}\right\Vert}_{2}^{2}$ \ \par 
and\ \par 
\quad \quad \quad $\displaystyle \forall p,s\in \lbrack 1,\infty \rbrack ,\ \frac{1}{s}=\frac{1}{p}+\frac{1}{q},\
 SH(p,s):\ \ \ \ \ {\left\Vert{k_{a}}\right\Vert}_{s'}\lesssim
 {\left\Vert{k_{a}}\right\Vert}_{p'}{\left\Vert{k_{a}}\right\Vert}_{q'},$
 \ \par 
hence here, with the right values of  $p,s\ :$ \ \par 
\quad \quad \quad $\displaystyle \ {\left\Vert{k_{a}}\right\Vert}_{2}\lesssim {\left\Vert{k_{a}}\right\Vert}_{p'}{\left\Vert{k_{a}}\right\Vert}_{s'}\leq
 {\left\Vert{k_{a}}\right\Vert}_{p'}{\times}\frac{{\left\Vert{k_{a}}\right\Vert}_{2}^{2}}{{\left\Vert{k_{a}}\right\Vert}_{s}},$
 \ \par 
the last inequality by  $SH(s),$  hence\ \par 
\quad \quad \quad $\displaystyle 1\lesssim \frac{{\left\Vert{k_{a}}\right\Vert}_{p'}{\left\Vert{k_{a}}\right\Vert}_{2}}{{\left\Vert{k_{a}}\right\Vert}_{s}}={\left\langle{\tilde
 \rho _{a},\ \frac{k_{a}}{{\left\Vert{k_{a}}\right\Vert}_{2}}}\right\rangle}.$
 \ \par 
So we have that  $S$  is dual bounded in  $H^{2}(\Omega )$  and
 by theorem~\ref{ConvInt1} we have that if  $\forall q<2$  then
  $S$  is  $H^{q}(\Omega )$  interpolating with the linear extension
 property.  $\blacksquare $ \ \par 

\begin{Rmrq}
The slight improvement from theorem~\ref{ConvInt1} done here
 relies only on the structural hypotheses, so it is in fact true
 in the abstract setting of uniform algebras.
\end{Rmrq}
\ \par 
\vfill\eject\ \par 

\section{Potential.~\label{AG1}}
\quad \quad 	Let us recall quickly how  Green formula gives us the Blaschke
 condition~\cite{zeroSkoda} :\ \par 
\quad \quad \quad $\displaystyle \ln  \left\vert{u(p)}\right\vert =\int_{\partial
 \Omega }{\ln  \left\vert{u(\zeta )}\right\vert P(p,\zeta )\,d\sigma
 (\zeta )}+\int_{\Omega }{\Delta \ln  \left\vert{u(z)}\right\vert
 G(p,z)\,dm(z)},$ \ \par 
where  $p\in \Omega ,\ G(p,z)$  is the Green kernel of  $\Omega
 $  with pole at  $p$  and  $P(p,\zeta )$  is the Poisson kernel
 of  $\Omega $  still with pole at  $p.$ \ \par 
\quad Let  $p\in \Omega $  fixed such that  $u(p)\neq 0$  and we suppose
  $u$  normalized to have\ \par 
\begin{equation} u(p)=1\Rightarrow \ln  \left\vert{u(p)}\right\vert
 =0.\end{equation}\ \par 
\quad Taking the positive Green function (minus the usual one) we have
  $G\geq 0,\ 0\leq P(p,\zeta )\leq {\left\Vert{P(p,\ \cdot )}\right\Vert}_{\infty
 },$  and we get\ \par 
\quad $\displaystyle \ \int_{\Omega }{\Delta \ln  \left\vert{u(z)}\right\vert
 G(p,z)\,dm(z)}=\int_{\partial \Omega }{\ln  ^{+}\left\vert{u}\right\vert
 P(p,\zeta )\,d\sigma }-\int_{\partial \Omega }{\ln  ^{-}\left\vert{u}\right\vert
 P(p,\zeta )\,d\sigma },$ \ \par 
\quad \quad \quad $\displaystyle \ \int_{\Omega }{\Delta \ln  \left\vert{u(z)}\right\vert
 G(p,z)\,dm(z)}\leq \int_{\partial \Omega }{\ln  ^{+}\left\vert{u(\zeta
 )}\right\vert P(p,\zeta )\,d\sigma }\leq $ \ \par 
\quad \quad \quad \quad \quad \quad \quad \quad \quad \quad \quad \quad \quad \quad \quad \quad \quad \quad \quad \quad $\leq {\left\Vert{P(p,\ \cdot )}\right\Vert}_{\infty }\int_{\partial
 \Omega }{\ln  ^{+}\left\vert{u(\zeta )}\right\vert \,d\sigma
 (\zeta )}\leq {\left\Vert{P(p,\ \cdot )}\right\Vert}_{\infty
 }{\left\Vert{u}\right\Vert}_{{\mathcal{N}}}.$ \ \par 
But  $\Delta \ln  \left\vert{u(z)}\right\vert =Tr\Theta ,$  the
 trace of  $\Theta ,$  so\ \par 
\quad \quad \quad \quad 	\begin{equation}  \ \int_{\Omega }{G(p,z)\mathrm{T}\mathrm{r}\Theta
 (z)dm(z)}\leq {\left\Vert{P(p,\ \cdot )}\right\Vert}_{\infty
 }\int_{\partial \Omega }{\ln  ^{+}\left\vert{u(\zeta )}\right\vert
 \,d\sigma (\zeta )}.\label{P1}\end{equation}\ \par 
\ \par 
\quad We have the known estimates(~\cite{Bidaut00} Prop 2.1).\ \par 

\begin{Prps}
~\label{P0}Let  $\Omega :=\lbrace x\in {\mathbb{R}}^{N}::\rho
 (x)<0\rbrace $  be a bounded domain of class  ${\mathcal{C}}^{2}$
  in   ${\mathbb{R}}^{N},$  defined by the function  $\rho $
  and  $a\in \Omega $  then there are constants  $c,\ c_{1},c_{2},$
  depending only on the regularity of  $\rho $  up to second
 order, such that, with  $P$  the Poisson kernel of  $\Omega
 ,$  with  $\displaystyle d(x)$  the distance from  $x$  to 
 $\displaystyle \partial \Omega ,$ \par 
\quad \quad \quad \quad \quad  $\displaystyle \forall (x,\zeta )\in \Omega {\times}\partial
 \Omega ,\ c_{1}\frac{d(x)}{\left\vert{\zeta -x}\right\vert ^{N}}\leq
 P(x,\zeta )\leq c_{2}\frac{d(x)}{\left\vert{\zeta -x}\right\vert ^{N}}.$ \par 
For the Green function  $G(x,z)$  of  $\Omega $  we have,\par 
\quad \quad \quad \quad \quad  $\displaystyle \forall (x,z)\in \Omega {\times}\Omega ,\ G(x,z)\geq
 c\frac{d(z)d(x)}{\left\vert{z-x}\right\vert ^{N}}.$ 
\end{Prps}
\ \par 
\quad \quad 	Using proposition~\ref{P0} we get\ \par 

\begin{Thrm}
Let  $\Omega :=\lbrace z\in {\mathbb{C}}^{n}::\rho (x)<0\rbrace
 $  be a bounded domain of class  ${\mathcal{C}}^{2}$  in  ${\mathbb{C}}^{n}\
 ;$  let  $\displaystyle G(p,z)$  be the positive Green function
 (minus the usual one) with pole  $\displaystyle p\in \Omega
 $  and  $u$  be a holomorphic function in  $\displaystyle \Omega
 $  such that  $\displaystyle u(p)=1,$  then we have\par 
\quad \quad \quad \quad \quad  $\displaystyle \ \int_{\Omega }{d(z)\mathrm{T}\mathrm{r}\Theta
 (z)dm(z)}\leq C\frac{R^{2n}}{r^{2n}}\int_{\partial \Omega }{\ln
  ^{+}\left\vert{u(\zeta )}\right\vert \,d\sigma (\zeta )}$ \par 
where  $r$  is the radius of the biggest ball  $\displaystyle
 B(p,r)$  centered at  $p$  and contained in  $\displaystyle
 \Omega $  and  $R$  is the radius of the smallest ball  $\displaystyle
 B(p,R)$  centered at  $p$  and containing  $\displaystyle \Omega
 .$  The constant  $C$  depends only on the regularity of  $\rho
 $  up to second order.
\end{Thrm}
\quad \quad 	Proof.\ \par 
We have by proposition~\ref{P0},  $\displaystyle \ {\left\Vert{P(p,\cdot
 )}\right\Vert}_{\infty }\leq c_{2}d(p)^{-2n+1}$  and  $\displaystyle
 G(p,z)\geq c\frac{d(z)d(p)}{\left\vert{z-p}\right\vert ^{2n}}$
  and, because  $\displaystyle d(p)=r$  and  $\displaystyle \
 \left\vert{z-p}\right\vert \leq R,$  we get  $\displaystyle
 \ {\left\Vert{P(p,\cdot )}\right\Vert}_{\infty }\leq c_{2}r^{-2n+1},\
 G(p,z)\geq c\frac{d(z)r}{R^{2n}}$  so, putting this in~(\ref{P1}),
 we get\ \par 
\quad \quad \quad \quad \quad  $\displaystyle c\frac{r}{R^{2n}}\int_{\Omega }{d(z)\mathrm{T}\mathrm{r}\Theta
 (z)dm(z)}\leq \frac{c_{2}}{r^{2n-1}}\int_{\partial \Omega }{\ln
  ^{+}\left\vert{u(\zeta )}\right\vert \,d\sigma (\zeta )},$ \ \par 
which gives the theorem with  $\displaystyle C:=\frac{c_{2}}{c}.$
   $\blacksquare $ \ \par 
\quad \quad 	Setting  $\displaystyle \ {\left\Vert{u}\right\Vert}_{{\mathcal{N}}(\Omega
 )}:=\int_{\partial \Omega }{\ln  ^{+}\left\vert{u(\zeta )}\right\vert
 \,d\sigma (\zeta )},$  the Nevanlinna norm of  $u,$  this prove
 that the zero set of a function in the Nevanlinna class verifies
 the Blaschke condition.\ \par 
\quad \quad 	Unfortunately the domains  $\displaystyle \Omega _{a}$  we are
 interested in have not the euclidean ball property that  $\displaystyle
 \ \frac{R}{r}\leq \gamma $   with a  $\gamma $  independent
 of  $a$  ; in fact they have it but for complex planes slices
 of  $\displaystyle \Omega _{a}$  with  $\displaystyle r,R$ 
 depending on the slices but still with  $\displaystyle \ \frac{R}{r}\leq
 \gamma ,\ \gamma $  independent of the slice, i.e. they have
 this type of property but for "ellipsoid" instead of balls.
 This is why the proofs are a little bit more involved.\ \par 

\subsection{Complex potential theory.}
\quad \quad 	In this section we shall use the notations  $\displaystyle dm:=d\sigma
 _{2n}$  for the Lebesgue measure in  ${\mathbb{C}}^{n}$  and
  $\displaystyle d\sigma $  for   $\displaystyle d\sigma _{2n-1}.$ \ \par 
\quad \quad 	Let  $\Omega $  be a domain in  ${\mathbb{C}}^{n}={\mathbb{R}}^{2n},$
  and  $u\in {\mathcal{N}}(\Omega )$  a holomorphic function
 in the Nevanlinna class of  $\Omega .$  With  $\Theta :=\Delta
 \ln \left\vert{u}\right\vert $  and  $\rho $  be a defining
 function for  $\Omega ,$  we have the lemma, application of
 the Green formula,\ \par 

\begin{Lmm}
~\label{1_P0}We have, with  $\eta $  the outward normal to  $\displaystyle
 \partial \Omega ,$ \par 
\quad \quad \quad \quad \quad  $\displaystyle \ \int_{\Omega }{(-\rho )\mathrm{T}\mathrm{r}\Theta
 dm}=\int_{\partial \Omega }{\ln \left\vert{u}\right\vert \frac{\partial
 \rho }{\partial \eta }d\sigma }-\int_{\Omega }{\ln \left\vert{u}\right\vert
 \Delta \rho dm}.$ 
\end{Lmm}
\quad \quad 	Proof.\ \par 
We have, by the Green formula,\ \par 
\quad \quad \quad \quad \quad  $\displaystyle \ \int_{\Omega }{\rho \Delta vdm}-\int_{\Omega
 }{v\Delta \rho dm}=\int_{\partial \Omega }{\rho \frac{\partial
 v}{\partial \eta }d\sigma }-\int_{\partial \Omega }{v\frac{\partial
 \rho }{\partial \eta }d\sigma }\ ;$ \ \par 
but  $\displaystyle \rho =0$  on  $\displaystyle \partial \Omega
 $  and changing sign, we get\ \par 
\quad \quad \quad \quad \quad  $\displaystyle \ \int_{\Omega }{(-\rho )\Delta vdm}=-\int_{\Omega
 }{v\Delta \rho dm}+\int_{\partial \Omega }{v\frac{\partial \rho
 }{\partial \eta }d\sigma }.$ \ \par 
\quad \quad 	Now setting  $\displaystyle v=\ln \left\vert{u}\right\vert $
  and approximating  $\displaystyle \ln \left\vert{u}\right\vert
 $  by smooth functions as usual, we get the lemma.  $\blacksquare $ \ \par 
\ \par 
\quad \quad 	The aim is to prove, under some circumstances, that we have\ \par 
\quad \quad \quad \quad \quad  $\displaystyle \ \int_{\Omega }{(-\rho )\mathrm{T}\mathrm{r}}\Theta
 \leq C\int_{\partial \Omega }{\ln ^{+}\left\vert{u}\right\vert
 d\sigma },$ \ \par 
with a {\sl good control} on  $C.$ \ \par 

\begin{Dfnt}
Let  ${\mathbb{S}}$  be the unit sphere in  ${\mathbb{C}}^{n}$
  and  $\Omega $  a domain in  ${\mathbb{C}}^{n},\ 0\in \Omega
 \ ;$  we shall say that  $\Omega $  si  ${\mathcal{C}}^{1}$
  {\bf starlike relatively to}  $0$  if  $\partial \Omega $ 
 admits a spherical parametrization, i.e. there is a function
  $R(\zeta )\in {\mathcal{C}}^{1}({\mathbb{S}}),\ R(\zeta )>0,$
  such that :\par 
\quad \quad \quad \quad \quad  $\partial \Omega =\lbrace z\in {\mathbb{C}}^{n}::\exists \zeta
 \in {\mathbb{S}},\ z=R(\zeta )\zeta \rbrace .$ 
\end{Dfnt}
This implies that  $\displaystyle \Omega =\lbrace z=tR(\zeta
 )\zeta ,\ \zeta \in {\mathbb{S}},\ t\in \lbrack 0,\ 1\lbrack \rbrace .$ \ \par 
\quad \quad 	Let  $\zeta \in {\mathbb{S}}$  and define  $\partial \Omega
 _{\zeta }$  to be the complex plane slice through  $\zeta \ :$ \ \par 
\quad \quad \quad \quad \quad  $\partial \Omega _{\zeta }:=\lbrace R(e^{i\theta }\zeta )e^{i\theta
 }\zeta ,\ \theta \in \lbrack 0,2\pi \rbrack \rbrace .$ \ \par 
The Lebesgue measure  $\displaystyle d\sigma _{\partial \Omega
 _{\zeta }}(\eta ),\ \eta =R(e^{i\theta }\zeta )e^{i\theta }\zeta
 ,$  on  $\partial \Omega _{\zeta }$  and  $d\theta $  on  $\lbrack
 0,2\pi \rbrack $  are related by\ \par 
\quad \quad \quad \quad \quad  $\displaystyle d\sigma _{\partial \Omega _{\zeta }}(\eta )={\sqrt{\left\vert{U_{\zeta
 }'(\theta )}\right\vert ^{2}+U_{\zeta }(\theta )^{2}}}d\theta ,$ \ \par 
where  $\displaystyle U_{\zeta }(\theta ):=R(e^{i\theta }\zeta ).$ \ \par 
Of course if  $\eta =e^{i\varphi }\zeta $  then  $\partial \Omega
 _{\eta }=\partial \Omega _{\zeta }$  and the measure is the same.\ \par 
\quad \quad 	We set  $\displaystyle \ {\left\Vert{U_{\zeta }'}\right\Vert}_{\infty
 }:=\sup _{\ \theta \in \lbrack 0,2\pi \rbrack }\left\vert{U_{\zeta
 }'(\theta )}\right\vert .$  We shall use the notations\ \par 
\quad \quad \quad \quad \quad  $\forall \zeta \in {\mathbb{S}},\ d_{\zeta }(0)=\inf _{\theta
 \in \lbrack 0,2\pi \rbrack }R(e^{i\theta }\zeta )\ ;\ d_{\zeta
 max}(0)=\sup _{\theta \in \lbrack 0,2\pi \rbrack }R(e^{i\theta
 }\zeta ).$ \ \par 
Now we have.\ \par 

\begin{Lmm}
~\label{Green50}Let  ${\mathbb{S}}$  be the unit sphere of  ${\mathbb{C}}^{n}$
  and  $\Omega $  a domain in  ${\mathbb{C}}^{n},\ 0\in \Omega
 $  such  $\Omega $  is  ${\mathcal{C}}^{1}$  starlike with respect
 to  $0.$  Then\par 
\quad \quad \quad \quad \quad \quad  $\displaystyle \ \int_{\Omega }{f(z)dm(z)}=c_{n}\frac{1}{2\pi
 }\int_{0}^{1}{(\int_{{\mathbb{S}}}{(\int_{\partial \Omega _{\zeta
 }}{f(t\eta )J_{\zeta }(\eta )t^{2n-1}d\sigma _{\zeta }(\eta
 )})d\sigma _{{\mathbb{S}}}(\zeta )})dt},$ \par 
with  $\displaystyle \ \frac{d_{\zeta }(0)^{2n}}{\ {\sqrt{{\left\Vert{U_{\zeta
 }'}\right\Vert}_{\infty }^{2}+d_{\zeta max}(0)^{2}}}}\leq J_{\zeta
 }(\eta )\leq \frac{d_{\zeta max}(0)^{2n}}{d_{\zeta }(0)}.$ 
\end{Lmm}
\quad \quad 	Proof.\ \par 
Integrating in spherical coordinates we get, with  $c_{n}=2nv_{n}/s_{n}$
  where  $v_{n}$  is the volume of the unit ball in  ${\mathbb{C}}^{n},\
 s_{n}$  the area of the unit sphere in  ${\mathbb{C}}^{n},$ \ \par 
\quad \quad \quad \quad \quad  $\displaystyle \ I:=\int_{\Omega }{f(z)dm(z)}=c_{n}\int_{{\mathbb{S}}}{\lbrace
 \int_{0}^{R(\zeta )}{r^{2n-1}f(r\zeta )dr}\rbrace d\sigma _{2n-1}(\zeta
 )}.$ \ \par 
Set  $\displaystyle t=\frac{r}{R(\zeta )}\Rightarrow dr=R(\zeta )dt$  and\ \par 
\quad \quad \quad \quad  $\displaystyle \forall \zeta \in {\mathbb{S}},\ \int_{0}^{R(\zeta
 )}{r^{2n-1}f(r\zeta )dr}=\int_{0}^{1}{R(\zeta )^{2n}f(tR(\zeta
 )\zeta )t^{2n-1}dt}.$ \ \par 
\quad \quad \quad \quad \quad  $\displaystyle I=c_{n}\frac{1}{2\pi }\int_{{\mathbb{S}}{\times}\lbrack
 0,1\rbrack {\times}\lbrack 0,2\pi \rbrack }{R(e^{i\theta }\zeta
 )^{2n}f(tR(e^{i\theta }\zeta )e^{i\theta }\zeta )t^{2n-1}dtd\theta
 d\sigma _{2n-1}(\zeta )}.$ \ \par 
Now we fix  $\zeta \in {\mathbb{S}},$  we get\ \par 
\quad \quad \quad \quad \quad  $\displaystyle \ \int_{\lbrack 0,1\rbrack {\times}\lbrack 0,2\pi
 \rbrack }{R(e^{i\theta }\zeta )^{2n}f(tR(e^{i\theta }\zeta )e^{i\theta
 }\zeta )t^{2n-1}dtd\theta }.$ \ \par 
Set  $\zeta \in {\mathbb{S}},\ \forall \theta \in \lbrack 0,2\pi
 \rbrack ,\ \eta =R(e^{i\theta }\zeta )e^{i\theta }\zeta \in
 \partial \Omega _{\zeta }$  and  $U_{\zeta }(\theta ):=R(e^{i\theta
 }\zeta )$  then we have\ \par 
\quad \quad \quad \quad \quad  $\partial \Omega _{\zeta }=\lbrace U_{\zeta }(\theta )e^{i\theta
 }\zeta ,\ \theta \in \lbrack 0,2\pi \rbrack \rbrace $ \ \par 
and\ \par 
\quad \quad \quad \quad \quad  $\displaystyle d\sigma _{\partial \Omega _{\zeta }}(\eta )={\sqrt{\left\vert{U_{\zeta
 }'(\theta )}\right\vert ^{2}+U_{\zeta }(\theta )^{2}}}d\theta ,$ \ \par 
so\ \par 
\quad \quad \quad \quad \quad  $\displaystyle \ \int_{\lbrack 0,2\pi \rbrack }{R(e^{i\theta
 }\zeta )^{2n}f(tR(e^{i\theta }\zeta )e^{i\theta }\zeta )d\theta
 }=\int_{\partial \Omega _{\zeta }}{f(t\eta )J_{\zeta }(\eta
 )d\sigma _{\zeta }(\eta )},$ \ \par 
where  $\displaystyle \ J_{\zeta }(\eta )=\frac{R(\eta )^{2n}}{D_{\zeta
 }(\eta )}$  and  $\displaystyle D_{\zeta }(\eta )={\sqrt{\left\vert{U_{\zeta
 }'(\theta )}\right\vert ^{2}+U_{\zeta }(\theta )^{2}}},$  expressed
 in  $\eta $  coordinates.\ \par 
\quad \quad 	So we have\ \par 
\quad \quad \quad \quad  $\displaystyle I=c_{n}\frac{1}{2\pi }\int_{0}^{1}{(\int_{{\mathbb{S}}}{(\int_{\partial
 \Omega _{\zeta }}{f(t\eta )J_{\zeta }(\eta )t^{2n-1}d\sigma
 _{\zeta }(\eta )})d\sigma _{{\mathbb{S}}}(\zeta )})dt}.$ \ \par 
Notice that in  $\partial \Omega _{\zeta }$  we have\ \par 
\quad \quad \quad \quad \quad  $\displaystyle d_{\zeta }(0)\leq R(\eta )\leq d_{\zeta max}(0)\
 ;\ d_{\zeta }(0)\leq D_{\zeta }(\eta )\leq {\sqrt{{\left\Vert{U_{\zeta
 }'}\right\Vert}_{\infty }^{2}+d_{\zeta max}(0)^{2}}}$ \ \par 
hence\ \par 
\quad \quad \quad \quad \quad  $\displaystyle \ \frac{d_{\zeta }(0)^{2n}}{\ {\sqrt{{\left\Vert{U_{\zeta
 }'}\right\Vert}_{\infty }^{2}+d_{\zeta max}(0)^{2}}}}\leq J_{\zeta
 }(\eta )\leq \frac{d_{\zeta max}(0)^{2n}}{d_{\zeta }(0)}.$ 
  $\blacksquare $ \ \par 

\begin{Dfnt}
~\label{CarlDomain31}The domain   $\Omega \subset {\mathbb{C}}^{n}$
  is  said to be  $\gamma $  {\bf balanced relatively to}  $0$  if :\par 
\quad \quad  $\bullet $   $\Omega $  is  ${\mathcal{C}}^{1}$  starlike with
 respect to  $0,$ \par 
\quad \quad  $\bullet $  all its slices  $\partial \Omega _{\zeta }$  through
 the origin verify\par 
\quad \quad \quad \quad \quad \quad \quad  $\forall \zeta \in {\mathbb{S}},\ d_{\zeta max}(0)\leq \gamma
 d_{\zeta }(0)\ ;$ 		 	 $\ {\left\Vert{U_{\zeta }'}\right\Vert}_{\infty
 }\leq \gamma d_{\zeta max}(0).$ 
\end{Dfnt}
\quad Set for any function  $\displaystyle v,\ v^{+}(z):=\max (v(z),0)\
 ;\ v^{-}(z):=-\max (-v(z),0).$ \ \par 
Then we have the lemmas.\ \par 

\begin{Lmm}
Suppose that  $v$  is a sub harmonic function in a  $\gamma $
  balanced domain  $D$  in  ${\mathbb{C}},$  such that  $v(0)=0,$  then\par 
\quad \quad \quad \quad \quad  $\displaystyle \ \int_{\partial D}{v^{-}(z)d\sigma (z)}\leq
 \gamma ^{2}\frac{c_{2}}{c_{1}}\int_{\partial D}{v^{+}(z)d\sigma (z)}.$ 
\end{Lmm}
\quad \quad 	Proof.\ \par 
Because  $v$  is sub harmonic we have\ \par 
\quad \quad \quad \quad \quad  $\displaystyle 0=v(0)\leq \int_{\partial D}{P(0,\zeta )v(\zeta
 )d\sigma (\zeta )}=\int_{\partial D}{P(0,\zeta )v^{+}(\zeta
 )d\sigma (\zeta )}-\int_{\partial D}{P(0,\zeta )v^{-}(\zeta
 )d\sigma (\zeta )},$ \ \par 
where  $P(0,\zeta )$  is the Poisson kernel of  $D$  for  $0\in D.$  So\ \par 
\quad \quad \quad \quad \quad  $\displaystyle \ \int_{\partial D}{P(0,\zeta )v^{-}(\zeta )d\sigma
 (\zeta )}\leq \int_{\partial D}{P(0,\zeta )v^{+}(\zeta )d\sigma
 (\zeta )}.$ \ \par 
Now we use the estimates in proposition~\ref{Green21}  $\displaystyle
 \ \frac{c_{1}d(0)}{d_{max}(0)^{2}}\leq \frac{c_{1}d(0)}{\left\vert{\zeta
 }\right\vert ^{2}}\leq P(0,\zeta )\leq \frac{c_{2}}{\left\vert{\zeta
 }\right\vert }\leq \frac{c_{2}}{d(0)}$  to get\ \par 
\quad \quad \quad \quad \quad  $\displaystyle \ \frac{c_{1}d(0)}{d_{max}(0)^{2}}\int_{\partial
 D}{v^{-}(\zeta )d\sigma (\zeta )}\leq \frac{c_{2}}{d(0)}\int_{\partial
 D}{v^{+}(\zeta )d\sigma (\zeta )},$ \ \par 
hence\ \par 
\quad \quad \quad \quad \quad  $\displaystyle \ \int_{\partial D}{v^{-}(\zeta )d\sigma (\zeta
 )}\leq \frac{c_{2}d_{max}(0)^{2}}{c_{1}d(0)^{2}}\int_{\partial
 D}{v^{+}(\zeta )d\sigma (\zeta )}\leq \gamma ^{2}\frac{c_{2}}{c_{1}}\int_{\partial
 D}{v^{+}(\zeta )d\sigma (\zeta )}.$   $\blacksquare $ \ \par 

\begin{Lmm}
~\label{Green51} Let  $\Omega \subset {\mathbb{C}}^{n}$  be a
  $\gamma $  balanced domain and let  $v$  be a pluri sub harmonic
 function in  $\Omega $  such that  $v(0)=0,$  then\par 
\quad \quad \quad \quad \quad  $\displaystyle \ \int_{\Omega }{\left\vert{v(z)}\right\vert
 dm(z)}\leq (1+2\gamma ^{2n+3}\frac{c_{2}}{c_{1}})\int_{\Omega
 }{v^{+}(z)dm(z)}$ 
\end{Lmm}
\quad \quad 	Proof.\ \par 
We shall use the decomposition of lemma~\ref{Green50}\ \par 
\quad \quad \quad \quad \quad  $\displaystyle \ \int_{\Omega }{v^{-}(z)dm(z)}=c_{n}\frac{1}{2\pi
 }\int_{0}^{1}{(\int_{{\mathbb{S}}}{(\int_{\partial \Omega _{\zeta
 }}{v^{-}(t\eta )J_{\zeta }(\eta )t^{2n-1}d\sigma _{\zeta }(\eta
 )})d\sigma _{{\mathbb{S}}}(\zeta )})dt}.$ \ \par 
But still by lemma~\ref{Green50} we have  $\displaystyle J_{\zeta
 }(\eta )\leq \frac{d_{\zeta max}(0)^{2n}}{d_{\zeta }(0)}$  hence\ \par 
\quad \quad \quad \quad \quad  $\displaystyle \ \int_{\partial \Omega _{\zeta }}{v^{-}(t\eta
 )J_{\zeta }(\eta )d\sigma _{\zeta }(\eta )}\leq \frac{d_{\zeta
 max}(0)^{2n}}{d_{\zeta }(0)}\int_{\partial \Omega _{\zeta }}{v^{-}(t\eta
 )d\sigma _{\zeta }(\eta )}.$ \ \par 
Doing the same we get\ \par 
\quad \quad \quad \quad \quad  $\displaystyle \ \int_{\partial \Omega _{\zeta }}{v^{+}(t\eta
 )J_{\zeta }(\eta )d\sigma _{\zeta }(\eta )}\geq \frac{d_{\zeta
 }(0)^{2n}}{\ {\sqrt{{\left\Vert{U_{\zeta }'}\right\Vert}_{\infty
 }^{2}+d_{\zeta max}(0)^{2}}}}\int_{\partial \Omega _{\zeta }}{v^{+}(t\eta
 )d\sigma _{\zeta }(\eta )}.$ \ \par 
\quad \quad 	Set\ \par 
\quad \quad \quad \quad \quad  $\displaystyle A:=\frac{d_{\zeta max}(0)^{2n}}{d_{\zeta }(0)}\
 ;\ B:=\frac{d_{\zeta }(0)^{2n}}{\ {\sqrt{{\left\Vert{U_{\zeta
 }'}\right\Vert}_{\infty }^{2}+d_{\zeta max}(0)^{2}}}}$ \ \par 
then\ \par 
\quad \quad \quad \quad \quad  $\displaystyle \ \int_{\partial \Omega _{\zeta }}{v^{-}(t\eta
 )J_{\zeta }(\eta )d\sigma _{\zeta }(\eta )}\leq A\int_{\partial
 \Omega _{\zeta }}{v^{-}(t\eta )d\sigma _{\zeta }(\eta )}$ \ \par 
and\ \par 
\quad \quad \quad \quad \quad  $\displaystyle \ \int_{\partial \Omega _{\zeta }}{v^{+}(t\eta
 )J_{\zeta }(\eta )d\sigma _{\zeta }(\eta )}\geq B\int_{\partial
 \Omega _{\zeta }}{v^{+}(t\eta )d\sigma _{\zeta }(\eta )}.$ \ \par 
\quad \quad 	But lemma~\ref{Green51} gives, because  $v$  being pluri sub
 harmonic in  $\Omega $  is sub harmonic in  $\Omega _{\zeta },$ \ \par 
\quad \quad \quad \quad \quad  $\displaystyle \ \int_{\partial \Omega _{\zeta }}{v^{-}(t\eta
 )d\sigma _{\zeta }(\eta )}\leq \gamma ^{2}\frac{c_{2}}{c_{1}}\int_{\partial
 \Omega _{\zeta }}{v^{+}(t\eta )d\sigma _{\zeta }(\eta )},$ \ \par 
so\ \par 
\quad \quad \quad \quad \quad  $\displaystyle \ \int_{\partial \Omega _{\zeta }}{v^{-}(t\eta
 )J_{\zeta }(\eta )d\sigma _{\zeta }(\eta )}\leq A\int_{\partial
 \Omega _{\zeta }}{v^{-}(t\eta )d\sigma _{\zeta }(\eta )}\leq
 A\gamma ^{2}\frac{c_{2}}{c_{1}}\int_{\partial \Omega _{\zeta
 }}{v^{+}(t\eta )d\sigma _{\zeta }(\eta )},$ \ \par 
hence continuing\ \par 
\quad \quad \quad \quad \quad  $\displaystyle \ \int_{\partial \Omega _{\zeta }}{v^{-}(t\eta
 )J_{\zeta }(\eta )d\sigma _{\zeta }(\eta )}\leq A\gamma ^{2}\frac{c_{2}}{c_{1}}\int_{\partial
 \Omega _{\zeta }}{v^{+}(t\eta )d\sigma _{\zeta }(\eta )}\leq $ \ \par 
\quad \quad \quad \quad \quad \quad \quad \quad \quad \quad  $\displaystyle \leq \frac{Ac_{2}}{Bc_{1}}\gamma ^{2}\int_{\partial
 \Omega _{\zeta }}{v^{+}(t\eta )J_{\zeta }(\eta )d\sigma _{\zeta
 }(\eta )}.$ \ \par 
So\ \par 
\quad \quad \quad \quad \quad  $\displaystyle \ \int_{\partial \Omega _{\zeta }}{v^{-}(t\eta
 )J_{\zeta }(\eta )d\sigma _{\zeta }(\eta )}\leq \frac{Ac_{2}}{Bc_{1}}\gamma
 ^{2}\int_{\partial \Omega _{\zeta }}{v^{+}(t\eta )J_{\zeta }(\eta
 )d\sigma _{\zeta }(\eta )}.$ \ \par 
Now we notice that\ \par 
\quad \quad \quad \quad \quad  $\displaystyle \ \frac{A}{B}=\frac{d_{\zeta max}(0)^{2n}\ {\sqrt{{\left\Vert{U_{\zeta
 }'}\right\Vert}_{\infty }^{2}+d_{\zeta max}(0)^{2}}}}{d_{\zeta
 }(0)^{2n+1}}\leq \gamma ^{2n+1}{\sqrt{1+\frac{\ {\left\Vert{U_{\zeta
 }'}\right\Vert}_{\infty }^{2}}{d_{\zeta max}(0)^{2}}}}\leq \gamma
 ^{2n+1}{\sqrt{\gamma ^{2}+1}}.$ \ \par 
Multiplying by  $t^{2n-1}$  and integrating on  ${\mathbb{S}}{\times}\lbrack
 0,1\rbrack $  give\ \par 
\quad \quad \quad \quad \quad  $\displaystyle \ \int_{\Omega }{v^{-}(z)dm(z)}\leq 2\gamma ^{2n+3}\frac{c_{2}}{c_{1}}\int_{\Omega
 }{v^{+}(z)dm(z)}\ ;$ \ \par 
but  $\displaystyle \ \left\vert{v(z)}\right\vert =v^{+}(z)+v^{-}(z)$
  hence\ \par 
\quad \quad \quad \quad \quad  $\displaystyle \ \int_{\Omega }{\left\vert{v(z)}\right\vert
 dm(z)}\leq (1+2\gamma ^{2n+3}\frac{c_{2}}{c_{1}})\int_{\Omega
 }{v^{+}(z)dm(z)}.$   $\blacksquare $ \ \par 

\begin{Lmm}
~\label{aspcGlo77}Let  $\Omega $  be a domain in  ${\mathbb{C}}^{n}$
  of class  ${\mathcal{C}}^{2},$  if  $v$  is a  positive sub
 harmonic function in  $\Omega ,$  then\par 
\quad \quad \quad \quad \quad  $\displaystyle \ \int_{\Omega }{v(z)dm(z)}\leq 2c_{2}\mathrm{d}\mathrm{i}\mathrm{a}\mathrm{m}(\Omega
 )\int_{\partial \Omega }{v(\zeta )d\sigma (\zeta )}.$ 
\end{Lmm}
\quad \quad 	Proof.\ \par 
Let  $P(z,\zeta )$  the Poisson kernel we have, by proposition~\ref{P0}, \ \par 
\quad \quad \quad \quad \quad  $\displaystyle \forall (z,\zeta )\in \Omega {\times}\partial
 \Omega ,\ c_{1}\frac{d(z)}{\left\vert{\zeta -z}\right\vert ^{2n}}\leq
 P(z,\zeta )\leq c_{2}\frac{d(z)}{\left\vert{\zeta -z}\right\vert
 ^{2n}},$ \ \par 
so, because  $d(x)\leq \left\vert{\zeta -x}\right\vert $  we get\ \par 
\quad \quad \quad \quad \quad  $\displaystyle \forall (z,\zeta )\in \Omega {\times}\partial
 \Omega ,\ P(z,\zeta )\leq c_{2}\frac{1}{\left\vert{\zeta -z}\right\vert
 ^{2n-1}}.$ \ \par 
Hence\ \par 
\quad \quad \quad  $\displaystyle \forall \zeta \in \partial \Omega ,\ \int_{\Omega
 }{P(z,\zeta )dm(z)}\leq c_{2}\int_{\Omega }{\frac{dm(z)}{\left\vert{\zeta
 -z}\right\vert ^{2n-1}}}\leq c_{2}\int_{B(0,\ \mathrm{d}\mathrm{i}\mathrm{a}\mathrm{m}(\Omega
 ))}{\frac{dm(z)}{\left\vert{\zeta -z}\right\vert ^{2n-1}}}\leq
 2c_{2}\mathrm{d}\mathrm{i}\mathrm{a}\mathrm{m}(\Omega ).$ \ \par 
\quad \quad 	Because  $v$  is sub harmonic we have\ \par 
\quad \quad \quad \quad \quad  $\displaystyle v(z)\leq \int_{\partial \Omega }{P(z,\zeta )v(\zeta
 )d\sigma (\zeta )},$ \ \par 
so, by Fubini Tonnelli, everything being positive,\ \par 
\quad \quad \quad \quad \quad  $\displaystyle \ \int_{\Omega }{v(z)dm(z)}\leq \int_{\Omega
 {\times}\partial \Omega }{P(z,\zeta )v(\zeta )dm(z)d\sigma (\zeta
 )}\leq 2c_{2}\mathrm{d}\mathrm{i}\mathrm{a}\mathrm{m}(\Omega
 )\int_{\partial \Omega }{v(\zeta )d\sigma (\zeta )}.$   $\blacksquare $ \ \par 

\begin{Prps}
~\label{aspcGlo78}Let  $\Omega $  be a domain in  ${\mathbb{C}}^{n}$
  of class  ${\mathcal{C}}^{2},\ \gamma $  balanced relatively
 to  $0\in \Omega \ ;$  if  $v$  is pluri sub harmonic in  $\Omega
 $  and  $v(0)=0$  then\par 
\quad \quad \quad \quad \quad  $\displaystyle \ \int_{\Omega }{\left\vert{v(z)}\right\vert
 dm(z)}\leq 2c_{2}\mathrm{d}\mathrm{i}\mathrm{a}\mathrm{m}(\Omega
 )(2\gamma ^{2n+3}\frac{c_{2}}{c_{1}}+1)\int_{\partial \Omega
 }{v^{+}(\zeta )d\sigma (\zeta )}.$ 
\end{Prps}
\quad \quad 	Proof.\ \par 
We apply successively lemma~\ref{Green51} and lemma~\ref{aspcGlo77},
 which can be done because  $v^{+}$  is still pluri sub harmonic
 in  $\Omega .$   $\blacksquare $ \ \par 

\begin{Thrm}
~\label{8_CarlDomain3}Let  $\Omega $  be a domain in  ${\mathbb{C}}^{n}$
  of class  ${\mathcal{C}}^{2},\ \gamma $  balanced relatively
 to  $0\in \Omega \ ;$  if  $u$  is holomorphic in  $\Omega $
  and  $\ \left\vert{u(0)}\right\vert =1$  then, with  $X:=u^{-1}(0),$ \par 
\quad \quad \quad \quad \quad  $\ {\left\Vert{\Theta _{X}}\right\Vert}_{B}:=\int_{\Omega }{d(z)\mathrm{T}\mathrm{r}}\Theta
 \leq C\int_{\partial \Omega }{\ln ^{+}\left\vert{u}\right\vert
 d\sigma (\zeta )}=:C{\left\Vert{u}\right\Vert}_{{\mathcal{N}}(\Omega )},$ \par 
with a constant  $\displaystyle C$  depending only on the constant
  $\gamma $  and the derivatives of  $\rho $  up to order  $\displaystyle 2.$ 
\end{Thrm}
\quad \quad 	Proof.\ \par 
By lemma~\ref{1_P0} we have\ \par 
\quad \quad \quad \quad \quad  $\displaystyle \ \int_{\Omega }{(-\rho )\mathrm{T}\mathrm{r}}\Theta
 =\int_{\partial \Omega }{\ln \left\vert{u}\right\vert \frac{\partial
 \rho }{\partial \eta }d\sigma }-\int_{\Omega }{\ln \left\vert{u}\right\vert
 \Delta \rho dm}.$ \ \par 
The function  $\displaystyle \ln \left\vert{u}\right\vert $ 
 is pluri sub harmonic in  $\Omega $  hence we can apply to it
 proposition~\ref{aspcGlo78} :\ \par 
\quad \quad \quad \quad \quad  $\displaystyle \ \int_{\Omega }{\left\vert{\ln \left\vert{u(z)}\right\vert
 }\right\vert dm(z)}\leq 2c_{2}\mathrm{d}\mathrm{i}\mathrm{a}\mathrm{m}(\Omega
 )(2\gamma ^{2n+3}\frac{c_{2}}{c_{1}}+1)\int_{\partial \Omega
 }{\ln ^{+}\left\vert{u(\zeta )}\right\vert d\sigma (\zeta )}.$ \ \par 
so we get, because  $\displaystyle 0<\frac{\partial \rho }{\partial
 \eta },$ \ \par 
\quad \quad \quad \quad \quad  $\displaystyle \ \int_{\Omega }{(-\rho )\mathrm{T}\mathrm{r}}\Theta
 =\int_{\partial \Omega }{\ln \left\vert{u}\right\vert \frac{\partial
 \rho }{\partial \eta }d\sigma }-\int_{\Omega }{\ln \left\vert{u}\right\vert
 \Delta \rho dm}\leq {\left\Vert{\frac{\partial \rho }{\partial
 \eta }}\right\Vert}_{\infty }\int_{\partial \Omega }{\ln ^{+}\left\vert{u}\right\vert
 d\sigma }+{\left\Vert{\Delta \rho }\right\Vert}_{\infty }\int_{\Omega
 }{\left\vert{\ln \left\vert{u(z)}\right\vert }\right\vert dm(z)}\leq $ \ \par 
\quad \quad \quad \quad \quad \quad \quad \quad \quad \quad  $\displaystyle \leq A\int_{\partial \Omega }{\ln ^{+}\left\vert{u(\zeta
 )}\right\vert d\sigma (\zeta )},$ \ \par 
with  $\displaystyle A:={\left\Vert{\frac{\partial \rho }{\partial
 \eta }}\right\Vert}_{\infty }+{\left\Vert{\Delta \rho }\right\Vert}_{\infty
 }(2c_{2}\mathrm{d}\mathrm{i}\mathrm{a}\mathrm{m}(\Omega )(2\gamma
 ^{2n+3}\frac{c_{2}}{c_{1}}+1)).$ \ \par 
\quad \quad 	But  $\displaystyle \ \frac{\partial \rho }{\partial \eta }(z)d(z)\simeq
 (-\rho (z))$  so, with  $\displaystyle M:={\left\Vert{\frac{1}{\partial
 \rho /\partial \eta }}\right\Vert}_{\infty },$  we get\ \par 
\quad \quad \quad \quad \quad  $\displaystyle \ \int_{\Omega }{d(z)\mathrm{T}\mathrm{r}}\Theta
 \leq M\int_{\Omega }{(-\rho )\mathrm{T}\mathrm{r}}\Theta .$ \ \par 
This proves the theorem with  $\displaystyle C=MA,$  a constant
 depending only on  $\gamma $  and the derivatives of  $\rho
 $  up to order  $\displaystyle 2.$   $\blacksquare $ \ \par 

\begin{Rmrq}
~\label{1_P1}This theorem will be applied to the domains  $\displaystyle
 \Omega _{a}$  built in section~\ref{6_CarlDomain33} and for
 these domains the derivatives of the defining function  $\displaystyle
 \rho _{a}$  are controlled by the derivatives of the global
 function  $\rho $  ; also the derivative  $\displaystyle \ \frac{\partial
 \rho _{a}}{\partial \eta }$  is bounded below uniformly independently
 of  $a$  by  $\displaystyle \ \frac{\partial \rho }{\partial
 \eta }$  so the constant  $C$  of theorem~\ref{8_CarlDomain3}
 is  independent of  $\displaystyle a.$ 
\end{Rmrq}
\ \par 

\bibliographystyle{/usr/local/texlive/2013/texmf-dist/bibtex/bst/base/plain}

\end{document}